\documentclass[11pt]{amsart}

\usepackage{amsmath, amssymb, amsthm}
\usepackage{mathrsfs}
\usepackage[margin=1in]{geometry}
\usepackage{xcolor}
\usepackage[expansion=false]{microtype}
\usepackage{enumitem}
\usepackage{booktabs}
\usepackage{array}
\usepackage{tikz}
\usetikzlibrary{arrows.meta,calc}
\usepackage{xurl}
\usepackage{hyperref}
\hypersetup{pdftitle={The Singular LeBrun Transform: Complex Conformal and Contact Geometry},
  pdfauthor={Yun-Gang Ye},
  pdfsubject={November 13, 2026},
  pdfkeywords={LeBrun-Salamon conjecture, quaternion-Kahler, contact Fano manifold, conformal geometry, twistor}}

\newtheorem{theorem}{Theorem}[section]
\newtheorem{proposition}[theorem]{Proposition}
\newtheorem{lemma}[theorem]{Lemma}
\newtheorem{corollary}[theorem]{Corollary}

\makeatletter\let\c@equation\c@theorem\makeatother

\newtheorem*{theoremD}{Theorem D}
\newtheorem*{theoremH}{Theorem H}
\newtheorem*{theoremJ}{Theorem J}
\newtheorem*{theoremQ}{Theorem Q}
\theoremstyle{definition}
\newtheorem{definition}[theorem]{Definition}
\newtheorem{remark}[theorem]{Remark}
\newtheorem{question}[theorem]{Question}
\newtheorem{example}[theorem]{Example}

\newcommand{\MC}{M_{\mathbb{C}}}
\newcommand{\MCbar}{\overline{M}_{\mathbb{C}}}
\newcommand{\CC}{\mathbb{C}}
\newcommand{\RR}{\mathbb{R}}
\newcommand{\ZZ}{\mathbb{Z}}
\newcommand{\HH}{\mathbb{H}}
\newcommand{\PP}{\mathbb{P}}
\newcommand{\OO}{\mathcal{O}}
\newcommand{\Gr}{\mathrm{Gr}}
\newcommand{\Sp}{\mathrm{Sp}}

\DeclareMathOperator{\im}{im}

\makeatletter

\def\part{\@startsection{part}{0}%
  \z@{2\linespacing\@plus\linespacing}{1.2\linespacing}%
  {\normalfont\large\bfseries\centering\def\@secnumfont{\bfseries}}}
\makeatother

\title[The Singular LeBrun Transform]
{The Singular LeBrun Transform:\\ Complex Conformal and Contact Geometry}

\author{Yun-Gang Ye}
\address{Associate Research Professor, Duke University, NC 27705, USA}
\email{david.ye@duke.edu}

\date{November 13, 2026}

\begin{document}

\begin{abstract}
A holomorphic conformal structure is rigid in the presence of rational curves: a complex conformal manifold that is not conformally flat contains no immersed
rational null geodesic. So a natural compactification of such a manifold, on which its null geodesics close up, must degenerate along the boundary. We study this
phenomenon for the complex conformal manifold that LeBrun attached to a complex contact manifold $Z$ in inverting Salamon's twistor construction: the space
$\MC$ of balanced conics of $Z$ transverse to the contact distribution, with its holomorphic Riemannian metric of holonomy contained in $\Sp(n,\CC)\Sp(1,\CC)$.
Its compactification by the degenerate conics of $Z$ forces the metric to become singular at the boundary. We develop a framework, the \emph{singular LeBrun
transform}, in which these singularities are a source of information rather than an obstacle: the singularities of the metric and of its curvature are read off
from the degenerate rational curves of $Z$, locally through the contact neighbourhood of a curve and globally through the incidence variety.

As an application, we give an alternative proof of the LeBrun--Salamon conjecture that all positive quaternionic-K\"ahler manifolds are symmetric,
which was recently claimed by OpenAI \cite{OAI26} and by Hu, Liu and Wan \cite{HLW26}. The emphasis of the paper, however, is less on this result than on
the framework developed here, which should apply to other problems of the same kind; we conclude with two such directions, negative quaternionic-K\"ahler
manifolds and a singular Merkulov transform for Legendre moduli spaces.
\end{abstract}

\maketitle

\setcounter{tocdepth}{1}
\tableofcontents

\section{Introduction}
\label{sec:intro}

\subsection{Singular complex conformal geometry}
\label{sec:introidea}

This paper grew out of a project that began more than thirty years ago, in the early 1990s, as an attempt to understand complex contact geometry and complex
conformal geometry. A \emph{complex contact manifold} is a complex manifold of odd dimension with a maximally non-integrable holomorphic distribution of
corank one, and a \emph{complex conformal manifold} is a complex manifold $M$ (not necessarily compact) with a \emph{holomorphic conformal structure}: a holomorphic section of
$S^{2}TM^{*}\otimes N$, for some line bundle $N$ on $M$, that is nondegenerate at every point. The hyperquadric $\mathbb{Q}^{m}$, with the conformal
structure induced by its defining quadratic form, is an example. What fascinated the author then, and motivates this paper, are the intricate relations between the geodesics of conformal geometry and the extremal
rational curves of algebraic geometry.

A holomorphic Riemannian metric, or a holomorphic conformal structure, on a complex manifold is far more rigid than its real counterpart, and most rigid in the
presence of rational curves. In \cite{Ye94} the author showed, using Mori's theory of extremal rays, that a projective manifold of dimension $m\ge3$ with a
holomorphic conformal structure and canonical bundle that is not nef is the hyperquadric $\mathbb{Q}^{m}$: by the cone theorem it contains rational curves, a minimal
rational curve through a general point is a null geodesic, and a projective conformal manifold containing an immersed rational null geodesic is conformally flat. Belgun
\cite{Bel01} removed the compactness: a complex conformal manifold that contains an immersed rational null geodesic is conformally flat.

So a complex conformal manifold that is not conformally flat contains no immersed rational null geodesic. If it has a natural compactification on which its null geodesics close up,
its structure must degenerate along the boundary. Projectivity alone does not provide such curves; for the space of conics of this paper, the compactification and its
curve geometry do (\S\ref{sec:mustdegenerate}). The point of view of this paper is that this degeneration is a source of
information rather than an obstacle: the singularities that the metric and its curvature acquire at the boundary encode the global geometry of the metric
and of the manifold.

\emph{Singular conformal structures.} We make this precise as follows (Definition~\ref{def:singconf}). Let $\overline M\supset M$ be a compactification of a complex conformal manifold $M$, an
irreducible complex space of dimension $m$ with smooth locus $\overline M_{\mathrm{reg}}$; there may be more than one way to compactify $M$. A \emph{conformal structure} on $\overline M$ is a holomorphic section
\[
g\in H^{0}\big(\overline M_{\mathrm{reg}},\,S^{2}T\overline M^{*}\otimes N\big),
\]
for some line bundle $N$ on $\overline M_{\mathrm{reg}}$, that is nondegenerate over a dense open set: $\det g$, a section of
$K_{\overline M}^{2}\otimes N^{m}$ over $\overline M_{\mathrm{reg}}$, is not identically zero. Sections that differ by a meromorphic factor, a meromorphic section of $N'\otimes N^{-1}$, define the same structure, and we take $g$ vanishing on no hypersurface.
The structure is \emph{regular} at $p$ if $\overline M$ is smooth at $p$ and $g(p)$ is nondegenerate, and \emph{singular} on $\overline M$ if it is not regular
everywhere. At a point $p$ where it is not regular, one of two things happens:
\begin{enumerate}[leftmargin=2em]
\item[(i)] $\overline M$ is smooth at $p$ and $g(p)$ is \emph{degenerate}: $p$ lies on the hypersurface $\{\det g=0\}$;
\item[(ii)] $\overline M$ itself is \emph{singular} at $p$.
\end{enumerate}
Both are intrinsic to the structure. A metric in it, $h=g/\sigma$ with $\sigma$ a meromorphic section of $N$, has in addition \emph{poles}, along the zero divisor of
$\sigma$, and these depend on $\sigma$: multiplying $h$ by a power of a local equation of a hypersurface trades a pole along it for a holomorphic form, possibly
degenerate. We call $h$ a \emph{singular conformal metric} on $\overline M$ if it does not extend to a holomorphic nondegenerate metric on all of $\overline M$. In this
paper $\overline M$ is the compactified space of conics of a contact manifold and $h$ is LeBrun's metric $g^{\CC}=G/\sigma_{\theta}$, the quotient of a quadratic form
$G$ with values in a line bundle by a section $\sigma_{\theta}$ of that line bundle (\S\ref{sec:introcorr}, Proposition~\ref{prop:levwronsk}).

\subsection{The LeBrun transform}
\label{sec:introcorr}

The space we study comes from twistor theory. Salamon \cite{Salamon82} attached to a quaternionic-K\"ahler manifold $(M^{4n},g)$ of positive scalar curvature
its twistor space $q\colon Z\to M$, a compact Fano manifold of dimension $2n+1$ with a holomorphic contact structure and a real structure, whose fibres, the twistor
lines, are rational curves. The holonomy $Sp(n)Sp(1)$ corresponds to a local factorization $TM\otimes\CC\cong E^{*}\otimes H^{*}$, with $E$ and $H$ of ranks $2n$
and $2$; see \cite{Salamon82} for details. LeBrun \cite{LeBrun89} inverted this construction. Specifically, let $Z$ be any complex contact manifold of dimension $2n+1$, with
contact line bundle $L$. The smooth rational curves of $Z$ of $L$-degree two (\emph{conics}) that are transverse to the contact distribution and have normal bundle
$\OO(1)^{2n}$, which we call \emph{balanced transverse conics}, form a complex manifold $\MC$ of dimension $4n$ carrying a holomorphic Riemannian metric
$g^{\CC}$, Einstein, with holonomy contained in $\Sp(n,\CC)\Sp(1,\CC)$. A point $z\in Z$ appears as the submanifold $\Sigma_{z}\subset\MC$ of the conics
through $z$, totally null and totally geodesic, a quaternionic $\alpha$-plane. If $Z$ has a real structure $\sigma$ fixing a balanced transverse conic, the
$\sigma$-invariant conics form a totally real slice of $\MC$ on which $g^{\CC}$ restricts to a real metric; for the twistor space of $M$ the twistor lines are
$\sigma$-invariant, and $M$ is recovered from $\MC$ as a component of this real slice, with $g^{\CC}|_{M}=g$; likewise the holomorphic curvature $R^{\CC}$
of $g^{\CC}$ restricts to the curvature of $g$, $R^{\CC}|_{M}=R$.

The curvature of $g^{\CC}$ is
\[
R^{\CC}=\rho_{0}R_{0}+R_{E}^{\CC},
\]
where $\rho_{0}\ne0$ is constant, $R_{0}$ is the curvature of quaternionic projective space, and $R_{E}^{\CC}$, a section of $S^{4}E_{\CC}$ for the factorization
$T\MC=E^{*}_{\CC}\otimes H^{*}_{\CC}$, is the \emph{quaternionic Weyl curvature}; on the real slice, $R_{E}^{\CC}|_{M}=R_{E}$ is the quaternionic Weyl curvature
of $g$. We call the correspondence through the incidence variety
\[
\MC\ \xleftarrow{\ \pi\ }\ F=\{([C],z):z\in C\}\ \xrightarrow{\ p\ }\ Z,
\]
with the transformation of tensors it induces, the \emph{LeBrun transform}.

One advantage of working with $R_{E}^{\CC}$ rather than $R_{E}$ is that $R_{E}^{\CC}$ is automatically holomorphic, so the tools of complex geometry apply to it.
On $M$ itself, $\nabla R_{E}=0$ would already imply the LeBrun--Salamon conjecture, since it makes $M$ locally symmetric. But by Salamon's description of the
complex structure of $Z$ \cite{Salamon82}, $\nabla R_{E}=0$ is equivalent to $q^{*}R_{E}$ being a holomorphic section of $q^{*}S^{4}E$ over $Z$, which is not
given in advance. On $\MC$ holomorphicity is given from the start, and the question becomes whether $R_{E}^{\CC}$ extends over the compactification
$\MCbar$ described below.

When $Z$ is compact and $L$ is ample, $(\MC,g^{\CC})$ is exactly the situation of \S\ref{sec:introidea}. The metric is nowhere conformally flat for $n\ge2$, and $\MC$ contains no
compact curve at all (Propositions~\ref{prop:notflat} and~\ref{prop:nocurves}). Its natural compactification is the closure $\MCbar$ of a component of $\MC$ in the
Hilbert scheme of $Z$, a component of the Hilbert scheme whose general member is a balanced transverse conic, and it is projective. So the structure degenerates:
$g^{\CC}$ is a singular conformal metric on $\MCbar$. The boundary consists of the degenerate conics of $Z$: pairs of lines, conics tangent to the contact
distribution (\emph{integral} conics), double lines and cuspidal curves; along it the metric acquires poles and degenerates, and $\MCbar$ may be singular. The
remaining points of $\MCbar\smallsetminus\MC$ form the locus $J$ of \emph{jumping} conics, smooth transverse conics whose normal bundle is not $\OO(1)^{2n}$;
along $J$ the factorization of $T\MC$ degenerates. The \emph{singular} LeBrun transform extends the
correspondence to these curves, and turns the singularities of $(g^{\CC},R_{E}^{\CC})$ into the geometry of the degenerate rational curves of the contact manifold
$Z$. It has a local side, where the singularity at a degenerate conic is read off from the contact neighbourhood of the curve in $Z$, the curvature at a conic
being the invariant of order three of its neighbourhood (Theorem~\ref{thm:thirdorder}); and a global side, where the compactness of $\MCbar$ acts, through the
incidence variety, on the geometry of $Z$.

\subsection{Main results}
\label{sec:intromain}

The main objective of this paper is to establish a framework for the correspondence between complex conformal geometry and complex contact geometry. The
correspondence was found by LeBrun for rational curves on a contact manifold \cite{LeBrun89}, and extended by Merkulov to compact Legendre submanifolds
\cite{Mer94,Mer97}. It lets one study algebraic geometry from the point of view of holomorphic Riemannian geometry, the complexification of classical
Riemannian geometry. Its most informative part is the boundary, where the conformal structure degenerates: there the metric and its curvature record the
algebraic geometry of the contact manifold, through its lines, its degenerate conics and the jumping of normal bundles. How the metric and its curvature
degenerate at the boundary is the main concern of this paper, and we call the correspondence together with its boundary the \emph{singular LeBrun
transform}. The proof of the LeBrun--Salamon conjecture below is one application of the framework, and its first test; Merkulov's setting, where the same philosophy leads to a
\emph{singular Merkulov transform}, is taken up in \S\ref{sec:boundaryconf}.

For quaternionic-K\"ahler geometry, the holonomy of LeBrun's metric $g^{\CC}$ lies in $\Sp(n,\CC)\Sp(1,\CC)$, the complexification of $Sp(n)Sp(1)$, the
quaternionic-K\"ahler entry of the metric part of Berger's list of holonomies of torsion-free connections \cite{Ber55}. The non-metric part of the list, in the
real and in the holomorphic setting, turned out to be incomplete: Bryant found the first exotic holonomies by twistor methods \cite{Bry91}; Merkulov showed that
every torsion-free connection with irreducibly acting holonomy can in principle be constructed on a moduli space of Legendre submanifolds \cite{Mer95,Mer97}; and
Merkulov and Schwachh\"ofer completed the classification \cite{MS99}; see Bryant's survey \cite{Bry00}. The link between this conformal geometry and
Fano contact manifolds is the conjecture of LeBrun and Salamon \cite{LS94}, which the framework proves as an application.

For a complex simple Lie algebra $\mathfrak g$, the \emph{adjoint variety} of $\mathfrak g$ is the closed orbit of the adjoint group in $\PP(\mathfrak g)$, the orbit
of the line spanned by a highest root vector. It is a homogeneous Fano contact manifold, and every homogeneous Fano contact manifold is of this form
\cite{Wolf65,Bea98}. The Wolf spaces \cite{Wolf65} are the compact quaternionic-K\"ahler symmetric spaces, one for each compact simple Lie group, and the twistor space of
a Wolf space is the adjoint variety of the complexified Lie algebra. LeBrun and Salamon conjectured that there are no other compact quaternionic-K\"ahler manifolds of positive scalar curvature, and Beauville
\cite{Bea98}, more generally, that every Fano contact manifold is homogeneous. The first conjecture was known for $n\le4$ \cite{Hit81,PS91,BWW22}; a recent preprint claims an extension to $n\le14$, that is to real dimension at most $56$
\cite{DGLS26}. By
\cite{LS94} the twistor space has $\mathrm{Pic}\,Z=\ZZ L$ unless $M$ is $\HH P^{n}$ or $\Gr_{2}(\CC^{n+2})$. In real dimension four, ``quaternionic-K\"ahler'' means
self-dual Einstein, and the identification of a positive quaternionic-K\"ahler metric with a Wolf-space metric is up to homothety.

In our approach both conjectures rest on one statement about the transform: \emph{for a compact contact Fano manifold $Z$ of dimension $2n+1\ge7$ with
$\mathrm{Pic}\,Z=\ZZ L$, and a component $\MCbar$ of the space of conics whose general member is a balanced transverse conic, the quaternionic Weyl curvature
$R_{E}^{\CC}$ is parallel}. Such components exist for every such $Z$: the nodal conic of two general lines through
a general point is balanced, and its general smoothing is a balanced transverse conic (Corollary~\ref{cor:balexist}). The statement is proved in two steps, which isolate the one place where the jumping conics enter. For a general
point $z\in Z$, let $\Gamma_{z}$ be the normalized graph compactification of the tangent-direction map $[C]\mapsto[T_{z}C]$ on the conics through $z$
(\S\ref{sec:fibre}), and $\widetilde\mu_{z}\colon\Gamma_{z}\to\PP(T_{z}Z)$ the induced proper morphism. A \emph{jumping divisor through $z$} is a prime divisor of $\Gamma_{z}$ whose general member is a jumping conic.

\begin{theoremD}[Theorem~\ref{thm:LSJ}]
Let $Z$ be a compact contact Fano manifold of dimension $2n+1$, $n\ge3$, with $\mathrm{Pic}\,Z=\ZZ L$, and $\MCbar$ a component of its Hilbert scheme of conics
whose general member is a balanced transverse conic. Suppose that the jumping conics are harmless: for a general $z$, no jumping divisor through $z$ is mapped by
$\widetilde\mu_{z}$ onto a divisor of $\PP(T_{z}Z)$. This holds in particular if no conic through a general point is jumping. Then $R_{E}^{\CC}$ is parallel.
\end{theoremD}

We borrowed two ideas from the OpenAI preprint \cite{OAI26}, the two-point vector lift and its measurement along test curves, to help remove the
condition in the above theorem for $R_{E}^{\CC}$ to be parallel (Propositions~\ref{prop:tplifts} and~\ref{prop:tpconst}).

\begin{theoremJ}[Theorem~\ref{thm:tpD}, Lemma~\ref{lem:twopoint} and Proposition~\ref{prop:cfromJ}]
Let $Z$ and $\MCbar$ be as in Theorem~D, with $n\ge2$. Then exactly one conic of $\MCbar$ passes through two general points of $Z$. Consequently, for a general
$z$, the map $\widetilde\mu_{z}$ is birational and contracts every jumping divisor through $z$, so the jumping conics are harmless.
\end{theoremJ}

Theorem~D is where the curvature is extended across the boundary (Parts~\ref{part:local} and~\ref{part:global}). Theorem~J is a count of conics, and is the
only place where the jumping conics have to be dealt with. Once $R_{E}^{\CC}$ is parallel, the rest is short (\S\ref{sec:introstrategy}). It gives Beauville's
conjecture for $\mathrm{Pic}\,Z=\ZZ L$,

\begin{theoremH}[Corollary~\ref{cor:beauville}]
Every compact contact Fano manifold $Z$ of dimension $2n+1\ge7$ with $\mathrm{Pic}\,Z=\ZZ L$ is homogeneous: it is the adjoint variety of a complex simple Lie
algebra.
\end{theoremH}

\noindent and, through the real structure of the twistor space, the conjecture of LeBrun and Salamon:

\begin{theoremQ}[Theorem~\ref{thm:main}]
Every compact quaternionic-K\"ahler manifold of positive scalar curvature is a Wolf space.
\end{theoremQ}

\subsection{The strategy of the proofs}
\label{sec:introstrategy}

\emph{The idea.} The application asks whether the complex conformal manifold $(\MC,g^{\CC})$ is locally symmetric, that is, whether the quaternionic Weyl
curvature is parallel:
\[
\nabla_{X}R_{E}^{\CC}=0\qquad\text{for every }X\in T\MC.
\]
Once this holds, the local symmetry of $\MC$ is carried by the transform to $Z$, the homogeneity extends to all of $Z$, and $Z$ is an adjoint variety. The proof
that $R_{E}^{\CC}$ is parallel rests on two observations.
\begin{itemize}[leftmargin=2em]
\item \emph{It suffices to look along one leaf.} For $z\in Z$, the balanced transverse conics through $z$ form a submanifold $\Sigma_{z}\subset\MC$, totally null and totally
geodesic. By the leaf principle (Proposition~\ref{prop:leaf}), $R_{E}^{\CC}$ is parallel if and only if, for a general point $z$, the curvature read at $z$ is
constant on $\Sigma_{z}$. Read at $z$, the curvature takes its values in one fixed vector space: it is a holomorphic map, the \emph{curvature quartic},
\[
W_{z}\colon\Sigma_{z}\longrightarrow S^{4}L^{\perp*}_{z}\otimes L_{z}^{2}\cong S^{4}\CC^{2n},
\]
into an affine space.
\item \emph{Compactness makes it constant.} $\Sigma_{z}$ is not compact, but its normalized graph compactification $\Gamma_{z}$ (\S\ref{sec:fibre}) is a projective variety. A holomorphic
map from a compact connected variety into an affine space is constant. So if $W_{z}$ extends holomorphically to the compactification, it is constant on
$\Sigma_{z}$, for a general $z$, and $R_{E}^{\CC}$ is parallel.
\end{itemize}
Everything else concerns this extension, that is, the behaviour of the curvature at the degenerate conics through $z$.

\emph{Where the extension is needed.} One need not extend $W_{z}$ across all of $\Gamma_{z}$. The tangent direction at $z$, $[C]\mapsto[T_{z}C]$, extends to a proper,
generically finite map $\widetilde\mu_{z}\colon\Gamma_{z}\to\PP(T_{z}Z)\cong\PP^{2n}$, a local isomorphism along $\Sigma_{z}$. A holomorphic map from $\PP^{2n}$ minus a set
of codimension two into a vector space is constant, by Hartogs. So $W_{z}$ is constant as soon as it has no pole along any divisor of $\Gamma_{z}$ whose image under
$\widetilde\mu_{z}$ is still a divisor. The proof is the bookkeeping of these divisors: for each kind, \emph{either the curvature extends across it, or
$\widetilde\mu_{z}$ contracts it, or it does not pass through a general point}.

\emph{The divisors that survive.} There are the jumping conics and three kinds of boundary point; every other degenerate conic through $z$ spans a divisor that
$\widetilde\mu_{z}$ contracts.
\begin{itemize}[leftmargin=2em]
\item \emph{Jumping conics.} There $\widetilde\mu_{z}$ is not a local isomorphism, so it contracts them once it is birational. Its degree is at most the number
$d_{2}$ of conics through two general points, and $d_{2}=1$ (Theorem~J, \S\ref{sec:twopointD}).
\item \emph{Pairs of lines with the node at $z$.} The curvature quartic, a function on the blown-up universal curve, extends across the generic point of the
exceptional divisor over them (\S\ref{sec:dilation}); on $\MCbar$ it becomes, after a double cover, a statement in a frame obtained by evaluating at a point of the neck (Corollary~\ref{cor:Rnodal}). The idea is a symmetry: the
boundary points are limits of interior points under a group of Heisenberg dilations, and every deformation of the contact neighbourhood of a nodal conic has
non-negative weight. The germ is compared with a finite normal form through the polar part of the curvature, which only a finite-order neighbourhood of
the conic sees.
\item \emph{Integral conics with normal bundle $\OO(1)^{2n}$.} The metric has a pole, but the curvature extends, and vanishes (\S\ref{sec:integralII}).
\item \emph{Conics with a cusp at $z$.} They do not pass through a general point (\S\ref{sec:cuspsII}). {The argument uses only the Levi form and the
$L$-degree: a cusp moving in a dominating family would fold its curve two-to-one onto its image by an explicit involution.}
\end{itemize}
So compactness acts on the curvature through three local problems, at the nodes, the integral conics and the cusps, and one global count. Table~\ref{tab:dictionary}
records how each kind of degenerate conic is disposed of. The three local problems concern the transform itself rather than the conjecture: they describe how the
metric and its curvature degenerate at each kind of boundary point, and they do not see the sign of the scalar curvature (\S\ref{sec:negQK}).

\emph{From parallel curvature to homogeneity.} Since $\rho_{0}$ is constant and $R_{0}$ is parallel, $g^{\CC}$ is locally symmetric, and $\rho_{0}\ne0$ identifies
the local model with the complexification of a Wolf space. On a twistor space the real slice $M$ is then a compact, simply connected, locally symmetric space,
hence a Wolf space. For a general $Z$ the compactness of $Z$ plays the role that geodesic completeness would play for $\MC$. The transform carries the local
isometries of $g^{\CC}$ to a developing map from the union of the balanced transverse conics to an adjoint variety, whose complement in $Z$ has codimension two because every
divisor meets every conic. The infinitesimal contact automorphisms of the adjoint variety pull back by it and extend by Hartogs, and the orbit of a general point
is a compact homogeneous space, hence all of $Z$ (Theorem~\ref{thm:homog}).

\subsection{Origins of this work and the use of AI}
\label{sec:introhistory}

In the early 1990s, using Mori's theory of extremal rays, the author showed that a complex projective contact threefold is either $\PP^{3}$ or the projectivized tangent bundle $\PP(T_{S})$ of a smooth projective
surface $S$ \cite{Ye94b}, and set out to extend this to Fano contact manifolds of higher dimension, which was later called Beauville's conjecture. Since a Fano contact
manifold is a twistor space exactly when it carries a K\"ahler--Einstein metric \cite{LeB95}, it contains the conjecture of LeBrun and Salamon, and is equivalent to it
for Fano contact manifolds with a K\"ahler--Einstein metric.

The idea of the present paper took shape in 1993, in a series of manuscripts \cite{Ye93,Ye93b}, shared with a few colleagues at the time (see \cite{LM01,LM07} for example) but not
published, that also contain several partial results. It was
first presented in a talk more than thirty years later, in April 2026, at Chongqing University of Technology, at the invitation of Kefeng Liu,
one of the authors of \cite{HLW26}; the author thanks him for his hospitality. In preparing that talk, the author reread these manuscripts.
Seeing many of his 1993 partial results confirmed in the subsequent literature was an encouragement to return to a problem
that had never ceased to interest him, and to take up the work where he had left it.

Since then, the ideas of the manuscripts have been developed into the present paper with the help of AI.
Claude (Anthropic) has been the main AI tool for the project, for reasons of cost and speed mostly Claude Opus 5.5 and occasionally Claude Fable 5.1, while
OpenAI's models have served as a referee for parts of the paper. The author first had Claude transcribe the manuscripts into \LaTeX. He then prompted it with ideas from them and worked with it interactively, one question and one
response at a time. No AI agent was run autonomously: the model served as an assistant, chiefly for computations, for concrete examples, and for writing out the
details of proof ideas. The author will describe this work and his experience with AI in a separate paper.

The author believes that the main value of this paper to the mathematical community lies less in the theorems themselves than in the paths it opens. The
correspondence between conformal geometry, complex holonomy groups and special algebraic varieties is rich; \S\ref{sec:future} begins to
explore it. He hopes to bring some intuition to this corner of mathematics, which is now being inspired by AI. His own experience is that AI, used
appropriately, can enhance and inspire mathematical research.

\subsection{Relation to other approaches}
\label{sec:introrelation}

Appendix~\ref{app:survey} summarizes related work since the early 1990s. Two proofs have recently been claimed, both for every contact Fano manifold, that is, for
Beauville's conjecture: by OpenAI \cite{OAI26} and by Hu, Liu and Wan \cite{HLW26}. Both aim at the space $H^{0}(Z,L)$, the Lie algebra of infinitesimal contact automorphisms
\cite[Prop.~2.1]{LeB95}, \cite[\S1]{Bea98}: once the contact vector fields span $T_{z}Z$ at every point, the contact automorphism group is transitive and $Z$ is an
adjoint variety.

The two arguments start in the same way, from the two-point tensor $(z-w)^{2}\partial_{z}\otimes\partial_{w}$ on the space of conics, contracted with the contact
form at one point. The OpenAI preprint constructs from it a symmetric section of $L\boxtimes L$ on $Z\times Z$ and extends it across the boundary of the
compactified space of conics; Hu, Liu and Wan extend the same tensor with its first point fixed. From there the two proofs diverge.
\begin{itemize}[leftmargin=2em]
\item In the OpenAI preprint the section $B\in H^{0}(Z\times Z,L\boxtimes L)$ has symbol $\theta^{2}$ along the diagonal. The diagonal jets of $B$ show that $H^{0}(Z,L)$
generates the first jets of $L$, and the contact moment map identifies $Z$ with an adjoint variety.
\item In the preprint of Hu, Liu and Wan the extended tensor, with its first point fixed, gives tangent sections on incidence varieties, hence local vector fields
that preserve the variety of minimal rational tangents of the contact lines; these show that this variety is irreducible and nonlinear, and a descent theorem of
Hwang, cited there, makes the vector fields global, so that they generate $TZ$.
\end{itemize}
In both, the object extended is a tensor with a closed formula, and the metric of $M$ is recovered only at the end. The present approach, which goes back to
\cite{Ye93}, extends instead the curvature $R_{E}$: a quartic rather than a quadratic form, whose behaviour at the boundary is more delicate, but from which
$\nabla R_{E}^{\CC}=0$ follows directly. This reflects the point of view of the paper: the behaviour of the curvature at the boundary, recorded in the dictionary of
Table~\ref{tab:dictionary}, is the subject, and the conjecture is derived from it. No vector fields are constructed in proving that the curvature is parallel: on the real slice Cartan's theorem and the
compactness of $M$ give Theorem~Q, and for a general $Z$ the global contact vector fields are obtained afterwards, pulled back from the adjoint variety by a
developing map (\S\ref{sec:qk}). On an adjoint variety the two
objects extended are the quartic and the quadratic invariants of $\mathfrak g$: Freudenthal's quartic (Remark~\ref{rem:Wolfconverse}) and the Killing form.

\subsection{Organization}

The paper has four parts: the first three develop the singular LeBrun transform, and the fourth applies it.

\emph{Part~\ref{part:transform}, the LeBrun transform,} sets up the objects. For a contact manifold with balanced transverse conics it recalls LeBrun's
holomorphic conformal structure on the space of conics, compactifies that space, and shows that the structure must degenerate at the boundary
(\S\S\ref{sec:lebrun}--\ref{sec:mustdegenerate}). It then introduces the tools that carry the curvature from the space of conics to $Z$: lines and the variety of
minimal rational tangents, and the incidence variety with its tangent-direction map $\widetilde\mu_{z}$, through which the curvature at a point $z$ becomes a
function $W_{z}$ on a compact space $\Gamma_{z}$ (\S\S\ref{sec:lines}--\ref{sec:incidence}).

\emph{Part~\ref{part:local}, the local theory,} studies the transform near each kind of degenerate conic. The curvature is a third-order invariant of the
contact neighbourhood of a conic (\S\ref{sec:formal}). With this, $W_{z}$ is shown to stay regular at the generic nodal conics (\S\ref{sec:dilation}) and at the
balanced integral conics (\S\ref{sec:integralII}), and cusps are excluded at a general point (\S\ref{sec:cuspsII}). Section~\ref{sec:dictionaryII} collects
the outcome in a dictionary between degenerate conics and singularities of the transform.

\emph{Part~\ref{part:global}, the global theory,} turns local regularity into a global statement. A Hartogs argument through $\widetilde\mu_{z}$ shows that
$W_{z}$ is constant on each leaf $\Sigma_{z}$, the \emph{leaf constancy} of the curvature, once the jumping conics are harmless (\S\ref{sec:descent}); they are, because exactly one conic passes
through two general points (\S\ref{sec:twopointD}).

\emph{Part~\ref{part:real}, consequences,} draws them from the parallel curvature: $Z$ is homogeneous, and Theorem~Q follows (\S\ref{sec:qk}). On an
adjoint variety the curvature quartic is Freudenthal's quartic invariant, and its zero locus is the tangent variety, and the dual variety, of the VMRT (Corollary~\ref{cor:curvhyp}); this links the curvature of the Wolf spaces to the Legendrian and Severi varieties of
Freudenthal's magic square (Remark~\ref{rem:magicsquare}). Section~\ref{sec:future} proposes two directions for future research, negative quaternionic-K\"ahler manifolds and a singular Merkulov transform for
Legendre moduli spaces, and Appendix~\ref{app:survey} summarizes related work.

\emph{The logic of the application.} The proof of Theorem~Q has two independent inputs, the regularity at the boundary and the uniqueness of the conic
through two general points:
\begin{center}
\begin{tikzpicture}[x=1cm,y=1cm,every node/.style={align=center}]
\node at (0,0) {Regularity at nodes and integral conics, no cusps at a general point\\
(Theorems~\ref{thm:Eregular}, \ref{thm:cuspgeneral} and Proposition~\ref{prop:intregular})\\[2pt]
$+$\quad uniqueness of the conic through two general points (Theorem~\ref{thm:tpD})};
\node at (0,-1.25) {$\Big\Downarrow$};
\node[anchor=west,font=\small] at (0.4,-1.25) {Theorem~\ref{thm:Ereduction}, Proposition~\ref{prop:cfromJ}};
\node at (0,-2.0) {Leaf constancy of the curvature quartic (Theorem~\ref{thm:LSJ})};
\node at (0,-2.7) {$\Big\Downarrow$};
\node at (0,-3.4) {$\nabla R_{E}^{\CC}=0$};
\node at (0,-4.1) {$\Big\Downarrow$};
\node[anchor=west,font=\small] at (0.4,-4.1) {Corollary~\ref{cor:endgame}};
\node at (0,-4.8) {LeBrun--Salamon conjecture (Theorem~Q)};
\end{tikzpicture}
\end{center}
Section~\ref{sec:twopointD} uses neither the curvature nor Theorem~\ref{thm:LSJ}, and Section~\ref{sec:future} and Appendix~\ref{app:survey}
are not used in the proof.

\emph{Conventions.} All varieties are over $\CC$. For a vector space or bundle $V$, $\PP(V)$ is the space of lines in $V$; so $\PP(T_{z}Z)$ is the space of
tangent directions at $z$. ``General'' means in a dense Zariski-open subset, which may shrink finitely many times in the course of an argument.

\part{The LeBrun Transform}
\label{part:transform}

\section{LeBrun's complex conformal structure on the space of conics}
\label{sec:lebrun}

\subsection{Contact manifolds and balanced transverse conics}
\label{sec:contactZ}

A \emph{complex contact manifold} is a complex manifold $Z$ of dimension $2n+1$ with a corank-one subbundle $L^{\perp}\subset TZ$, with quotient line bundle
$L=TZ/L^{\perp}$ and contact form $\theta\in H^{0}(Z,\Omega^{1}_{Z}\otimes L)$, $\ker\theta=L^{\perp}$, such that $d\theta$ restricted to $L^{\perp}$ is a nondegenerate
$L$-valued $2$-form, the \emph{Levi form}, written $\omega$. (On $L^{\perp}$ the form $d\theta$ does not depend on a local trivialization of $L$.) One has
$K_{Z}=L^{-(n+1)}$. A \emph{balanced transverse conic} is a smooth rational curve $C\subset Z$ with $L\cdot C=2$, transverse to $L^{\perp}$, with $N_{C/Z}\cong\OO(1)^{\oplus2n}$.
\begin{itemize}[leftmargin=2em]
\item \emph{Projective space.} Let $V$ be a symplectic vector space of dimension $2n+2$ with form $\Omega$. On $\PP^{2n+1}=\PP(V)$ the form $\theta_{v}=\Omega(v,dv)$
defines a contact structure with $L=\OO(2)$. A line $\PP(W)$ is a balanced transverse conic exactly when $\Omega|_{W}\ne0$, and the balanced transverse conics form the open subset of
$\Gr_{2}(V)$ of symplectic planes, the complexification of $\HH P^{n}$, whose quaternionic Weyl curvature vanishes.
\item \emph{Twistor spaces.} If $Z$ is the twistor space of a quaternionic-K\"ahler manifold $M$ of positive scalar curvature, $Z$ is a compact Fano contact
manifold and the twistor lines are balanced transverse conics (\S\ref{sec:twistor}).
\item \emph{Adjoint varieties.} For a simple complex Lie algebra $\mathfrak g$ the closed orbit of the adjoint group in $\PP(\mathfrak g)$ is a homogeneous Fano contact
manifold, with $\mathrm{Pic}\,Z=\ZZ L$ unless $\mathfrak g$ is of type $A$ or $C$; it is the twistor space of the corresponding Wolf space \cite{Wolf65}.
\item \emph{Germs.} LeBrun's construction is local along the curves, so it applies to a germ of a contact manifold along a balanced transverse conic, for instance to the twistor
space of a quaternionic-K\"ahler germ.
\end{itemize}
For a compact contact Fano manifold, $L$ is ample. Where the Picard group matters we assume
\[
\mathrm{Pic}\,Z=\ZZ L ;
\]
then $Z\not\cong\PP^{2n+1}$, and $Z$ is covered by lines (\S\ref{sec:vmrt}).

\subsection{Lines, conics and the space of conics}
\label{sec:conics}
\label{sec:component}

A \emph{line} is a rational curve $\ell\subset Z$ with $L\cdot\ell=1$, a \emph{conic} a connected curve of
arithmetic genus zero with $L\cdot C=2$, i.e.\ a subscheme with Hilbert polynomial $2t+1$ with respect to $L$.

\begin{lemma}[All or nothing]
\label{lem:allornothing}
Let $C\subset Z$ be a smooth conic. The composite $TC\to TZ|_{C}\to L|_{C}$ lies in
$\mathrm{Hom}(\OO(2),\OO(2))=\CC$. Hence $C$ is either transverse to $L^{\perp}$ at every point, or \emph{integral},
$\theta|_{C}\equiv0$. Every line is integral, since $\mathrm{Hom}(\OO(2),\OO(1))=0$.
\end{lemma}

For a transverse conic, $TZ|_{C}=TC\oplus L^{\perp}|_{C}$, so $N_{C/Z}\cong L^{\perp}|_{C}$. Since $L|_{C}\cong\OO(2)$ has
the square root $\OO(1)$ on $C$ we put
\[
B|_{C}:=L^{\perp}|_{C}\otimes L^{-1/2}|_{C},\qquad N_{C/Z}=B|_{C}\otimes\OO(1).
\]
The Levi form makes $B|_{C}$ a symplectic bundle (with values in $L\otimes L^{-1}=\OO$) of degree zero, so
\begin{equation}
B|_{C}\cong\bigoplus_{i=1}^{n}\big(\OO(k_{i})\oplus\OO(-k_{i})\big),\qquad k_{i}\ge0 .\label{eq:Bsplit}
\end{equation}
The conic is \emph{balanced} if $B|_{C}$ is trivial, that is $N_{C/Z}\cong\OO(1)^{\oplus2n}$.

\emph{The space of conics.} The local constructions below apply to complex contact manifolds and their analytic germs along curves; whenever Hilbert
schemes, projective compactifications or rational tensors are used, $Z$ is projective with $L$ ample, as it is for a contact Fano manifold. For such $Z$,
$\MCbar$ denotes an irreducible component of the Hilbert scheme of $Z$ (Hilbert polynomial $2t+1$ with respect to $L$) whose general member is a balanced transverse conic; such components exist when $Z$ is Fano with
$\mathrm{Pic}\,Z=\ZZ L$ and $n\ge3$ (Corollary~\ref{cor:balexist}). If $Z$ is a twistor space we take the component containing the twistor lines. Let $\MC\subset\MCbar$ be the open subset of balanced transverse conics. The degenerate members of $\MCbar$, those that are not smooth transverse
conics, form the \emph{boundary}; the other members of $\MCbar\smallsetminus\MC$ are the jumping conics (\S\ref{sec:jumpdef}).

\begin{proposition}[The general conic]
\label{prop:genconic}
Let $[C]\in\MC$ and $z\in C$.
\begin{enumerate}
\item $H^{1}(N_{C/Z})=0$, so $\MCbar$ is smooth of dimension $h^{0}(\OO(1)^{2n})=4n$ at $[C]$, and $\MC$ is a dense open subset of $\MCbar$.
\item $H^{1}(N_{C/Z}(-z))=0$ and $H^{0}(N_{C/Z}(-2z))=0$. So the conics of $\MCbar$ through $z$ form, near $[C]$, a smooth family of dimension $2n$,
and the map sending a conic through $z$ to its tangent direction at $z$ is a local isomorphism at $[C]$ onto an open subset of $\PP(T_{z}Z)$.
\item For $w\in C\smallsetminus\{z\}$, the conics through $z$ near $C$ sweep out a neighbourhood of $w$. Consequently, for a general $z\in Z$ the conics of
$\MCbar$ through $z$ cover $Z$, and through two general points of $Z$ pass finitely many conics of $\MCbar$, all of them smooth, transverse and balanced.
\end{enumerate}
\end{proposition}

\begin{proof}
(1) $N\cong\OO(1)^{2n}$ has $H^{1}=0$, and balancedness is open. (2) $N(-z)\cong\OO^{2n}$ and $N(-2z)\cong\OO(-1)^{2n}$. The tangent space to the conics
through $z$ is $H^{0}(N(-z))$ and the obstructions lie in $H^{1}(N(-z))=0$; the differential of the tangent-direction map at $[C]$ is the derivative at $z$
of a section of $N(-z)$, whose kernel $H^{0}(N(-2z))$ vanishes (Lemma~\ref{lem:dmu} below). (3) $N(-z)\cong\OO^{2n}$ is globally generated, so evaluation
$H^{0}(N(-z))\to N_{w}$ is onto, and the incidence variety of conics through $z$ with a marked point maps submersively to $Z$ near $([C],w)$. For a general $z$
the general conic through $z$ is balanced, so its image is dense, hence all of $Z$, the family being proper; the map from the conics through $z$ with a
marked point to $Z$ is then generically finite between varieties of dimension $2n+1$.
\end{proof}

\emph{Degenerate conics.} The members of $\MCbar$ that are not smooth transverse balanced conics are of two kinds. The boundary members are
the \emph{integral} conics, on which $\theta$ vanishes (Lemma~\ref{lem:allornothing}), the \emph{nodal conics}, two lines meeting transversally in one
point, and further degenerations (double lines, images of non-embedded maps), described in \S\ref{sec:P2}. The closure in $\MCbar$ of the nodal conics
whose branches have nef normal bundles is the \emph{nodal divisor} $D_{\mathrm{nod}}$, each of whose components is a divisor by
Lemma~\ref{lem:nodalsmooth} below. The second kind are the smooth transverse conics that are not balanced, the \emph{jumping
conics} of \S\ref{sec:jumpdef}. The integral and nodal conics are detected by the transversality section $\sigma_{\theta}$ of the following theorem.

\begin{theorem}[The transversality section]
\label{thm:thetasection}
Let $U^{\circ}\subset\MCbar$ be the open subset of conics that are smooth or nodal (two smooth rational curves meeting
transversally in one point), $\pi\colon\mathcal C\to U^{\circ}$ the universal curve and $ev\colon\mathcal C\to Z$ the evaluation map.
\begin{enumerate}
\item $\Lambda:=\pi_{*}(\omega_{\mathcal C/U^{\circ}}\otimes ev^{*}L)$ is a line bundle on $U^{\circ}$; on every fibre $C$,
$h^{0}(\omega_{C}\otimes L|_{C})=1$ and $h^{1}=0$.
\item The contact form induces a section $\sigma_{\theta}\in H^{0}(U^{\circ},\Lambda)$, the push-forward of the image of $ev^{*}\theta$
under $ev^{*}\Omega^{1}_{Z}\otimes ev^{*}L\to\Omega^{1}_{\mathcal C/U^{\circ}}\otimes ev^{*}L\to\omega_{\mathcal C/U^{\circ}}\otimes ev^{*}L$.
\item $\sigma_{\theta}([C])=0$ if and only if $\theta|_{C}\equiv0$. On a smooth conic the value of $\sigma_{\theta}$ is the scalar of
Lemma~\ref{lem:allornothing}; on a nodal conic whose branches have $L$-degree one, $\sigma_{\theta}$ vanishes.
\end{enumerate}
\end{theorem}

\begin{proof}
(1) On a smooth conic $\omega_{C}\otimes L|_{C}=\OO$. On a nodal conic $C=C_{1}\cup C_{2}$ with node $p$,
$\omega_{C}|_{C_{i}}=\omega_{C_{i}}(p)=\OO(-1)$ and $L|_{C_{i}}=\OO(1)$, so $F=\omega_{C}\otimes L|_{C}$ is trivial on each
component, and from $0\to F\to F|_{C_{1}}\oplus F|_{C_{2}}\to F_{p}\to0$ its sections are pairs of constants agreeing at
$p$: $h^{0}=1$. By Serre duality on the Gorenstein curve $C$, $h^{1}(F)=h^{0}(L^{-1}|_{C})=0$. Grauert's theorem gives
local freeness with base change. (2) is the displayed composite. (3) By base change the value of $\sigma_{\theta}$ at
$[C]$ is the restriction of $\theta$ to $C$, an element of $H^{0}(\omega_{C}\otimes L|_{C})$; on the smooth locus of $C$ it
is the composite $TC\to L|_{C}$. On a nodal conic it vanishes on each branch, every line being integral, and a section of
a line bundle vanishing off a finite set vanishes.
\end{proof}

\subsection{Jumping conics}
\label{sec:jumpdef}

A conic of $\MCbar$ is the image of a map $f\colon\PP^{1}\to Z$ with $\deg f^{*}L=2$ (an embedding, a map birational onto a singular image, or a double cover
of a line), or a nodal or non-reduced limit of such images. On a smooth transverse conic $C=f(\PP^{1})$, by \eqref{eq:Bsplit},
\begin{equation}
f^{*}TZ=TC\oplus L^{\perp}|_{C}\cong\OO(2)\oplus\bigoplus_{i=1}^{n}\big(\OO(1+k_{i})\oplus\OO(1-k_{i})\big).\label{eq:fTZ}
\end{equation}
For the general conic (Proposition~\ref{prop:genconic}) all $k_{i}=0$ and $f^{*}TZ\cong\OO(2)\oplus\OO(1)^{2n}$.

\begin{definition}[The jumping invariant]
\label{def:jinv}
For a map $f\colon\PP^{1}\to Z$ with $\deg f^{*}L=2$ put
\[
j(f):=h^{0}\big(\PP^{1},f^{*}TZ\otimes f^{*}L^{-1}\big)-1=h^{0}\big(\PP^{1},N_{f}\otimes f^{*}L^{-1}\big),\qquad N_{f}:=f^{*}TZ/T\PP^{1},
\]
the two expressions agreeing since $H^{1}(T\PP^{1}\otimes f^{*}L^{-1})=H^{1}(\OO)=0$. For a member $C$ of $\MCbar$ which is a local complete intersection in $Z$
put $j(C):=h^{0}(C,N_{C/Z}\otimes L^{-1}|_{C})$, the normal sheaf $N_{C/Z}$ being then locally free; for a smooth conic $j(C)=j(f)$. Smooth conics and
reduced pairs of lines meeting in a node are local complete intersections, and these are the only members at which $j(C)$ is used. Not every member of $\MCbar$ is
one: the Hilbert limit of irreducible conics acquiring a node or a cusp carries an embedded point (Remark~\ref{rem:nodalharmless}(2)) and is not
Cohen--Macaulay. Being a local complete intersection is an open condition in flat families, and over the open set $\MCbar^{\,\mathrm{lci}}$ of such members the
normal sheaf of the universal curve is locally free and commutes with base change; so $j$ is upper semicontinuous on $\MCbar^{\,\mathrm{lci}}$, and $j(f)$ is
upper semicontinuous on $\mathrm{Hom}_{2}(\PP^{1},Z)$. On a smooth transverse conic $j(C)=h^{0}(B|_{C}(-1))=\sum_{k_{i}>0}k_{i}$.
\end{definition}

\begin{definition}[Jumping conics]
\label{def:jumping}
A \emph{jumping conic} is a smooth transverse conic with $j(C)\ge1$, equivalently $N_{C/Z}\not\cong\OO(1)^{\oplus2n}$, equivalently $B|_{C}$ is not trivial. The
set $J\subset\MCbar\smallsetminus\MC$ of jumping conics is the \emph{jumping locus}; $\MC\cup J$ is the open subset of smooth transverse conics of $\MCbar$, and $J$
is closed in it.
\end{definition}

A nodal conic has $j=a$, the number of summands $\OO_{C}(1,1)$ of its normal bundle, and $j\ge1$ exactly when it is unbalanced
(Proposition~\ref{prop:balanced}, Remark~\ref{rem:jinv}). Integral conics may have $j\ge1$ even when $J=\emptyset$, as on the Wolf spaces of $E_{6}$ and
$\Gr_{4}(\RR^{n+4})$ (Example~\ref{ex:intWolf}).

LeBrun's structure (\S\ref{sec:lebrunstructure}) is built from $H^{0}(C,B|_{C})$ and the evaluation maps $\mathrm{ev}_{z}\colon H^{0}(C,B|_{C})\to B_{z}$, which
are isomorphisms exactly when $B|_{C}$ is trivial; so it lives on $\MC$. In the proof the jumping conics through a general point $z$ never have to
be confronted: at a jumping conic the tangent-direction map $\widetilde\mu_{z}$ is not a local isomorphism (Lemma~\ref{lem:dmu}), so once $\widetilde\mu_{z}$ is
birational they form no divisor with divisorial image (Proposition~\ref{prop:cfromJ}). Finer properties of $J$ are not needed.

\subsection{LeBrun's structure}
\label{sec:lebrunstructure}

We recall LeBrun's theorem \cite{LeBrun89} in the cotangent convention of Salamon,
$T^{*}\MC\cong E_{\CC}\otimes H_{\CC}$; LeBrun's bundle $E$ is our $E_{\CC}^{*}$. At a conic $[C]\in\MC$ put
\[
E^{*}_{\CC}|_{[C]}:=H^{0}(C,B|_{C})\cong\CC^{2n},\qquad H^{*}_{\CC}|_{[C]}:=H^{0}(C,L^{1/2}|_{C})\cong\CC^{2},
\]
so that $T_{[C]}\MC=H^{0}(C,N_{C/Z})=E^{*}_{\CC}\otimes H^{*}_{\CC}$. The Levi form induces a symplectic form $\omega_{E}$ on
$E^{*}_{\CC}$: for $s,t\in H^{0}(B|_{C})$ the function $d\theta(s,t)$ is a section of $\OO_{C}$, hence a constant. In particular,
for every $z\in C$ the evaluation map
\[
\mathrm{ev}_{z}\colon E^{*}_{\CC}|_{[C]}=H^{0}(C,B|_{C})\longrightarrow B_{z}
\]
is a symplectic isomorphism, $\omega_{E}(s,t)=\omega_{z}(s(z),t(z))$. For $z\in C$ let $h_{z}\subset H^{*}_{\CC}|_{[C]}$ be the line
of sections of $L^{1/2}|_{C}$ vanishing at $z$.

\begin{theorem}[LeBrun \cite{LeBrun89}]
\label{thm:lebrun}
Let $(Z,\theta)$ be a complex contact manifold of dimension $2n+1$ (not necessarily compact) and let $\mathcal M$ be a
connected family of transverse conics with trivial $B$, open in the space of such curves. Then:
\begin{enumerate}
\item $\mathcal M$ carries a holomorphic Riemannian metric $g^{\CC}$ whose conformal class is the unique one invariant under
$\Sp(n,\CC)\Sp(1,\CC)$, with null cone containing the decomposable tensors $e\otimes h_{z}$ (the first-order deformations of $C$ fixing a point), and whose scale is fixed, up to a
constant factor, by the contact form; it is Einstein. Its Levi-Civita connection preserves $E_{\CC}$, $H_{\CC}$ and symplectic
forms on both, so $\mathrm{Hol}(g^{\CC})\subset\Sp(n,\CC)\Sp(1,\CC)$. If $Z$ is the twistor space of $M$ and $\mathcal M=\MC$,
the twistor lines form a totally real submanifold $M\subset\MC$ on which $g^{\CC}$ restricts, for a suitable choice of
the constant, to the quaternionic-K\"ahler metric.
\item For $z\in Z$, the conics of $\mathcal M$ through $z$ form a submanifold $\Sigma_{z}$ of dimension $2n$, the \emph{leaf} of $z$,
with $T_{[C]}\Sigma_{z}=H^{0}(N_{C/Z}(-z))=E^{*}_{\CC}\otimes h_{z}$; it is totally null and totally geodesic. The connection of
$E^{*}_{\CC}$ restricted to $\Sigma_{z}$ is flat, and its parallel sections are the sections with constant value under
$\mathrm{ev}_{z}$: parallel transport along a path in $\Sigma_{z}$ from $[C]$ to $[C']$ is $\mathrm{ev}_{z}^{-1}\circ\mathrm{ev}_{z}$.
\item The curvature of $g^{\CC}$ is $\rho_{0}R_{0}+R_{E}^{\CC}$ with $\rho_{0}$ constant, $R_{0}$ the curvature tensor of quaternionic
projective space and $R_{E}^{\CC}\in H^{0}(\mathcal M,S^{4}E_{\CC})$ the quaternionic Weyl curvature; in the twistor case $R_{E}^{\CC}$
restricts on $M$ to the quaternionic Weyl curvature $R_{E}$ of $g$, the trace-free quartic component of the curvature of the $\Sp(n)$-factor.
\end{enumerate}
\end{theorem}

LeBrun's construction is local along the curves, and we use it for germs of contact manifolds as well (\S\ref{sec:dilation}). Part (2) is the quaternionic form of the basic property of Penrose's $\alpha$-surfaces: the connection on $E^{*}_{\CC}$ is the one for which a section is parallel along $\Sigma_{z}$ exactly when its
value at the common point $z$ is constant. The flatness itself is visible from (3): for $X=e\otimes h_{z}$, $Y=e'\otimes h_{z}$ every
term of the curvature of $E_{\CC}$ evaluated on $(X,Y)$ carries the factor $\omega_{H}(h_{z},h_{z})=0$. Note that $\mathrm{ev}_{z}$
takes values in $B_{z}=L^{\perp}_{z}\otimes L^{-1/2}_{z}$, which requires a choice of square root of the line $L_{z}$. More generally, $E_{\CC}$ and $H_{\CC}$ are
\emph{local} symplectic factors: no global lift of the structure group $(\Sp(n,\CC)\times\Sp(1,\CC))/\langle(-I,-I)\rangle$ to $\Sp(n,\CC)\times\Sp(1,\CC)$ is
assumed, and when $\mathrm{Pic}\,Z=\ZZ L$ the line bundle $L$ has no square root on $Z$. The associated bundles of representations on which $(-I,-I)$ acts trivially
are globally defined; in particular $T^{*}\MC$, $S^{4}E_{\CC}$ and $S^{5}E_{\CC}\otimes H_{\CC}$ are, and every tensor below has even weight in $L^{1/2}$, so no
choice is needed for it.
For the scale in (1) see \cite[\S1, Prop.~1.2 and the proof of Prop.~3.3]{LeBrun89}. Part (2) is read off from LeBrun's definitions
\cite[pp.~356--357]{LeBrun89}, in his notation $p\colon F\to Y$ (the space of conics) and $q\colon F\to Z$. On p.~356, $\nabla_{E}$ is defined as the push-forward
by $p$ of the differentiation $d_{\mathcal F}$ of sections of $q^{*}(D\otimes L^{-1/2})$ up the fibres of $q$, where $\Omega^{1}_{\mathcal F}=(\ker q_{*})^{*}$, through the exact
sequence $0\to q^{*}\OO(D^{*})\to p^{*}\Omega^{1}\to\Omega^{1}_{\mathcal F}\to0$. On p.~357 LeBrun proves that $\nabla_{E}\otimes1+1\otimes\nabla_{H}$ is torsion-free,
hence is the Levi-Civita connection of $g=\epsilon_{E}\otimes\epsilon_{H}$, which is $g^{\CC}$ up to a constant factor (\S\ref{sec:levwronsk}). A section $e$ of
$E^{*}_{\CC}$ (LeBrun's $E$) is the section $([C'],z')\mapsto e([C'])(z')$ of $q^{*}(D\otimes L^{-1/2})$. On the fibre $q^{-1}(z)$, which $p$ maps onto
$\Sigma_{z}$, this is $[C']\mapsto\mathrm{ev}_{z}(e([C']))$, with values in the fixed space $B_{z}$. Hence $\mathrm{ev}_{z}(\nabla_{v}e)=v(\mathrm{ev}_{z}\circ e)$ for
$v$ tangent to $\Sigma_{z}$: the Levi-Civita connection itself, and not merely a connection gauge-equivalent to it, has parallel transport
$\mathrm{ev}_{z}^{-1}\circ\mathrm{ev}_{z}$ along $\Sigma_{z}$. This also fixes the splitting of the Levi-Civita connection of $TY=E\otimes H$ into factors, which a
priori is determined only up to $(\nabla_{E}+\beta,\nabla_{H}-\beta)$: since $\omega_{E}(e,e')=d\theta_{z}(\mathrm{ev}_{z}e,\mathrm{ev}_{z}e')$ with $d\theta_{z}$ fixed
on $B_{z}$, the identity gives $v\,\omega_{E}(e,e')=\omega_{E}(\nabla_{v}e,e')+\omega_{E}(e,\nabla_{v}e')$ for $v$ tangent to $\Sigma_{z}$, and these $v$, $z\in C$, span
$T_{[C]}Y$; so $\nabla_{E}\omega_{E}=0$, the normalization of Theorem~\ref{thm:lebrun}(1). The parallel transport of $\nabla_{H}$ lifts paths tangent to
$\tilde D=q_{*}^{-1}(D)\supset\ker q_{*}$ \cite[p.~357]{LeBrun89}; it fixes $z$ and so preserves $h_{z}$. Hence $E^{*}_{\CC}\otimes h_{z}$ is parallel along $\Sigma_{z}$, which is therefore totally
geodesic, the connection being torsion-free.

\subsection{The metric as a Levi--Wronskian quotient}
\label{sec:levwronsk}

\begin{proposition}[The Levi--Wronskian formula]
\label{prop:levwronsk}
Let $C=h(\PP^{1})$ be a conic of $\MC$ and $\Lambda_{[C]}=H^{0}(\omega_{C}\otimes L|_{C})$.
\begin{enumerate}
\item Let $\nabla$ be a local connection on $L^{\perp}|_{C}=E^{*}_{\CC}\otimes L^{1/2}|_{C}$ of the form
$\nabla(e\otimes\sigma)=e\otimes\nabla\sigma$. For $s,t\in H^{0}(L^{\perp}|_{C})$ put
\[
G(s,t):=d\theta(s,\nabla t)-d\theta(\nabla s,t).
\]
Then $G(s,t)\in\Lambda_{[C]}$, it is independent of $\nabla$, symmetric, and $G=\omega_{E}\otimes\omega_{H}$, where
$\omega_{H}(\sigma,\tau)=\sigma\nabla\tau-\tau\nabla\sigma$ is the $\Lambda$-valued Wronskian on $H^{*}_{\CC}$.
\item LeBrun's metric is, up to a universal constant,
\[
g^{\CC}=\frac{G}{\sigma_{\theta}},
\]
the quotient of two elements of the line $\Lambda_{[C]}$.
\item ($\alpha$-plane form.) For $z\in C$, $s\in H^{0}(L^{\perp}|_{C}(-z))$ and $t\in H^{0}(L^{\perp}|_{C})$,
\[
g^{\CC}(s,t)=-\frac{d\theta\big(\nabla_{z}s,\,t(z)\big)}{\theta(dh_{z})},
\]
where $\nabla_{z}s\in L^{\perp}_{z}\otimes T^{*}_{z}C$ is the intrinsic derivative of $s$ at its zero.
\end{enumerate}
\end{proposition}

\begin{proof}
(1) For $s=e\otimes\sigma$, $t=e'\otimes\tau$ one has $d\theta(s,t)=\omega_{E}(e,e')\,\sigma\tau$, so
$G(s,t)=\omega_{E}(e,e')(\sigma\nabla\tau-\tau\nabla\sigma)$. The Wronskian is a section of $L|_{C}\otimes\Omega^{1}_{C}$ and does not
depend on $\nabla$: replacing $\nabla$ by $\nabla+a$ changes the two terms by $a\sigma\tau$ each. Symmetry follows from the
antisymmetry of $d\theta$ and of the Wronskian. (2) By Theorem~\ref{thm:lebrun}(1) the conformal class of $g^{\CC}$ is that of
$\omega_{E}\otimes\omega_{H}$, the unique (up to scale) $\Sp(n,\CC)\Sp(1,\CC)$-invariant symmetric form on $E^{*}_{\CC}\otimes H^{*}_{\CC}$ (its null cone
contains the decomposable tensors, and equals the set of them only for $n=1$); $\omega_{H}$ takes values in the line $\Lambda_{[C]}$, and on $\MC$ the section $\sigma_{\theta}$
trivializes $\Lambda$ (Theorem~\ref{thm:thetasection}(3)). In LeBrun's construction the Einstein representative of the
conformal class is the one determined by the contact form \cite{LeBrun89}, which is the trivialization by $\sigma_{\theta}$.
(3) Evaluate the constant $G(s,t)$ at $z$: since $s(z)=0$ the term $d\theta(s,\nabla t)$ vanishes there, $\nabla s(z)$ is
intrinsic, and $\sigma_{\theta}$ at $z$ is the scalar $\theta(dh_{z})$.
\end{proof}

In \cite{LeBrun89}, $g:=\epsilon_{E}\otimes\epsilon_{H}$, with $\epsilon_{E}$ the pull-back of $d\theta$ (\S1) and $\epsilon_{H}$ the Wronskian
$\Lambda^{2}H\to\tilde L$ divided by the canonical section $\tilde\theta$ of $\tilde L$ obtained by evaluating $\theta$ on the fibres (proof of Prop.~3.3), which is
$\sigma_{\theta}$; by Prop.~1.2 and Thm.~1.3 this $g$, up to a universal constant, is Einstein with holonomy in $\Sp(n,\CC)\Sp(1,\CC)$ and restricts to the
quaternionic-K\"ahler metric on the real slice. So the Einstein representative is $G/\sigma_{\theta}$ up to a universal constant.

\begin{corollary}[Rationality]
\label{cor:rational}
Let $Z$ be projective with $L$ ample. Then $g^{\CC}$ is the restriction to $\MC$ of a rational symmetric form on $T\MCbar$, and its Levi-Civita connection and
curvature are rational tensors on $\MCbar$.
\end{corollary}

\begin{proof}
On the dense open set $\MC$ the form $G$ is given by an algebraic expression in $\theta$, $d\theta$ and the first-order
deformations of the universal family, and $\sigma_{\theta}$ is an algebraic section of $\Lambda$; so $G/\sigma_{\theta}$ is a rational symmetric
form on $T\MCbar$. Christoffel symbols and curvature are polynomial in $g^{\CC}$, $(g^{\CC})^{-1}$ and their derivatives.
\end{proof}

\subsection{The curvature quartic at a point of a conic}
\label{sec:Wz}

The bundle $S^{4}E_{\CC}$ has no natural extension across the boundary, because $E_{\CC}$ involves $L^{1/2}$, which does not
exist on the universal curve over the nodal divisor. The following lattice avoids square roots. Let
\[
\overline F=\{([C],z):z\in C\}\subset\MCbar\times Z,\qquad\bar\pi\colon\overline F\to\MCbar,\quad\bar p\colon\overline F\to Z,
\]
be the universal conic, and
\[
\mathcal S:=\bar\pi_{*}\,\bar p^{\,*}\big(S^{4}L^{\perp*}\otimes L^{2}\big).
\]
On $\MC$, where $B|_{C}\cong\OO^{2n}$, an element $w\in H^{0}(C,S^{4}B^{*}|_{C})=H^{0}(C,S^{4}L^{\perp*}\otimes L^{2}|_{C})$
pairs with sections of $B|_{C}$ to a constant function; this gives an isomorphism
\[
\varpi\colon\mathcal S|_{\MC}\xrightarrow{\ \cong\ }S^{4}E_{\CC},\qquad\varpi(w)(\sigma^{4})=w(\sigma,\sigma,\sigma,\sigma),\ \sigma\in H^{0}(B|_{C}).
\]
Put $W^{\CC}:=\varpi^{-1}(R_{E}^{\CC})$. Evaluating at a point of the conic we obtain, for $([C],z)\in\overline F$ with
$[C]\in\MC$,
\[
W^{\CC}([C],z)\in V_{z}:=\big(S^{4}L^{\perp*}\otimes L^{2}\big)_{z},\qquad W^{\CC}([C],z)=(\mathrm{ev}_{z}^{-1})^{*}R_{E}^{\CC}([C]),
\]
a vector in a space that depends only on $z$. The curvature $R_{E}^{\CC}$ is the $S^{4}E_{\CC}$-component of the rational curvature
tensor of Corollary~\ref{cor:rational}, and $\mathrm{ev}_{z}$ is algebraic; so $W^{\CC}$ is a rational section of the square-root-free
bundle $\bar p^{\,*}(S^{4}L^{\perp*}\otimes L^{2})$ over $\overline F$ (no extension of $S^{4}E_{\CC}$ itself is needed), and for each $z$
we obtain a meromorphic map
\[
W_{z}\colon\overline F_{z}=\bar p^{-1}(z)\dashrightarrow V_{z},
\]
holomorphic on the transverse non-jumping conics through $z$. We call $W_{z}$ the \emph{curvature quartic at $z$}. (At the smooth integral conics with
normal bundle $\OO(1)^{2n}$, where $B|_{C}$ is no longer trivial, $W^{\CC}$ is also holomorphic: there the frame $N_{C/Z}\otimes L^{-1/2}|_{C}$ replaces
$B|_{C}$, see Proposition~\ref{prop:intregular}.)

\subsection{The leaf principle}
\label{sec:leaf}

\begin{proposition}[The leaf principle]
\label{prop:leaf}
The following are equivalent:
\begin{enumerate}
\item $\nabla R_{E}^{\CC}=0$ on $\MC$;
\item for every $z\in Z$, $W_{z}$ is locally constant on $\Sigma_{z}$;
\item there is a dense Zariski-open $V\subset Z$ such that for every $z\in V$, $W_{z}$ is locally constant on $\Sigma_{z}$.
\end{enumerate}
\end{proposition}

\begin{proof}
Fix $z$ and $[C]\in\Sigma_{z}$. By Theorem~\ref{thm:lebrun}(2), $\mathrm{ev}_{z}$ is parallel along $\Sigma_{z}$, so for a tangent
vector $X\in T_{[C]}\Sigma_{z}=E^{*}_{\CC}\otimes h_{z}$,
\[
X\big(W^{\CC}(\cdot,z)\big)=(\mathrm{ev}_{z}^{-1})^{*}\big(\nabla_{X}R_{E}^{\CC}\big).
\]
Hence (1)$\Rightarrow$(2)$\Rightarrow$(3). Assume (3). The conics meeting $Z\smallsetminus V$ in at most finitely many points form
a dense open subset of $\MC$ (indeed every conic not contained in $Z\smallsetminus V$). For such $[C]$ and every $z\in C\cap V$, $\nabla_{X}R_{E}^{\CC}([C])=0$ for $X\in E^{*}_{\CC}\otimes h_{z}$. As
$z$ runs over the infinite set $C\cap V$, the lines $h_{z}\subset H^{*}_{\CC}|_{[C]}\cong\CC^{2}$ are pairwise distinct ($h_{z}$
determines $z$ as the zero of its sections), so they span $H^{*}_{\CC}$ and the subspaces $E^{*}_{\CC}\otimes h_{z}$ span
$T_{[C]}\MC$. So $\nabla R_{E}^{\CC}$ vanishes on a dense open subset of $\MC$, which is smooth and irreducible, hence
connected, and $\nabla R_{E}^{\CC}$ vanishes identically.
\end{proof}

\begin{lemma}[The derivative of the quaternionic Weyl curvature]
\label{lem:S5}
On $\MC$, $\nabla R_{E}^{\CC}$ is a section of $S^{5}E_{\CC}\otimes H_{\CC}\subset T^{*}\MC\otimes S^{4}E_{\CC}$: for $e,e_{1},\dots,e_{4}\in E^{*}_{\CC}$
and $k\in H^{*}_{\CC}$ at a point of $\MC$, the number $(\nabla_{e\otimes k}R_{E}^{\CC})(e_{1},e_{2},e_{3},e_{4})$ is symmetric in $(e,e_{1},\dots,e_{4})$.
Consequently, in the formula $X(W^{\CC}(\cdot,z))=(\mathrm{ev}_{z}^{-1})^{*}(\nabla_{X}R_{E}^{\CC})$ of the proof of Proposition~\ref{prop:leaf}, for
$X=e\otimes k\in T_{[C]}\Sigma_{z}$, the right side is the contraction with $\mathrm{ev}_{z}(e)$ of a totally symmetric $5$-form on $B_{z}$.
\end{lemma}

\begin{proof}
By Theorem~\ref{thm:lebrun}(3) the curvature is $\rho_{0}R_{0}+R_{E}^{\CC}$ with $\rho_{0}$ constant, and $R_{0}$ is built from $g^{\CC}$, $\omega_{E}$ and $\omega_{H}$,
which are parallel; so the second Bianchi identity holds for $R_{E}^{\CC}$ alone. The $2$-forms of $R_{E}^{\CC}$ lie in $S^{2}E\otimes\Lambda^{2}H$ \cite{Salamon82}, so on
decomposable vectors $X_{i}=e_{i}\otimes k_{i}$,
\[
R_{E}^{\CC}(X_{1},X_{2},X_{3},X_{4})=R_{E}^{\CC}(e_{1},e_{2},e_{3},e_{4})\,\omega_{H}(k_{1},k_{2})\,\omega_{H}(k_{3},k_{4}).
\]
Put $T(e,k;e_{1},\dots,e_{4})=(\nabla_{e\otimes k}R_{E}^{\CC})(e_{1},\dots,e_{4})$, symmetric in $e_{1},\dots,e_{4}$. In the second Bianchi identity
$\sum_{\mathrm{cyc}(1,2,3)}(\nabla_{X_{1}}R_{E}^{\CC})(X_{2},X_{3},X_{4},X_{5})=0$ take $k_{2}=k_{3}=k$, and $k_{1},k_{4},k_{5}$ with $\omega_{H}(k,k_{1})\ne0\ne\omega_{H}(k_{4},k_{5})$.
The term with $\nabla_{X_{1}}$ carries $\omega_{H}(k_{2},k_{3})=0$, and the other two give
\[
T(e_{2},k;e_{3},e_{1},e_{4},e_{5})\,\omega_{H}(k,k_{1})+T(e_{3},k;e_{1},e_{2},e_{4},e_{5})\,\omega_{H}(k_{1},k)=0 ,
\]
that is, $T(e_{2},k;e_{1},e_{3},e_{4},e_{5})=T(e_{3},k;e_{1},e_{2},e_{4},e_{5})$. So exchanging the first argument of $T$ with any other leaves it unchanged, and with the
symmetry in the last four arguments $T$ is totally symmetric. The last assertion follows, since $\mathrm{ev}_{z}$ is an isomorphism.
\end{proof}

For $n=1$, Lemma~\ref{lem:S5} is the Bianchi identity $\nabla^{A}{}_{A'}\Psi_{ABCD}=0$ of a self-dual Einstein metric. A numerical check, carried out by machine and verified by the author, confirms
that the kernel of the second Bianchi map on $(E\otimes H)\otimes S^{4}E$ is exactly $S^{5}E\otimes H$ for $n=1,2$; the proof of Lemma~\ref{lem:S5} does not
depend on it.

\begin{theorem}[The extension principle]
\label{thm:extension}
Let $z\in Z$ be general and let $\overline F_{z}=\bar p^{\,-1}(z)\subset\MCbar$ be the compact variety of conics through $z$; its
irreducible components are $2n$-dimensional and contain transverse non-jumping conics (Theorem~\ref{thm:irred}). The curvature
quartic at $z$ is a meromorphic map
\[
W_{z}\colon\overline F_{z}\dashrightarrow V_{z}=\big(S^{4}L^{\perp*}\otimes L^{2}\big)_{z}\cong\CC^{N},\qquad
N=\dim S^{4}(L^{\perp}_{z})^{*}=\binom{2n+3}{4},
\]
holomorphic on the dense open set $F^{\mathrm{ne}}_{z}$ of transverse non-jumping conics, with values in a fixed complex vector
space depending only on $z$.
\begin{enumerate}
\item If $W_{z}$ extends to a holomorphic map on $\overline F_{z}$ for every $z$ in a dense Zariski-open
subset of $Z$, then $R_{E}^{\CC}$ is parallel; if $Z$ is a twistor space, $M$ is a Wolf space (Corollary~\ref{cor:endgame}).
\item Let $\mu_{z}\colon\overline F_{z}\dashrightarrow\PP(T_{z}Z)$, $[C]\mapsto[T_{z}C]$, be the tangent-direction map, defined on the
conics smooth at $z$, and let $\Gamma_{z}$ be the normalization of the closure of its graph, on which $\mu_{z}$ becomes a proper,
generically finite morphism $\widetilde\mu_{z}$ (\S\ref{sec:fibre}). Then it suffices that for such $z$ the map $W_{z}$ be regular at
the generic point of every prime divisor $D\subset\Gamma_{z}$ whose image $\widetilde\mu_{z}(D)$ is a divisor of $\PP(T_{z}Z)$. No
condition is needed along the divisors contracted by $\widetilde\mu_{z}$, nor along subsets of codimension at least two. Among the
divisors with divisorial image is the exceptional divisor over the conics with a node at $z$, where $\mu_{z}$ is undefined and whose limit
directions fill the contact hyperplane (Lemma~\ref{lem:muz}(3)); regularity there is Theorem~\ref{thm:Eregular}.
\end{enumerate}
\end{theorem}

\begin{proof}
(1) Each component is compact and irreducible, so each of the $N$ coordinates $\ell\circ W_{z}$ of a holomorphic map into $\CC^{N}$ is a
holomorphic function on a compact irreducible variety, hence constant on it. Thus $W_{z}$ is locally constant on $\Sigma_{z}$, and Proposition~\ref{prop:leaf}
applies.
(2) Each $\ell\circ W_{z}$ is a meromorphic function on the normal variety $\Gamma_{z}$, so its polar set is a divisor, and by hypothesis it is
bounded near the generic point of every prime divisor with divisorial image. Lemma~\ref{lem:hartogsmu}, applied to
$\widetilde\mu_{z}$, shows that such a function is constant: its elementary symmetric functions over the fibres of $\widetilde\mu_{z}$ are rational
functions on $\PP(T_{z}Z)$ with no polar divisor, hence constants. This is Theorem~\ref{thm:Ereduction}, proved in \S\ref{sec:hartogs}.
\end{proof}

\begin{remark}[Leaf constancy, and the role of the boundary]
\label{rem:descent}
Proposition~\ref{prop:leaf} says that $R_{E}^{\CC}$ is parallel if and only if the section $W^{\CC}$ of
$\bar p^{\,*}(S^{4}L^{\perp*}\otimes L^{2})$ is locally constant along the leaves; when the conics through a general point form an
irreducible family (Theorem~\ref{thm:irred}) this means that $W^{\CC}$ is pulled back from a section $s$ of
$S^{4}L^{\perp*}\otimes L^{2}$ on $Z$, whose value at $z$ is, when $Z$ is a twistor space, the curvature of $M$ at $q(z)$ read on the twistor line through $z$. Equivalently $R_{E}^{\CC}$
comes from a section on $Z$; we call this the \emph{leaf constancy} of the curvature. In \S\ref{sec:descent} it will be obtained in two steps:
first the symmetric functions of $W_{z}$ over the fibres of the tangent-direction map are functions on $\PP(T_{z}Z)$, where a Hartogs
argument makes them constant, and then constancy on each component of $\Gamma_{z}$ follows. In a non-compact
quaternionic-K\"ahler germ the curvature quartic at $z$ need not be constant, so constancy must come from the compactness of $Z$,
which acts on $W_{z}$ only through the boundary of the compactified space of conics through $z$. The parallel with
\cite{Ye94} is that there, too, the curvature of a holomorphic conformal structure was forced to vanish by a global
rational curve; here the rational curves are the conics, and the conformal structure is allowed to degenerate at the
boundary of their moduli.
\end{remark}

\subsection{The LeBrun transform and its compactification}
\label{sec:transform}

The correspondence space of the transform is the universal conic $\overline F\to\MCbar$ with its evaluation $\bar p\colon\overline F\to Z$ (\S\ref{sec:Wz}), or, in the language
of stable maps, the forgetful map $\overline{M}_{0,1}(Z,2\ell)\to\overline{M}_{0,0}(Z,2\ell)$ with the evaluation $\overline{M}_{0,1}(Z,2\ell)\to Z$, restricted to the component
of $\MCbar$. The transform carries a holomorphic bundle $V$ on $Z$ to $\bar\pi_{*}\bar p^{*}V$ on $\MCbar$; locally on $Z$, the bundle $E^{*}_{\CC}=\pi_{*}p^{*}B$ with its connection is the case
$V=B=L^{\perp}\otimes L^{-1/2}$, a twisted Ward transform (Theorem~\ref{thm:lebrun}(2)), while the global, square-root-free statement is
$\bar\pi_{*}\bar p^{*}\big(S^{4}L^{\perp*}\otimes L^{2}\big)\cong S^{4}E_{\CC}$ on $\MC$ (\S\ref{sec:Wz}). In the opposite direction a tensor on $\MC$ comes from $Z$ exactly when, read through the
incidence, it is constant along the $\alpha$-planes; for the curvature this is the leaf principle (Proposition~\ref{prop:leaf}): \emph{the quaternionic Weyl curvature of $\MC$ is parallel if and only if it is, locally on $Z$, the transform of a section of
$S^{4}L^{\perp*}\otimes L^{2}$}; on the real slice, \emph{$M$ is symmetric if and only if this holds}. In this language Theorem~\ref{thm:LSJ} says that compactness of $Z$ forces this
leaf constancy once the transform is understood at the boundary and two general points lie on only one conic.

The two compactifications are complementary. In the Hilbert scheme the boundary carries embedded points at nodes and cusps, and the flat limit
gives sharp control of limit directions (Lemmas~\ref{lem:P2proved} and~\ref{lem:ribbons}); in the space of stable maps every boundary point has the deformation--obstruction theory of stable maps, which at nodal sources also includes the
deformations of the domain and the smoothings of the node, and the strata are: reducible domains (line pairs), double covers of lines, maps tangent to the contact distribution ($\sigma_{\theta}=0$), maps that
are not embeddings, and maps with unbalanced $f^{*}TZ$ (jumping). The singular LeBrun transform is the study of the transform along these strata.

\begin{definition}[Singular conformal structures and metrics]
\label{def:singconf}
Let $\overline M$ be an irreducible complex space of dimension $m$ with smooth locus $\overline M_{\mathrm{reg}}$.
\begin{enumerate}
\item A \emph{conformal structure} on $\overline M$ is a holomorphic section $g\in H^{0}(\overline M_{\mathrm{reg}},S^{2}T\overline M^{*}\otimes N)$, for some line bundle $N$ on
$\overline M_{\mathrm{reg}}$, that is nondegenerate over a dense open subset: in local
coordinates $x^{1},\dots,x^{m}$ and a local frame of $N$,
$g=\sum g_{ij}\,dx^{i}dx^{j}$ with $g_{ij}$ holomorphic and $\det(g_{ij})\not\equiv0$; equivalently, $\det g\in H^{0}(\overline M_{\mathrm{reg}},K_{\overline M}^{2}\otimes N^{m})$ is not
identically zero. Sections $g$ of $S^{2}T^{*}\otimes N$ and $g'$ of $S^{2}T^{*}\otimes N'$ define the same structure if $g'=(\mathrm{id}\otimes\phi)(g)$ for a nonzero
meromorphic homomorphism $\phi\colon N\dashrightarrow N'$, that is, a meromorphic section of $N'\otimes N^{-1}$; locally $g'=fg$ with $f$ meromorphic, but the local
multipliers need not come from a global function. Equivalently, the structure is the line spanned by $g$ in $S^{2}T^{*}$ at the points where $g\ne0$.
\item If $g$ vanishes to order $k_{i}$ along the prime divisors $D_{i}$, it is a section of $S^{2}T^{*}\otimes N(-\sum k_{i}D_{i})$ that vanishes on no hypersurface.
This \emph{reduced} representative is unique up to isomorphism, and we use it from now on.
\item The \emph{regular locus} $M^{\circ}$ is the open set of points $p\in\overline M_{\mathrm{reg}}$ at which $g(p)$ is nondegenerate; its complement
$\{\det g=0\}$ in $\overline M_{\mathrm{reg}}$ is empty or a hypersurface. If $M^{\circ}\ne\overline M$, the structure is \emph{singular}, and a point
$p\in\overline M\smallsetminus M^{\circ}$ is a \emph{degeneracy} if $p\in\overline M_{\mathrm{reg}}$ and a \emph{singular point} if $p\notin\overline M_{\mathrm{reg}}$.
\item A \emph{metric} in the structure is $h=g/\sigma$, with $\sigma\not\equiv0$ a meromorphic section of $N$; every meromorphic section $h$ of $S^{2}T^{*}$ over
$\overline M_{\mathrm{reg}}$ with $\det h\not\equiv0$ is of this form, after clearing denominators and reducing. A point $p\in\overline M_{\mathrm{reg}}$ is a
\emph{pole} of $h$ if some coefficient $h_{ij}$ is not holomorphic at $p$, and a \emph{degeneracy} of $h$ if the $h_{ij}$ are holomorphic at $p$ and
$\det(h_{ij}(p))=0$. We call $h$ a \emph{singular conformal metric} on $\overline M$ if it does not extend to a holomorphic nondegenerate metric on all of
$\overline M$, that is, if it has a pole or a degeneracy or $\overline M$ is singular.
\end{enumerate}
The degeneracies and singular points of the structure are intrinsic; the poles of a metric are not. Multiplying $h$ by $f^{k}$, with $f$ a local equation of a
hypersurface, trades a pole of order $k$ along it for a holomorphic form, possibly degenerate there. The structure is regular at a smooth point $p$ exactly when
some local representative is holomorphic and nondegenerate at $p$: the coefficients of the reduced $g$ have no common factor in the factorial ring $\mathcal O_{p}$, so
any representative holomorphic at $p$ is a holomorphic multiple of $g$. In particular the degeneracies of the structure lie among the poles and degeneracies of any
metric in it, but need not fill them: for $n=1$ LeBrun's metric has a pole along the balanced integral divisor, on which the reduced representative is
nondegenerate (Remark~\ref{rem:intflat}(1)).
\end{definition}

\begin{remark}[The singularities of LeBrun's metric]
\label{rem:singLeBrun}
LeBrun's metric $g^{\CC}=G/\sigma_{\theta}$ (Proposition~\ref{prop:levwronsk}) is a singular conformal metric on $\MCbar$, and its singularities lie over $\MCbar\smallsetminus\MC$: the jumping conics and the
boundary.
\begin{itemize}[leftmargin=2em]
\item \emph{Poles.} Along a balanced integral divisor, where $\sigma_{\theta}$ vanishes to order one, $g^{\CC}$ has a pole of order two on a subspace of rank four and
of order one on the rest (Proposition~\ref{prop:intregular}(2)).
\item \emph{Degeneracies.} With $s$ a local equation of that divisor, $s^{2}g^{\CC}$ is holomorphic, and for $n\ge2$ degenerate of rank four along it: the divisor is a
complex conformal infinity (Remark~\ref{rem:intflat}). Since $s^{2}g^{\CC}$ vanishes on no hypersurface, it is the reduced representative near the divisor,
so the conformal structure is singular there for $n\ge2$, while for $n=1$ it is regular.
\item \emph{Degeneration of the factorization.} Along $J$ the factorization $T\MC=E^{*}_{\CC}\otimes H^{*}_{\CC}$ breaks down: $\mathrm{ev}_{z}$ and
$\omega_{E}$ degenerate (\S\ref{sec:jumpdef}). Along the nodal divisor $\sigma_{\theta}$ vanishes to order one at the line pairs whose branch tangents $u_{1},u_{2}$
satisfy $d\theta(u_{1},u_{2})\neq0$ (Proposition~\ref{prop:thetamult}), and the balanced smooth-conic factorization of $T\MC$ is no longer available canonically at the nodal members
(Table~\ref{tab:dictionary}).
\item \emph{Singular points.} $\MCbar$ is smooth at every member $C$ with $H^{1}(N_{C/Z})=0$, in particular at the balanced conics and at the nodal conics with nef
branches (Lemma~\ref{lem:nodalsmooth}); it may be singular at the other degenerate members. Where this matters the paper works on the normalization, on the space of
stable maps, or on the blown-up families $\Gamma_{z}$ of conics through a point.
\end{itemize}
The singular LeBrun transform reads these singularities, and those of $R_{E}^{\CC}$, off the degenerate conics of $Z$ (Table~\ref{tab:dictionary}).
\end{remark}

\section{Why the conformal structure must degenerate}
\label{sec:mustdegenerate}

The aim of this paper is to relate the differential geometry of the complex conformal manifold $(\MC,g^{\CC})$ to the algebraic geometry of the
contact manifold $Z$, and the two meet at the boundary. The metric $g^{\CC}=G/\sigma_{\theta}$ is defined on the balanced transverse conics. On $\MCbar\smallsetminus\MC$, at the jumping
conics and at the boundary, formed by the integral, nodal and other degenerate conics of $Z$, it acquires poles and degenerates (Definition~\ref{def:singconf}, Remark~\ref{rem:singLeBrun}), and the way it does so is governed by the algebraic geometry of those
curves: their normal bundles, their contact neighbourhoods, and the lines and the variety of minimal rational tangents of $Z$. This section sets the stage by showing that
the degeneration is forced. Three facts about $(\MC,g^{\CC})$ show it (\S\ref{sec:introidea}). For $n\ge2$ the Weyl tensor of $g^{\CC}$ vanishes nowhere, so by
Belgun's theorem no immersed rational null geodesic lies where $g^{\CC}$ is defined (Proposition~\ref{prop:notflat}). Every complete curve of conics through a point of $\MC$
meets the integral locus, where the denominator $\sigma_{\theta}$ of $g^{\CC}$ vanishes (Proposition~\ref{prop:nocurves}). And the null geodesics of the $\alpha$-planes,
closed up, run into the boundary (Proposition~\ref{prop:projflat}, Remark~\ref{rem:liouville}). So the global geometry of $(\MC,g^{\CC})$ has to be sought at its singular
locus. The rest of the paper studies the singularities of $(g^{\CC},R_{E}^{\CC})$ there: conic by conic in Part~\ref{part:local}, and globally, through the
tangent-direction map, in Part~\ref{part:global}.

\begin{proposition}[Nowhere conformally flat]
\label{prop:notflat}
Let $n\ge2$. The Weyl tensor of $g^{\CC}$ vanishes at no point of $\MC$. Consequently no immersed rational curve in $\MC$ is a null
geodesic of $g^{\CC}$.
\end{proposition}

\begin{proof}
By Theorem~\ref{thm:lebrun}(3) the curvature of $g^{\CC}$ is $\rho_{0}R_{0}+R_{E}^{\CC}$ with $\rho_{0}$ constant; $\rho_{0}\ne0$ by Theorem~\ref{thm:thirdorder}(4) (on a real slice it is a positive multiple of
the scalar curvature of $g$). The tensor $R_{E}^{\CC}\in S^{4}E_{\CC}$ has vanishing Ricci contraction \cite{Salamon82}, so it equals its own Weyl
part. The tensor $R_{0}$ is the curvature of $\HH P^{n}$, an invariant of the structure group $\Sp(n,\CC)\Sp(1,\CC)$; its Weyl part $W(R_{0})$ is again invariant,
and it is nonzero for $n\ge2$, because an Einstein metric with vanishing Weyl tensor has constant sectional curvature, while the sectional curvature of $\HH P^{n}$,
$n\ge2$, takes all values in an interval $[\kappa,4\kappa]$. So at every point the Weyl tensor of $g^{\CC}$ is $\rho_{0}W(R_{0})+R_{E}^{\CC}$, the sum of a nonzero
invariant and an element of $S^{4}E_{\CC}$, an irreducible representation of $\Sp(n,\CC)$ without invariants; the two lie in different isotypic components of the
space of Weyl tensors, so the sum is nonzero. The last assertion follows from Belgun's theorem \cite[Theorem~4]{Bel01}: a complex conformal manifold containing an
immersed rational null geodesic is conformally flat.
\end{proof}

\begin{proposition}[No compact curves]
\label{prop:nocurves}
Let $Z$ be a compact complex contact manifold whose contact line bundle $L$ is ample, and let $B$ be an irreducible complete curve of conics in $Z$ (in the Hilbert
scheme) whose general member is a smooth transverse conic. Then some member of $B$ is integral: its support is a smooth integral conic, a line pair, a line, or
the image of an integral map. In particular $\MC$ contains no compact curve, and every complete curve in $\MCbar$ through a point of $\MC$ meets the integral locus, the
members with integral support; on the open set $U^{\circ}$ of smooth and nodal members (Theorem~\ref{thm:thetasection}) this is the zero set of the transversality section
$\sigma_{\theta}$, the denominator of $g^{\CC}=G/\sigma_{\theta}$ (Proposition~\ref{prop:levwronsk}).
\end{proposition}

\begin{proof}
Suppose no member of $B$ is integral. After normalizing $B$ and a finite base change, the family of conics over $B$ is the family of images of a family of
genus-zero stable maps of $L$-degree two over a smooth complete curve $B'$, with the same cycles (the space of stable maps is proper). A stable map of this kind with
reducible domain maps onto two lines or onto one line, since a contracted component would need three special points; a double cover of a line has integral
image. Since lines are integral (Lemma~\ref{lem:allornothing}), every member is a map $f_{b}\colon\PP^{1}\to Z$ birational onto its image with
$\theta\circ df_{b}\ne0$; this is a constant in $H^{0}(K_{\PP^{1}}\otimes f_{b}^{*}L)=\CC$ (Lemma~\ref{lem:allornothing}), so $f_{b}$ is an immersion transverse to
$L^{\perp}$. So the domains form a smooth family of rational curves $\pi\colon S\to B'$, a $\PP^{1}$-bundle, with evaluation $F\colon S\to Z$. The section
$\sigma=\theta\circ dF|_{T_{S/B'}}$ of $\omega_{S/B'}\otimes F^{*}L$ restricts to a nonzero constant on every fibre, so it vanishes nowhere and
$F^{*}L\cong\omega_{S/B'}^{-1}$. For a $\PP^{1}$-bundle over a curve of genus $g$, $(\omega_{S/B'})^{2}=0$: $K_{S}^{2}=8(1-g)$ and $K_{S}\cdot\phi=-2$ for a fibre $\phi$, so
$(K_{S}-\pi^{*}K_{B'})^{2}=8(1-g)+4(2g-2)=0$. Hence $(F^{*}L)^{2}=0$. But the general members are pairwise distinct smooth conics, so $F(S)$ is a
surface, $F$ is generically finite onto it, and $(F^{*}L)^{2}=\deg(S\to F(S))\,(L|_{F(S)})^{2}>0$ since $L$ is ample. This contradiction proves the first assertion,
and the rest follows from Theorem~\ref{thm:thetasection}(3).
\end{proof}

Proposition~\ref{prop:nocurves} says that $\MC$ behaves like a
Stein manifold, or a Grauert tube with centre $M$, whose boundary at infinity is the integral divisor. Along a balanced integral divisor the metric has a pole of order two on a subspace of rank four and of order one
on the rest (Proposition~\ref{prop:intregular}(2), Remark~\ref{rem:intflat}(1)), so for $n\ge2$ the conformal structure does not extend there as a nondegenerate one.

\begin{proposition}[The $\alpha$-planes are projectively flat]
\label{prop:projflat}
Let $z\in Z$ and $\Sigma_{z}=\{[C]\in\MC:z\in C\}$. A null geodesic of $g^{\CC}$ with decomposable tangent $e\otimes h$, $h\in h_{z}$, lies in
$\Sigma_{z}$, and the geodesics of $\Sigma_{z}$ are exactly these. The tangent-direction map $\mu_{z}\colon\Sigma_{z}\to\CC^{2n}=\PP(T_{z}Z)\smallsetminus\PP(L^{\perp}_{z})$,
$[C]\mapsto[T_{z}C]$, is a local isomorphism which maps every geodesic of $\Sigma_{z}$ into an affine line, and every connected component of $\mu_{z}^{-1}(\ell)$,
for an affine line $\ell$, is a geodesic. So $\Sigma_{z}$ carries a flat projective structure with developing map $\mu_{z}$.
\end{proposition}

\begin{proof}
The inputs are Theorem~\ref{thm:lebrun}: (I1) $T\MC=E^{*}_{\CC}\otimes H^{*}_{\CC}$, with the Levi-Civita connection induced from connections $\nabla_{E},\nabla_{H}$ on the
factors; (I2) $\Sigma_{z}$ is totally geodesic with $T\Sigma_{z}=E^{*}_{\CC}\otimes h_{z}$, and a section $e$ of $E^{*}_{\CC}$ along a curve in $\Sigma_{z}$ is
$\nabla_{E}$-parallel if and only if $\mathrm{ev}_{z}(e)\in B_{z}$ is constant; and (I3) for a transverse conic, $d\mu_{z}$ at $[C]$ is $s\mapsto s'(z)$, which for
$s=e\otimes h$, $h$ spanning $h_{z}$, is $\mathrm{ev}_{z}(e)\otimes h'(z)$; so $d\mu_{z}$ is an isomorphism off $J$ (Lemma~\ref{lem:dmu}).

\emph{Step 1: $h_{z}$ is parallel along $\Sigma_{z}$.} Since $T\Sigma_{z}=E^{*}_{\CC}\otimes h_{z}$ is preserved by $\nabla$ along $\Sigma_{z}$ and
$\nabla_{X}(e\otimes h)=\nabla_{E,X}e\otimes h+e\otimes\nabla_{H,X}h$, we get $\nabla_{H,X}h\in h_{z}$: for a local nowhere vanishing section $h$ of $h_{z}$ there is a
$1$-form $\lambda$ on $\Sigma_{z}$ with $\nabla_{H,X}h=\lambda(X)h$.

\emph{Step 2: along a geodesic the $E$-factor is parallel up to scale.} Let $\gamma$ be a geodesic with $\gamma(0)\in\Sigma_{z}$ and $\gamma'(0)\in T\Sigma_{z}$. It stays
in $\Sigma_{z}$, and $\gamma'(t)=e(t)\otimes h(t)$ with $h(t)=h|_{\gamma}$. The geodesic equation reads $(\nabla_{E,\gamma'}e+\lambda(\gamma')e)\otimes h=0$, so
$\exp\big(\int_{0}^{t}\lambda(\gamma')\big)e(t)$ is $\nabla_{E}$-parallel, and by (I2) $\mathrm{ev}_{z}(e(t))=f(t)\,v$ with $f(t)\ne0$ and a fixed vector $v\in B_{z}$.

\emph{Step 3: the image is an affine line.} Write $\mu_{z}(\gamma(t))=[u(t)]$ with $u(t)\in T_{z}Z\smallsetminus L^{\perp}_{z}$. By (I3) the lift can be chosen with
$u'(t)=c(t)\,v+a(t)\,u(t)$, $c(t)\ne0$, where $v$ is read in $L^{\perp}_{z}$. The right side is tangent to the plane $P=\langle u(0),v\rangle$ whenever $u(t)\in P$, so
$u(t)\in P$ for all $t$; $\PP(P)$ meets the hyperplane at infinity only in $[v]$, so $\mu_{z}\circ\gamma$ is an immersion into the affine line $\PP(P)\smallsetminus\{[v]\}$.

\emph{Step 4: the converse.} Let $\ell$ be an affine line with point at infinity $[v]$ and $[C]\in\mu_{z}^{-1}(\ell)$; the geodesic with initial velocity
$\mathrm{ev}_{z}^{-1}(v)\otimes h$ maps into $\ell$ by Steps~2--3. Since $\mu_{z}$ is a local isomorphism, $\mu_{z}^{-1}(\ell)$ is a smooth curve which is locally a
geodesic at each of its points, so its connected components are geodesics. A geodesic contained in $\Sigma_{z}$ has its tangent in the totally null $T\Sigma_{z}$, so it
is null and decomposable, and conversely.
\end{proof}

In the developing coordinate the curvature quartic obeys the symmetric-derivative law (Proposition~\ref{prop:symlaw}): locally $W_{z}=\partial^{4}\Upsilon$, and leaf constancy is the
statement that $\Upsilon$ can be taken to be a polynomial of degree at most four.

\begin{remark}[Leaf constancy as a Liouville theorem]
\label{rem:liouville}
Let $z$ be general and $\ell\subset\CC^{2n}$ a general affine line, with closure $\bar\ell\subset\PP(T_{z}Z)$. The closure in $\Gamma_{z}$ (\S\ref{sec:fibre}) of
$\mu_{z}^{-1}(\ell)$ is the compact curve $\widetilde\mu_{z}^{-1}(\bar\ell)$, the union of the closures of null geodesics of the $\alpha$-plane $\Sigma_{z}$, and it reaches
the singular locus: over the point at infinity $\bar\ell\cap\PP(L^{\perp}_{z})$ it meets the boundary divisors dominating the contact hyperplane, in $\deg\widetilde\mu_{z}$ points,
counted with multiplicity, and over $\ell$ it can meet
the jumping conics and the boundary points with transverse limit direction. If the curvature quartic $W_{z}$ is regular at all these points, it is a holomorphic map
from a compact curve to the vector space $V_{z}$, hence constant on it; and the general lines through a point fill $\CC^{2n}$. This is the analogue of the theorems of
\cite{Ye94} and \cite{Bel01}, in which a compact rational null geodesic forces conformal flatness: here closed null geodesics of the $\alpha$-planes, which must meet the
singular locus, force local symmetry as soon as the curvature is regular where they meet it. The Hartogs argument of \S\ref{sec:hartogs} is the form of this argument
that sees only the divisors of $\Gamma_{z}$ with divisorial image, which are exactly the divisors a general line meets.
\end{remark}

\section{Lines, nodal conics and the variety of minimal rational tangents}
\label{sec:lines}

Throughout this section $\mathrm{Pic}\,Z=\ZZ L$, except in Lemmas~\ref{lem:branchres}--\ref{prop:balcrit} and
Proposition~\ref{prop:thetamult}, which hold for any contact Fano manifold.

\subsection{Lines and the variety of minimal rational tangents}
\label{sec:vmrt}

Lines exist on $Z$ and cover it: since $\mathrm{Pic}\,Z=\ZZ L$ gives $b_{2}(Z)=1$ and $Z\not\cong\PP^{2n+1}$, there is an irreducible component of the
space of rational curves whose members have $L$-degree one and cover $Z$ \cite[Sect.~2.3]{Keb01}, \cite[Sect.~1.1]{Keb05}. For a general point $z$: every
line through $z$ is free, i.e.\ $TZ|_{\ell}$ is nef \cite[II.3.11]{Kol96} (for curves of bounded degree, as here and for the conics below, the exceptional set
is a finite union of proper closed subsets, so ``general'' suffices); and by \cite[Thm.~1.1(1),(3),(4)]{Keb05}, every line through $z$ is smooth, two distinct
lines through $z$ have distinct tangent directions at $z$ and meet only at $z$, and the lines through $z$ sweep out a cone over the variety of
their tangent directions \cite[\S6.1]{Keb05}. \cite[Thm.~1.1(2)]{Keb05} also asserts that the lines through a general $z$ of one
component form an irreducible family. Its proof rests on \cite[Prop.~3.2]{Keb05}, which has a gap \cite[Rem.~3.2]{Buc10}, and the statement is
listed as open in \cite[\S6]{BM19}. It is not used in this paper: the cone structure and the other two statements are proved in \cite[\S\S4--6]{Keb05}
without it, and the variety of minimal rational tangents is treated throughout as possibly reducible (Proposition~\ref{prop:secantgeneral}).

A line $\ell$ is \emph{standard} if $TZ|_{\ell}\cong\OO(2)\oplus\OO(1)^{n-1}\oplus\OO^{n+1}$. Every line through a general point is
standard \cite[IV.2.9]{Kol96}.

\begin{lemma}[Splitting along a line]
\label{lem:linesplit}
Let $\ell$ be a standard line. Then
\[
TZ|_{\ell}\cong\OO(2)\oplus\OO(1)^{n-1}\oplus\OO^{n+1},\qquad N_{\ell/Z}\cong\OO(1)^{n-1}\oplus\OO^{n+1},
\]
\[
L^{\perp}|_{\ell}\cong\OO(2)\oplus\OO(1)^{n-1}\oplus\OO^{n-1}\oplus\OO(-1),
\]
with $T\ell$ the summand $\OO(2)$ of $L^{\perp}|_{\ell}$ and $\OO(-1)$ paired with it by the Levi form.
\end{lemma}

\begin{proof}
The first two splittings are the definition and its quotient by $T\ell$. The line is integral (Lemma~\ref{lem:allornothing}), so $T\ell\subset L^{\perp}|_{\ell}$, and $L^{\perp}|_{\ell}$ has
degree $(n+1)-1=n$. The Levi form makes $L^{\perp}|_{\ell}$ self-dual up to $L|_{\ell}=\OO(1)$, so its degrees pair as $(a,1-a)$. A
summand $\OO(k)$, $k\ge2$, of the subbundle $L^{\perp}|_{\ell}\subset TZ|_{\ell}$ maps injectively to the single summand $\OO(2)$ of
$TZ|_{\ell}$, so $T\ell$ is the only one. Its partner is $\OO(-1)$, and the remaining $2n-2$ summands have degrees $a,1-a\le1$, hence
pair as $(1,0)$.
\end{proof}

For a general $z$ the tangent directions of the lines through $z$ form the \emph{variety of minimal rational tangents}
(VMRT) $X_{z}\subset\PP(L^{\perp}_{z})\cong\PP^{2n-1}$, with affine cone $\hat X_{z}\subset L^{\perp}_{z}$. At $z$ the fibre of the
positive part $\OO(1)^{n-1}\subset N_{\ell/Z}$ is $\hat T_{[\ell]}\hat X_{z}/T_{z}\ell$, where $\hat T_{[\ell]}\hat X_{z}$ is the
affine tangent space of the cone; dually, the trivial summand $\OO^{n+1}\subset N^{*}_{\ell}$ annihilates the positive part of
$N_{\ell}$.

\begin{proposition}[Secant non-defectiveness]
\label{prop:secantgeneral}
Let $z\in Z$ be general, $X_{z}=X_{1}\sqcup\dots\sqcup X_{k}$ the decomposition into irreducible components, $\hat X_{i}$ the
affine cones, and $\omega$ the Levi form on $L^{\perp}_{z}$.
\begin{enumerate}
\item $X_{z}$ is nondegenerate in $\PP(L^{\perp}_{z})$, and its components are smooth, pairwise disjoint and Legendrian of
dimension $n-1$.
\item $\omega(\hat X_{i},\hat X_{j})\not\equiv0$ for some $i,j$; and $\omega(\hat X_{i},\hat X_{j})\equiv0$ only if $i=j$ and $X_{i}$
is a linear $\PP^{n-1}$.
\item If $\omega(\hat X_{i},\hat X_{j})\not\equiv0$, the join of $X_{i}$ and $X_{j}$ is $\PP(L^{\perp}_{z})$. In particular
$X_{z}$ is not secant-defective.
\end{enumerate}
\end{proposition}

\begin{proof}
(1) Nondegeneracy is Kebekus's spanning theorem \cite[Thm.~4.4]{Keb01}. \cite[Thm.~3.1(iii)]{Buc10} asserts that each
component is nondegenerate, by an extension of the proof of \cite[Thm.~4.4]{Keb01} that does not cover components permuted by monodromy. That
assertion is not used here; only the nondegeneracy of the union $X_{z}$ is. Smoothness and the Legendrian property, without
irreducibility, are \cite[Thm.~3.1]{Buc10}; smoothness means that the components are disjoint.

(2) If $\omega$ vanished on all pairs, the span of $\hat X_{z}$ would be isotropic, of dimension $\le n<2n$, contradicting (1).
Let $S_{i}$ be the span of $\hat X_{i}$. It contains a Lagrangian tangent space, so $S_{i}^{\perp}\subset S_{i}$. If
$\omega(\hat X_{i},\hat X_{j})\equiv0$ then $S_{j}\subset S_{i}^{\perp}\subset S_{i}$ and likewise $S_{i}\subset S_{j}$; so
$S=S_{i}=S_{j}$ is isotropic and contains a Lagrangian, hence is Lagrangian, and $X_{i}$, $X_{j}$ are irreducible of dimension
$n-1$ in $\PP(S)\cong\PP^{n-1}$: $X_{i}=X_{j}=\PP(S)$ and $i=j$.

(3) By Terracini's lemma, $\dim\hat J(X_{i},X_{j})=2n-\dim(\hat T_{a}\cap\hat T_{b})$ for general $a\in\hat X_{i}$, $b\in\hat X_{j}$. We
show that the addition map $\hat X_{i}\times\hat X_{j}\to L^{\perp}_{z}$ is dominant (the Legendrian secant lemma
\cite[proof of Prop.~17(2)]{LM07}, \cite[Prop.~3.5]{Buc10}). Otherwise take a curve $(a(t),b(t))$ in a positive-dimensional fibre,
$a+b=s$, through a general pair with $\omega(a,b)\neq0$. Then $\omega(a,b)=\omega(a,s)$ and
\[
\tfrac{d}{dt}\,\omega(a,s)=\omega(a',s)=\omega(a',a)+\omega(a',b)=\omega(a',a)-\omega(b',b)=0,
\]
because $a'=-b'$, the radial vector $a$ and the velocity $a'$ both lie in the Lagrangian space $\hat T_{a}$, and likewise $b,b'\in\hat T_{b}$.
So $\omega(a(t),b(t))=c\neq0$ is constant. On the
completed curve $a$ has a pole of some order $k\ge1$ (the fibre is positive-dimensional and affine); then $p=t^{k}a$ and $q=-t^{k}b$
are regular with $p(0)=q(0)\neq0$ and $p-q=t^{k}s$. By disjointness of the components, $i=j$. For a smooth Lagrangian cone,
$\omega(r(x),r(y))$ vanishes to third order along the diagonal, so $\omega(p,q)=O(t^{3k})$; but $\omega(p,q)=-t^{2k}c$, a
contradiction.
\end{proof}

\subsection{Nodal conics}
\label{sec:nodal}

Let $C=C_{1}\cup C_{2}$ be two lines meeting transversally at one point $p$, $N:=N_{C/Z}$, and $N_{p}$ its fibre at $p$.

\begin{lemma}[Nodal conics are smooth points]
\label{lem:nodalsmooth}
If $N_{C_{i}/Z}$ is nef for $i=1,2$, then $H^{1}(C,N)=0$, $N$ is globally generated at $p$, and $\MCbar$ is smooth of dimension
$4n$ at $[C]$; near $[C]$ the nodal conics form a smooth hypersurface whose reduced equation is the smoothing parameter of the
node.
\end{lemma}

\begin{proof}
By the proof of Lemma~\ref{lem:branchres} below, $N|_{C_{i}}$ is an upward elementary modification of $N_{C_{i}/Z}$, hence nef and
globally generated. In the gluing sequence $0\to N\to N|_{C_{1}}\oplus N|_{C_{2}}\to N_{p}\to0$ the middle term surjects onto $N_{p}$
on global sections, so $H^{1}(N)=0$. Hence $H^{0}(N)\to\mathcal Ext^{1}(\Omega_{C},\OO_{C})_{p}\cong\CC$ is onto, and the Hilbert
scheme is smooth at $[C]$ and submersive onto the versal deformation $xy=t$ of the node.
\end{proof}

From now on in this subsection the branches $C_{1},C_{2}$ are standard lines.

\begin{lemma}[Branch restriction]
\label{lem:branchres}
Let $\hat T_{i}\subset T_{p}Z$ be the fibre at $p$ of the positive part $\OO(2)\oplus\OO(1)^{n-1}$ of $TZ|_{C_{i}}$ (for a line
through a general point, the affine tangent space of the VMRT cone at $[C_{i}]$), and $P_{i}=\hat T_{i}/T_{p}C_{i}$. Let
$N^{0}_{p}\subset N_{p}$ be the image of $T_{p}Z\to N_{p}$; it is a hyperplane, $N^{0}_{p}\cong T_{p}Z/(T_{p}C_{1}+T_{p}C_{2})$, and
$N_{p}/N^{0}_{p}\cong T_{p}C_{1}\otimes T_{p}C_{2}$ is the smoothing direction. Call the node \emph{general} if $T_{p}C_{3-i}\notin\hat T_{i}$
for $i=1,2$. If the node is general, then $N|_{C_{i}}\cong\OO(1)^{n}\oplus\OO^{n}$, and the fibre $V_{i}\subset N_{p}$ of its positive
part satisfies
\[
V_{i}\cap N^{0}_{p}=\bar P_{i},\qquad V_{i}\not\subset N^{0}_{p},
\]
where $\bar P_{i}\cong P_{i}$ is the image of $\hat T_{i}$ in $N^{0}_{p}$, of dimension $n-1$. Moreover $\hat T_{i}\subset L^{\perp}_{p}$.
\end{lemma}

\begin{proof}
In coordinates $(x,y,w_{1},\dots,w_{2n-1})$ at $p$ with $C_{1}=\{y=w=0\}$, $C_{2}=\{x=w=0\}$, the ideals are $I_{C}=(xy,w)$ and
$I_{C_{1}}=(y,w)$. The map $N^{*}_{C/Z}|_{C_{1}}\to N^{*}_{C_{1}/Z}$ sends $xy\mapsto x\,dy$ and the $w$'s isomorphically. Dually,
$N_{C_{1}/Z}\hookrightarrow N|_{C_{1}}$ is an isomorphism off $p$ with cokernel of length one at $p$ in the $T_{p}C_{2}$-direction: $N|_{C_{1}}$ is
the elementary modification of $N_{C_{1}/Z}$ raising that direction by $\OO(p)$. At a general node the raised direction lies in a
trivial summand, turning one $\OO$ into $\OO(1)$: $N|_{C_{1}}\cong\OO(1)^{n}\oplus\OO^{n}$. Fix $q_{1}\in C_{1}\smallsetminus p$. Evaluation
at $p$ maps $H^{0}(N|_{C_{1}}(-q_{1}))$ isomorphically onto $V_{1}$. A section $s$ of $N|_{C_{1}}$ lies in $N_{C_{1}/Z}$ iff $s(xy)$
vanishes at $p$, i.e.\ iff $s(p)\in N^{0}_{p}$; so $V_{1}\cap N^{0}_{p}$ is the image of $H^{0}(N_{C_{1}/Z}(-q_{1}))$, the sections of the
positive part, whose values fill $P_{1}$, and $(N_{C_{1}/Z})_{p}\to N_{p}$ has kernel $T_{p}C_{2}\bmod T_{p}C_{1}\notin P_{1}$. The
remaining section has a simple pole in the $T_{p}C_{2}$-direction, so its value is not in $N^{0}_{p}$. Finally
$\hat T_{i}\subset L^{\perp}_{p}$ because $L^{\perp}|_{C_{i}}$ contains the positive part $\OO(2)\oplus\OO(1)^{n-1}$
(Lemma~\ref{lem:linesplit}).
\end{proof}

\begin{proposition}[Balanced decomposition]
\label{prop:balanced}
At a general node, with $q_{i}\in C_{i}\smallsetminus p$,
\[
a:=\dim(V_{1}\cap V_{2})=h^{0}\big(C,N(-q_{1}-q_{2})\big),\qquad
N\cong\OO_{C}(1,1)^{a}\oplus\OO_{C}(1,0)^{n-a}\oplus\OO_{C}(0,1)^{n-a}\oplus\OO_{C}^{a},
\]
where $\OO_{C}(a_{1},a_{2})$ restricts to $\OO(a_{i})$ on $C_{i}$. The nodal conic is \emph{balanced} if $a=0$:
$N\cong\OO_{C}(1,0)^{n}\oplus\OO_{C}(0,1)^{n}$; this is a Zariski-open condition.
\end{proposition}

\begin{proof}
For a bundle $F$ on $C$, $0\to F\to\nu_{*}\nu^{*}F\to F_{p}\to0$ with $\nu$ the normalization. For $F=N(-q_{1}-q_{2})$,
$N|_{C_{i}}(-q_{i})\cong\OO^{n}\oplus\OO(-1)^{n}$ has $H^{1}=0$ and $H^{0}\cong V_{i}$ by evaluation at $p$; so
$0\to H^{0}(F)\to V_{1}\oplus V_{2}\to N_{p}\to H^{1}(F)\to0$ with the difference map, and $h^{0}(F)=\dim(V_{1}\cap V_{2})$. For the
splitting choose a basis of $N_{p}$ adapted to both flags: $e_{1..a}$ a basis of $V_{1}\cap V_{2}$, $e_{a+1..n}$ completing to $V_{1}$,
$f_{1..n-a}$ completing to $V_{2}$, $g_{1..a}$ completing to $N_{p}$; take $W_{1}=\langle f,g\rangle$, $W_{2}=\langle e_{a+1..n},g\rangle$ as the
trivial complements. Each basis vector then spans a sub-line-bundle of bidegree $(1,1)$, $(1,0)$, $(0,1)$ or $(0,0)$. Openness is upper
semicontinuity of $h^{0}$.
\end{proof}

\begin{proposition}[Balancedness in terms of the VMRT]
\label{prop:balcrit}
Let $C=C_{1}\cup C_{2}$ be two standard lines meeting transversally at $p$, $u_{i}$ spanning $T_{p}C_{i}$. If the node is
general, then $a\ge\dim(\hat T_{1}\cap\hat T_{2})$, and
\begin{enumerate}
\item if $d\theta(u_{1},u_{2})\neq0$, then $a=\dim(\hat T_{1}\cap\hat T_{2})$, and $V_{1}\cap V_{2}=\bar P_{1}\cap\bar P_{2}$;
\item if $d\theta(u_{1},u_{2})=0$, then $a\ge1$.
\end{enumerate}
Hence $C$ is balanced if and only if $\hat T_{1}\cap\hat T_{2}=0$ and $d\theta(u_{1},u_{2})\neq0$.
\end{proposition}

\begin{proof}
Use the coordinates of Lemma~\ref{lem:branchres} with $\partial_{x}=u_{1}$, $\partial_{y}=u_{2}$ at $p$, trivialize $L$ and write
$\theta=\theta_{x}dx+\theta_{y}dy+\sum\theta_{j}dw_{j}$, so $D:=d\theta(u_{1},u_{2})=(\partial_{x}\theta_{y}-\partial_{y}\theta_{x})(p)$. Let
$\varepsilon\in N_{p}$ be given by $\varepsilon(xy)=1$, $\varepsilon(w_{j})=0$, and $\varepsilon^{*}$ the form with kernel $N^{0}_{p}$ and
$\varepsilon^{*}(\varepsilon)=1$.

\emph{Step 1.} Lines are integral, so $\theta_{x}$ vanishes on $C_{1}$, and $\theta_{y}|_{C_{1}}=x\,b_{1}$. For a local section $s$ of
$N|_{C_{1}}$ put $\theta_{1}(s)=b_{1}s(xy)+\sum\theta_{j}s(w_{j})$, the value of $\theta$ on the meromorphic field along $C_{1}$ with
$y$-component $s(xy)/x$ representing $s$. So $\theta_{1}\in H^{0}(C_{1},N^{*}|_{C_{1}}\otimes L)$ extends $\theta\colon N_{C_{1}/Z}\to L|_{C_{1}}$.
Define $\theta_{2}$ symmetrically. At $p$ both restrict to $\theta_{p}$ on $N^{0}_{p}$, and $\theta_{1}(p)-\theta_{2}(p)=D\,\varepsilon^{*}$.

\emph{Step 2.} Let $t_{1}\in H^{0}(N_{C_{1}/Z})$ with $t_{1}(p)=u_{2}\bmod u_{1}$, $\bar x$ a rational function on $C_{1}$ with divisor
$p-q_{1}$ and $d\bar x(u_{1})=1$, and $s_{1}=t_{1}/\bar x\in H^{0}(N|_{C_{1}}(-q_{1}))$; then $\varepsilon^{*}(s_{1}(p))=1$ and
$V_{1}=\bar P_{1}\oplus\CC s_{1}(p)$. The operator $V\mapsto(\iota_{V}d\theta+d(\theta(V)))|_{TC_{1}}$ is well defined on sections of
$N_{C_{1}/Z}$ with values in $K_{C_{1}}\otimes L\cong\OO(-1)$, so it kills $t_{1}$; at $p$, where $\theta(t_{1})=0$, this reads
$\partial_{x}\theta(t_{1})(p)=d\theta(u_{1},u_{2})$. Hence $\theta_{1}(p)(s_{1}(p))=D$, and symmetrically $\theta_{2}(p)(s_{2}(p))=-D$.

\emph{Step 3.} $\bar P_{1},\bar P_{2}$ lie in the image of $L^{\perp}_{p}$, where $\theta_{1}(p)=\theta_{2}(p)=0$; and
$\theta_{1}(p)(s_{2}(p))=-D+D=0$. So $\theta_{1}(p)$ annihilates $V_{2}$, $\theta_{2}(p)$ annihilates $V_{1}$, $\theta_{1}(p)$ vanishes on
$\bar P_{1}$ and equals $D$ on $s_{1}(p)$.

\emph{Step 4.} At a general node $\hat T_{1}\cap\hat T_{2}$ meets $\langle u_{1},u_{2}\rangle$ only in $0$ and maps isomorphically onto
$\bar P_{1}\cap\bar P_{2}\subset V_{1}\cap V_{2}$; so $a\ge\dim(\hat T_{1}\cap\hat T_{2})$. If $D\neq0$, the kernel of $\theta_{1}(p)$ on $V_{1}$
is $\bar P_{1}$, so $V_{1}\cap V_{2}\subset\bar P_{1}$, and likewise $\subset\bar P_{2}$: this is (1). If $D=0$, the nonzero form $\theta_{1}(p)$
annihilates $V_{1}+V_{2}$, so $\dim(V_{1}+V_{2})\le2n-1$: (2). At a non-general node some $u_{3-i}$ lies in $\hat T_{1}\cap\hat T_{2}$ and
$N|_{C_{i}}\not\cong\OO(1)^{n}\oplus\OO^{n}$, so both sides of the equivalence fail.
\end{proof}

\begin{corollary}[The general nodal conic is balanced]
\label{cor:balgeneric}
Let $z$ be general, $[C_{1}]\in X_{i}$, $[C_{2}]\in X_{j}$ general with $\omega(\hat X_{i},\hat X_{j})\not\equiv0$, and $n\ge3$. Then
$C_{1}\cap C_{2}=\{z\}$ and $C_{1}\cup C_{2}$ is balanced. Consequently the general member of a component of the nodal divisor whose
nodes fill $Z$ is balanced, unless all its branch pairs over a general point lie in $X_{i}\times X_{i}$ with $X_{i}$ linear.
\end{corollary}

\begin{proof}
Only finitely many lines join two points, so the lines through $z$ meeting $C_{1}$ again form a family of dimension $\le1<n-1=\dim X_{j}$.
Lines through a general point are standard and $\hat T_{i}$ is the affine tangent space of the cone at $[C_{i}]$. By
Proposition~\ref{prop:secantgeneral}(3) and Terracini's lemma, $\hat T_{1}\cap\hat T_{2}=0$ for a general pair, and $\omega(u_{1},u_{2})\neq0$
on a dense open subset of $\hat X_{i}\times\hat X_{j}$; apply Proposition~\ref{prop:balcrit}. If the nodes of a component of the nodal
divisor fill $Z$, its fibre over a general point has dimension $(4n-1)-(2n+1)=2n-2=\dim(X_{i}\times X_{j})$, so its branch pairs fill
components $X_{i}\times X_{j}$; the exception is Proposition~\ref{prop:secantgeneral}(2).
\end{proof}

\begin{proposition}[The transversality section at a node]
\label{prop:thetamult}
Let $C_{0}=\ell_{1}\cup\ell_{2}$ be as in Lemma~\ref{lem:nodalsmooth}, with node $z$ and branch tangents $u_{i}$. Near $[C_{0}]$,
$\sigma_{\theta}=t\cdot h$, where $t$ is the smoothing parameter and, at every nodal conic $[C']$ near $[C_{0}]$ with branch tangents
$u'_{1},u'_{2}$, $h([C'])=c\,d\theta(u'_{1},u'_{2})$ for a unit $c$. In particular, if $d\theta(u_{1},u_{2})\neq0$, then $h$ is a unit at $[C_{0}]$ and $\sigma_{\theta}$ has a simple zero along the nodal divisor; if
$d\theta(u_{1},u_{2})=0$, then $\sigma_{\theta}$ vanishes to order at least two at $[C_{0}]$ in the sense of the maximal ideal, which need not mean divisibility
by $t^{2}$.
\end{proposition}

\begin{proof}
Take a smoothing with smooth total surface, coordinates $(x,y)$ at the node, $C_{t}=\{xy=t\}$, evaluation $\varphi$, and
$\varphi^{*}\theta=a\,dx+b\,dy$ in a trivialization of $L$. Lines are integral, so $a=y\,a_{1}$, $b=x\,b_{1}$. The field
$v=x\partial_{x}-y\partial_{y}$ is tangent to $C_{t}$ and $\eta=dx/x$ generates $\omega_{C_{t}}$ near the node with $\eta(v)=1$; so
$\sigma_{\theta}|_{C_{t}}=\theta(v)\eta=(xa-yb)\eta=t(a_{1}-b_{1})\eta$, and $(a_{1}-b_{1})(0)=-d\theta(u_{1},u_{2})$. Since $t$ is a reduced
equation of the nodal locus (Lemma~\ref{lem:nodalsmooth}), $\sigma_{\theta}=t\cdot h$ with $h$ as stated.
\end{proof}

\begin{remark}[The jumping invariant at the boundary]
\label{rem:jinv}
On a nodal conic $L|_{C}\cong\OO_{C}(1,1)\cong\OO_{C}(q_{1}+q_{2})$, so $j(C)=h^{0}(N(-q_{1}-q_{2}))=a$: a balanced nodal conic has $j=0$, the
value of a non-jumping transverse conic, and at a non-general node $j\ge1$. Since $\{j\ge1\}$ is closed in the open set $\MCbar^{\,\mathrm{lci}}$
(Definition~\ref{def:jinv}), which contains the smooth and the nodal conics, the closure of $J$ meets the
nodal divisor only in the unbalanced locus, where the jumping summands of $N$ are the $\OO_{C}(1,1)^{a}$, with fibre $V_{1}\cap V_{2}$ at the
node.
\end{remark}

\begin{corollary}[Balanced transverse conics exist]
\label{cor:balexist}
Let $Z$ be a compact contact Fano manifold with $\mathrm{Pic}\,Z=\ZZ L$ and $n\ge3$. Then some irreducible component of the Hilbert scheme of conics of $Z$
has a smooth transverse balanced general member.
\end{corollary}

\begin{proof}
By Proposition~\ref{prop:secantgeneral}(2) there are components $X_{i},X_{j}$ of the VMRT at a general $z$ with $\omega(\hat X_{i},\hat X_{j})\not\equiv0$, and by
Corollary~\ref{cor:balgeneric} the nodal conic $C_{0}$ of two general lines through $z$ in these components is balanced, with $d\theta(u_{1},u_{2})\neq0$. By
Lemma~\ref{lem:nodalsmooth} the Hilbert scheme is smooth of dimension $4n$ at $[C_{0}]$ and its general member near $[C_{0}]$ is a smooth conic. Its jumping invariant
vanishes, since $j(C_{0})=0$ (Remark~\ref{rem:jinv}) and $j$ is upper semicontinuous; and it is transverse, since $\sigma_{\theta}=t\cdot h$ with $h([C_{0}])\neq0$
(Proposition~\ref{prop:thetamult}) and Theorem~\ref{thm:thetasection}(3). A smooth transverse conic with $j=0$ is balanced.
\end{proof}


\section{The incidence variety and the tangent-direction map}
\label{sec:incidence}

\subsection{The tangent-direction map and its resolution}
\label{sec:resolution}

Let $\overline F\subset\MCbar\times Z$ be the universal conic of \S\ref{sec:Wz}, $F=\bar\pi^{-1}(\MC\cup J)$, $\overline F_{z}=\bar p^{-1}(z)$ the
conics through $z$, and $F^{\circ}_{z}=\overline F_{z}\cap F$ the transverse conics through $z$. The map $\bar p$ is proper, $\overline F$
is irreducible of dimension $4n+1$, and the \emph{tangent-direction map}
\[
\mu\colon\overline F\dashrightarrow\PP(TZ),\qquad([C],z)\mapsto[T_{z}C],
\]
is a morphism where $C$ is smooth at $z$. It is undefined on the \emph{node section} $\mathcal B=\{([A\cup A'],p)\}\subset\overline F$ of the nodal
conics, $p$ the node.

Near a nodal conic with nef branches, $\overline F$ is smooth (Lemma~\ref{lem:nodalsmooth}) and has the local form $\{(x,y,\vec w)\}$ with
$\bar\pi(x,y,\vec w)=(t,\vec w)=(xy,\vec w)$, branches $\{y=0\}$, $\{x=0\}$ and $\mathcal B=\{x=y=0\}$, smooth of codimension two. The vertical
field $V=x\partial_{x}-y\partial_{y}$ is tangent to the fibres. Write $\alpha(\vec w)=\partial_{x}\bar p$, $\beta(\vec w)=\partial_{y}\bar p$ at $x=y=0$ for
the branch tangents; they lie in $L^{\perp}$ (Lemma~\ref{lem:allornothing}). Then $d\bar p(V)=x\alpha-y\beta+(x^{2}\bar p_{xx}-y^{2}\bar p_{yy})+O(3)$,
the cross term cancelling. Let $\widetilde F=\mathrm{Bl}_{\mathcal B}\overline F$ with exceptional divisor $E$; in the chart $x=ay$, $E=\{y=0\}$,
$t=ay^{2}$ and
\[
d\bar p(V)=y\big[a\alpha-\beta+y(a^{2}\bar p_{xx}-\bar p_{yy})+O(y^{2})\big].
\]
In particular
\begin{equation}
\operatorname{ord}_{E}t=2,\qquad\operatorname{ord}_{E}x=\operatorname{ord}_{E}y=1 .\label{eq:4-1}
\end{equation}

\begin{proposition}[Resolution at nodes]
\label{prop:resolution}
Near a nodal conic with independent branch tangents $\alpha,\beta$, $\mu$ extends to a morphism $\widetilde\mu$ on $\widetilde F$, and the
exceptional curve over a node maps by $a\mapsto[a\alpha-\beta]$ isomorphically onto the secant line
$\langle[\alpha],[\beta]\rangle\subset\PP(L^{\perp}_{p})$.
\end{proposition}

\begin{proof}
Projectivizing the displayed expression removes the factor $y$; on $E$ the value is $[a\alpha-\beta]$, defined for all $a$ when $\alpha,\beta$ are
independent. The other chart $y=bx$ gives $[\alpha-b\beta]$.
\end{proof}

\subsection{The conics through a general point}
\label{sec:fibre}

Let $z$ be general. The conics with a node at $z$ form $\mathcal B_{z}=\overline F_{z}\cap\mathcal B$, of codimension two in $\overline F_{z}$ where
non-empty; by generic smoothness applied to $\bar p$ and $\bar p|_{\mathcal B}$, $\overline F_{z}$ is smooth along $\mathcal B_{z}$ and meets
$\mathcal B$ transversally there, so $\widetilde F_{z}:=\widetilde F\cap\bar p^{-1}(z)=\mathrm{Bl}_{\mathcal B_{z}}\overline F_{z}$, with exceptional divisor
$E_{z}=E\cap\widetilde F_{z}$. The induced
\[
\widetilde\mu_{z}\colon\widetilde F_{z}\dashrightarrow\PP(T_{z}Z)\cong\PP^{2n}
\]
is a rational map. Let $\Gamma_{z}$ be the normalization of the closure of its graph, so that
\[
\widetilde\mu_{z}\colon\Gamma_{z}\longrightarrow\PP(T_{z}Z)
\]
is a proper morphism. We denote again by $F^{\circ}_{z}\subset\Gamma_{z}$ the preimage of the transverse conics, and by $F^{\mathrm{ne}}_{z}\subset F^{\circ}_{z}$
that of the transverse non-jumping conics; on $F^{\circ}_{z}$, $\widetilde\mu_{z}=\mu$. Write $\PP(L^{\perp}_{z})\subset\PP(T_{z}Z)$ for the
\emph{contact hyperplane} and $\CC^{2n}:=\PP(T_{z}Z)\smallsetminus\PP(L^{\perp}_{z})$. A transverse conic through $z$ has $[T_{z}C]\in\CC^{2n}$. By
\S\ref{sec:Wz}, $W_{z}$ is a meromorphic map $\Gamma_{z}\dashrightarrow V_{z}$, holomorphic on $F^{\mathrm{ne}}_{z}$. By Theorem~\ref{thm:irred},
whose proof is independent of what precedes it, $\Gamma_{z}$ has pure dimension $2n$ and every irreducible component of it contains transverse non-jumping conics
(there is only one component when $Z$ is a twistor space or $d_{2}=1$); so $F^{\mathrm{ne}}_{z}$ is a dense open subset.

\begin{lemma}[The differential of $\mu_{z}$]
\label{lem:dmu}
Let $C$ be a conic smooth at $z$ with $H^{1}(N_{C/Z}(-z))=0$. Then $\overline F_{z}$ is smooth at $[C]$ with tangent space
$H^{0}(N_{C/Z}(-z))$, and
\[
d\mu_{z}\colon H^{0}(N_{C/Z}(-z))\to N_{C/Z,z}\otimes T^{*}_{z}C=T_{[T_{z}C]}\PP(T_{z}Z),\qquad s\mapsto s'(z),
\]
the derivative of a section at its zero. Its kernel is $H^{0}(N_{C/Z}(-2z))$. For a transverse conic, $d\mu_{z}$ is identified with
$\mathrm{ev}_{z}\colon H^{0}(B|_{C})\to B_{z}$ (up to the isomorphisms $N(-z)\cong B$, $N_{z}\otimes T^{*}_{z}C\cong B_{z}$), and its kernel
has dimension $j(C)$.
\end{lemma}

\begin{proof}
The first statements are the deformation theory of curves through a fixed point, and the description of the differential is the
first-order variation of the tangent line. For a transverse conic, $N=B(1)$, so $N(-z)\cong B$ and $N(-2z)\cong B(-1)$, whose $h^{0}$ is
$j(C)$; writing a section of $N(-z)$ as $b\otimes\sigma_{z}$ with $\sigma_{z}$ a section of $\OO(1)$ vanishing at $z$, its derivative at $z$ is
$b(z)\otimes\sigma_{z}'(z)$.
\end{proof}

\begin{lemma}[The tangent-direction map on a fibre]
\label{lem:muz}
Let $z$ be general.
\begin{enumerate}
\item $\widetilde\mu_{z}$ is surjective and generically finite, and it is a local isomorphism into $\CC^{2n}$ at every point of $F^{\mathrm{ne}}_{z}$.
\item The nodal conics $A\cup A'$ with $z\in A$ away from the node map into $X_{z}\subset\PP(T_{z}Z)$, of codimension $n+1$; we write
$D'_{z}$ for their closure.
\item Where the branch tangents are independent, $E_{z}$ maps onto secant lines of $X_{z}$. If a component of $D_{\mathrm{nod}}$ has nodes filling
$Z$, its conics with node at $z$ form a family of dimension $2n-2$, a union of components $X_{i}\times X_{j}$ of the space of branch pairs; a
component of $E_{z}$ whose branch pairs come from $X_{i}\times X_{j}$ with $\omega(\hat X_{i},\hat X_{j})\not\equiv0$ dominates $\PP(L^{\perp}_{z})$.
\end{enumerate}
\end{lemma}

\begin{proof}
(1) On $F^{\mathrm{ne}}_{z}$, $d\mu_{z}=\mathrm{ev}_{z}$ is an isomorphism by Lemma~\ref{lem:dmu}. Hence $\widetilde\mu_{z}$ is generically finite;
it is proper and every component of $\Gamma_{z}$ has dimension $2n$ and meets $F^{\mathrm{ne}}_{z}$, so it is surjective on each component. (2) The tangent direction at $z$ is that of the line $A$.
(3) Proposition~\ref{prop:resolution}; the dimension count is that of Corollary~\ref{cor:balgeneric}, and the union of the secant lines over
$X_{i}\times X_{j}$ is the join, which is $\PP(L^{\perp}_{z})$ by Proposition~\ref{prop:secantgeneral}(3).
\end{proof}

\emph{The affine structure of $\CC^{2n}$, and the symmetric-derivative law.} The complement $\CC^{2n}=\PP(T_{z}Z)\smallsetminus\PP(L^{\perp}_{z})$ is an affine space over
$L^{\perp}_{z}\otimes L_{z}^{-1}$: for a line $\ell\not\subset L^{\perp}_{z}$, $T_{z}Z=\ell\oplus L^{\perp}_{z}$, and $T_{[\ell]}\PP(T_{z}Z)=\mathrm{Hom}(\ell,T_{z}Z/\ell)\cong
\mathrm{Hom}(\ell,L^{\perp}_{z})\cong L^{\perp}_{z}\otimes L_{z}^{-1}$ through $\theta\colon\ell\cong L_{z}$; these identifications are the translations of the affine space
of the lines spanned by $r+\xi$, $\xi\in L^{\perp}_{z}$, for a fixed $r$ with $\theta(r)\ne0$. A choice of $L^{1/2}_{z}$ identifies $L^{\perp}_{z}\otimes L_{z}^{-1}$ with
$B_{z}\otimes L_{z}^{-1/2}$, so the tangent spaces of $\CC^{2n}$ are identified with $B_{z}$ up to one constant factor, and the curvature quartic takes values in
$V_{z}=S^{4}B_{z}^{*}$. We write $\xi$ for an affine coordinate and $\partial$ for the flat derivative.

\begin{proposition}[The symmetric-derivative law]
\label{prop:symlaw}
Let $z\in Z$ and let $\Sigma'\subset\Sigma_{z}$ be an open set on which $\mu_{z}$ is injective, and regard $W_{z}|_{\Sigma'}$ as a function of $\xi\in\mu_{z}(\Sigma')\subset\CC^{2n}$ with
values in $S^{4}B_{z}^{*}$. Then $(\partial_{b}W_{z})(b_{1},\dots,b_{4})$ is symmetric in $(b,b_{1},\dots,b_{4})\in B_{z}^{5}$. Equivalently, near every point of
$\mu_{z}(\Sigma')$, $W_{z}=\partial^{4}\Upsilon$ for a holomorphic function $\Upsilon$, unique up to a polynomial of degree at most three in $\xi$.
\end{proposition}

\begin{proof}
Let $[C]\in\Sigma'$ and $X=e\otimes\sigma\in T_{[C]}\Sigma_{z}$, with $\sigma$ spanning $h_{z}$. By Lemma~\ref{lem:dmu}, $d\mu_{z}(X)=\mathrm{ev}_{z}(e)\otimes\sigma'(z)$, which under
the identifications above is $\vartheta\,\mathrm{ev}_{z}(e)$ for a scalar $\vartheta=\vartheta([C],\sigma)\ne0$; both $\vartheta$ and $\nabla_{X}R_{E}^{\CC}$ are linear in $\sigma$. By
Lemma~\ref{lem:S5}, $X(W_{z})=(\mathrm{ev}_{z}^{-1})^{*}(\nabla_{X}R_{E}^{\CC})=\mathcal T(\mathrm{ev}_{z}(e),\cdot,\cdot,\cdot,\cdot)$ with $\mathcal T$ totally symmetric, so
$\partial_{\mathrm{ev}_{z}(e)}W_{z}=\vartheta^{-1}\mathcal T(\mathrm{ev}_{z}(e),\cdot,\cdot,\cdot,\cdot)$, and $\mathrm{ev}_{z}$ is onto. For the potential, the symmetry says that for fixed
$b_{2},b_{3},b_{4}$ the $1$-form $b\mapsto W_{z}(b,b_{2},b_{3},b_{4})$ is closed; a primitive taken linearly in $W_{z}$ (by radial integration) is an $S^{3}B_{z}^{*}$-valued
function $\Upsilon_{3}$ with $W_{z}=\partial\Upsilon_{3}$ and $\partial\Upsilon_{3}$ symmetric, and three more steps give $\Upsilon$. The kernel of $\partial^{4}$ consists of the
polynomials of degree at most three.
\end{proof}

\begin{corollary}[Boundary values and leading poles]
\label{cor:symlaw}
Let $z\in Z$.
\begin{enumerate}
\item Let $y_{0}$ be a point of $\Gamma_{z}$ near which $\Gamma_{z}$ is smooth, $\widetilde\mu_{z}$ is a local isomorphism onto a neighbourhood of a point of
$\PP(L^{\perp}_{z})$, the points off $Y:=\widetilde\mu_{z}^{-1}(\PP(L^{\perp}_{z}))$ are transverse non-jumping conics, and $W_{z}$ is holomorphic. Regard
$W_{\infty}:=W_{z}|_{Y}$, through $\widetilde\mu_{z}$, as a function on an open cone of $L^{\perp}_{z}\otimes L_{z}^{-1}\cong B_{z}$, homogeneous of degree $0$. Then
\[
W_{\infty}=\tfrac{1}{24}\,\partial^{4}w_{\infty},\qquad w_{\infty}(\xi)=W_{\infty}(\xi)(\xi,\xi,\xi,\xi).
\]
\item Let $Y$ be a prime divisor of $\Gamma_{z}$, at whose generic point $\widetilde\mu_{z}$ is finite, with $\widetilde\mu_{z}(Y)\not\subset\PP(L^{\perp}_{z})$, and suppose that
$W_{z}$ has a pole of order $p\ge1$ along $Y$: $W_{z}=y^{-p}(W_{(0)}+O(y))$ with $y$ a local equation of $Y$ and $W_{(0)}|_{Y}\not\equiv0$. Let $t$ be a local equation of the
hypersurface $\widetilde\mu_{z}(Y)$. Then $W_{(0)}|_{Y}=c_{Y}\,(dt)^{4}$ for a function $c_{Y}$ on $Y$, where $dt$ is read in $T^{*}\CC^{2n}\cong B_{z}^{*}$.
\end{enumerate}
\end{corollary}

\begin{proof}
(1) In projective coordinates near the image point, with $\PP(L^{\perp}_{z})=\{x_{0}=0\}$, $W_{z}$ is holomorphic in $(x_{0},u)$, and its expansion in powers of $x_{0}$ is an
expansion $W_{z}=\sum_{j\ge0}W_{j}$ on a conic neighbourhood of infinity, with $W_{j}$ homogeneous of degree $-j$ in $\xi$ and $W_{0}=W_{\infty}$. As $\partial$ lowers the
degree by one, every $W_{j}$ satisfies the law of Proposition~\ref{prop:symlaw}, in particular $F:=W_{\infty}$. By Euler's identity $\xi^{m}\partial_{m}F=0$, and with the
symmetry $\partial_{m}(\xi^{i}F_{ijkl})=F_{mjkl}+\xi^{i}\partial_{i}F_{mjkl}=F_{mjkl}$; so $F=\partial F_{1}$ with $F_{1}=\iota_{\xi}F$, homogeneous of degree one, whose
derivative $F$ is symmetric. In the same way (Euler's identity now contributing $1$, $2$, $3$) $F_{1}=\frac12\partial(\iota_{\xi}F_{1})$,
$\iota_{\xi}F_{1}=\frac13\partial(\iota_{\xi}^{2}F_{1})$ and $\iota_{\xi}^{2}F_{1}=\frac14\partial w_{\infty}$, which gives the formula.
(2) At a generic point of $Y$, $Y\to\widetilde\mu_{z}(Y)$ is a local isomorphism and $\widetilde\mu_{z}$ is $(y,x')\mapsto(y^{e},x')=(t,x')$ for some $e\ge1$; nearby, the
complement of $Y$ consists of transverse non-jumping conics. On a branch, the most singular part of $\partial W_{z}$ is $-\frac{p}{e}\,y^{-p-e}\,dt\otimes W_{(0)}$; the derivatives
in $x'$, and those of the higher terms, are less singular. So $dt\otimes W_{(0)}$ is symmetric by Proposition~\ref{prop:symlaw}. If $Q\in S^{4}B_{z}^{*}$ and $dt\otimes Q$ is
symmetric, then for $b\in\ker dt$ and all $b_{1}$, the symmetry gives $0=(dt\otimes Q)(b,b_{1},\dots)=(dt\otimes Q)(b_{1},b,\dots)=dt(b_{1})\,Q(b,\dots)$, so $Q$ vanishes
on $\ker dt$ and $Q\in\CC\,(dt)^{4}$.
\end{proof}

Part (1) applies to the components of $E_{z}$ dominating $\PP(L^{\perp}_{z})$ at whose generic point $\Gamma_{z}\cong\widetilde F_{z}$ and $\widetilde\mu_{z}$ has
multiplicity one (the nearby conics are then smooth, transverse and non-jumping by Proposition~\ref{prop:thetamult} and Remark~\ref{rem:jinv}), and to the
balanced integral conics (Corollary~\ref{cor:intvanish}).

\subsection{Limit directions}
\label{sec:P2}

To control the boundary points of $\Gamma_{z}$ we use stable maps. Let $\overline{\mathcal M}_{0,1}(Z,2)$ be the space of stable maps of genus zero
with one marked point whose image has $L$-degree two \cite{FP97}; it is proper. Let $y\in\Gamma_{z}\smallsetminus F^{\circ}_{z}$. Since $F^{\mathrm{ne}}_{z}$
is dense open in $\Gamma_{z}$ (Theorem~\ref{thm:irred}(1)), there is an analytic arc $\gamma\colon\Delta\to\Gamma_{z}$ with $\gamma(0)=y$ and $\gamma(t)\in F^{\mathrm{ne}}_{z}$ for
$t\neq0$. Write $\gamma(t)=([C_{t}],z)$ and parametrize $C_{t}=f_{t}(\PP^{1})$, $f_{t}(p_{t})=z$. After a finite base change the family
$(f_{t},p_{t})$ extends over $\Delta$ to a family of stable maps; call its central member $(f_{0},p_{0})$, the \emph{stable limit} of $\gamma$. The cycle
of $f_{0}$ is that of the Hilbert point under $y$, and $\widetilde\mu_{z}(y)=\lim_{t\to0}[df_{t}(p_{t})]$. Lemma~\ref{lem:allornothing} applies to
maps: for $f\colon\PP^{1}\to Z$ with $f^{*}L\cong\OO(2)$, $\theta\circ df\in H^{0}(K_{\PP^{1}}\otimes f^{*}L)=\CC$, so $f$ is either integral or a
transverse immersion.

\begin{lemma}[Limit directions]
\label{lem:P2proved}
Let $z$ be general and $\gamma$ an arc as above with stable limit $(f_{0},p_{0})$. Then $\widetilde\mu_{z}(\gamma(0))\in\PP(L^{\perp}_{z})$, except in the
following two cases:
\begin{itemize}[leftmargin=2.5em]
\item[(i)] the domain of $f_{0}$ is irreducible and $f_{0}$ is a transverse immersion that is not an embedding; then
$\widetilde\mu_{z}(\gamma(0))=[df_{0}(p_{0})]\in\CC^{2n}$;
\item[(ii)] the domain of $f_{0}$ is irreducible, $f_{0}$ is integral and birational onto its image, and $df_{0}(p_{0})=0$.
\end{itemize}
\end{lemma}

\begin{proof}
\emph{The domain of $f_{0}$ is reducible.} A genus-zero stable map of degree two with one marked point and reducible domain has two components
of degree one, mapped onto lines $\ell_{1},\ell_{2}$, and possibly a contracted component, carrying $p_{0}$, between them; no other contracted
component is stable. If $p_{0}$ is a smooth point of a non-contracted component, the limit is the tangent direction at $z$ of a line through $z$,
which is integral (lines through the general point $z$ are smooth, so the map from the component is an isomorphism onto the line and its
differential does not vanish). Otherwise $p_{0}$ lies on the contracted component or at the node $q$, and both lines pass through $z=f_{0}(q)$.
Forgetting the marked point and contracting the contracted component, the family of domains near the node becomes $\{xy=t^{k}\}$, on which the
maps $f_{t}$ are defined, with $p_{t}$ tending to the node. If the tangent directions $\alpha,\beta$ of the lines at $z$ are independent, then
$df_{t}(x\partial_{x}-y\partial_{y})=x\alpha-y\beta+O\big((|x|+|y|)(|x|+|y|+|t|)\big)$ at $p_{t}$, so the limit lies in $\langle\alpha,\beta\rangle\subset L^{\perp}_{z}$. If
they are dependent, then $\ell_{1}=\ell_{2}$, because distinct lines through a general point have distinct tangent directions (\S\ref{sec:vmrt},
\cite[Thm.~1.1]{Keb05}),
and we are in the non-reduced case below.

\emph{The domain is irreducible and $df_{0}(p_{0})\neq0$.} The limit is $[df_{0}(p_{0})]$, by continuity of the derivative along the fibres of the
family of maps. If $f_{0}$ is integral, it lies in $L^{\perp}_{z}$. If $f_{0}$ is transverse and an embedding, the Hilbert limit of the $C_{t}$ is
$f_{0}(\PP^{1})$, a smooth transverse conic, so $\gamma(0)\in F^{\circ}_{z}$, which is excluded; otherwise we are in case (i).

\emph{The domain is irreducible and $df_{0}(p_{0})=0$.} Then $f_{0}$ is not a transverse immersion, so it is integral; it is either birational onto
its image, case (ii), or a double cover of a line $\ell$ through $z$.

\emph{The non-reduced case.} If $f_{0}$ is a double cover of a line or the two lines coincide, the Hilbert limit $C_{0}$ of the $C_{t}$ is a
subscheme with support a line $\ell$ through $z$, standard since $z$ is general. If $g_{i}(x,t)$ generate the ideal of the flat family near $z$,
then $dg_{i}(z,t)$ vanishes on $T_{z}C_{t}$, so every limit of tangent directions lies in the Zariski tangent space $T_{z}C_{0}$; we show
$T_{z}C_{0}\subset L^{\perp}_{z}$. Let $I_{\ell}\supset I$ be the ideals of $\ell$ and $C_{0}$. From $0\to I_{\ell}/I\to\OO_{C_{0}}\to\OO_{\ell}\to0$ and
$\chi(\OO_{C_{0}})=1$ we get $\chi(I_{\ell}/I)=0$. The sheaf $I_{\ell}/I$ is supported on $\ell$ with generic rank one; since $\OO_{C_{0}}$ has length
two at the generic point of $\ell$, $I_{\ell}^{2}\subset I$ generically, so $(I+I_{\ell}^{2})/I$ is torsion, and the torsion length $e$ below counts both torsion pieces; its quotient
$I_{\ell}/(I+I_{\ell}^{2})$ is a quotient of $N^{*}_{\ell}\cong\OO(-1)^{n-1}\oplus\OO^{n+1}$ (Lemma~\ref{lem:linesplit}), and the rest is torsion. Write
$I_{\ell}/I$ as an extension of a torsion sheaf of length $e\ge0$ by a quotient line bundle $\OO(b)$ of $N^{*}_{\ell}$. A nonzero map
$\OO(-1)^{n-1}\oplus\OO^{n+1}\to\OO(b)$ forces $b\ge-1$, and $\chi=(b+1)+e=0$ forces $b=-1$, $e=0$. So $C_{0}$ is a ribbon with $I_{\ell}/I\cong\OO(-1)$;
the surjection $N^{*}_{\ell}\to\OO(-1)$ vanishes on $\OO^{n+1}$, so its kernel $K$ contains $\OO^{n+1}$. Hence $T_{z}C_{0}=T_{z}\ell\oplus\lambda$, where
$\lambda\subset N_{\ell,z}$ is the annihilator of $K_{z}$; since $\OO^{n+1}$ annihilates the positive part $\OO(1)^{n-1}$ of $N_{\ell}$, whose fibre at
$z$ is $\hat T_{[\ell]}\hat X_{z}/T_{z}\ell$, we get $T_{z}C_{0}\subset\hat T_{[\ell]}\hat X_{z}\subset L^{\perp}_{z}$, the last inclusion because $X_{z}$ is
Legendrian (Proposition~\ref{prop:secantgeneral}(1)).
\end{proof}

The two exceptional cases are controlled as follows.

\begin{lemma}[Non-embedded transverse conics]
\label{lem:nonembedded}
In the space of maps modulo $\mathrm{Aut}(\PP^{1})$, the non-jumping transverse immersions $f\colon\PP^{1}\to Z$ of $L$-degree two that are not
injective form a locus of dimension at most $2n+1$.
Consequently, for a general $z$, the directions $[df_{0}(p_{0})]$ arising in case (i) of Lemma~\ref{lem:P2proved} with $f_{0}$ non-jumping form a
set $R_{z}$ whose closure has dimension at most one.
\end{lemma}

\begin{proof}
At such an $f$, $N_{f}=f^{*}TZ/T\PP^{1}\cong f^{*}L^{\perp}\cong\OO(1)^{2n}$, so the space of maps modulo reparametrization is smooth of dimension
$h^{0}(N_{f})=4n$ near $f$. Consider triples $(f,s,s')$ with $s\neq s'\in\PP^{1}$, a manifold of dimension $4n+2$, and
$\Phi(f,s,s')=(f(s),f(s'))\in Z\times Z$. Since a section of $\OO(1)^{2n}$ is determined by its values at two points, $H^{0}(N_{f})\to N_{f,s}\oplus N_{f,s'}$
is an isomorphism, and with the tangent directions $df(s)$, $df(s')$ this shows that $d\Phi$ is an isomorphism. So $\Phi^{-1}(\Delta_{Z})$ has dimension
$2n+1$ near each of its points, and its image in the space of maps, the non-injective ones, has dimension $\le2n+1$ near each of its points,
hence everywhere. The pairs $(f,p)$ with $f$ in this family form a set of
dimension $\le2n+2$; its fibre over a general $z$ under $(f,p)\mapsto f(p)$ has dimension $\le1$, and so does its image under $(f,p)\mapsto[df(p)]$.
\end{proof}

For the jumping ones, only a bound on dimension is available, and it suffices.

\begin{lemma}[Non-injective transverse conics through a general point]
\label{lem:nodalbound}
Let $z$ be general. The non-injective transverse immersions $f\colon\PP^{1}\to Z$ of $L$-degree two with $z\in f(\PP^{1})$ and $j(f):=h^{0}(N_{f}(-2))=j$ form a
family of dimension $\le j+1$. Consequently $\dim S^{(i)}_{z}\le n+1$ (Definition~\ref{def:Siz}), and for $n\ge3$ the image under $\widetilde\mu_{z}$ of the case-(i) points of $\Gamma_{z}$ has
codimension $\ge n-1\ge2$ in $\PP(T_{z}Z)$; in particular a divisor of $\Gamma_{z}$ consisting of such points is contracted.
\end{lemma}

\begin{proof}
Every conic through the general point $z$ is free \cite[II.3.11]{Kol96}; since the maps have bounded degree, the exceptional set in
\cite[II.3.11]{Kol96} is a finite union of proper closed subsets, so ``general'' rather than ``very general'' suffices, for all such maps at once. Freeness is used only
to get $j\le n$. So $f^{*}T_{Z}$ is nef, and for a transverse immersion
$N_{f}=f^{*}T_{Z}/T_{\PP^{1}}\cong f^{*}L^{\perp}$ is nef. The $L$-valued symplectic form on $L^{\perp}$ gives $N_{f}\cong N_{f}^{*}\otimes f^{*}L=N_{f}^{*}(2)$, so the splitting
type is symmetric about $1$ and, being nef, lies in $\{0,1,2\}$: $N_{f}\cong\OO(2)^{j}\oplus\OO(1)^{2n-2j}\oplus\OO^{j}$ with $j=h^{0}(N_{f}(-2))\le n$.
Such an $f$ is birational onto its image (otherwise the image would be a line, which is integral), so it has finitely many double points.

Fix $w\in Z$ and consider the triples $(f,s,s')$, $s\neq s'$, with $f(s)=f(s')=w$, modulo $\mathrm{Aut}(\PP^{1})$. Let $\mathcal N_{j}$ denote the set of
triples $(f,s,s')$ with $s\neq s'$, $f(s)=f(s')$ and $j(f)=j$, fibred over $w=f(s)$ by evaluation; $(f,s,s')\mapsto f$ is finite-to-one. The Zariski tangent space at $(f,s,s')$
of the fibre over $(w,w)$ of $\mathrm{Hom}_{2}(\PP^{1},Z)\times(\PP^{1}\times\PP^{1}\smallsetminus\Delta)\to Z\times Z$ (here $T_{f}\mathrm{Hom}=H^{0}(f^{*}T_{Z})$, so no smoothness
is needed), before dividing by $\mathrm{Aut}$, is $\{(\xi,a,b):\xi(s)+a\,df(s)=0,\ \xi(s')+b\,df(s')=0\}$ with $\xi\in H^{0}(f^{*}T_{Z})$, $a\in T_{s}\PP^{1}$, $b\in T_{s'}\PP^{1}$, i.e.\ the
sections $\xi$ with $\xi(s)\in\im df(s)$ and $\xi(s')\in\im df(s')$. Its image in $H^{0}(N_{f})$ is $H^{0}(N_{f}(-s-s'))=H^{0}(N_{f}(-2))$ (a class with
$\bar\xi(s)=\bar\xi(s')=0$ lifts to $\xi\in H^{0}(f^{*}T_{Z})$ because $H^{1}(T_{\PP^{1}})=0$, and the lift automatically satisfies the two conditions), and
the kernel is $H^{0}(T_{\PP^{1}})\cdot df$, the tangent space of the $\mathrm{Aut}$-orbit. So the nodal maps with node at $w$ form a family of dimension $\le j$
near $f$; the bound is local at $f$ with $j=j(f)$, and holds on the whole stratum $\{j(f)=j\}$ by upper semicontinuity of $j$. Hence $\mathcal N_{j}$ has
dimension $\le(2n+1)+j$. The pointed family
$\{(f,p):f\in\mathcal N_{j},\ p\in\PP^{1}\}$ has dimension $\le2n+2+j$ and maps to $Z$ by $(f,p)\mapsto f(p)$; its fibre over a general $z$ is empty or of
dimension $\le2n+2+j-(2n+1)=j+1$. By Lemma~\ref{lem:P2proved}(i), every point of $\Gamma_{z}$ whose stable limit is $(f,p)$ has image $[df(p)]$, whatever the
number of such points (in the Hilbert model there may be a family of them, from the embedded point at the node, Remark~\ref{rem:nodalharmless}(2)); so
$S^{(i)}_{z}$ has dimension $\le\max_{j}(j+1)\le n+1$. For $n\ge3$ this is $<2n-1$, the dimension of a divisor of $\Gamma_{z}$, so such a divisor is
contracted.
\end{proof}

Lemma~\ref{lem:nonembedded} is essentially the case $j=0$. The two-point evaluation is a local isomorphism only for $j=0$, but the Zariski-tangent-space bound
above does not need it.

\begin{definition}[Transverse branch directions]
\label{def:Siz}
For a general $z$ let $S^{(i)}_{z}\subset\CC^{2n}$ be the set of directions $[df(p)]$, $f(p)=z$, of the non-injective transverse immersions $f$ of
$L$-degree two through $z$, jumping or not; so $R_{z}\subset S^{(i)}_{z}$, and $\dim\overline{S^{(i)}_{z}}\le n+1$ by Lemma~\ref{lem:nodalbound}.
\end{definition}

\begin{remark}[Why irreducible nodal conics are harmless]
\label{rem:nodalharmless}
(1) \emph{Curvature.} LeBrun's structure at $[f]$ uses $T_{[f]}=H^{0}(N_{f})$ with $N_{f}\cong\OO(1)^{2n}$ and the contact form, not the injectivity of $f$.
So at a non-jumping nodal conic, as a point of the stable-map space, $g^{\CC}$ is non-degenerate and $W^{\CC}$ is holomorphic; there is nothing to extend.
At a jumping nodal conic the question is the same as at a jumping smooth conic, and by Lemma~\ref{lem:nodalbound} it never arises along a divisor of
$\Gamma_{z}$ with divisorial image.
(2) \emph{The Hilbert picture.} Over a nodal $\Sigma$ the Hilbert scheme carries an embedded point at the node $w$ with a direction $v\in\PP(T_{w}Z/T_{w}\Sigma)$,
and these points may form a boundary divisor of $\overline F_{z}$; but $\widetilde\mu_{z}$ sends all of them to the branch tangents of $\Sigma$ at $z$, so such a divisor
is contracted (its image has dimension $\le n+1$) and is ignored by Theorem~\ref{thm:Ereduction}.
(3) \emph{Lower bounds.} Every component of the locus of free irreducible nodal conics has dimension $\ge2n+1$ (a node is $2n-1$ conditions on the $4n$-dimensional family), so irreducible nodal conics, when they
exist and sweep out $Z$, pass through a general point in a family of dimension $\ge1$, never finitely many, and their nodes can fill $Z$. (Here, as in this whole remark,
``nodal'' means the image of a non-injective transverse immersion, not a line pair.)
\end{remark}

The boundary points of $\Gamma_{z}$ whose limit direction lies in the contact hyperplane include two further kinds: the ribbons of Lemma~\ref{lem:P2proved} and the images of integral maps that are not embeddings and have no cusp over $z$. Both are contracted,
with no hypothesis.

\begin{lemma}[Ribbons are contracted]
\label{lem:ribbons}
Let $z$ be general and $y\in\Gamma_{z}$ a point whose image in $\overline F_{z}$ is supported on a line $\ell\ni z$. Then this image is a ribbon $C_{0}$, and
$\widetilde\mu_{z}(y)$ lies in the line $\PP(T_{z}C_{0})$, a tangent line of $X_{z}$ at $[\ell]$. The ribbons through $z$ form a family of dimension $\le2n-3$, so the
points of $\Gamma_{z}$ over ribbons have image of dimension $\le2n-2$, and every prime divisor of $\Gamma_{z}$ whose generic point lies over a ribbon, including any
component of $E_{z}$ lying over ribbons, is contracted by $\widetilde\mu_{z}$.
\end{lemma}

\begin{proof}
By the non-reduced case of the proof of Lemma~\ref{lem:P2proved}, $C_{0}$ is a ribbon with $I_{\ell}/I\cong\OO(-1)$ and no embedded points, given by a
surjection $N^{*}_{\ell}=\OO(-1)^{n-1}\oplus\OO^{n+1}\to\OO(-1)$. It kills $\OO^{n+1}$, so it is a nonzero element of $\mathrm{Hom}(\OO(-1)^{n-1},\OO(-1))=\CC^{n-1}$, and an
embedded ribbon in a smooth ambient space is determined by its conormal surjection. So the ribbons on $\ell$ form at most a $\PP^{n-2}$, independent of the point $z$ of
$\ell$, and with the $(n-1)$-dimensional family of lines through $z$ the ribbons through $z$ form a family of dimension $\le2n-3$. Locally
$I=(w_{2},\dots,w_{2n},w_{1}^{2})$, where $dw_{2},\dots,dw_{2n}$ span the kernel $K$ of the surjection, a subbundle; so $T_{x}C_{0}=T_{x}\ell\oplus K^{\perp}_{x}$ is
two-dimensional at every $x\in\ell$, and $K^{\perp}\cong\OO(1)$ lies in the positive part $\OO(1)^{n-1}$ of $N_{\ell}$, whose fibre at $z$ is
$\hat T_{[\ell]}\hat X_{z}/T_{z}\ell$. Hence $\PP(T_{z}C_{0})$ is a tangent line of $X_{z}$ at $[\ell]$, inside $\PP(L^{\perp}_{z})$. Since $F^{\mathrm{ne}}_{z}$ is dense in $\Gamma_{z}$
(Theorem~\ref{thm:irred}), $y$ is the end of an arc from $F^{\mathrm{ne}}_{z}$; along it the Hilbert limit is the image of $y$ under the morphism
$\Gamma_{z}\to\widetilde F_{z}\to\overline F_{z}$, so the flat-limit argument of Lemma~\ref{lem:P2proved} puts $\widetilde\mu_{z}(y)$ in $T_{z}C_{0}$, however singular
$\widetilde F_{z}$ may be over the ribbons. The image of the points over ribbons lies in a union of lines indexed by a family of dimension $\le2n-3$.
\end{proof}

\begin{lemma}[Non-embedded integral conics are contracted]
\label{lem:intnonembedded}
Let $z$ be general and $n\ge2$. Let $D\subset\Gamma_{z}$ be a prime divisor whose generic point lies over a member of $\MCbar$ supported on $f(\PP^{1})$, with $f$
integral, birational onto its image and not an embedding, and with no cusp of $f$ over $z$. Then $\dim\widetilde\mu_{z}(D)\le2n-2$: $D$ is contracted.
\end{lemma}

\begin{proof}
\emph{Setup.} The generic point $y$ of $D$ is the end of an arc from $F^{\mathrm{ne}}_{z}$, whose stable limit $(f_{0},p_{0})$ has as cycle the reduced,
irreducible support $\Sigma$; so its domain is irreducible and $f_{0}$ is the normalization of $\Sigma$. The smooth conics of $\MCbar$ lie in the component $\mathcal M$ of
$\overline{\mathcal M}_{0,0}(Z,2)$ through its balanced transverse conics, irreducible of dimension $4n$; so $f_{0}\in\mathcal M$. On the open set of maps with irreducible domain,
$\theta\circ df$ is a section of the line bundle with fibre $H^{0}(K_{\PP^{1}}\otimes f^{*}L)=\CC$, not identically zero, so every component of the integral locus has
dimension $\le4n-1$. A birational $f$ has no automorphisms, and $T_{f}\mathcal M\subset H^{0}(N_{f})$, $N_{f}=f^{*}T_{Z}/T_{\PP^{1}}$, torsion at cusps included.

\emph{Singular point at $z$.} If $z=f_{0}(p_{0})$ is a singular point of $\Sigma$ (a node, tacnode or multiple point; $df_{0}(p_{0})\neq0$), the limit directions are branch
tangents by Lemma~\ref{lem:P2proved}. The pairs $(f,w)$ with $f$ integral and singular at $w$ form a family of dimension $\le4n-1$, whose fibre over a general $w=z$ has dimension
$\le2n-2$; so the image has dimension $\le2n-2$.

\emph{Smooth point at $z$.} Otherwise $f_{0}^{-1}(z)=\{p_{0}\}$ and $f_{0}$ is immersive there, so near $p_{0}$ the stable family is an embedding, the flat family near
$z$ is smooth, and $\widetilde\mu_{z}(y)=[df_{0}(p_{0})]$ depends only on $(f_{0},p_{0})$. Hence $\widetilde\mu_{z}(D)$ lies in the closure of $\mu_{z}(\mathcal D_{z})$, where
$\mathcal D$ runs over the finitely many families of non-embedded integral maps and $\mathcal D_{z}=\{(f,p):f\in\mathcal D,\ f(p)=z\}$. If $\dim\mathcal D\le4n-2$, then
$\dim\mathcal D_{z}\le2n-2$. So we may assume that $\mathcal D$ is a component of the integral locus, of dimension $4n-1$, with non-embedded generic member $f$, and put
$T=T_{f}\mathcal D$, of dimension $4n-1$. For general $z$ the generic point of each $(2n-1)$-dimensional component of $\mathcal D_{z}$ avoids $\mathrm{Sing}\,\mathcal D$, the points
$p$ over singular points of $f(\PP^{1})$, the cusp points, and the proper closed subset over which the double-point incidence fails to be generically smooth.

\emph{First-order integral deformations.} For $\xi\in T$ put $h=\theta(\xi)\in H^{0}(f^{*}L)=H^{0}(\OO(2))$. The derivative of $\theta\circ df$ at $\xi$ is the
section $d\theta(\tilde\xi,df)+d(\theta(\tilde\xi))$ of $K_{\PP^{1}}\otimes f^{*}L\cong\OO$, a constant, which vanishes. At immersed points this says
$\xi\equiv X_{h}\bmod TC^{\perp}$, i.e.\ the image of $\xi$ in $TZ|_{C}/TC^{\perp}\cong J^{1}(L|_{C})$ is $j^{1}h$ (Lemma~\ref{lem:intnormal}(1)); so if $h=0$, then $\xi$
is a section of $W=TC^{\perp}/TC$, of rank $2n-2$, near every immersed point. At a cusp $s_{0}$, where $df(s_{0})=0$, the constant equals $h'(s_{0})$, so
$h'(s_{0})=0$.

\emph{The singularity costs one condition on $h$.} At a double point $f(s)=f(s')=w$, generic smoothness of the double-point incidence over $\mathcal D$ shows that
every $\xi\in T$ satisfies $\tilde\xi(s)+a\,f'(s)=\tilde\xi(s')+b\,f'(s')$ in $T_{w}Z$ for some $a,b$. Since $f'(s),f'(s')\in L^{\perp}_{w}$, applying $\theta_{w}$ gives
$h(s)=h(s')$ in $L_{w}$, a nonzero linear condition on $H^{0}(\OO(2))$. At a cusp the condition is $h'(s_{0})=0$. Either way $\dim h(T)\le2$, so
$T_{W}:=T\cap\ker h$ has dimension $\ge4n-3$.

\emph{Contraction.} At a generic point $(f,p)$ of $\mathcal D_{z}$, with $p$ immersed, $\mathcal D_{z}$ is smooth with tangent space
$T(-p)=\{\xi\in T:\xi(p)=0\}$. The jets $J^{1}_{p}W$ form a space of dimension $4n-4<4n-3$, so some nonzero $\xi\in T_{W}$ vanishes to second order at $p$. It lies in $T(-p)$,
and $d\mu_{z}(\xi)=\xi'(p)=0$ (Lemma~\ref{lem:dmu}). (It may be the torsion section $f'/(s-s_{0})$ at a cusp $s_{0}\neq p$; it is still a nonzero tangent vector killed by
$d\mu_{z}$.) So the generic rank of $d\mu_{z}$ on the reduced $\mathcal D_{z}$ is $<2n-1$, and $\dim\mu_{z}(\mathcal D_{z})\le2n-2$.
\end{proof}

\begin{remark}[Checks]
\label{rem:intnonembeddedchecks}
For an embedded smooth integral conic there is no condition on $h$: $h(T)$ is three-dimensional, $T_{W}=H^{0}(W)$, and when $W\cong\OO(1)^{2n-2}$ the map
$j^{1}_{p}$ is an isomorphism, so no kernel is forced; this is Lemma~\ref{lem:intcontract}(1). For $n=1$, $W=0$ and the count shows that double points of Legendrian conics have
codimension $\ge1$, as in Legendrian deformation theory. The proof also shows that such a $\mathcal D$ can exist only if, for immersed $f$, $h^{0}(W)\ge4n-3$, that is, only if $W$ contains
$\OO(a)\oplus\OO(2-a)$ with $a\ge4$. Non-birational integral maps, the double covers of lines, give ribbons (Lemma~\ref{lem:ribbons}).
\end{remark}

The remaining case, a cusp at $z$, does not occur for a general $z$: this is the following condition, which holds on every projective contact manifold
(Theorem~\ref{thm:cuspgeneral}).

\begin{itemize}[leftmargin=2.5em]
\item[\textbf{(K)}] \emph{For a general $z\in Z$, no member of $\MCbar$ has as support (reduced underlying curve) the image $f(\PP^{1})$ of a map $f$ that
is integral, birational onto its image, with $df(p)=0$ at a point $p$ over $z$.}
\end{itemize}

\begin{definition}[The transverse boundary directions]
\label{def:Sz}
For a general $z$ let
\[
S_{z}:=\widetilde\mu_{z}\big(\Gamma_{z}\smallsetminus F^{\circ}_{z}\big)\cap\CC^{2n}\subset\PP(T_{z}Z),
\]
the limiting tangent directions at $z$ of boundary points of $\Gamma_{z}$ that are not in the contact hyperplane. By Lemma~\ref{lem:P2proved},
$S_{z}\subseteq S^{(i)}_{z}\cup S^{(ii)}_{z}$, where $S^{(i)}_{z}$ consists of the branch tangents $[df(p)]$, $f(p)=z$, of the non-injective transverse immersions
through $z$, and $S^{(ii)}_{z}$ of the transverse limit directions at cuspidal integral conics with cusp at $z$. The set $R_{z}$ of Lemma~\ref{lem:nonembedded}
is the non-jumping part of $S^{(i)}_{z}$ (Definition~\ref{def:Siz}); condition (K) (Theorem~\ref{thm:cuspgeneral}) gives $S_{z}\subseteq S^{(i)}_{z}$.
\end{definition}

\emph{Why (K) is expected, and why it holds.} A cusp $df(p)=0$ is $2n+1$ conditions on the $(4n+1)$-dimensional space of pointed maps, so the cuspidal stratum of $\MCbar$
has codimension at most $2n$ wherever it is non-empty; if it has this maximal codimension, the cusps sweep out a subset of $Z$ of dimension $\le2n$, and (K)
holds. This is the expected behaviour, but with no margin; Theorem~\ref{thm:cuspgeneral} proves it on every projective contact
manifold, {by an involution argument that uses only the Levi form and the bound on the $L$-degree}. Every cuspidal conic is integral, since $\theta\circ df$ is a constant (Lemma~\ref{lem:allornothing}) that vanishes at the cusp; and on a
contact Fano manifold on which $L$ is very ample there are no irreducible singular conics at all, an irreducible curve of degree two in projective space being a
smooth plane conic.

\begin{proposition}[The boundary at infinity]
\label{prop:boundaryinfinity}
Let $n\ge3$ and let $z$ be general. Then
\[
\widetilde\mu_{z}\big(\Gamma_{z}\smallsetminus F^{\circ}_{z}\big)\subset\PP(L^{\perp}_{z})\cup\overline{S^{(i)}_{z}},\qquad\dim\overline{S^{(i)}_{z}}\le n+1\le2n-2 .
\]
Consequently, for every closed $S\subset\PP(T_{z}Z)$ containing
$\overline{S^{(i)}_{z}}$, the map $F^{\circ}_{z}\smallsetminus\widetilde\mu_{z}^{-1}(S)\to\CC^{2n}\smallsetminus S$
is proper.
\end{proposition}

\begin{proof}
By (K) (Theorem~\ref{thm:cuspgeneral}), case (ii) of Lemma~\ref{lem:P2proved} does not occur for general $z$, and case (i) lands in $S^{(i)}_{z}$, of dimension $\le n+1$ by
Lemma~\ref{lem:nodalbound}. Every point of $\Gamma_{z}\smallsetminus F^{\circ}_{z}$ is the end of an arc as above. For the last assertion,
$\widetilde\mu_{z}$ is proper and the preimage of $\CC^{2n}\smallsetminus S$ is contained in $F^{\circ}_{z}$.
\end{proof}

Theorem~\ref{thm:LSJ} and Theorem~Q therefore carry no hypothesis on cusps.

\part{Local Theory: The Neighbourhood of a Conic}
\label{part:local}

In this part the questions are local along a conic: how the contact neighbourhood of a conic $C\subset Z$ determines the germ of the transform at $[C]$, and how the
transform degenerates when $C$ does. Section~\ref{sec:formal} develops contact neighbourhoods and proves Theorem~\ref{thm:thirdorder}, that the curvature of $g^{\CC}$ at a balanced transverse conic is the
invariant of its neighbourhood of order three. The following sections treat the degenerate conics in turn: the nodal conics, across whose exceptional divisor the curvature
extends regularly (\S\ref{sec:dilation}, Theorem~\ref{thm:Eregular}); the integral conics, along which the conformal structure goes to infinity while the curvature extends, and vanishes on the balanced ones
(\S\ref{sec:integralII}); and the cusps, which do not occur at a general point
(\S\ref{sec:cuspsII}). Section~\ref{sec:dictionaryII} collects the results. Where the arguments use a general point $z$ and the space $\Gamma_{z}$, $Z$ is compact.

\section{Contact neighbourhoods of conics}
\label{sec:formal}

This section is the local half of the transform: it compares the contact neighbourhood of a conic in $Z$ with the conformal geometry of $\MC$ at the corresponding
point. It develops normal forms for a contact form near a curve (\S\ref{sec:normalform}), the graded automorphism sheaves of the formal neighbourhood of a curve
(\S\ref{sec:autsheaves}), the rigidity of the neighbourhood of a balanced transverse conic and the identification of its invariant of order three with the curvature
(\S\ref{sec:thirdorder}, Theorem~\ref{thm:thirdorder}), and the cohomology of the graded automorphism sheaves on a balanced nodal conic
(\S\ref{sec:nodalgroups}). The normal forms, the automorphism sheaves and the nodal groups are the tools of the dilation argument of \S\ref{sec:dilation}; the
nonvanishing of $\rho_{0}$ (Theorem~\ref{thm:thirdorder}(4)) is used in \S\ref{sec:mustdegenerate} and \S\ref{sec:homog}; and Theorem~\ref{thm:thirdorder} is the first entry of the
dictionary of \S\ref{sec:dictionaryII}, where the behaviour of the neighbourhoods at the degenerate conics is collected. Here $Z$ is any contact manifold of dimension
$2n+1$.

\subsection{Normal forms}
\label{sec:normalform}

\begin{proposition}[Local normal forms]
\label{prop:normalform}
Let $\theta$ be a holomorphic contact form near $p\in Z$.
\begin{enumerate}
\item If $A$, $A'$ are smooth curve germs through $p$ tangent to $\ker\theta$ with $d\theta(T_{p}A,T_{p}A')\neq0$, there are holomorphic
coordinates $(a,Q,b,c,P)$, $b,c\in\CC^{n-1}$, centred at $p$, with
\[
\theta=dP+Q\,da+\textstyle\sum_{i<n}c_{i}\,db_{i},\qquad A=\{Q=b=c=P=0\},\quad A'=\{a=b=c=P=0\}.
\]
\item A single smooth curve germ tangent to $\ker\theta$ is the $a$-axis in coordinates of the same form.
\item If $C$ is a smooth curve germ transverse to $\ker\theta$, there are a unit $f$ and coordinates with $f\theta=dz-\sum_{i}p_{i}\,dq_{i}$ and
$C$ the $z$-axis.
\end{enumerate}
\end{proposition}

\begin{proof}
By Darboux's theorem $\theta=dz+\lambda_{0}$ with $\lambda_{0}$ a $1$-form on $(\CC^{2n},0)$, $\omega_{0}=d\lambda_{0}$ symplectic. A curve tangent to
$\ker\theta$ is transverse to the Reeb field $\partial_{z}$, so it projects immersively and is the horizontal lift of its projection. If $\phi$ is a
symplectomorphism germ, $\phi^{*}\lambda_{0}-\lambda_{0}=dF$ with $F(0)=0$, and $\Phi(m,z)=(\phi(m),z-F(m))$ preserves $\theta$ and carries horizontal lifts
to horizontal lifts (likewise between two primitives $\lambda,\lambda'$). So it suffices to move the projections by a symplectomorphism.

(1) The projections $\gamma_{A},\gamma_{A'}$ satisfy $\omega_{0}(T_{0}\gamma_{A},T_{0}\gamma_{A'})\neq0$. They lie on a smooth surface germ $\Sigma$ with
$T_{0}\Sigma=V:=T_{0}\gamma_{A}+T_{0}\gamma_{A'}$: write $\CC^{2n}=V\oplus V^{\omega}$; the projections of the two curves to $V$ are transverse smooth germs,
so after a holomorphic change of coordinates $(x,y)$ on $V$ they are the two axes, and the curves are graphs $x\mapsto g_{A}(x)$, $y\mapsto g_{A'}(y)$ over them,
with values in $V^{\omega}$; take $\Sigma$ the graph of $g_{A}(x)+g_{A'}(y)$, which contains both since $g_{A}(0)=g_{A'}(0)=0$. $\Sigma$ is symplectic, so by the relative Darboux theorem (Moser's proof
applies to holomorphic germs) there are Darboux coordinates $(a,Q,b,c)$ with $\Sigma=\{b=c=0\}$. On $\Sigma$ choose coordinates $(u,v)$ with
$\gamma_{A}=\{v=0\}$, $\gamma_{A'}=\{u=0\}$, $\omega_{0}|_{\Sigma}=f\,du\wedge dv$, $f(0)\neq0$; then $a'=-u$, $Q'=\int_{0}^{v}f(u,\tau)\,d\tau$ satisfy
$dQ'\wedge da'=\omega_{0}|_{\Sigma}$, and replacing $(a,Q)$ by $(a',Q')$ keeps Darboux coordinates. Lift with $\lambda'=Q\,da+\sum c_{i}\,db_{i}$ and put
$P=z$; as $\lambda'$ vanishes on both axes, their horizontal lifts through $p$ are the axes. (2) follows from (1) with an auxiliary $A'$. (3)
Choose a section $h$ of $L$ whose contact field $X_{h}$ is tangent to $C$ and nonvanishing (possible since $X_{h}|_{C}$ is determined by
$h|_{C}=\theta(v)$ and by derivatives of $h$ in the directions of $\ker\theta$, which may be prescribed). Then $X_{h}$ is the Reeb field of
$\theta/h$; flow-box coordinates give $\theta/h=dz+\lambda$ with $C=\{m=0\}$, and Darboux on the quotient finishes.
\end{proof}

\subsection{Automorphism sheaves}
\label{sec:autsheaves}

For a local section $h$ of $L$ let $X_{h}$ be the contact vector field with $\theta(X_{h})=h$. Let $C\subset Z$ be a curve with ideal $I$, a local
complete intersection, and let $\mathcal A_{j}$ be the sheaf on $C$ of germs of Hamiltonians $h$ with $X_{h}\in I^{j}T_{Z}$; the flows of such
fields are the contact automorphisms of the formal neighbourhood that are the identity on $C^{(j-1)}=(C,\OO_{Z}/I^{j})$. Put
\[
\mathcal G_{j}=\mathcal A_{j}/\mathcal A_{j+1},\qquad S_{j}=I^{j+1}L/I^{j+2}L=L\otimes S^{j+1}N^{*},
\]
the graded automorphism sheaf and the \emph{naive} sheaf.

\begin{lemma}[The basic filtration]
\label{lem:basicfilt}
For $j\ge1$, $\mathcal A_{j+1}\subset I^{j+1}L\subset\mathcal A_{j}$, and each $\mathcal A_{j}$ is an $\OO_{Z}$-module. With
$\mathcal P_{j}=\mathcal A_{j}/I^{j+1}L$, an $\OO_{C}$-submodule of $S_{j-1}$, there is an exact sequence
\[
0\to S_{j}/\mathcal P_{j+1}\to\mathcal G_{j}\to\mathcal P_{j}\to0 .
\]
\end{lemma}

\begin{proof}
$h=\theta(X_{h})$, so $X_{h}\in I^{j}T$ forces $h\in I^{j}L$; $X_{h}$ is a first-order operator in $h$, so $h\in I^{j+1}L$ gives $X_{h}\in I^{j}T$;
$X_{fh}=fX_{h}+h\,\sigma(df)$ with $\sigma$ the symbol gives $f\mathcal A_{j}\subset\mathcal A_{j}$; $\mathcal P_{j}$ is killed by $I$; and the
sequence is $0\to I^{j+1}L/\mathcal A_{j+1}\to\mathcal A_{j}/\mathcal A_{j+1}\to\mathcal A_{j}/I^{j+1}L\to0$.
\end{proof}

\begin{proposition}[Transverse and integral curves]
\label{prop:contactauto}
\begin{enumerate}
\item If $C$ is transverse to $L^{\perp}$, then $X_{h}\in I^{j}T_{Z}\iff h\in I^{j+1}L$; so $\mathcal P_{j}=0$ and $\mathcal G_{j}=S_{j}$. The contact
neighbourhoods of order $k$ inducing a given one of order $k-1$ form a torsor under $H^{1}(C,L\otimes S^{k+1}N^{*})$, for $k\ge2$.
\item If $C$ is integral (tangent to the contact distribution), in coordinates $\theta=dz-\sum p_{i}dq_{i}$ with $C$ the $q_{1}$-axis, $X_{h}\in I^{j}T_{Z}\iff h\in c(q_{1})z^{j}+I^{j+1}L$;
so $\mathcal P_{j}\cong L^{1-j}$, spanned by $z^{j}$, where $L^{-1}\subset N^{*}$ is the conormal line of $\theta$.
\end{enumerate}
\end{proposition}

\begin{proof}
(1) In the coordinates of Proposition~\ref{prop:normalform}(3), $I=(q,p)$ and $X_{h}=\sum(-h_{p_{i}}\partial_{q_{i}}+(h_{q_{i}}+p_{i}h_{z})\partial_{p_{i}})
+(h-\sum p_{i}h_{p_{i}})\partial_{z}$. Expanding $h$ in homogeneous parts in $(q,p)$, the conditions kill every part of degree $\le j$. The torsor
statement is Griffiths' gluing argument: two neighbourhoods inducing the same one of order $k-1$ are locally isomorphic over it, and differences
of local isomorphisms form a cocycle in $\mathcal G_{k}$, abelian for $k\ge2$. (2) The same computation with $I=(z,p,q_{2},\dots,q_{n})$ leaves
exactly $c(q_{1})z^{j}$ in degree $j$.
\end{proof}

\subsection{The curvature is the third-order invariant of the neighbourhood}
\label{sec:thirdorder}

Let $C$ be a transverse conic and $I$ its ideal. The \emph{contact neighbourhood of order $k$} of $C$ is the formal neighbourhood $(C,\OO_{Z}/I^{k+1})$ with the
contact structure, up to contact isomorphisms defined modulo the automorphisms that induce the identity on it; in terms of \S\ref{sec:autsheaves}, a germ of
$Z$ along $C$ is glued from Darboux charts (Proposition~\ref{prop:normalform}(3)) by contact transformations, and its contact neighbourhood of order $k$ is the
class of the gluing modulo $\exp(\mathcal A_{k+1})$, since $X_{h}\in I^{k+1}T_{Z}$ for $h\in\mathcal A_{k+1}$. For a balanced transverse conic put $E:=E_{\CC}|_{[C]}$, so that
$B|_{C}\cong E^{*}\otimes\OO_{C}$, $N\cong E^{*}\otimes L^{1/2}|_{C}$ and
\begin{equation}
L\otimes S^{k+1}N^{*}\ \cong\ S^{k+1}E\otimes L^{(1-k)/2}|_{C}\ \cong\ S^{k+1}E\otimes\OO(1-k).\label{eq:Gk}
\end{equation}
Since $\theta$ restricts to an isomorphism $TC\to L|_{C}$, $L^{-1}|_{C}\cong K_{C}$ canonically, and $H^{1}(C,K_{C})=\CC$; so
\begin{equation}
H^{1}(C,L\otimes S^{4}N^{*})\ \cong\ S^{4}E\otimes H^{1}(C,K_{C})\ =\ S^{4}E\label{eq:H1S4}
\end{equation}
canonically, the space in which $R_{E}^{\CC}([C])$ lies. We compare with a line $\ell=\PP(W)\subset\PP^{2n+1}=\PP(V)$, $\Omega|_{W}\ne0$ (\S\ref{sec:contactZ}). The
subgroup $\Sp(W^{\perp})\cong\Sp(n,\CC)$ of $\Sp(V)$ fixes $\ell$, acts on its neighbourhood by contact transformations, and acts on $E^{*}_{\CC}|_{[\ell]}=H^{0}(\ell,B|_{\ell})\cong
W^{\perp}$ by the standard representation.

\begin{theorem}[The curvature is the third-order invariant]
\label{thm:thirdorder}
Let $(Z,\theta)$ be a complex contact manifold of dimension $2n+1$, $n\ge1$, and $C\subset Z$ a balanced transverse conic.
\begin{enumerate}
\item The contact neighbourhood of order two of $C$ is isomorphic to that of $\ell$.
\item The contact neighbourhoods of order three inducing a given one of order two form a torsor under $H^{1}(C,L\otimes S^{4}N^{*})\cong S^{4}E$.
\item Let $\iota$ be an isomorphism of the contact neighbourhoods of order two of $C$ and of $\ell$, and let $\kappa_{\iota}(C)\in S^{4}E$ be the difference between the
contact neighbourhood of order three of $C$ and the one induced by $\iota$ from $\ell$. Then
\[
R_{E}^{\CC}([C])\ =\ c_{n}\,\kappa_{\iota}(C)
\]
for a constant $c_{n}\ne0$ depending only on $n$, on the normalization of $g^{\CC}$ and on that of $H^{1}(C,K_{C})\cong\CC$ (by the trace). In particular $\kappa_{\iota}(C)$ does not depend on $\iota$, $R_{E}^{\CC}([C])$
depends only on the contact neighbourhood of order three of $C$, and $R_{E}^{\CC}([C])=0$ if and only if that neighbourhood is isomorphic to the one of a line of
$\PP^{2n+1}$ transverse to the contact structure.
\item For $g^{\CC}=G/\sigma_{\theta}$ (Proposition~\ref{prop:levwronsk}(2)), the constant $\rho_{0}$ of Theorem~\ref{thm:lebrun}(3) is the same for every $Z$: the value for
the complexification of $\HH P^{n}$. In particular $\rho_{0}\ne0$.
\end{enumerate}
\end{theorem}

\begin{proof}
(1) In Darboux charts along $C$ (Proposition~\ref{prop:normalform}(3)) the gluing of the neighbourhood of order one consists of the identification of $TC$ with $L|_{C}$
and of linear symplectic identifications of the normal spaces, that is, of the symplectic normal bundle $(N,\omega)$; the extension class of the first-order neighbourhood
lies in $H^{1}(C,TC\otimes N^{*})=H^{1}(\OO(1))^{2n}=0$. For $C$ and for $\ell$ these data are $\PP^{1}$, $\OO(2)$ and $E^{*}\otimes\OO(1)$ with a constant symplectic form
on $E^{*}$, so the neighbourhoods of order one agree. By Proposition~\ref{prop:contactauto}(1) the extensions of order two of a given contact neighbourhood of order one
are classified by $H^{1}(C,\mathcal G_{2})$, $\mathcal G_{k}=L\otimes S^{k+1}N^{*}$, and by \eqref{eq:Gk} $H^{1}(C,\mathcal G_{2})=S^{3}E\otimes H^{1}(\OO(-1))=0$.
(2) is Proposition~\ref{prop:contactauto}(1) with $k=3$ and \eqref{eq:H1S4}.

(3) \emph{Step 1: presentation.} Cover $C$ by two discs $C_{0},C_{\infty}$ and choose Darboux charts $U_{0},U_{\infty}$ along them (Proposition~\ref{prop:normalform}(3)), so
that the germ of $Z$ along $C$ is glued by a contact transformation $\Phi$ of the overlap $U_{0\infty}$. Using $\iota$ we take the same charts for $\ell$, glued by
$\Phi_{\ell}$, and may assume that $\Phi_{\ell}^{-1}\circ\Phi$ induces the identity on the neighbourhood of order two; modulo $\exp(\mathcal A_{4})$ it is then $\exp(X_{h})$ with
$h\in\Gamma(U_{0\infty},I^{4}L)=\Gamma(U_{0\infty},\mathcal A_{3})$ (Proposition~\ref{prop:contactauto}(1)), and the class of $h_{4}:=h\bmod I^{5}$ in
$H^{1}(C,\mathcal G_{3})$ is $\kappa_{\iota}(C)$.

\emph{Step 2: parametrization.} Parametrize the conics near $C$ by their first-order part $s\in H^{0}(N)$, $|s|=\varepsilon$, written in the two charts as graphs over
$C_{0}$, $C_{\infty}$; the members are the solutions of the gluing equation, which is solvable to all orders since $H^{1}(N)=0$, and we normalize the terms of order
$\ge2$ in $\varepsilon$ to lie in a fixed complement of $H^{0}(N)$ in the \v{C}ech complex, the same for $Z$ and for $\ell$. The metric at a conic is computed from its
first-order deformations, transported between the charts by the differential of the gluing, and from $\theta$ and $d\theta$ along it (Proposition~\ref{prop:levwronsk}).
A different normalization changes the parametrization by a diffeomorphism tangent to the identity, which does not affect the curvature at $[C]$.

\emph{Step 3: linearity.} Compare $Z$ with $\ell$ in the common charts. Since $X_{h}\in I^{3}T_{Z}$, a conic of $Z$ with first-order part $s$ differs from the conic of
$\PP^{2n+1}$ with the same first-order part by a displacement of order $\varepsilon^{3}$, whose leading term is the normalized solution of the linearized gluing equation
with right-hand side the cubic term $X_{h_{4}}(s,s,s)$; it is linear in $h_{4}$ and depends otherwise only on the neighbourhood of order one. The first-order
deformations of the conic with first-order part $s$, from which the metric is computed, are the derivatives in $s$ of this family, so they differ from those for
$\ell$ by $O(\varepsilon^{2})$, linearly in $h_{4}$ at leading order. The transition of normal
vectors between the charts changes by $d\exp(X_{h_{4}})=\mathrm{id}+O(\varepsilon^{2})$, again linearly in $h_{4}$ at leading order, and $\exp(X_{h})$ rescales $\theta$ by
a factor $1+O(\varepsilon^{4})$, since $\mathcal L_{X_{h}}\theta$ is $\theta$ times the Reeb derivative of $h\in I^{4}L$. Terms quadratic in $h_{4}$ (from composing two cubic fields) displace the conic at order $\varepsilon^{5}$ at the earliest, so they change the metric at
order $\ge\varepsilon^{4}$; likewise $h_{5}$ changes it at order $\ge\varepsilon^{3}$. So through order two only the term linear in $h_{4}$ enters. Everything else that enters the formula of Proposition~\ref{prop:levwronsk}
(the charts, $\theta$ and $d\theta$ in them, and the normalization of Step~2) is common to $Z$ and $\ell$, so there is no other term of order $\le2$. Hence $g^{\CC}_{Z}-g^{\CC}_{\ell}=O(\varepsilon^{2})$, with leading term quadratic in $s$,
linear in $h_{4}$, and depending otherwise only on the neighbourhood of order two; in particular the curvature at $[C]$ depends only on the neighbourhood of order
three. The first derivatives of $g^{\CC}_{Z}-g^{\CC}_{\ell}$ vanish at $[C]$, so the Christoffel symbols agree there and the difference of the curvatures at $[C]$ is the
second derivative of the leading term:
\[
\mathrm{Rm}_{Z}([C])-\mathrm{Rm}_{\ell}([\ell])\ =\ \Lambda(\kappa_{\iota}(C))
\]
for a linear map $\Lambda$ from $S^{4}E$ to the space of algebraic curvature tensors, which depends only on the neighbourhood of order two and is therefore the same for
every $Z$. (A coboundary $h$ changes the charts, not the germ, and does not change the curvature, so $\Lambda$ is defined on the class.) Every class is realized by a
germ: regluing the germ of $\PP^{2n+1}$ along $\ell$ by $\Phi_{\ell}\circ\exp(X_{h_{4}})$, with $h_{4}$ a cocycle of quartic polynomials in the Darboux coordinates, gives a
contact germ along a balanced transverse conic with class $[h_{4}]$. (The curvature quartic itself is written explicitly in a pair of charts, as a universal rational
function of the $4$-jet of the family of conics through two points, in Lemma~\ref{lem:findet}(1); we do not need the value of $c_{n}$.)

\emph{Step 4: Schur's lemma.} The group $\Sp(W^{\perp})$ acts on the neighbourhoods of $\ell$, fixing the neighbourhood of order three of $\ell$, and $\Lambda$ is equivariant.
By Theorem~\ref{thm:lebrun}(3) the curvature of $g^{\CC}$ at a balanced transverse conic of any contact germ is $\rho_{0}R_{0}+R_{E}^{\CC}$ with $R_{E}^{\CC}\in S^{4}E$, so $\Lambda$
takes values in $\CC R_{0}\oplus S^{4}E$, where $R_{0}$ is $\Sp(W^{\perp})$-invariant, and $S^{4}E$ is an irreducible representation of $\Sp(W^{\perp})$ without invariants. So
$\Lambda$ maps into $S^{4}E$ and is a multiple $c_{n}$ of the identity, and the scalar part of the curvature is not changed. Since $R_{E}^{\CC}([\ell])=0$, this gives $R_{E}^{\CC}([C])=c_{n}\kappa_{\iota}(C)$ and (4).

\emph{Step 5: $c_{n}\ne0$.} Otherwise $R_{E}^{\CC}$ would vanish for every contact manifold; but on the twistor space of a Wolf space other than $\HH P^{n}$ it restricts to the
quaternionic Weyl curvature of $M$, which is not zero. The remaining assertions of (3) follow, $R_{E}^{\CC}([C])$ being intrinsic. (Independently of (3), $\kappa_{\iota}$ does not depend on $\iota$: the Lie algebra of the automorphism group of the neighbourhood of order two of
$\ell$ is $H^{0}(TC)\oplus H^{0}(\mathcal G_{1})\oplus H^{0}(\mathcal G_{2})=\mathfrak{sl}(W)\oplus\mathfrak{sp}(W^{\perp})\oplus0$, realized by $\Sp(W)\times\Sp(W^{\perp})$, which
preserves the neighbourhood of order three of $\ell$.)
\end{proof}

\begin{corollary}[The curvature quartic as a neighbourhood invariant]
\label{cor:thirdorderWz}
For $z\in C$, the curvature quartic of \S\ref{sec:Wz} is $W^{\CC}([C],z)=c_{n}\,(\mathrm{ev}_{z}^{-1})^{*}\kappa_{\iota}(C)$, the value at $z$ of the class of the neighbourhood of
order three. In particular $R_{E}^{\CC}$ is parallel if and only if, for general $z$, these values are locally constant on the conics through $z$
(Proposition~\ref{prop:leaf}).
\end{corollary}

\begin{remark}[Higher orders]
\label{rem:higherorders}
By \eqref{eq:Gk}, for $k\ge3$ the neighbourhoods of order $k$ inducing a given one of order $k-1$ form a torsor under
$S^{k+1}E\otimes H^{1}(C,\OO(1-k))\cong S^{k+1}E\otimes S^{k-3}H$, with $H=H_{\CC}|_{[C]}$. For $k=4$ this is $S^{5}E\otimes H$, a summand of
$S^{4}E\otimes E\otimes H$, the space of $\nabla R_{E}^{\CC}$ at $[C]$. This suggests that the neighbourhood of order $k$ encodes the $(k-3)$-jet of the curvature. There is no canonical reference
neighbourhood of order four when $R_{E}^{\CC}([C])\ne0$; measured against a locally symmetric model with the same curvature at $[C]$, one expects the class of order four to
be the relevant part of $\nabla R_{E}^{\CC}([C])$. For $n=1$ and to first order, Theorem~\ref{thm:thirdorder}(3) is
the contour-integral description of the linearized self-dual curvature by fourth derivatives of a Hamiltonian in linearized twistor theory.
\end{remark}

\subsection{The nodal groups}
\label{sec:nodalgroups}

Let $C_{0}=A\cup A'$ be a nodal conic with node $p$ whose branches are standard lines, with $N=N_{C_{0}/Z}$ balanced; then
$d\theta(T_{p}A,T_{p}A')\neq0$ (Proposition~\ref{prop:balcrit}). Let $E^{A}_{+}\cong\OO(1)^{n}\subset N|_{A}$ be the positive part (the maximal
destabilizing subbundle), and $E^{A'}_{+}$ likewise. Since $A'$ is integral, $\theta$ defines $\theta_{A}\in H^{0}(A,N^{*}|_{A}\otimes L)$ (as in the proof
of Proposition~\ref{prop:balcrit}), and likewise $\theta_{A'}$.

\begin{proposition}[The naive groups]
\label{prop:formalnbhd}
$L\otimes S^{k}N^{*}=\bigoplus_{a+b=k}\OO_{C_{0}}(1-a,1-b)^{m_{a}m_{b}}$ with $m_{a}=\dim S^{a}\CC^{n}$, and
\[
h^{1}\big(C_{0},L\otimes S^{k}N^{*}\big)=\begin{cases}0,&k\le2,\\ 2\binom{n+2}{3},&k=3,\\
4\binom{n+3}{4}+2n\binom{n+2}{3}+\binom{n+1}{2}^{2},&k=4.\end{cases}
\]
\end{proposition}

\begin{proof}
$N\cong\OO_{C_{0}}(1,0)^{n}\oplus\OO_{C_{0}}(0,1)^{n}$ and $L|_{C_{0}}=\OO_{C_{0}}(1,1)$. On the nodal curve, $h^{1}(\OO_{C_{0}}(p,q))$ is $0$ for
$p,q\ge0$, $-q-1$ for $p\ge0>q$, and $-(p+q+1)$ for $p,q<0$.
\end{proof}

\begin{proposition}[Cross-annihilation]
\label{prop:Hholds}
$\theta_{A'}(p)$ annihilates $E^{A}_{+}(p)$ and $\theta_{A}(p)$ annihilates $E^{A'}_{+}(p)$.
\end{proposition}

\begin{proof}
This is Step~3 of the proof of Proposition~\ref{prop:balcrit}, with $V_{1}=E^{A}_{+}(p)$, $V_{2}=E^{A'}_{+}(p)$, $\theta_{1}=\theta_{A}$, $\theta_{2}=\theta_{A'}$.
\end{proof}

\begin{theorem}[The nodal groups are the naive ones]
\label{thm:naive}
For $j\ge1$:
\begin{enumerate}
\item $\mathcal P_{j}\cong\iota_{A*}L^{1-j}(-p)\oplus\iota_{A'*}L^{1-j}(-p)\cong\OO_{A}(-j)\oplus\OO_{A'}(-j)$;
\item $h^{1}(C_{0},\mathcal G_{j})=h^{1}(C_{0},S_{j})$ and $h^{0}(C_{0},\mathcal G_{j})=h^{0}(C_{0},S_{j})+2$.
\end{enumerate}
\end{theorem}

\begin{proof}
(1) Use the coordinates of Proposition~\ref{prop:normalform}(1), in which $I=(aQ,b,c,P)$. Away from $p$, Proposition~\ref{prop:contactauto}(2) gives
$\mathcal P_{j}|_{A\smallsetminus p}=L\otimes(\theta_{A}L^{-1})^{j}\cong L^{1-j}$. At $p$, a section of $\mathcal P_{j}$ restricts to $\varphi(a)P^{j}$ on $A$
and to $\psi(Q)(P+aQ)^{j}$ on $A'$ (on $A'$ the $\theta$-coordinate is $P+aQ$); these agree at $p$, and $P(p)$, $(P+aQ)(p)$ are independent in
$I/I^{2}\otimes k(p)$, so $\varphi(0)=\psi(0)=0$. Conversely $aP^{j}$ and $Q(P+aQ)^{j}$ lie in $\mathcal A_{j}$: for $u=P+aQ$ and $h=Qu^{j}$,
$h_{a}-Qh_{P}=0$, $h_{Q}=u^{j}+jaQu^{j-1}$, $h-Qh_{Q}-\sum c_{i}h_{c_{i}}=-jaQ\cdot Qu^{j-1}$, $c_{i}h_{P}=jc_{i}Qu^{j-1}$, all in $I^{j}$, while
$h_{b_{i}}=h_{c_{i}}=0$; $aP^{j}$ is similar. As $\mathcal A_{j}$ is an $\OO_{Z}$-module, this gives (1).

(2) $H^{0}(\mathcal P_{j})=0$ and $h^{1}(\mathcal P_{j})=2(j-1)$. Put $T=S_{j}/\mathcal P_{j+1}$. Lemma~\ref{lem:basicfilt} gives
$h^{1}(\mathcal G_{j})=h^{1}(T)+2(j-1)$, and $0\to\mathcal P_{j+1}\to S_{j}\to T\to0$ gives $h^{1}(S_{j})=h^{1}(T)+2j-r_{j}$ with
$r_{j}=\dim\ker(H^{1}(\mathcal P_{j+1})\to H^{1}(S_{j}))$. We show $r_{j}=2$. Restricting to $A$: the image of $\mathcal P^{A'}_{j+1}$ in $S_{j}|_{A}$ is
zero, and $\mathcal P^{A}_{j+1}=\Lambda^{A}_{j}(-p)\subset\Lambda^{A}_{j}\subset S_{j}|_{A}$ with $\Lambda^{A}_{j}=L\otimes(\theta_{A}L^{-1})^{j+1}$. The kernel of
$H^{1}(\Lambda^{A}_{j}(-p))\to H^{1}(\Lambda^{A}_{j})$ is spanned by $\epsilon_{A}=\delta(1)$, $1\in H^{0}(\CC_{p})$, and $H^{1}(\Lambda^{A}_{j})\to H^{1}(S_{j}|_{A})$
is injective because $\theta_{A}L^{-1}$ is a direct summand of $N^{*}|_{A}$: indeed $\ker\theta_{A}\subset N|_{A}$ is the elementary modification at $p$,
in the direction $T_{p}A'$, of $L^{\perp}|_{A}/TA\cong\OO(1)^{n-1}\oplus\OO^{n-1}\oplus\OO(-1)$ (Lemma~\ref{lem:linesplit}), and since
$d\theta(T_{p}A,T_{p}A')\neq0$ the direction $T_{p}A'$ has a nonzero $\OO(-1)$-component, so $\ker\theta_{A}$ has only $n-1$ summands $\OO(1)$; hence
$\theta_{A}|_{E^{A}_{+}}\colon\OO(1)^{n}\to\OO(1)$ is a nonzero constant surjection, which splits. So $\ker\subset\langle\epsilon_{A},\epsilon_{A'}\rangle$ and $r_{j}\le2$.

Now $\epsilon_{A}$ maps to zero in $H^{1}(S_{j})$ iff $e(p)\in\mathcal V_{A}+\mathcal V_{A'}$, where $e$ is a local generator of $\Lambda^{A}_{j}$ and $\mathcal V_{A}$, $\mathcal V_{A'}$ are the
images of $H^{0}(A,S_{j}|_{A})$, $H^{0}(A',S_{j}|_{A'})$ in $S_{j}(p)$ (\v{C}ech computation on the cover by a neighbourhood of $p$ and
$A\smallsetminus p$, $A'\smallsetminus p$). Write $N^{*}(p)=W_{A}\oplus W_{A'}$ with $W_{A}=E^{A}_{+}(p)^{\perp}$, $W_{A'}=E^{A'}_{+}(p)^{\perp}$, complementary for a
balanced conic. The sections of $S_{j}|_{A}$ not vanishing at $p$ come from $\OO(1)\otimes S^{j+1}\OO^{n}$ and
$\OO(1)\otimes S^{j}\OO^{n}\otimes\OO(-1)^{n}$, so $\mathcal V_{A}=L(p)\otimes\{$monomials with at most one factor from $W_{A'}\}$, and symmetrically. By
Proposition~\ref{prop:Hholds}, $\theta_{A}(p)\in W_{A'}$, so $e(p)\in L(p)\otimes S^{j+1}W_{A'}\subset\mathcal V_{A'}$. Likewise for $\epsilon_{A'}$. So $r_{j}=2$,
and (2) follows, together with $\chi(\mathcal G_{j})=\chi(S_{j})+2$.
\end{proof}

In particular, by Proposition~\ref{prop:formalnbhd}, $H^{1}(C_{0},\mathcal G_{1})=0$. The groups $H^{1}(C_{0},\mathcal G_{j})$, $j\ge2$, are the spaces
of higher-order deformations of the contact neighbourhood of $C_{0}$; they do not vanish, and the second-order neighbourhood of a nodal conic is
not rigid, unlike that of a balanced transverse conic (Theorem~\ref{thm:thirdorder}(1)).

\subsection{Symbols of contactomorphisms}
\label{sec:symbols}

The normalization of \S\ref{sec:dilation} improves a contact transition order by order without taking logarithms. It uses the following symbol, which must be
read off from the actual displacement of the map: an arbitrary vector field congruent to it modulo $I^{j+1}T$ will not do. (For $\theta=dP+Q\,da+\sum c_{i}\,db_{i}$
and $C$ the $a$-axis, the fields $0$ and $Q^{3}\partial_{P}$ are congruent modulo $I^{3}T$, but $\theta(Q^{3}\partial_{P})=Q^{3}$ has $X_{Q^{3}}=3Q^{2}\partial_{a}-2Q^{3}\partial_{P}\notin I^{3}T$.)

\begin{lemma}[Symbols of contactomorphisms]
\label{lem:symbol}
Let $C\subset Z$ be a smooth curve with ideal $I$, let $j\ge2$, and let $R$ be a holomorphic contactomorphism defined near an open subset of $C$ that is the identity on
the neighbourhood of order $j-1$: $R^{*}f-f\in I^{j}$ for every holomorphic $f$. In holomorphic coordinates $x=(x_{a})$ let
\[
V_{R}=\sum_{a}(x_{a}\circ R-x_{a})\,\frac{\partial}{\partial x_{a}},\qquad r_{R}=\theta(V_{R}),
\]
the displacement field of $R$ and its Hamiltonian.
\begin{enumerate}
\item $V_{R}-X_{r_{R}}\in I^{j+1}T$; in particular $r_{R}\in\mathcal A_{j}$.
\item The class $\sigma_{j}(R)\in\mathcal G_{j}$ of $r_{R}$, the \emph{symbol} of $R$, does not depend on the coordinates.
\item $\sigma_{j}(R\circ S)=\sigma_{j}(R)+\sigma_{j}(S)$ for two such maps; $\sigma_{j}(\exp X_{h})=[h]$ for $h\in\mathcal A_{j}$; and $\sigma_{j}(R)=0$ exactly when
$R^{*}f-f\in I^{j+1}$ for every $f$.
\end{enumerate}
\end{lemma}

\begin{proof}
Work locally, with $L$ trivialized, so that $R^{*}\theta=\lambda\theta$. The coefficients of $V_{R}$ lie in $I^{j}$.

(1) Taylor expansion gives $R^{*}\theta-\theta-\mathcal L_{V_{R}}\theta\in I^{2j-1}\Omega^{1}\subset I^{j+1}\Omega^{1}$: the nonlinear terms contain either two displacement factors,
of order $2j$, or a displacement and the derivative of one, of order at least $2j-1$. So $\mathcal L_{V_{R}}\theta\equiv(\lambda-1)\theta$ modulo $I^{j+1}\Omega^{1}$. Put
$W=V_{R}-X_{r_{R}}$. Then $\theta(W)=0$ and $\iota_{W}d\theta=\mathcal L_{V_{R}}\theta-\mathcal L_{X_{r_{R}}}\theta$ is a multiple of $\theta$ modulo $I^{j+1}\Omega^{1}$, so
$\iota_{W}d\theta$ restricted to $\ker\theta$ has coefficients in $I^{j+1}$. The Levi form is nondegenerate on $\ker\theta$, so $W\in I^{j+1}T$. Hence
$X_{r_{R}}=V_{R}-W\in I^{j}T$, that is, $r_{R}\in\mathcal A_{j}$.

(2) Under a change of coordinates $y=F(x)$, Taylor expansion shows that the two displacement fields, compared as vector fields, differ by an element of $I^{2j}T$. Their
Hamiltonians differ by an element of $I^{2j}L\subset I^{j+2}L\subset\mathcal A_{j+1}$ (Lemma~\ref{lem:basicfilt}).

(3) The displacement field of $R\circ S$ differs from $V_{R}+V_{S}$ by an element of $I^{2j-1}T\subset I^{j+1}T$, since $S^{*}g-g\in I^{2j-1}$ for $g\in I^{j}$. Applying (1) to
$R$, $S$ and $R\circ S$ gives the additivity. The time-one flow of $X_{h}$, $h\in\mathcal A_{j}$, has displacement $X_{h}$ modulo $I^{2j-1}T$, so $\sigma_{j}(\exp X_{h})=[h]$.
Finally, by (1), $V_{R}\in I^{j+1}T$ if and only if $X_{r_{R}}\in I^{j+1}T$, that is, $r_{R}\in\mathcal A_{j+1}$.
\end{proof}


\section{Regularity of the curvature along the node exceptional divisor}
\label{sec:dilation}

In a general quaternionic-K\"ahler germ the curvature quartic $W_{z}$ can have poles along the divisor $D'_{z}$ of nodal conics passing through
$z$ on a branch; this is a local phenomenon, which the Hartogs reduction discharges because $D'_{z}$ is contracted by $\widetilde\mu_{z}$
(Example~\ref{ex:bl0}). Along the exceptional divisor $E$ over the nodes, which maps onto the contact hyperplane, the curvature must be shown
to be regular. This section proves it for every germ along a balanced nodal conic, to all orders and for all $n$. The mechanism is a symmetry of
the homogeneous model rather than a computation: points of $E$ are dilation images of interior points, and every deformation of the germ has
non-negative weight under the dilation.

\subsection{Overview of the proof}
\label{sec:overview6}

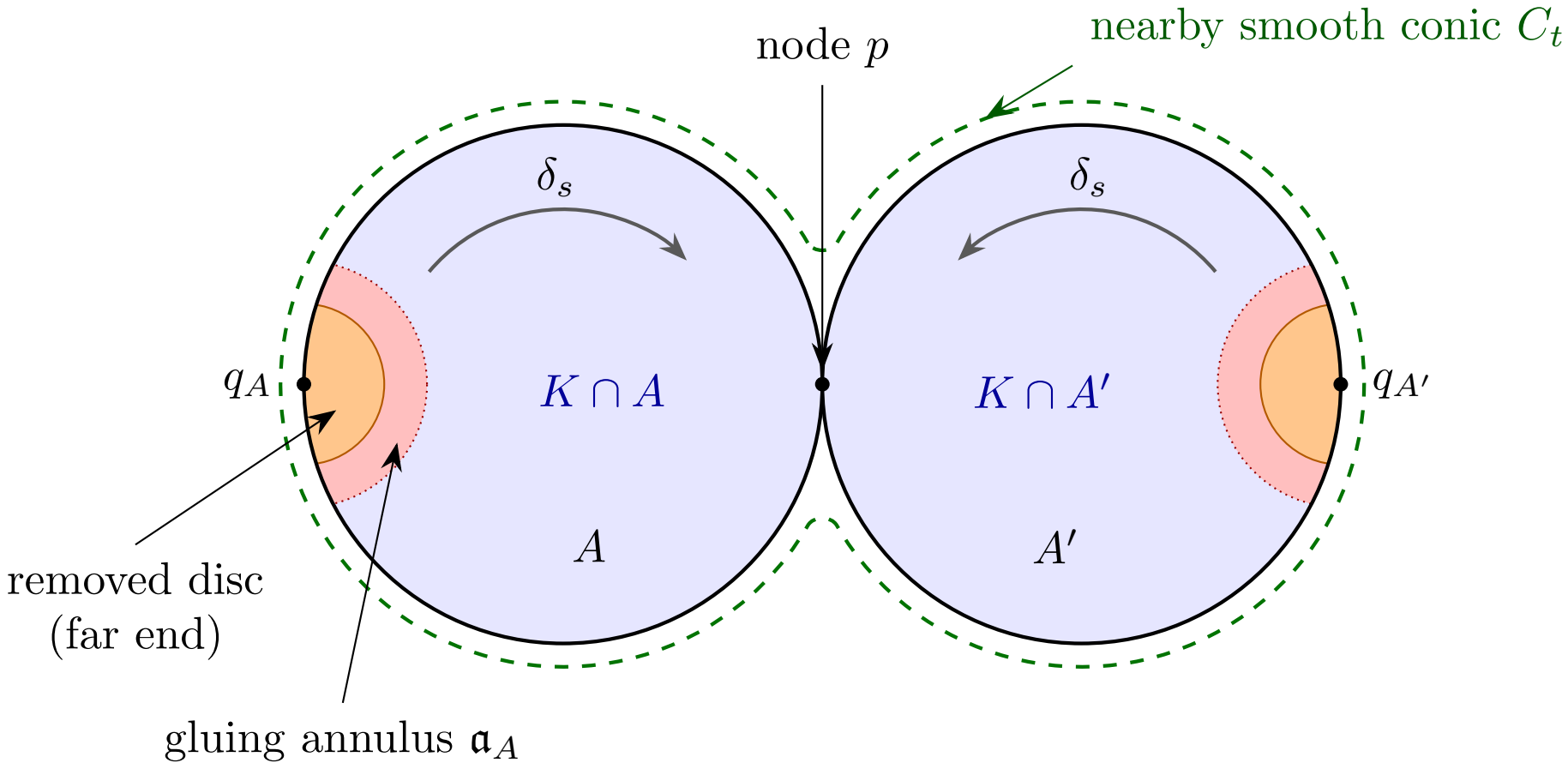
\begin{figure}[ht]
\centering
\begin{tikzpicture}[>={Stealth[length=2.2mm]},font=\small]
  \def\R{2.0}\def\rc{0.62}\def\ra{0.95}\def\ro{2.18}
  \coordinate (cA) at (-\R,0); \coordinate (cB) at (\R,0); \coordinate (p) at (0,0);
  \coordinate (qA) at (-2*\R,0); \coordinate (qB) at (2*\R,0);
  \foreach \c/\q in {cA/qA,cB/qB}{
    \fill[blue!10] (\c) circle (\R);
    \begin{scope}\clip (\c) circle (\R);
      \fill[red!25] (\q) circle (\ra);
      \fill[orange!45] (\q) circle (\rc);
      \draw[orange!70!black] (\q) circle (\rc);
      \draw[red!60!black,densely dotted] (\q) circle (\ra);
    \end{scope}
    \draw[thick] (\c) circle (\R);}
  \draw[dashed,thick,green!45!black]
    ($(cA)+(30:\ro)$) arc[start angle=30,end angle=330,radius=\ro]
    to[out=60,in=120] ($(cB)+(210:\ro)$)
    arc[start angle=210,end angle=510,radius=\ro]
    to[out=240,in=300] ($(cA)+(30:\ro)$);
  \foreach \pt in {p,qA,qB} \fill (\pt) circle (1.6pt);
  \node[left=3pt] at (qA) {$q_{A}$};
  \node[right=3pt] at (qB) {$q_{A'}$};
  \node (lp) at (0,2.6) {node $p$};
  \draw[->,thin] (lp.south) -- ($(p)+(0,0.1)$);
  \node at ($(cA)+(0.2,-1.25)$) {$A$};
  \node at ($(cB)+(-0.2,-1.25)$) {$A'$};
  \node[blue!60!black] at ($(cA)+(0.3,-0.05)$) {$K\cap A$};
  \node[blue!60!black] at ($(cB)+(-0.3,-0.05)$) {$K\cap A'$};
  \draw[->,thick,gray!70!black] ($(cA)+(140:1.35)$) arc[start angle=140,end angle=45,radius=1.35];
  \node at ($(cA)+(92:1.6)$) {$\delta_{s}$};
  \draw[->,thick,gray!70!black] ($(cB)+(40:1.35)$) arc[start angle=40,end angle=135,radius=1.35];
  \node at ($(cB)+(88:1.6)$) {$\delta_{s}$};
  \node[align=center] (lc) at (-5.3,-1.75) {removed disc\\(far end)};
  \draw[->,thin] (lc.north) -- ($(qA)+(0.25,-0.2)$);
  \node[align=center] (la) at (-3.7,-2.75) {gluing annulus $\mathfrak a_{A}$};
  \draw[->,thin] (la.north) -- ($(qA)+(0.72,-0.45)$);
  \node[align=center,green!35!black] (lt) at (3.9,2.75) {nearby smooth conic $C_{t}$};
  \draw[->,thin,green!35!black] (lt.south west) -- ($(cB)+(110:\ro)$);
\end{tikzpicture}
\caption{The nodal conic $C_{0}=A\cup A'$, its branches drawn schematically as spheres touching at the node $p$ (in $Z$ they meet transversally).
$K$ (blue) is $C_{0}$ minus small discs around the far ends $q_{A},q_{A'}$ (orange). A neighbourhood of $K$ in $Z$ is contactomorphic to one in the model,
and so are neighbourhoods of the orange discs; $Z$ is obtained by gluing these pieces on the annuli $\mathfrak a_{A},\mathfrak a_{A'}$ (red), where all the
information about $Z$ along $C_{0}$ sits. For $s\to0$ the dilation $\delta_{s}$ moves the points of $K$ towards $p$, and $\delta_{s}^{-1}$ moves them
towards the far ends. A nearby smooth conic (dashed) is a single sphere with a neck near $p$. It passes through the gluing regions, so no conic lies in a
neighbourhood of $K$.}
\label{fig:K}
\end{figure}

The object is $W^{\CC}$ on $\widetilde F$, a rational section of $\tilde p^{*}(S^{4}L^{\perp*}\otimes L^{2})$ (\S\ref{sec:Wz}), and the claim is that it has
no pole along $E$. The proof does not compute $W^{\CC}$ near $E$. It shows that $W^{\CC}$ has no pole along arcs that end on $E$, and meromorphy turns
this into regularity. In one line: approaching $E$ transversally is the same as dilating a fixed interior point; the dilation does not change
$W^{\CC}$ but changes the germ, the dilated germs of a finite normal form converge because every deformation of the germ has weight $\ge0$, and
the germ has the same polar part along $E$ as a finite normal form of high enough order.
\begin{enumerate}[leftmargin=2em]
\item \emph{Reduction.} If $W^{\CC}=y^{-k}h$ near a generic point of $E=\{y=0\}$ with $k\ge1$, it is unbounded along every arc meeting $E$ transversally
at a point where $h\neq0$. So it suffices to find such arcs, ending at a dense set of points of $E$, along which $W^{\CC}$ is bounded. Boundedness is
measured in a fixed frame near the node, since the bundle is pulled back from $Z$. Near $E$, $W^{\CC}$ is a quotient of holomorphic functions
(Lemma~\ref{lem:findet}(2)).
\item \emph{Transverse arcs are dilation orbits.} The Heisenberg dilation $\delta_{s}$ of the model scales the node coordinates by $s$ and the smoothing
parameter by $s^{2}$, so the orbit of an interior point of a smooth conic is a disc transverse to $E$ (\S\ref{sec:conv6}(3), Lemma~\ref{lem:dilE}).
\item \emph{Covariance.} $W^{\CC}$ has weight $0$ under $\delta_{s}$, so its value at the dilated point of $Z$ equals its value at the fixed interior point
for the dilated germ $\delta_{s}^{*}g$ (Lemma~\ref{lem:dilcov}).
\item \emph{Normal form.} Near $K$, the nodal conic with small discs around its far ends removed (Figure~\ref{fig:K}), $Z$ is contactomorphic to the
model, and so it is near the two discs. So $Z$ is these pieces glued on two annuli, by maps that are the identity to first order and whose higher-order
terms are classes of weight $\ge0$ (Lemmas~\ref{lem:firstorder} and~\ref{lem:normalN}, Proposition~\ref{prop:dilweights}). The first step is first-order rigidity: the first-order contact
neighbourhood of a balanced nodal conic is that of the model. At the node this is an orbit problem for the data that glue the two branch neighbourhoods,
which reduces to finite jets and is solved by a Mayer--Vietoris argument (Lemma~\ref{lem:firstorder}).
\item \emph{Convergence.} Stop the normal form at an order $N_{0}$; this is the \emph{finite model}. Conjugation by $\delta_{s}$ multiplies the weight-$w$
part of its gluing by $s^{w}$. As no weight is negative, the dilated finite models converge to the model reglued by the weight-$0$ parts, and $W^{\CC}$ at the
fixed point converges with them (Lemma~\ref{lem:jetconv}).
\item \emph{Comparison.} The germ and its finite model agree to order $N_{0}$ along $C_{0}$, and so, by a property of deformations over Artinian
rings, do their universal families of conics. $W^{\CC}$ is a quotient $\Pi/\Delta$ of universal polynomials in finitely many jets of the universal family, and the
order $d_{\Delta}$ of $\Delta$ at a general point of $E$ is fixed by the germ. So for $N_{0}\ge d_{\Delta}+4$ the germ and its finite model have the same polar part along matching arcs
(Lemma~\ref{lem:findet}(2),(3)). The comparison is formal: no estimate on shrinking domains is needed.
\item \emph{Conclusion.} The dilation arcs end at a dense set of points of $E$, and $W^{\CC}$ has no pole along them (Theorem~\ref{thm:Eregular}).
\end{enumerate}
The model $Z_{0}$ of \S\ref{sec:model} is the twistor space of $\Gr_{2}(\CC^{n+2})$. It serves as a reference germ, not as an assumption on $Z$: what the
proof uses is a Darboux theorem along $K$ rather than at a point, together with a dilation defined along the whole of $C_{0}$, and the model is the simplest
place where both are explicit. Every conic near $C_{0}$ passes through the gluing regions (Figure~\ref{fig:K}), so the identification near $K$ says nothing
directly about the metric near $[C_{0}]$; this is why the proof follows the gluing under the dilation.

\subsection{Conventions and the local picture at the node}
\label{sec:conv6}

Let $\sigma\colon\widetilde F=\mathrm{Bl}_{\mathcal B}\overline F\to\overline F$ be the blow-up of \S\ref{sec:resolution}, with exceptional divisor $E$,
and put
\[
\tilde\pi=\bar\pi\circ\sigma\colon\widetilde F\to\MCbar,\qquad\tilde p=\bar p\circ\sigma\colon\widetilde F\to Z .
\]
\begin{enumerate}[leftmargin=2em]
\item \emph{The local form.} Near a nodal conic $C_{0}=A\cup A'$ with nef branches, $\MCbar$ is smooth and maps submersively onto the versal deformation
$xy=t$ of the node (Lemma~\ref{lem:nodalsmooth}). So near the node $\overline F$ has coordinates $(x,y,\vec w)$, $\vec w\in\CC^{4n-1}$, with
$\bar\pi(x,y,\vec w)=(t,\vec w)=(xy,\vec w)$, and $t$ is a reduced equation of the nodal divisor $\mathcal D$ (\S\ref{sec:resolution}). The fibre of $\bar\pi$ over
$(t,\vec w)$ is the curve $\{xy=t\}$: for $t\neq0$ the \emph{neck} of the conic, an annulus with coordinate $x$; for $t=0$ the branches $A=\{y=0\}$ and $A'=\{x=0\}$.
Hence
\[
\bar\pi^{*}\mathcal D=D_{A}+D_{A'},\qquad D_{A}=\{y=0\},\quad D_{A'}=\{x=0\},
\]
where $D_{A}$ (resp.\ $D_{A'}$) consists of the pairs $([C_{0}],z)$ with $z\in A$ (resp.\ $z\in A'$). Both components are smooth and reduced, and they meet
transversally along the node section $\mathcal B=\{x=y=0\}$.
\item \emph{The blow-up.} In the chart $x=ay$ of $\widetilde F$, $E=\{y=0\}$ and the strict transform $\widetilde D_{A'}$ of $D_{A'}$ is $\{a=0\}$; in the chart
$y=bx$, $E=\{x=0\}$ and $\widetilde D_{A}=\{b=0\}$. The fibre of $E$ over a point of $\mathcal B$ is the line $\PP^{1}$ of \emph{smoothing ratios} $a=x/y$; it meets
$\widetilde D_{A'}$ at $a=0$ and $\widetilde D_{A}$ at $a=\infty$. Since $D_{A}$ and $D_{A'}$ are smooth and contain the centre, $\sigma^{*}D_{A}=\widetilde D_{A}+E$ and
$\sigma^{*}D_{A'}=\widetilde D_{A'}+E$, so
\begin{equation}
\tilde\pi^{*}\mathcal D=\widetilde D_{A}+\widetilde D_{A'}+2E,\qquad\text{in the chart }x=ay:\ t=a\,y^{2}.\label{eq:6-0}
\end{equation}
This is \eqref{eq:4-1}; the multiplicity $2$ is the multiplicity of $xy$ along $\mathcal B$.
\item \emph{The double cover.} By \eqref{eq:6-0} a disc transverse to $E$, say $y=s$ with $a\in\CC^{*}$ and $\vec w$ fixed, maps to $\MCbar$ by $s\mapsto(as^{2},\vec w)$, a
double cover branched along $\mathcal D$. Correspondingly a smooth conic $C_{t}$ near $C_{0}$ has two points with ratio $a$, namely $\pm(\sqrt{at},\sqrt{t/a})$, on
opposite sides of the neck, and the monodromy around $\mathcal D$ exchanges them. A neck point with a given ratio is therefore a single-valued function of $s$, $s^{2}=t$,
but not of $t$. This is why Corollary~\ref{cor:Rnodal} works on the double cover $t=s^{2}$, and why the dilation parameter $s$ of \S\ref{sec:model}, which scales the
node coordinates by $s$ and $t$ by $s^{2}$, is a coordinate transverse to $E$ (Lemma~\ref{lem:dilE}). The normalization of the base change of $\widetilde F$ by
$s\mapsto s^{2}=t$ is, locally, the product of $\CC^{4n-1}$ with the minimal resolution of the surface $xy=s^{2}$. Its central fibre is the reduced chain
$\widetilde A\cup R\cup\widetilde A'$, where $R$ is a $(-2)$-curve mapping two-to-one onto the fibre of $E$, branched at $a=0,\infty$; the two points of $R$ over $a$
are the limits of the two neck points with ratio $a$. In a surface slice $\vec w=\mathrm{const}$ of $\widetilde F$ itself, the fibre $\ell$ of $E$ is a $(-1)$-curve and
the central fibre is the non-reduced curve $\widetilde A+2\ell+\widetilde A'$.
\item \emph{What is extended.} The object studied on $\widetilde F$ is $W^{\CC}$, a rational section of $\tilde p^{*}(S^{4}L^{\perp*}\otimes L^{2})$. This bundle is
pulled back from $Z$ and needs no square root of $L$ (\S\ref{sec:Wz}). Over the transverse non-jumping conics the evaluation map
$\mathrm{ev}\colon\tilde\pi^{*}E^{*}_{\CC}\to\tilde p^{*}B$, $\mathrm{ev}([C],z)=\mathrm{ev}_{z}$, is an isomorphism, and
\[
W^{\CC}=(\mathrm{ev}^{-1})^{*}\,\tilde\pi^{*}R_{E}^{\CC},
\]
the pull-back of $R_{E}^{\CC}$ written in the moving frame $\mathrm{ev}^{-1}(b_{i}(z))$. Regularity along $E$ is a statement about $W^{\CC}$, not about
$\tilde\pi^{*}R_{E}^{\CC}$, which would require an extension of $E_{\CC}$ across $\mathcal D$. Unlike a pull-back, $W^{\CC}$ is not constant on the fibres of $\tilde\pi$:
along the fibre of $E$ its value can depend on the ratio $a$ (Theorem~\ref{thm:Eregular}), and along $\widetilde D_{A}$ and $\widetilde D_{A'}$ it may have poles
(Remark~\ref{rem:Dprime}). The translation back to $R_{E}^{\CC}$ on $\MCbar$ is Corollary~\ref{cor:Rnodal}.
\item \emph{Notation in the model.} In \S\S\ref{sec:model}--\ref{sec:Ethm} the model coordinates $a,Q$ play the roles of $x,y$, and $\tau$ that of $t$, so there the
smoothing ratio is written $a:Q$ (Lemma~\ref{lem:dilE}). From Corollary~\ref{cor:Rnodal} on we return to the notation $x,y,a$ of \S\ref{sec:resolution}.
\end{enumerate}

\subsection{The model and its dilation}
\label{sec:model}

Let $Z_{0}=F(1,n+1;\CC^{n+2})$, $n\ge1$, with the chart $x=(1,a,b_{1},\dots,b_{n-1},P)$, $\xi=(R,Q,c_{1},\dots,c_{n-1},1)$,
$R=-(aQ+\sum b_{i}c_{i}+P)$, in which $\theta=dP+Q\,da+\sum c_{i}\,db_{i}$. Its nodal conic $C_{0}=A\cup A'$ has $A$ the $a$-axis, $A'$ the $Q$-axis and
node $p$ at the origin; $q_{A}$, $q_{A'}$ denote the far ends of $A$ and $A'$. The conics $C_{\tau}=\{b=c=P=0,\ aQ=\tau\}$ smooth the node. The element
$\delta_{s}=\mathrm{diag}(1,s,\dots,s,s^{2})$ of $PGL_{n+2}$ acts on $Z_{0}$ by contactomorphisms; it fixes $p,q_{A},q_{A'}$, preserves $A$ and $A'$, and
$\delta_{s}(C_{\tau})=C_{s^{2}\tau}$. In the chart, $a,Q,b_{i},c_{i}$ have weight $1$, $P$ and $R$ weight $2$, and $\delta_{s}^{*}\theta=s^{2}\theta$.

\emph{Weights.} We use pull-back weights: $f\circ\delta_{s}=s^{w}f$ for a function, so $\partial_{a}$ has weight $-1$, and the frame $e_{L}$ of $L$ with
$\theta=\theta_{\mathrm{fn}}e_{L}$ has weight $-2$. For a Hamiltonian $h=h_{\mathrm{fn}}e_{L}$ the field $X_{h}$ has weight $\operatorname{wt}h_{\mathrm{fn}}-2$; this is the
\emph{weight} of $h$, and $\delta_{s}^{*}X_{h}=X_{s^{w}h}$ for $h$ of weight $w$. At $q_{A}$ the coordinates
\[
\alpha=x_{0}/x_{1},\quad\beta_{i}=x_{i+1}/x_{1},\quad\gamma=x_{n+1}/x_{1},\quad\rho=\xi_{0}/\xi_{n+1},\quad\chi_{i}=\xi_{i+1}/\xi_{n+1}
\]
have weights $-1,0,1,2,1$, and $L_{q_{A}}$ has weight $-1$ (at $p$ it has weight $-2$). In the node chart $\beta_{i}=b_{i}/a$. Symmetrically at $q_{A'}$,
where $\beta'_{i}=c_{i}/Q$ has weight $0$. The conormal coordinates of $A$ at $q_{A}$ are $\beta,\gamma,\rho,\chi$.

The universal family of conics near $C_{0}$ has the local form of \S\ref{sec:resolution} with $x\leftrightarrow a$, $y\leftrightarrow Q$,
$t\leftrightarrow\tau$; by \eqref{eq:4-1},
\begin{equation}
\operatorname{ord}_{E}\tau=2,\qquad\operatorname{ord}_{E}a=\operatorname{ord}_{E}Q=1 .\label{eq:6-1}
\end{equation}

\begin{lemma}[Covariance]
\label{lem:dilcov}
Let $g$ be a contact manifold germ along $C_{0}$, given by charts that are open pieces of $Z_{0}$ with contact transition maps, and let
$\delta_{s}^{*}g$ be the same germ with every chart precomposed with $\delta_{s}$. Then for every non-jumping transverse conic $C$ of $\delta_{s}^{*}g$ and $z\in C$,
\[
W^{\CC}_{g}([\delta_{s}C],\delta_{s}z)=W^{\CC}_{\delta_{s}^{*}g}([C],z),
\]
where $\delta_{s}C$ is the same curve in the charts of $g$, both sides being written in the frames $\partial_{a}-Q\partial_{P},\partial_{Q},
\partial_{b_{i}}-c_{i}\partial_{P},\partial_{c_{i}}$ of $L^{\perp}$ and $\theta$ of $L$.
\end{lemma}

\begin{proof}
$g$ and $\delta_{s}^{*}g$ are the same contact manifold in different charts, so they have the same conics, the same LeBrun structure and the same
$W^{\CC}$, a section of $S^{4}L^{\perp*}\otimes L^{2}$. Under $\delta_{s}$ the frame vectors of $L^{\perp}$ scale by $s$ and $\theta$ by $s^{2}$, so the
components change by $s^{4}\cdot s^{-4}=1$.
\end{proof}

\begin{lemma}[Points of $E$ are dilation images]
\label{lem:dilE}
Let $z'\in C_{1}$ with $a(z')=a'\neq0$, $Q(z')=1/a'$. The arc $s\mapsto(C_{s^{2}},\delta_{s}z')$ meets $E$ transversally at the point over
$([C_{0}],p)$ with smoothing ratio $a':Q(z')$. As $z'$ moves on $C_{1}$ these ratios fill $\CC^{*}\subset\PP^{1}$.
\end{lemma}

\begin{proof}
$\delta_{s}z'=(sa',\dots,s/a',\dots)$ lies on $C_{s^{2}}$. By \eqref{eq:6-1}, $\operatorname{ord}_{E}s=1$, and the ratio $a:Q$ tends to $a':Q(z')$.
\end{proof}

\subsection{Weights of the deformations of a nodal neighbourhood}
\label{sec:weights}

The graded automorphism sheaves $\mathcal G_{j}$ of \S\ref{sec:autsheaves} for $C_{0}\subset Z_{0}$ are $\CC^{*}$-equivariant.

\begin{proposition}[Weights of the nodal classes]
\label{prop:dilweights}
For all $n\ge1$ and $j\ge1$, every weight of $\CC^{*}$ on $H^{1}(C_{0},\mathcal G_{j})$ is $\ge0$. More precisely:
\begin{enumerate}
\item classes coming from $H^{1}(S_{j}|_{A})$ have weight $\ge$ (the least weight of the fibre $S_{j}(q_{A})$) $+1\ \ge0$, and likewise on $A'$;
\item classes supported at the node have weight $\ge j-1$, and classes from $\mathcal P_{j}$ have weight $\ge j$;
\item weight $0$ occurs only for $n\ge2$ and $j\ge2$, exactly on the classes $\alpha^{-1}\beta^{m}e'_{L}$, $|m|=j+1$, at $q_{A}$ ($e'_{L}$ the equivariant
frame of $L$ at $q_{A}$, of weight $-1$) and their mirrors at $q_{A'}$: $\binom{j+n-1}{n-2}$ on each branch.
\end{enumerate}
\end{proposition}

\begin{proof}
\emph{A bound on $\PP^{1}$.} Let $V$ be a $\CC^{*}$-equivariant bundle on $A\cong\PP^{1}$. A global section of weight $w$ of an equivariant bundle
expands near $q_{A}$ as $\sum c\,\alpha^{l}f_{i}$ in an equivariant frame; as $\alpha$ has weight $-1$, $w\le\max_{i}\operatorname{wt}f_{i}$. Equivariant Serre
duality, $H^{1}(V)\cong H^{0}(V^{*}\otimes K_{A})^{*}$ (the trace class $a^{-1}da\in H^{1}(K_{A})$ has weight $0$ and $K_{A}$ has weight $-1$ at $q_{A}$),
turns this into: every weight of $H^{1}(A,V)$ is at least $\min\operatorname{wt}V(q_{A})+1$.

\emph{The fibre at $q_{A}$.} $S_{j}=L\otimes S^{j+1}N^{*}$; at $q_{A}$, $N^{*}(q_{A})=N^{*}_{A}(q_{A})$ is spanned by $d\beta_{i},d\gamma,d\rho,d\chi_{i}$ of
weights $0,1,2,1$, and $L(q_{A})$ has weight $-1$. So the weights of $S_{j}(q_{A})$ are $\ge-1$, with equality exactly on the $\beta$-monomials of degree
$j+1$. This gives (1).

\emph{The node and $\mathcal P_{j}$.} By Mayer--Vietoris, $H^{1}(C_{0},S_{j})$ is an extension of $H^{1}(S_{j}|_{A})\oplus H^{1}(S_{j}|_{A'})$ by a quotient of
$S_{j}(p)$. At $p$ the conormal generators are $aQ,P$ (weight $2$) and $b_{i},c_{i}$ (weight $1$), and $L_{p}$ has weight $-2$, so $S_{j}(p)$ has weights
$\ge j-1$. By Theorem~\ref{thm:naive}(1), $\mathcal P_{j}\cong\iota_{A*}L^{1-j}(-p)\oplus\iota_{A'*}L^{1-j}(-p)$, with frame $aP^{j}e_{L}$ near $p$ on $A$
(weight $2j-1$) and degree $-j$, so the weights of its $H^{1}$ are $j,\dots,2j-2$. By Lemma~\ref{lem:basicfilt}, $H^{1}(\mathcal G_{j})$ is an extension of
$H^{1}(\mathcal P_{j})$ by $H^{1}(S_{j}/\mathcal P_{j+1})$, a quotient of $H^{1}(S_{j})$ (as $H^{0}(\mathcal P_{j})=0$ and $H^{2}=0$). This proves (2) and the first
statement.

\emph{Weight $0$.} Only $\beta$-monomials reach weight $-1$ at $q_{A}$. In $N^{*}|_{A}$ the span of the $d\beta_{i}$ ($=db_{i}$ near $p$) is an equivariant
summand $\OO(-1)^{n-1}$, so $\beta^{m}e_{L}$, $|m|=j+1$, spans a summand $\OO(-j)$ of $S_{j}|_{A}$ with weights $j-1$ at $p$ and $-1$ at $q_{A}$; its $H^{1}$
has weights $0,\dots,j-2$, the weight-$0$ class being $\alpha^{-1}\beta^{m}e'_{L}$. The node and $\mathcal P$ contributions have positive weight, so these
classes survive and form the whole weight-$0$ part. Their number is the number of $\beta$-monomials of degree $j+1$ in $n-1$ variables.
\end{proof}

\subsection{Charts adapted to the dilation}
\label{sec:charts}

\begin{lemma}[First-order rigidity]
\label{lem:firstorder}
Let $C_{0}=A\cup A'\subset(Z,\theta)$ be a nodal conic whose branches are standard lines, with $N_{C_{0}/Z}$ balanced, and let $\iota$ be an isomorphism of
$C_{0}$ onto the nodal conic of $Z_{0}$ with $\iota(p)=p$. Then there is an isomorphism $\OO_{Z}/I_{C_{0}}^{2}\to\OO_{Z_{0}}/I_{C_{0}}^{2}$ of the first-order
neighbourhoods over $\iota$ that is induced, near every point of $C_{0}$, by a contactomorphism: the first-order contact neighbourhood of $C_{0}$ (the germ modulo
$\mathcal A_{2}$) is that of the model.
\end{lemma}

\begin{proof}
\textbf{Step 1: rigidity of the branches.} \emph{Local triviality.} Near $p$ the germ is that of $Z_{0}$ by Proposition~\ref{prop:normalform}(1); by relative Moser in the $(a,Q)$-plane, with
conformal constant $\varphi_{A}'(0)\varphi_{A'}'(0)$, the normal form can realize any pair of branch reparametrizations fixing $p$. Near other points of
$A$ or $A'$ the germ is that of $Z_{0}$ by Proposition~\ref{prop:normalform}(2). So the first-order neighbourhood of $C_{0}$ is a torsor $\mathcal T$ under
the sheaf of groups $\mathcal A_{1}/\mathcal A_{2}$ of $Z_{0}$, which is not abelian.

\emph{The branches.} Let $\mathsf G^{A}$ be the sheaf of first-order contact automorphisms of the neighbourhood of $A$ that are the identity on $A$.
Put $SN_{A}=(TA)^{\perp}/TA$, where $(TA)^{\perp}\subset L^{\perp}|_{A}$ is the Levi-orthogonal of $TA$; it is an $L|_{A}$-valued symplectic bundle, and
$\mathsf G^{A}$ maps onto $\mathcal R=\mathrm{GSp}(SN_{A})$, its conformal-symplectic group with multiplier $L|_{A}$ (for $n=1$, $\mathbb G_{m}$, the
conformal factor of $\theta$). The image of $\mathcal T|_{A}$ in $H^{1}(A,\mathcal R)$ is the class of $SN_{A}$ as an $L|_{A}$-valued symplectic bundle. By
Lemma~\ref{lem:linesplit}, $SN_{A}\cong V_{1}\oplus V_{0}$ with $V_{1}=\OO(1)^{n-1}$, $V_{0}=\OO^{n-1}$, and multiplier $\OO(1)$. Any $\OO(1)$-valued
symplectic form $\omega$ on $V_{1}\oplus V_{0}$ is equivalent to that of $Z_{0}$:
\begin{itemize}[leftmargin=2em]
\item $\omega$ vanishes on $\Lambda^{2}V_{1}$, as $\Lambda^{2}V_{1}^{*}\otimes\OO(1)\cong\OO(-1)^{\binom{n-1}{2}}$ has no sections;
\item it pairs $V_{1}$ with $V_{0}$ by a constant matrix $C$, a section of $V_{1}^{*}\otimes V_{0}^{*}\otimes\OO(1)\cong\OO^{(n-1)^{2}}$, invertible since $\omega$ is
nondegenerate;
\item on $\Lambda^{2}V_{0}$ it is a skew matrix $S$ of linear forms, and the automorphism of $V_{1}\oplus V_{0}$ that fixes $V_{1}$ and sends $v\in V_{0}$ to
$v+\tfrac12C^{-T}Sv$ (a map $V_{0}\to V_{1}$ being a matrix of linear forms) turns $S$ into $S-S=0$;
\item a constant change of basis of $V_{0}$ then makes $C$ the identity.
\end{itemize}
So the class is that of $Z_{0}$.

The kernel $\mathcal U$ of $\mathsf G^{A}\to\mathcal R$ is unipotent. Its Lie algebra $\mathfrak u$ is the kernel of
$(L\otimes S^{2}N_{A}^{*})/L^{-1}\to L\otimes S^{2}SN_{A}^{*}$, that is $\big(L\otimes\mathrm{Ann}(SN_{A})\cdot N_{A}^{*}\big)/L^{-1}$, where
$\mathrm{Ann}(SN_{A})=(N_{A}/SN_{A})^{*}$. Here $N_{A}/SN_{A}=TZ|_{A}/(TA)^{\perp}$ is an extension of $L|_{A}\cong\OO(1)$ by
$L^{\perp}/(TA)^{\perp}\cong T^{*}A\otimes L\cong\OO(-1)$; it is non-split, because otherwise this quotient of the nef bundle $TZ|_{A}$ would have a summand $\OO(-1)$ (cf.\ \cite{Ye93}, where the contact sequence is shown to be non-split on every curve of
positive $L$-degree). So $N_{A}/SN_{A}\cong\OO^{2}$ and
$\mathrm{Ann}(SN_{A})\cong\OO^{2}$. So $\mathfrak u$ is a quotient of $L\otimes\mathrm{Ann}(SN_{A})\otimes N_{A}^{*}\cong\OO^{2n-2}\oplus\OO(1)^{2n+2}$. It is also a direct summand
of $(L\otimes S^{2}N_{A}^{*})/L^{-1}$, since $N_{A}^{*}\to SN_{A}^{*}$ splits ($\mathrm{Ext}^{1}=H^{1}(SN_{A})^{2}=0$). Filter $\mathfrak u$ by contact weight
($\mathrm{Ann}$ having weight $2$ and $SN_{A}^{*}$ weight $1$): the graded pieces are $\OO^{2n-2}\oplus\OO(1)^{2n-2}$ and a quotient of $\OO(1)^{3}$, all
of non-negative degree. By d\'evissage along this filtration, $H^{1}(A,\mathcal U)$ is trivial (a point), so the fibre of $H^{1}(A,\mathsf G^{A})\to
H^{1}(A,\mathcal R)$ over the class of $Z_{0}$ is a point and $\mathcal T|_{A}$ is trivial; likewise $\mathcal T|_{A'}$. Since $\mathcal A_{1}/\mathcal A_{2}$
is not abelian (with $\theta=dP+Q\,da$, $[X_{P},X_{Q^{2}}]=X_{Q^{2}}$), the vanishing of $H^{1}(\mathcal G_{1})$ alone would not suffice; this is why the
branches are treated through $\mathcal R$ and $\mathcal U$.

\textbf{Step 2: the matching data and their finite jets.} \emph{Matching at the node.} In the coordinates of Proposition~\ref{prop:normalform}(1), $I_{A}=(Q,b,c,P)$, $I_{A'}=(a,b,c,P)$ and $I_{C_{0}}=(aQ,b,c,P)$. A
monomial $a^{i}Q^{j}m$, with $m$ a monomial of degree $k$ in $b,c,P$, lies in $I_{A}^{2}\cap I_{A'}^{2}$ iff $j+k\ge2$ and $i+k\ge2$, and in $I_{C_{0}}^{2}$ iff
$k\ge2$, or $k=1$ and $i,j\ge1$, or $k=0$ and $i,j\ge2$; these conditions agree, so $I_{A}^{2}\cap I_{A'}^{2}=I_{C_{0}}^{2}$ and
\[
\OO/I_{C_{0}}^{2}=\OO/I_{A}^{2}\times_{B}\OO/I_{A'}^{2},\qquad B=\OO/(I_{A}^{2}+I_{A'}^{2}),
\]
near $p$. ($B$ has basis $1,a,Q,b_{i},c_{i},P,aQ$: the normal order along $C_{0}$ is not the Taylor order at $p$, since $aQ\in I_{C_{0}}$ has Taylor degree $2$.)
Let $\mathsf G^{A}$ be the sheaf on $A$ of automorphisms of $\OO_{Z_{0}}/I_{A}^{2}$ that are, near each point, induced by a contactomorphism fixing $A$
pointwise; its Lie algebra is $\mathfrak g^{A}=\mathcal A^{A}_{1}/\mathcal A^{A}_{2}$, acting by the derivations $D_{h}$ induced by the $X_{h}$
(\S\ref{sec:autsheaves}). Define $\mathsf G^{A'}$, $\mathfrak g^{A'}$ likewise, and let $\mathcal K_{p}$ be the group of germs at $p$ of automorphisms of
$\OO_{Z_{0}}/I_{C_{0}}^{2}$ induced by contactomorphisms fixing $C_{0}$ pointwise; its Lie algebra is $\mathfrak k_{p}=(\mathcal A_{1}/\mathcal A_{2})_{p}=\mathcal G_{1,p}$ of
$C_{0}$. By the fibre product, restriction $\rho=(\rho_{A},\rho_{A'})\colon\mathcal K_{p}\to\mathsf G^{A}_{p}\times\mathsf G^{A'}_{p}$ is injective, and likewise on Lie
algebras, since $X_{h}\in I_{A}^{2}T\cap I_{A'}^{2}T=I_{C_{0}}^{2}T$ forces $h\in\mathcal A_{2}$.

\emph{The gluing data.} By the branch step, and composing with an automorphism of $Z_{0}$ preserving $A_{0}$ (the stabilizer of the corresponding line and
hyperplane acts on $A_{0}$ through $PGL_{2}$), there is an isomorphism $\varphi_{A}$ from the first-order neighbourhood of $A_{0}$ onto that of $A$ over $\iota^{-1}$,
induced near each point by a contactomorphism; likewise $\varphi_{A'}$. By the local triviality, there is a contactomorphism $\varphi_{p}$ from a neighbourhood of $p$ in
$Z_{0}$ onto one in $Z$ carrying $C_{0}$ to $C_{0}$ over $\iota^{-1}$. Put
\[
g_{A}=\varphi_{A}^{-1}\varphi_{p}\in\mathsf G^{A}_{p},\qquad g_{A'}=\varphi_{A'}^{-1}\varphi_{p}\in\mathsf G^{A'}_{p},
\]
restrictions to the branch neighbourhoods near $p$. Let $G=H^{0}(A,\mathsf G^{A})\times H^{0}(A',\mathsf G^{A'})$, acting on germs at $p$ by left multiplication.
\emph{It suffices to show $(g_{A},g_{A'})\in G\cdot\rho(\mathcal K_{p})$.} Indeed, if $g_{A}=h_{A}\rho_{A}(k)$ and $g_{A'}=h_{A'}\rho_{A'}(k)$, then $\varphi_{A}h_{A}$ and
$\varphi_{A'}h_{A'}$ agree near $p$ with the first-order map of the contactomorphism $\varphi_{p}\circ\tilde k^{-1}$, $\tilde k$ a contactomorphism inducing $k$, on the
respective branch neighbourhoods; by the fibre product they glue to an isomorphism $\OO_{Z}/I_{C_{0}}^{2}\to\OO_{Z_{0}}/I_{C_{0}}^{2}$ over $\iota$, induced near
every point by a contactomorphism. Conversely, the class of $(g_{A},g_{A'})$ in $G\backslash(\mathsf G^{A}_{p}\times\mathsf G^{A'}_{p})/\rho(\mathcal K_{p})$ does not
depend on the choices of $\varphi_{A}$, $\varphi_{A'}$, $\varphi_{p}$, and determines the first-order neighbourhood: given $x=(g_{A},g_{A'})$, choose contactomorphisms $\psi_{A}$, $\psi_{A'}$ inducing
$g_{A}$, $g_{A'}$ on a neighbourhood $U$ of $p$ in $Z_{0}$, and glue $U$ to disjoint tubes around $A_{0}\smallsetminus D$ and $A'_{0}\smallsetminus D'$, with $D$, $D'$
small closed discs about $p$ in $U$, by $\psi_{A}$ and $\psi_{A'}$ on the overlaps; this is a contact manifold germ $Z_{x}$ along $C_{0}$, whose first-order neighbourhood is the one with data $x$. The model is $Z_{(1,1)}$.

\emph{Reduction to finite jets.} Let $\mathfrak m$ be the ideal of $p$. An element of $\mathsf G^{A}_{p}$ fixes $p$, so it acts on the finite-dimensional
algebra $\OO/(I_{A}^{2}+\mathfrak m^{M})$; let $J^{A}_{M}$ be the image and $\mathcal N^{A}_{M}$ the kernel. In the coordinates above, write
$y=(Q,b,c,P)$; an element of $\mathsf G^{A}_{p}$ is $a\mapsto a+\langle u(a),y\rangle$, $y\mapsto\Lambda(a)y$ modulo $I_{A}^{2}$, and it lies in $\mathcal N^{A}_{M}$
iff $u$ and $\Lambda-1$ have coefficients in $a^{M-1}\OO$. For $h=\kappa(a)P+q$, $q$ quadratic in $y$ with coefficients in $\OO_{A}$, one finds modulo $I_{A}^{2}T$
\[
X_{h}\equiv q_{Q}\,\partial_{a}+\textstyle\sum_{i}q_{c_{i}}\,\partial_{b_{i}}+\sum_{i}(\kappa c_{i}-q_{b_{i}})\,\partial_{c_{i}}+(\kappa Q-\kappa'P)\,\partial_{Q}+\kappa P\,\partial_{P}.
\]
Every quadratic monomial in $y$ other than $P^{2}\in\mathcal A^{A}_{2}$ contains $Q$, $b$ or $c$, so if the coefficients of $D_{h}$ lie in $a^{l}\OO$, then $h\in
a^{l}\mathcal A^{A}_{1}+\mathcal A^{A}_{2}$.

\emph{Claim.} For $M\ge4$, $\rho(\mathcal K_{p})\supset\mathcal N^{A}_{M}\times\mathcal N^{A'}_{M}$, and $\rho(\mathfrak k_{p})\supset\mathrm{Lie}\,\mathcal N^{A}_{M}\oplus
\mathrm{Lie}\,\mathcal N^{A'}_{M}$.

(i) Let $h\in a^{3}\mathcal A^{A}_{1}$. As $\mathcal A^{A}_{1}$ is an $\OO$-module and $a\in I_{A'}$, $h\in\mathcal A^{A}_{1}\cap I_{A'}^{3}L$, so $X_{h}\in I_{A}T\cap I_{A'}^{2}T$
(Lemma~\ref{lem:basicfilt}); hence $X_{h}\in I_{C_{0}}T$, $h\in\mathfrak k_{p}$, and $\rho(h)=([h],0)$. The same holds for a time-dependent $h_{t}\in a^{3}\mathcal A^{A}_{1}$:
its flow up to time $1$, defined near $p$ because $X_{h_{t}}$ vanishes on $C_{0}$, fixes $C_{0}$ pointwise, is the identity on the first-order neighbourhood of $A'$, and
induces on that of $A$ the solution of $\dot g_{t}=D_{h_{t}}\circ g_{t}$, $g_{0}=1$, a linear equation, pointwise in $a$, for the coefficients $u$, $\Lambda$ of $g_{t}$. This gives the Lie algebra statement, since $\mathrm{Lie}\,\mathcal N^{A}_{M}\subset
a^{M-1}\mathfrak g^{A}$.

(ii) Let $\delta\in\mathcal N^{A}_{M}$ be induced by a contactomorphism $\varphi$ fixing $A$ pointwise. The weighted dilation
$\mu_{t}(a,Q,b,c,P)=(a,t^{2}Q,tb,tc,t^{2}P)$ satisfies $\mu_{t}^{*}\theta=t^{2}\theta$ and fixes $A$ pointwise. Since $\varphi$ preserves the Levi-orthogonal
$(TA)^{\perp}=\{dQ=dP=0\}$ along $A$, every monomial in the Taylor expansion of $\varphi^{*}y_{i}$ has weight at least that of $y_{i}$ (weights $2,1,1,2$ for $Q,b,c,P$);
so $\varphi_{t}=\mu_{t}^{-1}\varphi\mu_{t}$, whose components are $\varphi_{t}^{*}y_{i}=t^{-w_{i}}\mu_{t}^{*}\varphi^{*}y_{i}$, is a family of contactomorphisms fixing $A$
pointwise, holomorphic in $t$ for $|t|\le1$, with $\varphi_{1}=\varphi$ and $\varphi_{0}$ weighted-homogeneous. Its first-order maps $\delta_{t}$ stay in $\mathcal N^{A}_{M}$, as $\mu_{t}$ does not involve $a$. Comparing coefficients in
$\varphi_{0}^{*}\theta=\lambda\theta$, $\lambda=\varphi^{*}\theta/\theta|_{A}$, the map $\varphi_{0}$ is $a\mapsto a$, $(b,c)\mapsto S(a)(b,c)$, $Q\mapsto\lambda Q-\lambda'P+q_{2}$,
$P\mapsto\lambda P+p_{2}$, with $q_{2},p_{2}$ quadratic in $(b,c)$ and $S(a)$ conformally symplectic with multiplier $\lambda$ (for $n=1$ read $CSp(0)$ as $\CC^{*}$ and $S$ as $\lambda$); so $\delta_{0}$ is determined by $S$, and
$S\equiv1$ modulo $a^{M-1}$. Since $\delta_{0}\in\mathcal N^{A}_{M}$, also the coefficient $\lambda'$ of $P$ in the image of $Q$ lies in $a^{M-1}\OO$. The path
$S_{t}=\exp(t\log S)$, $t\in[0,1]$, has multiplier $\lambda_{t}=\exp(t\log\lambda)$, and $\lambda_{t}'=t\lambda_{t}\lambda'/\lambda\in a^{M-1}\OO$; with
$S_{t}-1\in a^{M-1}\OO$ this keeps the path in $\mathcal N^{A}_{M}$. It has velocity $\log S$, a germ of $\mathfrak{csp}(2n-2)$ with coefficients in $a^{M-1}\OO$,
which is $D_{h}$ for some $h=\kappa P+q_{2}$ with coefficients in $a^{M-1}\OO$. Concatenating, $\delta$ is the end point of a piecewise smooth path in
$\mathcal N^{A}_{M}$ from $1$ whose right velocities $\dot\delta_{t}\delta_{t}^{-1}$ are derivations with coefficients in $a^{M-1}\OO$; on the second piece they are the
first-order maps of the contact fields $\dot\varphi_{t}\circ\varphi_{t}^{-1}$, which vanish on $A$. By the computation above they are $D_{h_{t}}$ with
$h_{t}\in a^{M-1}\mathcal A^{A}_{1}$, piecewise smooth in $t$; as $M-1\ge3$, (i) gives $k\in\mathcal K_{p}$ with $\rho(k)=(\delta,1)$. Likewise on $A'$.

\textbf{Step 3: the orbit argument.} Fix $M=4$ and let $\ell\colon\mathsf G^{A}_{p}\times\mathsf G^{A'}_{p}\to\mathbb G=J^{A}_{4}\times J^{A'}_{4}$ be the jet map. Let $\mathfrak j^{A}\subset\mathrm{Der}_{\CC}\big(\OO/(I_{A}^{2}+\mathfrak m^{4})\big)$ be the image of $\mathfrak g^{A}_{p}$, a complex Lie
subalgebra, and $H^{A}$ its connected analytic subgroup in the automorphism group of this finite-dimensional algebra. The $4$-jet of the flow of a Hamiltonian
representative of $v\in\mathfrak j^{A}$ is $\exp v$, so $H^{A}\subset J^{A}_{4}$. Conversely, the construction of (ii) applies to any $\delta\in\mathsf G^{A}_{p}$,
joining $S$ to $1$ through $S(sa)$, $s\in[0,1]$, and a path in $\mathrm{CSp}(n-1,\CC)$, each piece realized by flows of contact fields vanishing on $A$; so every element of
$J^{A}_{4}$ is joined to $1$ by a piecewise smooth path whose right logarithmic derivatives lie in $\mathfrak j^{A}$, and such a path stays in $H^{A}$. Hence
$J^{A}_{4}=H^{A}$, a connected complex Lie group with Lie algebra $\mathfrak j^{A}$ (no closedness is needed). Products of the flows of Hamiltonian
representatives of a basis of $\mathfrak j^{A}$, with small complex times, give local holomorphic charts of $J^{A}_{4}$, in which the first-order neighbourhood
data vary holomorphically. Likewise on $A'$, so $\mathbb G$ is a connected complex Lie group. By the claim, the
double coset of $(g_{A},g_{A'})$ is determined by $x=\ell(g_{A},g_{A'})$, and $(g_{A},g_{A'})\in G\cdot\rho(\mathcal K_{p})$ iff $x\in\ell(G)\,\ell(\rho(\mathcal K_{p}))$. Call
$x\in\mathbb G$ \emph{admissible} if $C_{0}$ is balanced in $Z_{x}$. The normal bundle of $C_{0}$ in $Z_{x}$ depends only on the first-order neighbourhood: its
restriction to $A$ is the elementary modification at $p$ of $N_{A_{0}}$ in the direction that $g_{A}$ assigns to $T_{p}A'$, and likewise on $A'$. For every $x$ the node is
general in the sense of Lemma~\ref{lem:branchres}: near $p$, $Z_{x}$ is the model, so $T_{p}A$ and $T_{p}A'$ are not Levi-orthogonal, while the positive part of
$TZ_{x}|_{A}$, which contains $T_{p}A$ and lies in $L^{\perp}$, is Levi-isotropic (its summands have degree $\ge1$ and $L|_{A}=\OO(1)$, Lemma~\ref{lem:linesplit}), so it does not contain $T_{p}A'$, and symmetrically. By
Lemma~\ref{lem:branchres} the restriction to each branch is $\OO(1)^{n}\oplus\OO^{n}$, so the fibres $E^{A}_{+}(p)$, $E^{A'}_{+}(p)$ of the positive parts depend
holomorphically on $x$, and admissibility, $E^{A}_{+}(p)\cap E^{A'}_{+}(p)=0$, is an open determinantal condition: locally on $\mathbb G$ it is the non-vanishing of
the determinant of holomorphic frames of the two spaces. The non-admissible set is thus a closed complex-analytic subset of $\mathbb G$, proper since $1$ is
admissible, and the admissible set $\Omega$, its complement in a connected complex manifold, is connected, and it is invariant under $\ell(G)$ and $\ell(\rho(\mathcal K_{p}))$, which do not change the first-order neighbourhood
of $Z_{x}$.

Let $x=\ell(g)\in\Omega$. The infinitesimal first-order automorphisms of $Z_{x}$ fixing $C_{0}$ form the sheaf $\mathcal G_{1}(Z_{x})=\mathcal A_{1}/\mathcal A_{2}$ of
\S\ref{sec:autsheaves}. Restricted to the branch neighbourhoods and transported by the branch charts $\varphi_{A},\varphi_{A'}$ of $Z_{x}$, which are first-order contact
isomorphisms and so respect $\mathcal A_{1}/\mathcal A_{2}$, it maps to $\nu_{*}(\mathfrak g^{A}\oplus\mathfrak g^{A'})$; its stalk at $p$ goes to $\mathrm{Ad}_{g}\rho(\mathfrak k_{p})$,
as $Z_{x}$ is the model near $p$ and $g_{A}=\varphi_{A}^{-1}\varphi_{p}$. By the fibre
product, there is an exact sequence of sheaves on $C_{0}$,
\[
0\to\mathcal G_{1}(Z_{x})\to\nu_{*}\big(\mathfrak g^{A}\oplus\mathfrak g^{A'}\big)\to\mathcal Q_{x}\to0,\qquad
\mathcal Q_{x}=\big(\mathfrak g^{A}_{p}\oplus\mathfrak g^{A'}_{p}\big)\big/\mathrm{Ad}_{g}\rho(\mathfrak k_{p})\ \text{at }p,
\]
with $\nu$ the normalization, so $H^{0}(A,\mathfrak g^{A})\oplus H^{0}(A',\mathfrak g^{A'})\to\mathcal Q_{x}\to H^{1}(C_{0},\mathcal G_{1}(Z_{x}))$ is exact. As $C_{0}$
is balanced in $Z_{x}$, Theorem~\ref{thm:naive}(2) and Proposition~\ref{prop:formalnbhd} (with $j=1$, $k=2$), applied to the germ $Z_{x}$, give
$H^{1}(C_{0},\mathcal G_{1}(Z_{x}))=0$. By the Lie algebra part of the claim, $\mathrm{Ad}_{g}\rho(\mathfrak k_{p})$ contains the kernel of $d\ell$, so
$\mathcal Q_{x}$ is the quotient of $T_{x}\mathbb G$ (right-translated to $\mathrm{Lie}\,\mathbb G$) by the tangent space of $x\,\ell(\rho(\mathcal K_{p}))$. Hence the map
$(\xi,\eta)\mapsto\ell\big(\exp\xi\cdot g\cdot\exp(-\rho(\eta))\big)$, with $\xi\in H^{0}(A,\mathfrak g^{A})\oplus H^{0}(A',\mathfrak g^{A'})$ and $\eta$ in a
finite-dimensional subspace of $\mathfrak k_{p}$ mapping onto $d\ell(\rho(\mathfrak k_{p}))$, is a submersion at $0$ (here $\exp\xi\in G$, and $\exp\eta\in\mathcal K_{p}$ is
the flow of $X_{\eta}$, which vanishes on $C_{0}$), so its image contains a neighbourhood of $x$: every class of $\Omega$ under $x\sim\ell(h)\,x\,\ell(\rho(k))^{-1}$ is
open. As
$\Omega$ is connected, it is a single class, that of $1$. The datum of $Z$ is admissible, so $(g_{A},g_{A'})\in G\cdot\rho(\mathcal K_{p})$.
\end{proof}

\begin{lemma}[Charts adapted to the dilation]
\label{lem:normalN}
In the situation of Lemma~\ref{lem:firstorder}, let $K\subset C_{0}$ be the complement of small open discs around $q_{A}$ and $q_{A'}$, and let $N_{0}\ge2$.
There are neighbourhoods $U_{K}$ of $K$ and $U_{A}$, $U_{A'}$ of closed discs around $q_{A}$, $q_{A'}$ in $Z$, with contactomorphisms
$\psi_{K},\psi_{A},\psi_{A'}$ from neighbourhoods of the corresponding pieces of the nodal conic of $Z_{0}$ inducing $\iota$, such that on the annuli
$\mathfrak a_{A}=\psi_{K}^{-1}(U_{K}\cap U_{A})$ and $\mathfrak a_{A'}$
\[
\psi_{K}^{-1}\psi_{A}=\exp X_{h_{2}}\circ\cdots\circ\exp X_{h_{N_{0}}}\circ R_{N_{0}},
\]
and likewise at $q_{A'}$, where each $h_{j}$ is a finite sum of $\CC^{*}$-homogeneous Laurent monomials in $\alpha$ of nonnegative weight, polynomial in the normal
coordinates, representing a class of $H^{1}(C_{0},\mathcal G_{j})$, whose weight-$0$ part is exactly $\sum_{|m|=j+1}c_{m}\alpha^{-1}\beta^{m}e'_{L}$; and
$R_{N_{0}}$ is a contactomorphism, holomorphic near the annulus, which is the identity on the neighbourhood of order $N_{0}$ of $A$: $R_{N_{0}}^{*}f-f\in I_{A}^{N_{0}+1}$
for every holomorphic $f$.
\end{lemma}

\begin{proof}
\emph{Charts.} $K$ is two closed discs glued at $p$; it is contractible, and $K$ and its model $K_{0}\subset Z_{0}$ have Stein neighbourhoods \cite{Siu76}.
Restricting the isomorphism of Lemma~\ref{lem:firstorder} to $K$ gives an isomorphism of the first-order neighbourhoods, $\OO_{Z}/I_{K}^{2}\to\OO_{Z_{0}}/I_{K_{0}}^{2}$ over
$\iota$, induced near each point by a contactomorphism. It is realized by a holomorphic map, although $K$, being singular at $p$, is not a holomorphic retract of any of its neighbourhoods. Embed a Stein
neighbourhood $U$ of $K$ as a closed submanifold of some $\CC^{N}$, and let $\rho\colon W\to U$ be a holomorphic retraction of a neighbourhood $W$ of $U$ in
$\CC^{N}$ \cite{DG60}; this uses only that $U$ is smooth. The coordinate functions of $\CC^{N}$, restricted to $U$, give sections of $\OO_{Z_{0}}/I_{K_{0}}^{2}$ on
a Stein neighbourhood $U_{0}$ of $K_{0}$, which lift to holomorphic functions on $U_{0}$ because $H^{1}(U_{0},I_{K_{0}}^{2})=0$ (Cartan's Theorem~B). This gives a
holomorphic $F\colon U_{0}\to\CC^{N}$ with the prescribed $1$-jet along $K_{0}$; after shrinking $U_{0}$, $F(U_{0})\subset W$, and $\Psi_{0}:=\rho\circ F$ has the same
$1$-jet along $K_{0}$, since $F$ and $\rho\circ F$ agree on $K_{0}$ and $d\rho$ is the identity on $TU$, which contains the image of $dF$ along $K_{0}$. So $\Psi_{0}$ is a
local biholomorphism along $K_{0}$ (its differential is the given isomorphism, at $p$ that of the normal-form chart), and on $K_{0}$ it is the bijection
$K_{0}\to K$ given by $\iota$. A local biholomorphism near a compact set on which it is injective is injective on a neighbourhood of it: otherwise there are
$x_{k}\neq y_{k}$ converging to $K_{0}$ with $\Psi_{0}(x_{k})=\Psi_{0}(y_{k})$, whose limits $x,y\in K_{0}$ satisfy $\Psi_{0}(x)=\Psi_{0}(y)$, so $x=y$, contradicting
injectivity near $x$. Hence $\Psi_{0}$ is a biholomorphism from a neighbourhood of $K_{0}$ onto a neighbourhood of $K$ inducing the given first-order isomorphism.

\emph{Second-order contact correction.} $\Psi_{0}$ need not be contact to second order, and composing it with contactomorphisms cannot make it so: a contactomorphism
preserving $K_{0}$ acts invertibly on the class of $\Psi_{0}^{*}\theta-u\theta_{0}$ in $I\Omega^{1}/(I^{2}\Omega^{1}+I\theta_{0})$, and changing $u$ does not alter
it. (For instance, in the coordinates of
Proposition~\ref{prop:normalform}(1), $\Psi_{0}=(a,Q,b,c,P+b_{1}^{2})$ induces the identity to first order, and $\Psi_{0}^{*}\theta_{0}-\theta_{0}=2b_{1}\,db_{1}$.) The
correction is made with an ordinary vector field. Write $I=I_{K_{0}}$, and trivialize $L$ near $K$ and near $K_{0}$ (on each branch, matched at $p$, and extended by
Theorem~B), so that $\theta$, $\theta_{0}$ are $1$-forms; let $R_{0}$ be the Reeb field of $\theta_{0}$, $\theta_{0}(R_{0})=1$.
\begin{itemize}[leftmargin=2em]
\item \emph{Local form.} Near each point $x$ of $K_{0}$ the first-order isomorphism is the $1$-jet of a local contactomorphism $\psi_{x}$, $\psi_{x}^{*}\theta=u_{x}\theta_{0}$
(by its construction, Lemma~\ref{lem:firstorder}). Then $\varphi_{x}=\psi_{x}^{-1}\circ\Psi_{0}$ is the identity modulo $I^{2}$, $\varphi_{x}=\mathrm{id}+\epsilon$ in
coordinates with $\epsilon\in I^{2}$. Expanding, $\varphi_{x}^{*}\theta_{0}\equiv\theta_{0}+d\big(\theta_{0}(\epsilon)\big)$ modulo $I^{2}\Omega^{1}$, so
\[
\Psi_{0}^{*}\theta\equiv u_{x}\theta_{0}+dh_{x}\pmod{I^{2}\Omega^{1}},\qquad h_{x}=u_{x}\,\theta_{0}(\epsilon)\in I^{2}.
\]
\item \emph{Globalization.} In particular $\Psi_{0}^{*}\theta\equiv u_{x}\theta_{0}$ modulo $I\Omega^{1}$, so the $u_{x}$ agree modulo $I$; let $u$ be a holomorphic function near
$K_{0}$ congruent to them (Theorem~B), nonvanishing after shrinking, and $\beta_{0}=\Psi_{0}^{*}\theta-u\theta_{0}$. The map $\bar d\colon I^{2}\to\mathcal Q=I\Omega^{1}/(I^{2}\Omega^{1}+I\theta_{0})$,
$h\mapsto dh$, is $\OO$-linear, since $d(fh)-f\,dh=h\,df\in I^{2}\Omega^{1}$. The class of $\beta_{0}$ in $\mathcal Q$ is locally $\bar d h_{x}$; as $\ker\bar d$ is coherent and
$H^{1}(U_{0},\ker\bar d)=0$ on a Stein neighbourhood $U_{0}$, it is $\bar d h$ for some $h\in H^{0}(U_{0},I^{2})$: $\beta_{0}=dh+v\theta_{0}+\beta'$ with $v\in I$,
$\beta'\in I^{2}\Omega^{1}$.
\item \emph{Correction.} Put $Y=-(h/u)R_{0}\in I^{2}T$ and replace $\Psi_{0}$ by $\Psi_{0}\circ\exp Y$. Since $Y$ vanishes on $K_{0}$, its flow exists up to time one near
$K_{0}$, and $\exp Y$ is the identity modulo $I^{2}$, so the first-order isomorphism is unchanged. For such a map, $(\exp Y)^{*}\omega\equiv\omega$ for
$\omega\in I\Omega^{1}$ and $(\exp Y)^{*}\theta_{0}\equiv\theta_{0}+d(\theta_{0}(Y))$, modulo $I^{2}\Omega^{1}$. Hence
\[
(\Psi_{0}\circ\exp Y)^{*}\theta\equiv u\theta_{0}+u\,d(-h/u)+dh+v\theta_{0}\equiv(u+v)\,\theta_{0}\pmod{I^{2}\Omega^{1}}.
\]
\end{itemize}
Replacing $u$ by $u+v$, now $\beta_{0}=\Psi_{0}^{*}\theta-u\theta_{0}\in I^{2}\Omega^{1}$, so $\theta_{t}=(1-t)u\theta_{0}+t\Psi_{0}^{*}\theta$ agree to second order along $K$ and are contact near $K$ for $t\in[0,1]$. Gray's argument \cite{Gra59} gives
$V_{t}$ with $V_{t}\in\ker\theta_{t}$ and $\iota_{V_{t}}d\theta_{t}=-\dot\theta_{t}$ on $\ker\theta_{t}$; this equation is pointwise linear in
$\dot\theta_{t}=\Psi_{0}^{*}\theta-u\theta_{0}=\beta_{0}\in I^{2}\Omega^{1}\otimes L$, so $V_{t}\in I^{2}T$. As $V_{t}$ vanishes on $K$, its flow is defined up to $t=1$
on a neighbourhood of the compact $K$, fixes $K$ pointwise and is the identity modulo $I^{2}$, so it carries $\theta_{0}$ to a multiple of $\Psi_{0}^{*}\theta$
without changing the first-order class. This gives $\psi_{K}$; the same argument on discs gives $\psi_{A}$, $\psi_{A'}$.

\emph{Symbols.} A contactomorphism $R$ near the annulus that is the identity on the neighbourhood of order $j-1$ of $A$ has a symbol $\sigma_{j}(R)\in\mathcal G_{j}$,
read off from its displacement field (Lemma~\ref{lem:symbol}); it is additive, $\sigma_{j}(\exp X_{h})=[h]$ for $h\in\mathcal A_{j}$, and it vanishes exactly when $R$ is the
identity on the neighbourhood of order $j$. No logarithm of $R$ is needed.

\emph{Normalization.} The traces on $C_{0}$ of $U_{K}$, $U_{A}$, $U_{A'}$ have Stein neighbourhoods and $\mathcal G_{j}$ is coherent, so by Leray's
theorem the cover computes $H^{1}(C_{0},\mathcal G_{j})$: a class is a pair of sections of $\mathcal G_{j}$ on the two annuli, modulo differences of
sections on the trace of $U_{K}$ and on the discs. Suppose
$\psi_{K}^{-1}\psi_{A}=\exp X_{h_{2}}\circ\cdots\circ\exp X_{h_{j-1}}\circ R$ with $R$ a contactomorphism near the annulus that is the identity on the neighbourhood of
order $j-1$ of $A$ (for $j=2$, $R=\psi_{K}^{-1}\psi_{A}$, both charts inducing the same first-order isomorphism). The symbol $\sigma_{j}(R)$, with its analogue at $q_{A'}$,
is such a $1$-cocycle. Its class is a finite sum of classes of pure weight, since $H^{1}(C_{0},\mathcal G_{j})$
is finite-dimensional and $\CC^{*}$-graded; and each class of pure weight $w$ has an \emph{equivariant representative}: in an equivariant frame of
$\mathcal G_{j}$ near $q_{A}$ (monomials in the normal coordinates times $e'_{L}$, by the grading), the Laurent expansion in $\alpha$ of a section on the
annulus, which is invariant under the unit circle, is its decomposition into weights, and its weight-$w$ part is a Laurent polynomial of weight $w$
representing the weight-$w$ component of the class. Let $h_{j}$ be the sum of these representatives. The class has no component of negative weight
(Proposition~\ref{prop:dilweights}), so $h_{j}$ has weights $\ge0$, and its weight-$0$ part lies in the span of the $\alpha^{-1}\beta^{m}e'_{L}$, which
span the weight-$0$ part of $H^{1}$ (Proposition~\ref{prop:dilweights}(3)). The symbol minus $h_{j}$ has class zero, so by Leray it is a coboundary
$c_{K}-c_{A}$, with $c_{K}$, $c_{A}$ sections of $\mathcal G_{j}$ on the traces of $U_{K}$ and $U_{A}$; no splitting of Laurent series at the node is needed. Lift $c_{K},c_{A}$ to Hamiltonians in $\mathcal A_{j}$ (possible on Stein
neighbourhoods, $\mathcal A_{j+1}$ being coherent by Lemma~\ref{lem:basicfilt}) and replace $\psi_{K}$ by $\psi_{K}\circ\exp X_{c_{K}}$, $\psi_{A}$ by
$\psi_{A}\circ\exp X_{c_{A}}$. Since $(\exp X_{h})^{*}X_{c}=X_{c'}$ with $c'-c\in\mathcal A_{j+1}$ for $h\in\mathcal A_{2}$, the new transition is
$\exp X_{h_{2}}\circ\cdots\circ\exp X_{h_{j}}\circ R'$ with $R'=\exp(-X_{h_{j}})\circ\exp(-X_{c'_{K}})\circ R\circ\exp X_{c_{A}}$, a contactomorphism with
$\sigma_{j}(R')=\sigma_{j}(R)-h_{j}-c_{K}+c_{A}=0$ by Lemma~\ref{lem:symbol}(3); so $R'$ is the identity on the neighbourhood of order $j$. Each step shrinks the neighbourhoods slightly; stop after
$N_{0}-1$ steps, and let $R_{N_{0}}$ be the last remainder.

\emph{The choices.} Each $h_{j}$ is determined by its class only up to the choice of representative, and the construction is inductive, so $h_{j}$ for
$j\le N_{0}$ does not depend on $N_{0}$. Only $h_{2},h_{3},h_{4}$ enter the value along $E$ (Lemma~\ref{lem:findet}(1), $k_{0}=4$), and we fix their
choices once and for all; all later statements refer to this choice.
\end{proof}

\subsection{Finite determinacy, polar parts, and the dilated finite models}
\label{sec:jetconv}

\begin{lemma}[Finite determinacy]
\label{lem:findet}

Let $(M,\theta)$ be a complex contact manifold of dimension $2n+1$, $C\subset M$ a transverse non-jumping conic, and $z\neq w$ points of $C$. Choose
holomorphic charts $x=(x_{0},x')$, $x'\in\CC^{2n}$, centred at $z$, and $\tilde x=(\tilde x_{0},\tilde x')$ centred at $w$, in which $C$ is a graph over the unit disc
of the $x_{0}$-axis, resp.\ of the $\tilde x_{0}$-axis, and put $T_{z}=\{x_{0}=0\}$, $T_{w}=\{\tilde x_{0}=0\}$. For $y=(z_{1},w_{1})\in T_{z}\times T_{w}$ near $(z,w)$
let $C(y)$ be the conic near $C$ through $z_{1}$ and $w_{1}$. It exists, is unique, and is transverse and non-jumping, because $H^{0}$ and $H^{1}$ of
$N_{C/M}(-z-w)\cong\OO(-1)^{2n}$ vanish. Near $z$ it is the graph $x'=\Phi(y;x_{0})$, $|x_{0}|<1$, of a holomorphic function $\Phi$, the \emph{two-point family}
of $(M,C;z,w)$. Fix holomorphic frames of $L^{\perp}$ and of $L$ near $z$, in which $W^{\CC}([C],z)$ has components.
\begin{enumerate}
\item The components of $W^{\CC}([C],z)$ are a universal rational function $\mathcal W$ of the $4$-jet of $\Phi$ at $((z,w);0)$, of the $3$-jet of $\theta$ at $z$
and of the values of the frames at $z$, all written in the chart. The denominator of $\mathcal W$ does not vanish on the data of a transverse non-jumping
conic. In particular $k_{0}=4$ will do below. The $4$-jet of $\Phi$ at $((z,w);0)$ is determined by the neighbourhood of order $4$ of $C$, so $W^{\CC}([C],z)$ depends
only on that neighbourhood and on the $3$-jet of $\theta$ at $z$.
\item \emph{At a balanced nodal conic.} Let $C_{0}=A\cup A'\subset M$ be a balanced nodal conic with standard branches and node $p$, let $r\in A$ and
$r'\in A'$ be points other than $p$, and $T\ni r$, $T'\ni r'$ small complex hypersurfaces transverse to $A$ and $A'$.
\begin{enumerate}[label=(\alph*)]
\item The conics of $M$ near $C_{0}$ meet $T$ and $T'$ once each, and $C\mapsto m(C)=(C\cap T,\,C\cap T')$ is a chart, centred at $m_{0}=(r,r')$, of the space of
conics of $M$ near $[C_{0}]$. We write $m=(u,u')$, $C(m)$ for the conic with coordinates $m$, and $Y\subset T\times T'$ for the domain of the chart. In particular the
two-point family $\Phi$ of $(M,C_{0};r,r')$, with $T_{r}=T$ and $T_{r'}=T'$, is defined near $(m_{0};0)$, although $C_{0}$ is singular.
\item Fix a holomorphic frame of $L$ near $r$. Let $U$ be a coordinate ball about $p$ carrying holomorphic frames of $L^{\perp}$ and of $L$, and let
$F=(F_{1},\dots,F_{2n})$ be holomorphic functions
on $Y\times U$ generating the ideal of the incidence $\overline F_{U}=\{(m,z):z\in C(m)\}$. There are universal polynomials, depending only on $n$, which,
evaluated on the derivatives of order $\le4$ of $\Phi$ at $(m;0)$ and of $F$ at $(m,z)$, and on those of the coefficients of $\theta$ and of the frames at $u$ and
at $z$, give holomorphic functions $\Pi$ and $\Delta$ on $Y\times U$ with the following properties: $\Delta$ does not vanish identically on $\overline F_{U}$, and
\[
W^{\CC}\big([C(m)],z\big)=\Pi(m,z)/\Delta(m,z),
\]
in the frames at $z$, whenever $C(m)$ is transverse and non-jumping, $z\in C(m)\cap U$ and $\Delta(m,z)\neq0$. In particular $W^{\CC}$ is meromorphic on
$\overline F_{U}$, and on its blow-up along the node section.
\end{enumerate}
\item \emph{Comparison of polar parts.} Let $(M,\theta)$ and $(M',\theta')$ be complex contact manifolds containing the same balanced nodal conic
$C_{0}$ with standard branches, and an open set $V\ni r,r',p$ on which they coincide as contact manifolds. Let $N\ge4$, and suppose that the neighbourhoods
of order $N$ of $C_{0}$, $(C_{0},\OO_{M}/I^{N+1})$ and $(C_{0},\OO_{M'}/I'^{N+1})$, are isomorphic by an isomorphism whose restriction over $V$ is the identity of
the whole thickening $\OO_{V}/I^{N+1}$, not only of $C_{0}\cap V$. Take
$T,T',U\subset V$ and the frames to be the same for $M$ and $M'$, mark the objects of (2) for $M'$ by a prime, and let $\mathfrak m$ be the ideal of $m_{0}$.
\begin{enumerate}[label=(\alph*)]
\item $\Phi-\Phi'\in\mathfrak m^{N+1}$, and given $F$ one may choose $F'$ with $F-F'\in\mathfrak m^{N+1}\OO_{Y\times U}$. Then $\Pi-\Pi'$ and $\Delta-\Delta'$
lie in $\mathfrak m^{N-3}\OO_{Y\times U}$.
\item Let $\gamma'(s)=(m(s),z'(s))$ be a holomorphic arc in $\overline F'_{U}$ with $\gamma'(0)=b_{0}=(m_{0},p)$, whose lift to the blow-up $\widetilde F'$ of
$\overline F'_{U}$ along the node section meets the exceptional divisor $E'$ transversally. Then there is a formal arc $\gamma(s)=(m(s),z(s))$ in $\overline
F_{U}$ with $z-z'\in s^{N}\CC[[s]]^{2n+1}$ (this part needs only $N\ge2$). Its lift to $\widetilde F$ meets $E$ transversally, at the same point: $E$ and $E'$ have the same fibre
$\PP\big(T_{b_{0}}\overline F_{U}/T_{b_{0}}\mathcal B\big)$ over $b_{0}$, where $\mathcal B$ is the node section, and $\gamma$, $\gamma'$ end at the same
point of it.
\item For such arcs, $\Pi(\gamma)-\Pi'(\gamma')$ and $\Delta(\gamma)-\Delta'(\gamma')$ lie in $s^{N-3}\CC[[s]]$ (componentwise for the vector $\Pi$). Let
$d_{\Delta}$ be the order of $\Delta(\gamma(s))$ in $s$, and suppose $N\ge d_{\Delta}+4$. Then $\Delta'(\gamma'(s))$ has order $d_{\Delta}$ as well, and if
$W^{\CC}_{M'}$ has no pole along $\gamma'$, that is $W'=\Pi'(\gamma')/\Delta'(\gamma')\in\CC[[s]]$, then
\[
W^{\CC}_{M}\big(\gamma(s)\big)-W'\in s^{N-3-d_{\Delta}}\,\CC[[s]],
\]
where $W^{\CC}_{M}(\gamma(s))$ means $\Pi(\gamma)/\Delta(\gamma)\in\CC((s))$. So $W^{\CC}_{M}$ has no pole along $\gamma$, and the two values at $s=0$ agree.
\end{enumerate}
\end{enumerate}

\end{lemma}

The point of (3) is that the loss $d_{\Delta}$ can be fixed before $N$ is chosen. For a fixed presentation $W^{\CC}=\Pi/\Delta$ write
$\sigma^{*}\Delta=y_{E}^{d}v$, $E=\{y_{E}=0\}$; at a point $e_{0}\in E$ with $v(e_{0})\neq0$ every arc meeting $E$ transversally at $e_{0}$ has $d_{\Delta}=d$. (Elsewhere
the order may depend on the arc: for $\Delta=y^{d}u$ the arcs $(y,u)=(s,s^{k})$ end at the same point with orders $d+k$.) So the loss can be paid for by taking $N$
large; the
comparison is formal, and no estimate on shrinking domains is needed.

\begin{proof}

(1) The map $y\mapsto[C(y)]$ is a local coordinate system on the space of conics near $[C]$. It is injective because a conic near $C$ meets $T_{z}$ and
$T_{w}$ once, and its differential is bijective because $H^{0}(N(-z-w))=0$. The conics $C(z_{1},\cdot)$ all pass through $z_{1}$, so at $[C(y)]$ the vectors
$\partial/\partial w_{1}^{a}$ span $T\Sigma_{z_{1}}=E^{*}_{\CC}\otimes h_{z_{1}}$; likewise the $\partial/\partial z_{1}^{b}$ span $E^{*}_{\CC}\otimes h_{w_{1}}$. Both
subspaces are totally null (Theorem~\ref{thm:lebrun}(2)), so $g^{\CC}=2\sum_{a,b}g_{ab}(y)\,dw_{1}^{a}\,dz_{1}^{b}$. The variation fields
$\partial_{w_{1}^{a}}\Phi$ vanish at $x_{0}=0$, because $\Phi(y;0)=x'(z_{1})$. Proposition~\ref{prop:levwronsk}(3), applied at the point $z_{1}$, gives, up to the
universal constant,
\[
g_{ab}(y)=-\frac{d\theta\big(\pi\,\partial_{x_{0}}\partial_{w_{1}^{a}}\Phi(y;0),\ \pi\,\partial_{z_{1}^{b}}\Phi(y;0)\big)}{\theta\big(\partial_{x_{0}}+\partial_{x_{0}}\Phi(y;0)\big)}.
\]
Here $\theta$ and $d\theta$ are taken at $z_{1}$, and $\pi$ is the projection of $T_{z_{1}}M$ onto $L^{\perp}_{z_{1}}$ along $T_{z_{1}}C(y)$, which is the
identification $N_{C(y)}\cong L^{\perp}$ at $z_{1}$; this is where transversality is used. So $g_{ab}$ is a rational function of the $2$-jet of $\Phi$ at $(y;0)$
and the $1$-jet of $\theta$ at $z_{1}$, with denominator $\theta(\partial_{x_{0}}+\partial_{x_{0}}\Phi)\neq0$. Moreover $\det(g_{ab})\neq0$, because $g^{\CC}$ is
nondegenerate at non-jumping conics. The curvature tensor $\mathrm{Rm}$ of $g^{\CC}$ at $[C]$ is polynomial in $g_{ab}$, $\det(g_{ab})^{-1}$ and the
$y$-derivatives of $g_{ab}$ of order $\le2$ at $(z,w)$. These involve the $4$-jet of $\Phi$ at $((z,w);0)$ and the $3$-jet of $\theta$ at $z$.

$W^{\CC}([C],z)$ is read off from $\mathrm{Rm}$ as follows. Let $\tau=\partial_{x_{0}}+\partial_{x_{0}}\Phi((z,w);0)$. For $v\in L^{\perp}_{z}$ let
$X_{w}(v)=\sum_{b}X^{b}\partial_{z_{1}^{b}}$ be the vector with $\pi(\sum_{b}X^{b}e_{b})=v$, and $X_{z}(v)=\sum_{a}X^{a}\partial_{w_{1}^{a}}$ the vector with
$\sum_{a}X^{a}\,\pi\,\partial_{x_{0}}\partial_{w_{1}^{a}}\Phi((z,w);0)=v$. These invert the isomorphisms $X\mapsto X(z)$ of $E^{*}_{\CC}\otimes h_{w}$ onto $L^{\perp}_{z}$
and $X\mapsto\nabla_{z}X(\partial_{x_{0}})$ of $E^{*}_{\CC}\otimes h_{z}$ onto $L^{\perp}_{z}$, so $X_{w}(v)=e\otimes h_{w}$ and $X_{z}(v)=e'\otimes h_{z}$ with $e'$
proportional to $e$. In the decomposition $\mathrm{Rm}=\rho_{0}R_{0}+R^{\CC}_{E}$ of Theorem~\ref{thm:lebrun}(3), every term of $R_{0}$ evaluated on
$(X_{z},X_{w},X_{z},X_{w})$ contains $\omega_{E}$ on two proportional vectors, so it vanishes. Hence, up to the universal constant of $g^{\CC}$ and the sign
convention for $\mathrm{Rm}$,
\[
W^{\CC}([C],z)(v,v,v,v)=\theta(\tau)^{2}\,\mathrm{Rm}\big(X_{z}(v),X_{w}(v),X_{z}(v),X_{w}(v)\big),
\]
and polarization gives all components. The denominators are $\theta(\tau)$, $\det(g_{ab})$, and the determinants of $\pi$ on $\{dx_{0}=0\}$ and of the vectors
$\pi\,\partial_{x_{0}}\partial_{w_{1}^{a}}\Phi((z,w);0)$. The last is nonzero because $H^{0}(N(-2z))=0$, and none of them vanishes at a transverse non-jumping conic.

Finally, the coefficients of the Taylor expansion of $\Phi$ at $((z,w);0)$ of total order $\le4$ are obtained by solving the equations of the conic through
$z_{1}$ and $w_{1}$ order by order, as in the proof of Theorem~\ref{thm:thirdorder}. At each order the equations involve the transition maps only modulo the
ideal of $C$ to the fifth power. So they are determined by the neighbourhood of order $4$.

(2)(a) Write $N_{C_{0}}=N_{C_{0}/M}$. By the proof of Proposition~\ref{prop:balanced}, $H^{0}$ and $H^{1}$ of $N_{C_{0}}(-r-r')$ vanish when $C_{0}$ is
balanced, so evaluation $H^{0}(N_{C_{0}})\to N_{C_{0},r}\oplus N_{C_{0},r'}$ is an isomorphism. By Lemma~\ref{lem:nodalsmooth}, whose proof is local along $C_{0}$, the space of conics of $M$ is smooth of
dimension $4n$ at $[C_{0}]$, with tangent space $H^{0}(N_{C_{0}})$. A conic near $C_{0}$ meets $T$ and $T'$ transversally, once each, so $m$ is holomorphic;
its differential at $[C_{0}]$ is the evaluation above, with $N_{C_{0},r}\cong T_{r}T$ and $N_{C_{0},r'}\cong T_{r'}T'$. So $m$ is a chart. Near $r$ the conics $C(m)$ are graphs
over the $x_{0}$-axis, which defines $\Phi$.

(2)(b) \emph{The metric.} Wherever $C(m)$ is transverse and non-jumping, the proof of (1) applies with $(z_{1},w_{1})=(u,u')$: it uses at $u$ only that $C(m)$ is
smooth and transverse there. So $g^{\CC}=2\sum g_{ab}\,du'^{a}du^{b}$, with $g_{ab}$ given by the formula for $g_{ab}$ in the proof of (1), and the components of $\mathrm{Rm}$ in the coordinates
$m$ are polynomials in the derivatives of order $\le4$ of $\Phi$ at $(m;0)$ and of order $\le3$ of $\theta$ at $u$, divided by powers of
$\theta(\partial_{x_{0}}+\partial_{x_{0}}\Phi(m;0))$, of the determinant of $\pi$ on $\{dx_{0}=0\}$, and of $\det(g_{ab})$ with these denominators cleared.

\emph{Values at $z$.} At a smooth point $z$ of $C(m)$ the differential $d_{z}F$ has rank $2n$. It identifies $N_{C(m),z}$ with $\CC^{2n}$, and the tangent line is
spanned by the vector $\tau_{z}$ whose components are the signed $2n$-minors of $d_{z}F$. The variation of the conics along $\partial/\partial m_{j}$ has normal value
$-\partial_{m_{j}}F(m,z)$ at $z$. Let $B_{z}$ be the matrix $d_{z}F(b_{1},\dots,b_{2n})$, $b_{i}$ the frame of $L^{\perp}$; its determinant is nonzero exactly when
$\theta(\tau_{z})\neq0$, and then $B_{z}^{-1}$ is the identification of $N_{C(m),z}$ with $L^{\perp}_{z}$ in the frame $b$. So the value at $z$ of a field
$X=\sum X^{j}\partial_{m_{j}}$ is $\mathrm{ev}_{z}(X)=-B_{z}^{-1}\sum_{j}X^{j}\,\partial_{m_{j}}F(m,z)$. If it vanishes, its derivative along $\tau_{z}$ is
$-B_{z}^{-1}\sum_{j}X^{j}\,\partial_{m_{j}}(d_{z}F)(\tau_{z})$, since at a zero the derivative of a section of $N_{C(m)}$ does not depend on the trivialization.

\emph{Reading off $W^{\CC}$.} The proof of (1) reads $W^{\CC}([C],z)$ off $\mathrm{Rm}$ by means of a second point $w\neq z$ of $C$, and uses of $w$ only that $X_{w}(v)$ is
a field vanishing at $w$ with value $v$ at $z$. We take $w=u$. The fields vanishing at $u$ are the combinations of the $\partial_{u'^{a}}$; if $E_{z}$ is the matrix of
their values at $z$, then $X_{u}(v)=\sum c^{a}\partial_{u'^{a}}$ with $c=E_{z}^{-1}v$. The fields vanishing at $z$ are
$X=\sum c^{b}\partial_{u^{b}}+\sum d^{a}\partial_{u'^{a}}$ with $d=-E_{z}^{-1}\mathrm{ev}_{z}(\sum c^{b}\partial_{u^{b}})$; their derivative along $\tau_{z}$ is a
linear function $H_{z}c$, and $X_{z}(v)$ is the field with $c=H_{z}^{-1}v$. Then, as in the proof of (1),
\[
W^{\CC}\big([C(m)],z\big)(v,v,v,v)=\theta(\tau_{z})^{2}\,\mathrm{Rm}\big(X_{z}(v),X_{u}(v),X_{z}(v),X_{u}(v)\big)
\]
up to the universal constant; the right side does not change when $\tau_{z}$ is rescaled. Polarizing, and clearing the denominators
$\theta(\partial_{x_{0}}+\partial_{x_{0}}\Phi)$, $\det\pi$, $\det(g_{ab})$, $\det B_{z}$, $\det E_{z}$ and $\det H_{z}$, gives $\Pi$ and $\Delta$; only derivatives
of order $\le4$ occur. Every smooth conic $C(m)$ is transverse and non-jumping, for $m$ near $m_{0}$: by Proposition~\ref{prop:thetamult} the transversality section is
$\sigma_{\theta}=t\cdot h$ with $h([C_{0}])=c\,d\theta(u_{1},u_{2})\neq0$, as $C_{0}$ is balanced (Proposition~\ref{prop:balcrit}), and $j(C_{0})=0$ with $j$
upper semicontinuous (Remark~\ref{rem:jinv}). So the formula holds wherever $t(m)\neq0$ and $\Delta(m,z)\neq0$. Where $C(m)$ is smooth and $z\neq u$, which
holds on $U$ if $U$ is small, none of these denominators vanishes: with $N_{C(m)}\cong\OO(1)^{2n}$, $E_{z}$ is
invertible because $H^{0}(N_{C(m)}(-u))\to N_{C(m),z}$ is an isomorphism, and $H_{z}$ because $H^{0}(N_{C(m)}(-u-z))=H^{0}(N_{C(m)}(-2z))=0$. These points are dense in
$\overline F_{U}$, which is smooth and connected near $(m_{0},p)$ (Lemma~\ref{lem:nodalsmooth}), so $\Delta\not\equiv0$ there.

(3)(a) Let $Y_{N}$ be the point $m_{0}$ with the structure sheaf $\OO_{Y}/\mathfrak m^{N+1}$. The universal conic of $M$ over $Y_{N}$ is a flat family with central
fibre $C_{0}$. On it $I\,\OO\subset\mathfrak m\,\OO$, as one sees modulo $\mathfrak m$, so $I^{N+1}\OO=0$: it is a subspace of $Y_{N}\times(C_{0},\OO_{M}/I^{N+1})$. Through the
isomorphism of neighbourhoods it becomes a flat family of conics of $M'$ over $Y_{N}$, whose intersections with $T$ and $T'$ are computed over $V$, where the
isomorphism is the identity; they are given by the coordinates $m$. By (2)(a) for $M'$ and the universal property of the Douady space of $M'$ \cite{Dou66}, the
family is the restriction to $Y_{N}$ of the universal conic of $M'$. So the two universal conics coincide over $Y_{N}$: over $V$,
$I_{\mathcal C}+\mathfrak m^{N+1}\OO=I_{\mathcal C'}+\mathfrak m^{N+1}\OO$. Near $r$ both are graphs, so $\Phi-\Phi'\in\mathfrak m^{N+1}$. Near $p$, after shrinking
$Y$ and $U$, write $F_{i}=F'_{i}+R_{i}$ with $F'_{i}\in I_{\mathcal C'}$ and $R_{i}\in\mathfrak m^{N+1}\OO$. The $F'_{i}$ generate $I_{\mathcal C'}$: they agree
with the $F_{i}$ over $m_{0}$, so $\{F'=0\}$ has fibre $C_{0}\cap U$ there; it is cut out by $2n$ equations in the $(6n+1)$-dimensional $Y\times U$ and its fibres
over $Y$ have dimension $\le1$ near $(m_{0},p)$, so it is a complete intersection, flat over $Y$ \cite[Thm.~23.1]{Mat86}; and $\OO/(F')\to\OO_{\mathcal C'}$ is a
surjection of flat $\OO_{Y}$-modules that is an isomorphism over $m_{0}$, hence an isomorphism by Nakayama's lemma. Finally, derivatives in $m$ of order
$\le j$ map $\mathfrak m^{N+1}\OO$ into $\mathfrak m^{N+1-j}\OO$, derivatives in $z$ preserve it, and $\theta$ and the frames are the same; as $\Pi$ and $\Delta$
involve derivatives of order $\le4$, $\Pi-\Pi'$ and $\Delta-\Delta'$ lie in $\mathfrak m^{N-3}\OO$.

(3)(b) By Lemma~\ref{lem:nodalsmooth}, $\overline F'_{U}$ is smooth at $b_{0}$, with coordinates $(x,y,\vec w)$ in which the projection to $Y$ is $(xy,\vec w)$ and
the node section is $\{x=y=0\}$. The ideal generated on $\overline F'_{U}$ by the $2n$-minors of $d_{z}F'$ is the first Fitting ideal of
$\Omega_{\overline F'_{U}/Y}$, which in these coordinates is $(x,y)$, the ideal of the node section; its pull-back to $\widetilde F'$ is $\OO(-E')$. So along
$\gamma'$ all these minors have order $\ge1$ and one of them, $\mu'$, has order $1$. Let $\mu$ be the corresponding minor of $d_{z}F$. By (a),
$\mu(m(s),z'(s))\equiv\mu'(\gamma'(s))$ modulo $s^{N+1}$, so it has order $1$. Fix the coordinate of $z$ not involved in $\mu$ to be that of $z'(s)$, and solve
$G(z):=F(m(s),z)=0$ for the other $2n$ coordinates. Here $G(z'(s))=(F-F')(\gamma'(s))\in s^{N+1}\CC[[s]]^{2n}$, and the Jacobian of $G$ at $z'(s)$ has
determinant of order $\delta=1$. \emph{Newton's lemma}: if $G(z')\in s^{K}$ and the Jacobian $J$ of $G$ at $z'$ has determinant of order $\delta$ with $K>2\delta$,
then $G$ has a zero $z$ with $z-z'\in s^{K-\delta}$. Indeed, put $z=z'+s^{K-\delta}\zeta$ and $J^{-1}=s^{-\delta}A$, with $A$ over $\CC[[s]]$; by Taylor's formula
the equation becomes $\zeta=-s^{-K}AG(z')-s^{K-2\delta}AH(\zeta)$ with $H$ a power series, which has a solution by successive approximation in the $s$-adic
topology, as $K-2\delta\ge1$. With $K=N+1$ and $\delta=1$ we obtain $\gamma$, with $z-z'\in s^{N}$.

The minors of $d_{z}F$ along $\gamma$ differ from those of $d_{z}F'$ along $\gamma'$ by elements of $s^{N}\CC[[s]]$, so they generate $(s)$: the lift of $\gamma$
meets $E$ transversally. At points over $m_{0}$ all derivatives of $F-F'$ of order $\le N$ vanish. As $N\ge2$, $\overline F_{U}$ and $\overline F'_{U}$ have the same
tangent space at $b_{0}$, and so have their node sections, which are cut out by $F$, resp.\ $F'$, and the $2n$-minors of $d_{z}F$, resp.\ $d_{z}F'$; and
$\dot\gamma(0)=\dot\gamma'(0)$. The lift of an arc that is transverse to a smooth centre ends at the class of its tangent vector in the projectivized normal
space. So the two lifts end at the same point.

(3)(c) $\Pi(\gamma)-\Pi'(\gamma')=(\Pi-\Pi')(\gamma')+\big(\Pi(\gamma)-\Pi(\gamma')\big)$. The first term lies in $s^{N-3}\CC[[s]]$ by (a), as $\mathfrak m$
pulls back into $(s)$ along $m(s)$, and the second in $s^{N}\CC[[s]]$ by (b); likewise for $\Delta$. As $d_{\Delta}<N-3$, $\Delta'(\gamma')$ has order $d_{\Delta}$ too,
and if $W'=\Pi'(\gamma')/\Delta'(\gamma')$ is a power series, then
\[
\frac{\Pi(\gamma)}{\Delta(\gamma)}-W'=\frac{\big(\Pi(\gamma)-\Pi'(\gamma')\big)-W'\big(\Delta(\gamma)-\Delta'(\gamma')\big)}{\Delta(\gamma)}\in s^{N-3-d_{\Delta}}\,\CC[[s]],
\]
since $\Pi'(\gamma')=W'\Delta'(\gamma')$. For $N\ge d_{\Delta}+4$ the exponent is $\ge1$.
\end{proof}

\begin{lemma}[The dilated finite models converge near $C_{1}$]
\label{lem:jetconv}
In the situation of Lemma~\ref{lem:normalN}, let $g$ be the germ of $Z$ in the charts $\psi$; let $g^{(N_{0})}$ be $Z_{0}$ reglued on $\{0<|\alpha|<\epsilon\}$ and
its mirror by $\exp X_{h_{2}}\circ\cdots\circ\exp X_{h_{N_{0}}}$, the \emph{finite model} of order $N_{0}$, and $g^{(N_{0})}_{\infty}$ be $Z_{0}$ reglued by the
weight-$0$ parts $h_{j,0}$, $j\le N_{0}$, alone. Fix $y=([C_{1}],z')$ as in Lemma~\ref{lem:dilE} and a second point $w'\in C_{1}$, both in the node chart.
Near the far end of $A$, $C_{1}=\{\beta=\gamma=\chi=0,\ \rho=-1\}$ (as $R=-aQ=-1$ on $C_{1}$), and the regluings, their inverses and their finite compositions
are defined on a common neighbourhood of the relevant part of $C_{1}$, holomorphically in $s$ for $|s|$ small (see the proof).
\begin{enumerate}
\item $\delta_{s}^{*}g^{(N_{0})}$ is $Z_{0}$ reglued on the \emph{fixed} region by $\prod_{j}\exp X_{\sum_{w}s^{w}h_{j,w}}$; it depends analytically on $s$ and equals
$g^{(N_{0})}_{\infty}$ at $s=0$.
\item $X_{h_{j,0}}\in\beta^{j}T$ with $j\ge2$, so the gluing of $g^{(N_{0})}_{\infty}$ is the identity modulo $I^{2}_{C_{1}}$; thus $g^{(N_{0})}_{\infty}$ agrees with $Z_{0}$
to first order along $C_{1}$, and to second order up to isomorphism, since $H^{1}(C_{1},L\otimes S^{3}N^{*})=0$. So $C_{1}$ is a conic of $g^{(N_{0})}_{\infty}$ with
normal bundle $\OO(1)^{2n}$, and $W^{\CC}_{g^{(N_{0})}_{\infty}}([C_{1}],z')$ is finite; for $N_{0}\ge4$ it is independent of $N_{0}$, by Lemma~\ref{lem:findet}(1): the $h_{j}$, $j\le N_{0}$, do not depend on
$N_{0}$, and the terms with $j\ge5$ do not change the neighbourhood of order $4$ of $C_{1}$.
\item For $|s|$ small, $M_{s}=\delta_{s}^{*}g^{(N_{0})}$ has a unique conic $C_{s}$ near $C_{1}$ through $z'$ and $w'$. It depends holomorphically on
$s$, with $C_{s}=C_{1}$ at $s=0$, and $W^{\CC}_{M_{s}}([C_{s}],z')$ is holomorphic in $s$, with value $W^{\CC}_{g^{(N_{0})}_{\infty}}([C_{1}],z')$ at $s=0$. For
$s\neq0$, $\delta_{s}C_{s}$ is a conic of $g^{(N_{0})}$ through $\delta_{s}z'$. The arc $\gamma'(s)=([\delta_{s}C_{s}],\delta_{s}z')$ in the incidence variety of
$g^{(N_{0})}$ extends holomorphically to $s=0$, with $\gamma'(0)=([C_{0}],p)$; its lift meets the exceptional divisor transversally at the point over $([C_{0}],p)$ with
smoothing ratio $a':Q(z')$; and
\[
W^{\CC}_{g^{(N_{0})}}\big(\gamma'(s)\big)=W^{\CC}_{M_{s}}([C_{s}],z')\qquad(s\neq0).
\]
\end{enumerate}
\end{lemma}

\begin{proof}
\emph{Common domains.} The $h_{j}$ are polynomial in $\rho$, so the fields $X_{h_{j,w}}$ are holomorphic on a neighbourhood of $C_{1}$ in the gluing annuli; this alone
does not make their time-one flows defined there. By (2) below, the weight-zero fields $X_{h_{j,0}}$, $j\ge2$, lie in $\beta^{j}T\subset I_{C_{1}}T$ and vanish along
$C_{1}$. Fix a relatively compact part $K$ of the gluing annuli. A holomorphic field vanishing on the compact set $K\cap C_{1}$ has its flow defined up to time one,
in both directions, on a neighbourhood of $K\cap C_{1}$, with values in a given neighbourhood. For a fixed $N_{0}$ only finitely many flows are composed, and the
positive-weight terms $s^{w}h_{j,w}$, $w\ge1$, are holomorphic perturbations that tend to zero with $s$. Shrinking the neighbourhood and the disc $|s|<s_{0}$
finitely many times, all the maps $\prod_{j}\exp X_{\sum_{w}s^{w}h_{j,w}}$, their inverses and the transition maps of the reglued manifold are defined on common
domains and depend holomorphically on $(s,\cdot)$.

(1) Each $h_{j,w}$, $w\ge0$, is defined on the whole punctured neighbourhood, and conjugation by $\delta_{s}$ replaces $X_{h_{j,w}}$ by $X_{s^{w}h_{j,w}}$. (2)
$X_{h}$ involves first derivatives of $h=\alpha^{-1}\beta^{m}e'_{L}$, $|m|=j+1$, so $X_{h}\in\beta^{j}T$ and $\exp X_{h}\equiv\mathrm{id}$ modulo $I^{2}_{C_{1}}$; the
second-order neighbourhoods then differ by a class in $H^{1}(C_{1},L\otimes S^{3}N^{*})=0$ (\S\ref{sec:autsheaves}).

(3) By (1), the $M_{s}$, $|s|<s_{0}$, are the fibres of a complex manifold $\mathcal M$ over a disc, obtained by gluing pieces of
$Z_{0}\times\{|s|<s_{0}\}$ by maps holomorphic in all variables; $\theta$ and the charts at $z'$ and $w'$, in the node chart away from the regluing, do not depend on
$s$. By (2), $C_{1}$ is a transverse non-jumping conic of $M_{0}$, so $H^{0}$ and $H^{1}$ of $N_{C_{1}}(-z'-w')\cong\OO(-1)^{2n}$ vanish, and by the implicit
function theorem for the relative Douady space of $\mathcal M$ over the disc the conics of $M_{s}$ near $C_{1}$ through $z'$ and $w'$ form a single holomorphic
family $C_{s}$, transverse and non-jumping for $s$ small. The two-point families of $(M_{s},C_{s};z',w')$ in the fixed charts depend holomorphically on $s$, so by
Lemma~\ref{lem:findet}(1) so does $W^{\CC}_{M_{s}}([C_{s}],z')$, which at $s=0$ is the value for $g^{(N_{0})}_{\infty}$.

For $s\neq0$, $\delta_{s}$ is the identity map between the charts of $M_{s}=\delta_{s}^{*}g^{(N_{0})}$ and of $g^{(N_{0})}$, a contactomorphism. So
$\delta_{s}C_{s}$ is a conic of $g^{(N_{0})}$, depending holomorphically on $s\neq0$, and $W^{\CC}_{g^{(N_{0})}}([\delta_{s}C_{s}],\delta_{s}z')=W^{\CC}_{M_{s}}([C_{s}],z')$
by Lemma~\ref{lem:dilcov}. As $C_{s}$ is $O(s)$-close to $C_{1}$, and $\delta_{s}$ multiplies the coordinates normal to $C_{1}$, in the node chart and at the far ends, by
non-negative powers of $s$, $\delta_{s}C_{s}$ is $O(s)$-close to $\delta_{s}C_{1}=C_{s^{2}}$, so it tends to $C_{0}$. Its far-point coordinates (Lemma~\ref{lem:findet}(2)(a)) are holomorphic for $s\neq0$ and
bounded, so they extend holomorphically to $s=0$, with value $m_{0}$; and $\delta_{s}z'=(sa',s/a',0,\dots,0)$ tends to $p$.

Near $p$, $\delta_{s}C_{s}$ differs from $C_{s^{2}}$ by $O(s^{2})$, because $C_{s}$ is $O(s)$-close to $C_{1}$ and every node coordinate has weight $\ge1$; so the
node moves by $O(s^{2})$, while the marked point is exactly $\delta_{s}z'$. In the coordinates $(x,y,\vec w)$ of the incidence variety near $([C_{0}],p)$
(\S\ref{sec:resolution}), $x$ is $a\circ\bar p$ minus the offset of the node: a correction that vanishes on the model family $C_{\tau}$ and is $O(s^{2})$ along
our arc, which is $O(s^{2})$-close to the model arc $(C_{s^{2}},\delta_{s}z')$; likewise for $y$ and $Q\circ\bar p$. Since $a(\delta_{s}z')=sa'$ and
$Q(\delta_{s}z')=s/a'$ exactly, $x=sa'+O(s^{2})$ and $y=s/a'+O(s^{2})$ along the arc. As $\operatorname{ord}_{E}s=1$, the arc meets the exceptional divisor
transversally at the point over $([C_{0}],p)$ with ratio $a':Q(z')$, as in Lemma~\ref{lem:dilE}.
\end{proof}

\subsection{Regularity along $E$ and across the nodal divisor}
\label{sec:Ethm}

\begin{theorem}[$W^{\CC}$ is regular along $E$]
\label{thm:Eregular}
Let $Z$ be a compact contact Fano manifold, and $\mathcal D$ a component of the nodal divisor of $\MCbar$ whose generic member $C_{0}=A\cup A'$ has standard
branches and balanced normal bundle. Then $W^{\CC}$ is regular at the generic point of the component of $E$ over $\mathcal D$. Its value at the point over
$([C_{0}],p)$ with smoothing ratio $a':Q(z')$, $Q(z')=1/a'$, is $W^{\CC}_{g^{(N_{0})}_{\infty}}([C_{1}],z')$ for any $N_{0}\ge4$, computed from the weight-$0$ part of
the germ of $Z$ along $C_{0}$.

If $\mathrm{Pic}\,Z=\ZZ L$ and $n\ge3$, these hypotheses hold for every component of the nodal divisor whose nodes fill $Z$ and whose branch pairs come
from components $X_{i},X_{j}$ of the VMRT with $\omega(\hat X_{i},\hat X_{j})\not\equiv0$; in particular hypothesis (b) of Definition~\ref{def:abc} holds.
\end{theorem}

\begin{proof}
\emph{Choices depending only on $Z$.} Fix a member $C_{*}$ of $\mathcal D$ with standard branches and balanced normal bundle, and choose $T$, $T'$, $U$,
the frames and $F$ as in Lemma~\ref{lem:findet}(2) for $(Z,C_{*})$. They serve for every member $C_{0}$ of $\mathcal D$ near $C_{*}$, with $r$, $r'$ the points
where $T$, $T'$ meet its branches. Let $\Pi$ and $\Delta$ be the functions of Lemma~\ref{lem:findet}(2)(b), and $E_{\mathcal D}$ the part over $Y\times U$ of the
component of $E$ over $\mathcal D$; it is irreducible after shrinking, being a $\PP^{1}$-bundle over the node section. With $y_{E}$ a local equation of $E$, write
$\sigma^{*}\Delta=y_{E}^{d_{\Delta}}\,v$ and $\sigma^{*}\Pi=y_{E}^{d_{\Pi}}v'$, where $v|_{E}$ and $v'|_{E}$ do not vanish identically, the order $d_{\Pi}$ of the vector
$\Pi$ being the least order of its components (if $\Pi\equiv0$ there is nothing to prove); let
$S=V(v)\cup V(v')\subset E_{\mathcal D}$ be the \emph{union} of the zero set of $v$ and the common zero set of the components of $v'$ on
$E_{\mathcal D}$, a proper analytic subset. The numbers $d_{\Delta}$, $d_{\Pi}$ and the set $S$ depend only on these choices.

\emph{The point.} The fibres of $E_{\mathcal D}$ over the node section are the lines of smoothing ratios, and only over a proper analytic subset of
$\mathcal D$ is the fibre contained in $S$. Fix a member $C_{0}$ of $\mathcal D$ near $C_{*}$ outside it, with standard branches and balanced normal bundle (open
conditions), and $a'\in\CC^{*}$ such that the point $e_{0}$ over $([C_{0}],p)$ with ratio $a':Q(z')$ is not in $S$.

\emph{The finite model.} Let $N_{0}\ge d_{\Delta}+4$, and fix the charts of Lemma~\ref{lem:normalN} along $C_{0}$ with this $N_{0}$, the discs about the far ends avoiding
$r$ and $r'$. After shrinking $T$, $T'$ and $U$ about $r$, $r'$ and $p$, they lie in the neighbourhood $U_{K}$ of $K$. The finite model $g^{(N_{0})}$ of
Lemma~\ref{lem:jetconv} and $Z$ both contain $U_{K}$, on which they coincide as contact manifolds, and their neighbourhoods of order $N_{0}$ along $C_{0}$ are
isomorphic by the identity of the charts, which over $U_{K}$ is the identity of the whole thickening, since the transitions $G\circ R_{N_{0}}$ and $G$ agree modulo
$I^{N_{0}+1}$ (Lemma~\ref{lem:normalN}). So
Lemma~\ref{lem:findet}(3) applies, with $M=Z$, $M'=g^{(N_{0})}$, $N=N_{0}$ and $V=U_{K}$.

\emph{Comparison.} Let $\gamma'$ be the arc of Lemma~\ref{lem:jetconv}(3). It ends transversally at the point of the exceptional divisor of $g^{(N_{0})}$ over
$([C_{0}],p)$ with ratio $a':Q(z')$, and $W^{\CC}_{g^{(N_{0})}}$ is holomorphic along it, with value $W^{\CC}_{g^{(N_{0})}_{\infty}}([C_{1}],z')$ at $s=0$ (holomorphy does not depend on the
holomorphic frames, and values are read in those of Lemma~\ref{lem:dilcov}).
Lemma~\ref{lem:findet}(3)(b) gives a formal arc $\gamma(s)=(m(s),z(s))$ in the incidence variety of $Z$ whose lift meets $E$ transversally, at the
corresponding point. That point is $e_{0}$: the computation at the end of the proof of Lemma~\ref{lem:jetconv}(3) applies to $\gamma$, since near $p$ the conic
$C(m(s))$ of $Z$ agrees with $\delta_{s}C_{s}$ modulo $s^{N_{0}+1}$ (Lemma~\ref{lem:findet}(3)(a)), hence is $O(s^{2})$-close to $C_{s^{2}}$, and
$z(s)=\delta_{s}z'+O(s^{N_{0}})$.

\emph{Regularity.} As $e_{0}\notin S$ and $\gamma$ meets $E$ transversally at $e_{0}$, $\Delta(\gamma(s))$ has order exactly $d_{\Delta}$ and $\Pi(\gamma(s))$
order exactly $d_{\Pi}$. By Lemma~\ref{lem:findet}(3)(c), $W^{\CC}_{Z}(\gamma(s))=\Pi(\gamma)/\Delta(\gamma)$ is a power series in $s$. Its order is
$d_{\Pi}-d_{\Delta}$, so $d_{\Pi}\ge d_{\Delta}$, and $\sigma^{*}W^{\CC}=y_{E}^{d_{\Pi}-d_{\Delta}}v'/v$ is holomorphic on $E_{\mathcal D}\smallsetminus S$, which contains
the generic point of the component of $E$.

\emph{The value.} By Lemma~\ref{lem:findet}(3)(c) again, the value of $W^{\CC}$ at $e_{0}$ is the constant term $W^{\CC}_{g^{(N_{0})}_{\infty}}([C_{1}],z')$. This holds for all $a'\in\CC^{*}$ outside a finite set, and, by continuity in $a'$, at every point of the
fibre over $([C_{0}],p)$ with $a'\in\CC^{*}$ at which $W^{\CC}$ is regular; by Lemma~\ref{lem:jetconv}(2) the value does not depend on $N_{0}\ge4$.

This argument is local along $C_{0}$: it proves the regularity along $E$ for every germ of a contact manifold along a balanced nodal conic with standard
branches, and does not use the rationality of $W^{\CC}$ (Corollary~\ref{cor:rational}).

For the last paragraph: by Corollary~\ref{cor:balgeneric} the generic member of such a component is balanced, and its branches, lines through a
general point, are standard. For general $z$ the generic points of the components of $E_{z}$ are generic points of $E$. A component of $E_{z}$
dominating $\PP(L^{\perp}_{z})$ comes from a component $X_{i}\times X_{j}$ of branch pairs that is not of the form $X_{i}\times X_{i}$ with $X_{i}$ linear,
since the secant lines of a linear $\PP^{n-1}$ stay inside it; so $\omega(\hat X_{i},\hat X_{j})\not\equiv0$ by Proposition~\ref{prop:secantgeneral}(2).
\end{proof}

\begin{remark}[What the value along $E$ depends on]
\label{rem:Elimit}
For $n=1$ there are no weight-$0$ classes, so $g^{(N_{0})}_{\infty}=Z_{0}$ and $W^{\CC}|_{E}$ is the value for the homogeneous model, for every germ. For
$n\ge2$ the weight-$0$ classes $\alpha^{-1}\beta^{m}e'_{L}$, $|m|=j+1\ge3$, are as many as the Taylor coefficients of order $j+1$ of a function of $n-1$
variables: the coefficients of a generating function of a Legendrian germ through $[A]$ tangent to the linear component of the model. So the values of
$W^{\CC}$ along $E$ are governed by higher-order data of the VMRT at the two branch directions. Only the weight-$0$ parts of $h_{2},h_{3},h_{4}$
(Lemma~\ref{lem:normalN}, with the choices fixed there), the classes with $|m|=3,4,5$, enter the value; the value itself does not depend on those
choices, since $W^{\CC}$ is regular along $E$ and therefore equal to the limit along any arc ending at the given point. Whatever these data are, along a component of $E_{z}$ at whose generic point
$\widetilde\mu_{z}$ is a local isomorphism the value obeys the boundary law of Corollary~\ref{cor:symlaw}(1): it is $\frac1{24}\partial^{4}$ of its own diagonal.
\end{remark}

\begin{remark}[Why the divisor $D'_{z}$ is different]
\label{rem:Dprime}
In a germ along a nodal conic nothing prevents $W^{\CC}$ from having a pole along the nodal divisor at points $([C],z)$ with $z$ on a branch away from
the node. Such points are
$\delta_{s}$-images of points escaping to the far end of $C_{1}$, where the regluing sits, so the dilation argument does not apply there, consistently
with Example~\ref{ex:bl0}: those points form the divisor $D'_{z}$, which $\widetilde\mu_{z}$ contracts into $X_{z}\subset\PP(T_{z}Z)$, and they impose no condition in the
Hartogs reduction.
\end{remark}

Theorem~\ref{thm:Eregular} can be read back on $\MCbar$, after a double cover: there the curvature extends across the nodal divisor, in the frame of $E_{\CC}$
obtained by evaluating at a point of the neck.

\begin{corollary}[The curvature across the nodal divisor]
\label{cor:Rnodal}
Let $\mathcal D$ be as in Theorem~\ref{thm:Eregular}, $C_{0}=A\cup A'$ a generic member with node $p$, and $t$ the smoothing parameter, so that $\mathcal D=\{t=0\}$
near $[C_{0}]$. On the double cover $t=s^{2}$ of a neighbourhood of $[C_{0}]$ choose node coordinates $xy=t$ and, for $a\in\CC^{*}$ generic, the point
$z_{a}(C)=(x,y)=(s\sqrt a,\,s/\sqrt a)$ of the neck of $C$; let $b_{1},\dots,b_{2n}$ be a holomorphic frame of $B=L^{\perp}\otimes L^{-1/2}$ near $p$. Then the frame
\[
e_{i}(C):=\mathrm{ev}_{z_{a}(C)}^{-1}\big(b_{i}(z_{a}(C))\big)
\]
extends the pull-back of $E^{*}_{\CC}$ across $s=0$; in it $\omega_{E}$ is holomorphic and nondegenerate, and $R_{E}^{\CC}$ has holomorphic components. So, after the base
change $t=s^{2}$ and in this frame, $R_{E}^{\CC}$ extends holomorphically across the generic point of $\mathcal D$. An extension on the $t$-base would need compatibility with
$s\mapsto-s$ and a lattice on the $t$-base; it is not needed, the Hartogs argument using the quartic $W^{\CC}$ on the blown-up incidence space.
\end{corollary}

\begin{proof}
By \S\ref{sec:Wz}, $R_{E}^{\CC}(e_{i},e_{j},e_{k},e_{l})=W^{\CC}\big([C],z_{a}(C)\big)(b_{i},b_{j},b_{k},b_{l})$ and
$\omega_{E}(e_{i},e_{j})=\omega_{z_{a}(C)}(b_{i},b_{j})$, the Levi form at $z_{a}(C)$, which tends to the nondegenerate Levi form at $p$. As $s\to0$ the points
$([C],z_{a}(C))$ tend to the point of $E$ over $([C_{0}],p)$ with ratio $a$, by the definition of the chart $x=ay$ of the blow-up (\S\ref{sec:resolution}). By Theorem~\ref{thm:Eregular}, $W^{\CC}$ is holomorphic near the generic
point of the component of $E$ over $\mathcal D$, so the components of $R_{E}^{\CC}$ are bounded near $s=0$, hence holomorphic there.
\end{proof}


\section{Integral conics: a complex conformal infinity}
\label{sec:integralII}
\label{sec:integral}

The smooth integral conics form a boundary divisor of $\MCbar$ along which LeBrun's metric $g^{\CC}=G/\sigma_{\theta}$ (Proposition~\ref{prop:levwronsk}) has its
denominator vanishing, so a priori the curvature quartic is not controlled there. They nevertheless impose no condition in the Hartogs reduction, for one of
two reasons and without any hypothesis: either $\widetilde\mu_{z}$ contracts them (Lemma~\ref{lem:intcontract}), or their normal bundle is $\OO(1)^{2n}$,
and then the curvature quartic is regular, indeed zero, along them, although $g^{\CC}$ has a pole of order two there (Proposition~\ref{prop:intregular},
Corollary~\ref{cor:intvanish}). Since
$d_{2}=1$ (Theorem~\ref{thm:tpD}), the second case does not occur, a posteriori, by Theorem~\ref{thm:LSJ} and Remark~\ref{rem:intflat}(3), unless $R_{E}^{\CC}\equiv0$; on the real slice
$R_{E}^{\CC}\equiv0$ would make $M=\HH P^{n}$, which is excluded by $\mathrm{Pic}\,Z=\ZZ L$. This section records what is proved, and what is
not.

\begin{lemma}[Integral conics through a general point]
\label{lem:intdivisor}
Let $D$ be an irreducible component of the closure in $\MCbar$ of the locus of smooth integral conics. Then $D$ is a divisor of $\MCbar$. For a general $z$, the
strict transform $D_{z}\subset\Gamma_{z}$ of the conics of $D$ through $z$ is empty or of dimension $2n-1$, and $\widetilde\mu_{z}(D_{z})\subset\PP(L^{\perp}_{z})$.
\end{lemma}

\begin{proof}
Smooth conics lie in the open set $U^{\circ}$ of Theorem~\ref{thm:thetasection}, and by part (3) of that theorem a smooth conic is integral if and only if
$\sigma_{\theta}$ vanishes at it. The section $\sigma_{\theta}$ of the line bundle $\Lambda$ on the irreducible $U^{\circ}$ is not identically zero, the general member
being transverse, so by Krull's principal ideal theorem every irreducible component of its zero locus has codimension one. The nodal conics form a closed subset
of $U^{\circ}$, so the smooth integral conics are open in the zero locus, and the components of their closure are components of the zero locus. Hence $D$ is a
divisor. Its incidence variety $\bar\pi^{-1}(D)$ has dimension $4n$; if it does not dominate $Z$, $D_{z}=\emptyset$ for general $z$, and otherwise its fibre over a
general $z$ has dimension $4n-(2n+1)=2n-1$. A smooth integral conic through $z$ is smooth at $z$, where $\mu$ is a morphism, so this fibre has a strict transform
in $\Gamma_{z}$ of the same dimension; its tangent direction at $z$ lies in $L^{\perp}_{z}$ by Lemma~\ref{lem:allornothing}.
\end{proof}

\begin{lemma}[Contraction criterion]
\label{lem:intcontract}
Let $z$ be general, $D'$ an irreducible component of $D_{z}$, and $C$ its generic member, $N=N_{C/Z}$. Then $N$ is nef of degree $2n$, $\overline F_{z}$ is
smooth at $[C]$, and:
\begin{enumerate}
\item if $N\cong\OO(1)^{2n}$, then $\widetilde\mu_{z}$ is a local isomorphism at $[C]$ and maps $D'$ onto $\PP(L^{\perp}_{z})$;
\item otherwise $h:=h^{0}(N(-2z))\ge1$ and $\dim\widetilde\mu_{z}(D')\le2n-2$, so that $D'$ is contracted. More precisely, the fibres of
$\widetilde\mu_{z}|_{D'}$ have dimension at least $h-1$, and exactly $h$ if the derivative of $\sigma_{\theta}$ at $[C]$ does not vanish on $H^{0}(N(-z))$;
this holds if $L^{\perp}|_{C}$ is nef, in particular if $h=1$.
\end{enumerate}
\end{lemma}

\begin{proof}
For a general $z$ every rational curve through $z$ of $L$-degree at most two is free \cite[II.3.11]{Kol96} (stated there for a very general point;
as these curves form a bounded family, finitely many components of the space of morphisms are involved, and a general $z$ suffices), so $TZ|_{C}$ is nef, hence so is its quotient $N$,
of degree $-K_{Z}\cdot C-2=2n$. Thus $H^{1}(N(-z))=0$, $\overline F_{z}$ is smooth at $[C]$ with tangent space $H^{0}(N(-z))$ (Lemma~\ref{lem:dmu}), and, $C$
being smooth at $z$, $\Gamma_{z}\to\overline F_{z}$ is an isomorphism near $[C]$. A nef bundle of rank and degree $2n$ on $\PP^{1}$ other than $\OO(1)^{2n}$ has a
summand of degree $\ge2$, so $h=0$ exactly in case~(1), and $h=1$ exactly when $N\cong\OO(2)\oplus\OO(1)^{2n-2}\oplus\OO$.

\emph{The derivative of $\sigma_{\theta}$.} At the integral conic $C$ the derivative of $\sigma_{\theta}$ is an intrinsic linear form
$d\sigma_{\theta}\colon H^{0}(N)\to\Lambda_{[C]}$. For $\xi\in H^{0}(N(-z))$ choose a lift $\tilde\xi$ to $TZ|_{C}$ near $z$ with $\tilde\xi(z)\in T_{z}C$. In a
trivialization of $L$, $d\sigma_{\theta}(\xi)=\big(\iota_{\tilde\xi}\,d\theta+d(\theta(\tilde\xi))\big)|_{TC}$, a constant section of
$\omega_{C}\otimes L|_{C}\cong\OO$. Since $\theta(TC)=0$, $\theta(\tilde\xi)=\theta(\xi)$ is the image of $\xi$ under $\theta\colon N\to L|_{C}$, a section of
$L|_{C}(-z)$, and $d\theta(\tilde\xi(z),T_{z}C)=0$; evaluating at $z$,
\[
d\sigma_{\theta}(\xi)=d\big(\theta(\xi)\big)_{z},
\]
the derivative at $z$ of a section vanishing there. In particular $H^{0}(N(-2z))\subset\ker d\sigma_{\theta}$. If $L^{\perp}|_{C}$ is nef, so is
$Q:=L^{\perp}|_{C}/TC$, so $H^{1}(Q(-z))=0$, and the sequence $0\to Q\to N\xrightarrow{\ \theta\ }L|_{C}\to0$ shows that
$H^{0}(N(-z))\to H^{0}(L|_{C}(-z))\cong H^{0}(\OO(1))$ is onto; a section $\xi$ whose image vanishes simply at $z$ has $d\sigma_{\theta}(\xi)\neq0$. If $h=1$,
then $TZ|_{C}\cong TC\oplus N$ (as $H^{1}(N^{*}(2))=0$) has no summand of degree $>2$; the Levi form pairs the summands of $L^{\perp}|_{C}$ in degrees $a$ and
$2-a$, and $L^{\perp}|_{C}\subset TZ|_{C}$ forces $a\le2$, so every $a\ge0$ and $L^{\perp}|_{C}$ is nef.

(1) If $h=0$, then $d\mu_{z}$ is injective on $H^{0}(N(-z))$ (Lemma~\ref{lem:dmu}), so $\widetilde\mu_{z}$ is a local isomorphism at $[C]$; it maps the
$(2n-1)$-dimensional $D'$ into $\PP(L^{\perp}_{z})$ (Lemma~\ref{lem:intdivisor}) with injective differential, hence onto it.
(2) $\ker d\mu_{z}=H^{0}(N(-2z))$ has dimension $h$, and $T_{[C]}D'$ is a hyperplane of $H^{0}(N(-z))$, $D'$ being smooth at its generic point; so
$\ker d\mu_{z}\cap T_{[C]}D'$ has dimension $\ge h-1$. If $d\sigma_{\theta}$ does not vanish on $H^{0}(N(-z))$, then $\sigma_{\theta}|_{\overline F_{z}}$
vanishes to order one at $[C]$, $T_{[C]}D'=\ker d\sigma_{\theta}\cap H^{0}(N(-z))\supset H^{0}(N(-2z))$, and the intersection has dimension $h$. In
characteristic zero the dimension of $\widetilde\mu_{z}(D')$ is the rank of $d(\widetilde\mu_{z}|_{D'})$ at the generic point: at most $(2n-1)-(h-1)\le2n-2$
if $h\ge2$, and $2n-2$ if $h=1$.
\end{proof}

\begin{lemma}[The normal bundle of a smooth integral conic]
\label{lem:intnormal}
Let $C$ be a smooth integral conic, $TC^{\perp}\subset L^{\perp}|_{C}$ the Levi-orthogonal of $TC$, and $W_{C}:=TC^{\perp}/TC$, a bundle of rank
$2n-2$ and degree $2n-2$ with a nondegenerate $L|_{C}$-valued symplectic form, so that $W_{C}\cong W_{C}^{*}\otimes L|_{C}$.
\begin{enumerate}
\item There is a canonical isomorphism $TZ|_{C}/TC^{\perp}\cong J^{1}(L|_{C})$ with the bundle of $1$-jets of sections of $L|_{C}$; in particular
$TZ|_{C}/TC^{\perp}\cong\OO(1)^{2}$.
\item If $W_{C}$ has no summand of negative degree, then $N_{C/Z}\cong W_{C}\oplus\OO(1)^{2}$. This holds if $L^{\perp}|_{C}$ is nef, in particular if
$h^{0}(N_{C/Z}(-2))\le1$.
\end{enumerate}
So, when $L^{\perp}|_{C}$ is nef, $C$ is balanced if and only if $W_{C}\cong\OO(1)^{2n-2}$, $H^{0}(N_{C/Z}(-2z))=H^{0}(W_{C}(-2z))$ for $z\in C$, and when
$h^{0}(N_{C/Z}(-2))=1$ the summand $\OO(2)$ of $N_{C/Z}$ lies in $W_{C}$, that is, in the image of $L^{\perp}|_{C}$.
\end{lemma}

\begin{proof}
(1) For $x\in C$ and a local section $h$ of $L$ near $x$ let $X_{h}$ be the infinitesimal contact automorphism with $\theta(X_{h})=h$; in a local
trivialization of $L$, $\iota_{X_{h}}d\theta=dh(R)\,\theta-dh$ with $R$ the Reeb field, so $d\theta(X_{h},T)=-dh(T)$ for $T\in T_{x}C\subset L^{\perp}_{x}$. If
$h(x)=0$ and $dh|_{T_{x}C}=0$, then $X_{h}(x)\in L^{\perp}_{x}$ and $\omega(X_{h}(x),T_{x}C)=0$, so $X_{h}(x)\in TC^{\perp}_{x}$. Hence
$j^{1}_{x}(h|_{C})\mapsto X_{h}(x)\bmod TC^{\perp}_{x}$ is well defined; it is linear and injective by the same two identities, between bundles of rank two, and
it does not depend on the trivialization, $h\mapsto X_{h}$ being the canonical correspondence between sections of $L$ and infinitesimal contact automorphisms.
The jet sequence $0\to K_{C}\otimes L|_{C}\to J^{1}(L|_{C})\to L|_{C}\to0$ is $0\to\OO\to J^{1}\to\OO(2)\to0$ with extension class the Atiyah class of
$L|_{C}$, which is nonzero as $\deg L|_{C}=2$; so $J^{1}(L|_{C})\cong\OO(1)^{2}$.
(2) $N_{C/Z}=TZ|_{C}/TC$ contains $W_{C}$ with quotient $TZ|_{C}/TC^{\perp}\cong\OO(1)^{2}$, and $\mathrm{Ext}^{1}(\OO(1)^{2},W_{C})=H^{1}(W_{C}(-1))^{2}=0$ when
$W_{C}$ has no negative summand. If $h^{0}(N_{C/Z}(-2))\le1$, every summand of $N_{C/Z}$ has degree $\le2$ and at most one has degree $2$; as the
$2n$ degrees sum to $2n$, the splitting type is $(1^{2n})$ or $(2,1^{2n-2},0)$, so $N_{C/Z}$ is nef, $TZ|_{C}\cong TC\oplus N_{C/Z}$ (as
$H^{1}(N_{C/Z}^{*}(2))=0$) has no summand of degree $>2$, and $L^{\perp}|_{C}$ is nef by the argument at the end of the second paragraph of the proof of
Lemma~\ref{lem:intcontract}. If $L^{\perp}|_{C}$ is nef, its degrees lie in $[0,2]$ by the self-duality $L^{\perp}\cong L^{\perp*}\otimes L$. The map
$L^{\perp}|_{C}\to TC^{*}\otimes L|_{C}\cong\OO$, $v\mapsto\omega(v,\cdot)|_{TC}$, is onto with kernel $TC^{\perp}$, and it vanishes on the summands of positive
degree; so $TC^{\perp}$ has no negative summand, nor has its quotient $W_{C}$. The final assertions follow, since $h^{0}(\OO(1)^{2}(-2z))=0$ and $W_{C}$ has
degrees in $[0,2]$ symmetric about $1$.
\end{proof}

For a conic on a Legendrian quadric $\Lambda$ (LeBrun's sky conics; \cite{LeBrun89}) one finds $W_{C}=N_{C/\Lambda}\oplus\OO^{n-1}\cong\OO(2)^{n-1}\oplus\OO^{n-1}$,
and for a contact line of $\PP^{2n+1}$, $W_{C}\cong\OO(1)^{2n-2}$. Lemma~\ref{lem:intnormal} is used in the proof of Lemma~\ref{lem:intnonembedded}.

\emph{The normal-bundle frame.} LeBrun's bundle $E^{*}_{\CC}=H^{0}(B|_{C})$ is built from $B|_{C}=L^{\perp}|_{C}\otimes L^{-1/2}|_{C}$, which is not
trivial at an integral conic even when $N_{C/Z}\cong\OO(1)^{2n}$: then $TZ|_{C}\cong TC\oplus N_{C/Z}$, and the argument at the end of the previous proof gives
$L^{\perp}|_{C}\cong\OO(2)\oplus\OO(1)^{2n-2}\oplus\OO$, with $TC$ the summand $\OO(2)$, since $L^{\perp}|_{C}/TC\subset N_{C/Z}$ has no summand of degree
$2$. On a transverse conic $N_{C/Z}\cong L^{\perp}|_{C}$ (\S\ref{sec:conics}), so one may equally use
\[
\check B|_{C}:=N_{C/Z}\otimes L^{-1/2}|_{C},\qquad\check E^{*}|_{[C]}:=H^{0}(C,\check B|_{C}),
\]
which is trivial on every conic with $N_{C/Z}\cong\OO(1)^{2n}$, transverse or not. On a transverse conic the two frames are related by the projection along
$TC$,
\[
\varphi_{C}\colon N_{C/Z}\longrightarrow L^{\perp}|_{C},\qquad\varphi_{C}(v)=v-\frac{\theta(v)}{\sigma_{\theta}},
\]
where for $v\in T_{x}Z$, $x\in C$, $\theta(v)/\sigma_{\theta}$ denotes the vector $u\in T_{x}C$ with $\sigma_{\theta}(u)=\theta(v)$; $\varphi_{C}$ is the
inverse of the composite $L^{\perp}|_{C}\to TZ|_{C}\to N_{C/Z}$, and it acquires a pole at the integral conics.

\begin{proposition}[Regularity along a balanced integral divisor]
\label{prop:intregular}
Let $n\ge2$, let $C_{0}$ be a smooth integral conic of $\MCbar$ with $N_{C_{0}/Z}\cong\OO(1)^{2n}$, $\mathcal D$ the component of the closure of the smooth
integral conics through $[C_{0}]$, and $\mathcal V$ a small neighbourhood of $[C_{0}]$ in $\MCbar$.
\begin{enumerate}
\item $\mathcal V$ is smooth of dimension $4n$, its members are smooth conics with normal bundle $\OO(1)^{2n}$, $\sigma_{\theta}$ vanishes to order one
along $\mathcal D$, $\{\sigma_{\theta}=0\}\cap\mathcal V=\mathcal D\cap\mathcal V$ is smooth, and $\mathcal V\smallsetminus\mathcal D\subset\MC$. The spaces $\check E^{*}$ and $H^{*}_{\CC}$ form holomorphic bundles on $\mathcal V$ with
$T\mathcal V=\check E^{*}\otimes H^{*}_{\CC}$, and on $\mathcal V\smallsetminus\mathcal D$ the projection $\varphi$ identifies $\check E^{*}$ with $E^{*}_{\CC}$, compatibly with
Theorem~\ref{thm:lebrun}.
\item Write $\sigma_{\theta}=s\lambda$ with $\lambda$ a generator of $\Lambda|_{\mathcal V}$. In the frame $\check E^{*}$ the forms $s\,\omega_{E}$ and $s^{2}g^{\CC}$
extend holomorphically over $\mathcal V$. At $[C]\in\mathcal D\cap\mathcal V$, $s\,\omega_{E}$ has rank two, and the coefficient of $s^{-2}$ in $g^{\CC}$ has rank four.
\item $R_{E}^{\CC}$ extends to a holomorphic section of $S^{4}\check E$ over $\mathcal V$, and $W^{\CC}$ is holomorphic on $\bar\pi^{-1}(\mathcal V)$. In particular, for a
general $z$, $W_{z}$ is regular at the generic point of every component of $D_{z}$ whose conics lie in $\mathcal D$.
\item For $([C],z)\in\bar\pi^{-1}(\mathcal D\cap\mathcal V)$, $W^{\CC}([C],z)\in S^{4}A_{[C],z}\otimes L^{2}_{z}\subset V_{z}$, where $A_{[C],z}\subset L^{\perp*}_{z}$ is the
line of forms vanishing on the Levi-orthogonal complement $(T_{z}C)^{\perp_{\omega}}\subset L^{\perp}_{z}$.
\end{enumerate}
\end{proposition}

\begin{proof}
(1) $H^{1}(N_{C_{0}/Z})=0$, so $\MCbar$ is smooth of dimension $4n$ at $[C_{0}]$, and being smooth with balanced normal bundle is an open condition. By
the paragraph preceding the proposition $L^{\perp}|_{C_{0}}$ is nef, so by the proof of Lemma~\ref{lem:intcontract} the derivative of $\sigma_{\theta}$ at
$[C_{0}]$ does not vanish (already on $H^{0}(N(-z))$ for any $z\in C_{0}$): $\{\sigma_{\theta}=0\}$ is a smooth hypersurface near $[C_{0}]$, and
$\sigma_{\theta}$ vanishes on it to order one. Off $\{\sigma_{\theta}=0\}$ the members of $\mathcal V$ are transverse
(Theorem~\ref{thm:thetasection}(3)) with trivial $B$, hence in $\MC$. Shrinking $\mathcal V$ to a polydisc, the universal conic $\bar\pi^{-1}(\mathcal V)\to\mathcal V$
is a holomorphically trivial $\PP^{1}$-bundle \cite{FG65}, so $\bar p^{\,*}L$ has a square root $\mathcal L^{1/2}$ of degree one on the fibres. Let
$\mathcal N$ be the relative normal bundle and $\check B=\mathcal N\otimes\mathcal L^{-1/2}$, trivial on every fibre; then $\check E^{*}=\bar\pi_{*}\check B$ and
$H^{*}_{\CC}=\bar\pi_{*}\mathcal L^{1/2}$ are vector bundles, $\bar\pi^{*}\check E^{*}\cong\check B$, and multiplication gives
$\check E^{*}\otimes H^{*}_{\CC}\cong\bar\pi_{*}\mathcal N=T\mathcal V$. On a transverse conic $\varphi_{C}$ is an isomorphism of bundles, so it induces
$\check E^{*}|_{[C]}\cong H^{0}(B|_{C})=E^{*}_{\CC}|_{[C]}$, compatibly with the factorizations of $T_{[C]}\MC=H^{0}(N_{C/Z})$.

(2) For every $[C]\in\mathcal V$, $\lambda([C])$ is a nonzero element of $H^{0}(\omega_{C}\otimes L|_{C})\cong H^{0}(\OO)$, hence vanishes nowhere on $C$. Let
$\tau_{x}\in T_{x}C\otimes L_{x}^{-1}$ be defined by $\lambda(\tau_{x})=1$; it depends holomorphically on $([C],x)$, and
$\varphi_{C}(v)=v-s^{-1}\theta(v)\,\tau_{x}$. In a local trivialization of $L$ near $x$, $\theta$ is a scalar $1$-form and $d\theta$ a $2$-form whose
restriction to $L^{\perp}$ is the Levi form. For $v,w\in T_{x}Z$, since $\theta(v)\tau_{x}$ and $\theta(w)\tau_{x}$ lie on the line $T_{x}C$,
\[
s^{2}\,d\theta\big(\varphi_{C}v,\varphi_{C}w\big)=d\theta\big(sv-\theta(v)\tau_{x},\,sw-\theta(w)\tau_{x}\big)
=s^{2}\,d\theta(v,w)+s\,(\theta\wedge\beta)(v,w),\qquad\beta:=d\theta(\cdot,\tau_{x}),
\]
the term $\theta(v)\theta(w)\,d\theta(\tau_{x},\tau_{x})$ being zero. So $s\,\omega_{E}=s\,d\theta+\theta\wedge\beta$ on $\check B$ is holomorphic. It is
well defined on $N_{x}=T_{x}Z/T_{x}C$: since $\theta(\tau_{x})=\sigma_{\theta}(\tau_{x})=s$, replacing $v$ by $v+c\tau_{x}$ changes
$(s\,d\theta+\theta\wedge\beta)(v,w)$ by $-cs\,\beta(w)+cs\,\beta(w)=0$. It is
independent of the trivialization where $s\neq0$, hence everywhere. For $e,e'\in\check E^{*}$ the function $(s\,\omega_{E})(e,e')$ on the universal conic is
constant on the fibres over $\mathcal V\smallsetminus\mathcal D$ (\S\ref{sec:lebrunstructure}), hence on all fibres, so $s\,\omega_{E}$ is a holomorphic $2$-form on
$\check E^{*}$. At $[C]\in\mathcal D$ it is $\theta\wedge\beta$, of rank two on $N_{x}$: $\theta$ vanishes on the hyperplane $L^{\perp}_{x}/T_{x}C$, on which
$\beta$ restricts to $\omega(\cdot,\tau_{x})\neq0$. Finally the Wronskian of a basis of $H^{0}(\OO(1))$ vanishes nowhere, so $\omega_{H}/\lambda$ is a
holomorphic symplectic form on $H^{*}_{\CC}$ over $\mathcal V$, and by Proposition~\ref{prop:levwronsk}
\[
g^{\CC}=\frac{\omega_{E}\otimes\omega_{H}}{\sigma_{\theta}}=\frac{(s\,\omega_{E})\otimes(\omega_{H}/\lambda)}{s^{2}},
\]
whose $s^{-2}$-coefficient on $\mathcal D$ is $(\theta\wedge\beta)\otimes(\omega_{H}/\lambda)$, of rank $2\cdot2=4$.

(3) The conics of $\mathcal V$ form a complete family of rational curves with normal bundle $\OO(1)^{2n}$, so $T\mathcal V=\check E^{*}\otimes H^{*}_{\CC}$ is a holomorphic
almost-Grassmannian structure on $\mathcal V$, determined by the curves alone, not by $\theta$. On $\mathcal V^{\circ}=\mathcal V\smallsetminus\mathcal D$ the Levi-Civita
connection of $g^{\CC}=\omega_{E}\otimes(\omega_{H}/\sigma_{\theta})$ is torsion-free and compatible (Theorem~\ref{thm:lebrun}(1)). Lemma~\ref{lem:curvcomp} below, with
$\omega'_{H}=\omega_{H}/\sigma_{\theta}$, gives a holomorphic tensor $\mathcal W$ on $\mathcal V$ with
\[
\omega_{E}\big(\mathcal W(e_{1}\otimes h_{1},e_{2}\otimes h_{2})\,e_{3},\,e_{4}\big)=c\,\frac{\omega_{H}}{\sigma_{\theta}}(h_{1},h_{2})\,R_{E}^{\CC}(e_{1},e_{2},e_{3},e_{4})\qquad\text{on }\mathcal V\smallsetminus\mathcal D.
\]
Substituting $\omega_{E}=(s\,\omega_{E})/s$ and $\omega_{H}/\sigma_{\theta}=(\omega_{H}/\lambda)/s$, the factors $s$ cancel:
\begin{equation}
R_{E}^{\CC}(e_{1},e_{2},e_{3},e_{4})=\frac{(s\,\omega_{E})\big(\mathcal W(e_{1}\otimes h_{1},e_{2}\otimes h_{2})\,e_{3},\,e_{4}\big)}{c\,(\omega_{H}/\lambda)(h_{1},h_{2})}\label{eq:7-1}
\end{equation}
for any $h_{1},h_{2}$ with $(\omega_{H}/\lambda)(h_{1},h_{2})\neq0$. The right side is holomorphic on $\mathcal V$, so $R_{E}^{\CC}$ extends holomorphically, as a section
of $S^{4}\check E$ by continuity. For the curvature quartic, let $\check{\mathrm{ev}}_{z}\colon\check E^{*}\to\check B_{z}=N_{C/Z,z}\otimes L^{-1/2}_{z}$ be
evaluation, an isomorphism for every $[C]\in\mathcal V$, and $\iota_{C,z}\colon L^{\perp}_{z}\to N_{C/Z,z}$ the map induced by $T_{z}Z\to T_{z}Z/T_{z}C$. On transverse
conics $\varphi_{C}^{-1}=\iota_{C}$, so $\check{\mathrm{ev}}_{z}=\iota_{C,z}\circ\mathrm{ev}_{z}$ under $\check E^{*}\cong E^{*}_{\CC}$, and
\[
W^{\CC}([C],z)=(\mathrm{ev}_{z}^{-1})^{*}R_{E}^{\CC}=\iota_{C,z}^{*}\big((\check{\mathrm{ev}}_{z}^{-1})^{*}R_{E}^{\CC}\big),
\]
in which every term is holomorphic in $([C],z)\in\bar\pi^{-1}(\mathcal V)$. For a general $z$ the generic point of a component of $D_{z}$ whose conics lie in
$\mathcal D$ lies in $\bar\pi^{-1}(\mathcal V)$ for a suitable $C_{0}$ (the members of $\mathcal D$ excluded by the genericity conditions form a proper closed
subset, whose incidence variety meets $\bar p^{-1}(z)$, for general $z$, in dimension $<2n-1$), and $\Gamma_{z}\cong\overline F_{z}$ there
(Lemma~\ref{lem:intcontract}).

(4) At $[C]\in\mathcal D$, $s\,\omega_{E}=\theta\wedge\beta$, and \eqref{eq:7-1} vanishes when $e_{4}$ lies in $\mathcal K:=\ker(\theta\wedge\beta)\subset\check E^{*}$, of
codimension two; by symmetry $R_{E}^{\CC}([C])\in S^{4}(\mathcal K^{\perp})$, $\mathcal K^{\perp}\subset\check E$ the plane spanned by $\theta$ and $\beta$ (evaluated at $z$).
Under $\iota_{C,z}^{*}$, $\theta\mapsto0$, because $\theta(L^{\perp}_{z})=0$, and $\beta\mapsto\omega(\cdot,\tau_{z})$, which spans $A_{[C],z}$ since $\tau_{z}$
spans $T_{z}C$.
\end{proof}

The comparison used in the proof of Proposition~\ref{prop:intregular}(3) is the following; no contact geometry enters in it. Let $\mathcal V$ be a
complex manifold with $T\mathcal V=\check E^{*}\otimes H^{*}$, for holomorphic bundles $\check E^{*}$, $H^{*}$ of ranks $2n\ge4$ and $2$ (a holomorphic
almost-Grassmannian structure of type $(2,2n)$ \cite{BE91}). Put $T=T\mathcal V$ and $\mathfrak g_{0}=\big(\mathfrak{gl}(\check E^{*})\oplus\mathfrak{gl}(H^{*})\big)/\CC$,
acting on $T$. A connection on $T$ is \emph{compatible} if it is induced by connections on $\check E^{*}$ and $H^{*}$; its curvature lies in $\Lambda^{2}T^{*}\otimes\mathfrak g_{0}$.
Under $G_{0}=GL(\check E^{*})\times GL(H^{*})$ the module $\Lambda^{2}T^{*}\otimes\mathfrak g_{0}$ contains exactly one copy of $\Lambda^{2}H\otimes(S^{3}\check E\otimes\check E^{*})_{0}$
(it can only come from $S^{2}\check E\otimes\Lambda^{2}H\otimes\mathfrak{sl}(\check E^{*})$); for $R\in\Lambda^{2}T^{*}\otimes\mathfrak g_{0}$ let $\mathcal H(R)$ be its
component there, the \emph{harmonic part} \cite[\S4.1.3]{CS09}.

\begin{lemma}[The curvature comparison]
\label{lem:curvcomp}
In this situation:
\begin{enumerate}
\item If the torsion of some compatible connection on an open set $\mathcal V^{\circ}$ vanishes, and $\mathcal V^{\circ}$ is dense, then every point of $\mathcal V$ has a
neighbourhood carrying a holomorphic torsion-free compatible connection.
\item For torsion-free compatible connections $\nabla,\nabla'$ on an open set, $\mathcal H(R^{\nabla'})=\mathcal H(R^{\nabla})$. So $\mathcal W:=\mathcal H(R^{\nabla})$ is a
holomorphic tensor on $\mathcal V$, the harmonic curvature of the structure.
\item Suppose that on $\mathcal V^{\circ}$ there are symplectic forms $\omega_{E}$ on $\check E^{*}$ and $\omega'_{H}$ on $H^{*}$ such that the Levi-Civita connection of
$g=\omega_{E}\otimes\omega'_{H}$ is compatible, with curvature $\rho_{0}R_{0}+R_{E}$ as in Theorem~\ref{thm:lebrun}(3), $R_{E}\in S^{4}\check E$. Then on $\mathcal V^{\circ}$, for
$e_{i}\in\check E^{*}$, $h_{1},h_{2}\in H^{*}$ and a universal constant $c\neq0$ fixed by the conventions of Theorem~\ref{thm:lebrun}(3),
\[
\omega_{E}\big(\mathcal W(e_{1}\otimes h_{1},e_{2}\otimes h_{2})\,e_{3},\,e_{4}\big)=c\,\omega'_{H}(h_{1},h_{2})\,R_{E}(e_{1},e_{2},e_{3},e_{4}).
\]
\end{enumerate}
\end{lemma}

\begin{proof}
(1) Near a point, the trivial connection $\nabla^{0}$ of local frames of $\check E^{*}$ and $H^{*}$ is compatible; let $T^{0}$ be its torsion. The other compatible
connections are $\nabla^{0}+A$, $A\in T^{*}\otimes\mathfrak g_{0}$, with torsion $T^{0}+\partial A$, where $\partial\colon T^{*}\otimes\mathfrak g_{0}\to\Lambda^{2}T^{*}\otimes T$
is the alternation, a bundle map of constant rank. On $\mathcal V^{\circ}$ the torsion of some compatible connection vanishes, so $T^{0}\in\operatorname{im}\partial$
there, hence everywhere, $\operatorname{im}\partial$ being a closed subbundle. A holomorphic right inverse $\sigma$ of $\partial$ onto its image gives the torsion-free
compatible connection $\nabla^{0}-\sigma(T^{0})$.

(2) $\nabla'=\nabla+A$ with $\partial A=0$. The kernel of $\partial$ is the first prolongation of $\mathfrak g_{0}$, which for type $(2,2n)$, $n\ge1$, is
$\mathfrak g_{1}\cong T^{*}$, $A(X)=[\Upsilon,X]$ \cite[\S4.1.3]{CS09}. Then $R^{\nabla'}-R^{\nabla}=d^{\nabla}A+[A\wedge A]$ is the image of $\nabla\Upsilon\otimes\cdot$ and of
$\Upsilon\otimes\Upsilon$ under $G_{0}$-equivariant maps defined on $T^{*}\otimes T^{*}$. As a $G_{0}$-module $T^{*}\otimes T^{*}=(\check E\otimes\check E)\otimes(H\otimes H)$ has
no constituent of type $\Lambda^{2}H\otimes(S^{3}\check E\otimes\check E^{*})_{0}$, whose $\check E$-part has a dual factor; so $\mathcal H(R^{\nabla'}-R^{\nabla})=0$. By (1)
and this, the local tensors $\mathcal H(R^{\nabla})$ glue to a holomorphic $\mathcal W$.

(3) Take $\nabla$ the Levi-Civita connection. Under $\Sp(n,\CC)\times\Sp(1,\CC)$, $\Lambda^{2}T^{*}\otimes\mathfrak g_{0}$ contains $S^{4}\check E$ (with trivial
$\Sp(1,\CC)$-part) exactly once, inside $S^{2}\check E\otimes\Lambda^{2}H\otimes S^{2}\check E$: the other constituents of $\Lambda^{2}T^{*}=S^{2}\check E\otimes\Lambda^{2}H\oplus
\Lambda^{2}\check E\otimes S^{2}H$ and of $\mathfrak g_{0}\cong S^{2}\check E\oplus\Lambda^{2}_{0}\check E\oplus\CC\oplus S^{2}H$ give, by highest weights, no $S^{4}\check E$ with trivial
$\Sp(1,\CC)$-part. The harmonic module contains such a copy, $S^{4}\check E\subset(S^{3}\check E\otimes\check E)_{0}\cong(S^{3}\check E\otimes\check E^{*})_{0}$ via
$\omega_{E}$. The curvature $R_{E}$ takes values in this copy: its $2$-forms lie in $S^{2}\check E\otimes\Lambda^{2}H$ and its values in $\mathfrak{sp}(\check E^{*})\cong
S^{2}\check E$ \cite{Salamon82}, and $R_{E}(e_{1}\otimes h_{1},e_{2}\otimes h_{2})=c\,\omega'_{H}(h_{1},h_{2})\,R_{E}(e_{1},e_{2},\cdot,\cdot)$, with the last two
slots read as an endomorphism through $\omega_{E}$. So $\mathcal H(R_{E})=R_{E}$. The term $\rho_{0}R_{0}$ is given by the same universal formula as the
curvature of the Levi-Civita connection of the symmetric space $\Sp(n+1,\CC)/\Sp(n,\CC)\times\Sp(1,\CC)$, an open subset of $\Gr_{2}(\CC^{2n+2})$ whose
almost-Grassmannian structure is flat; by (2), applied there to the Levi-Civita connection and to a flat connection of a Grassmannian chart, $\mathcal H(R_{0})=0$.
Hence $\mathcal W=\mathcal H(\rho_{0}R_{0}+R_{E})=R_{E}$, which is the displayed formula.
\end{proof}

\begin{corollary}[The curvature vanishes on the balanced integral divisor]
\label{cor:intvanish}
In the situation of Proposition~\ref{prop:intregular}, $R_{E}^{\CC}$ vanishes on $\mathcal D\cap\mathcal V$, and $W^{\CC}$ vanishes on $\bar\pi^{-1}(\mathcal D\cap\mathcal V)$.
\end{corollary}

\begin{proof}
Let $[C]\in\mathcal D\cap\mathcal V$ and $z\in C$. Since $N=N_{C/Z}\cong\OO(1)^{2n}$, $H^{1}(N(-z))=H^{0}(N(-2z))=0$, so $\overline F_{z}$ is smooth at $[C]$, it agrees
with $\Gamma_{z}$ there, and $\widetilde\mu_{z}$ is a local isomorphism there. This holds for every $z\in C$, not only for a general one: $\Gamma_{z}$ was
constructed in \S\ref{sec:fibre} for a general $z$, but the members of $\mathcal V$ are smooth conics, so near $[C]$ there is no node to blow up and
$\widetilde\mu_{z}$ is the morphism $\mu_{z}$ (Lemma~\ref{lem:dmu}); it maps the integral conics through $z$ near $[C]$ into $\PP(L^{\perp}_{z})$, the other conics near
$[C]$ through $z$ are transverse and non-jumping (Proposition~\ref{prop:intregular}(1)), and $W_{z}$ is holomorphic near $[C]$ (Proposition~\ref{prop:intregular}(3)). So
Corollary~\ref{cor:symlaw}(1) applies at $[C]$. By Proposition~\ref{prop:intregular}(4), $W_{\infty}(v)=\chi(v)\,\omega(\cdot,v)^{4}$ for $v$ spanning the tangent
line of the integral conic, so $w_{\infty}(\xi)=\chi(\xi)\,\omega(\xi,\xi)^{4}=0$ and $W_{\infty}=0$ (already $\iota_{\xi}W_{\infty}=0$). In particular $W^{\CC}([C],z)=0$, for every $z\in C$.

For $R_{E}^{\CC}$: by the proof of Proposition~\ref{prop:intregular}(4), $R_{E}^{\CC}([C])\in S^{4}(\mathcal K^{\perp})$, where the plane $\mathcal K^{\perp}\subset\check E|_{[C]}$
is spanned by $\theta_{x}$ and $\beta_{x}$ for every $x\in C$, and $W^{\CC}([C],z)$ is the image of $R_{E}^{\CC}([C])$ under the fourth power of a map
$\mathcal K^{\perp}\to L^{\perp*}_{z}$ of rank one with kernel $\CC\theta_{z}$. Read $R_{E}^{\CC}([C])$ as a binary quartic on $(\mathcal K^{\perp})^{*}$; then $W^{\CC}([C],z)=0$
says that it vanishes at the point annihilated by $\theta_{z}$. Here $\theta_{z}(e)$ is the value at $z$ of $\theta(e)\in H^{0}(C,L^{1/2}|_{C})\cong H^{0}(\OO(1))$, and
$e\mapsto\theta(e)$ is onto, because $\theta\colon N\to L|_{C}$ is a surjection $\OO(1)^{2n}\to\OO(2)$, whose $2n$ sections of $\OO(1)$ have no common zero and hence span $H^{0}(\OO(1))$; so $z\mapsto[\theta_{z}]$
is injective, and the quartic vanishes at infinitely many points. Hence $R_{E}^{\CC}([C])=0$.
\end{proof}

So the balanced integral divisor is a conformal infinity at which the quaternionic Weyl curvature extends and vanishes; what it carries about the curvature is in the normal
derivative, the term $W_{1}$ of the expansion in the proof of Corollary~\ref{cor:symlaw}(1). Nothing is claimed here for the unbalanced integral conics, which are contracted
(Lemma~\ref{lem:intcontract}(2)).

\begin{corollary}[Integral conics impose no condition]
\label{cor:intharmless}
Let $n\ge2$, $z$ general, and $D'$ an irreducible component of $D_{z}$. Then either $\widetilde\mu_{z}(D')$ has codimension at least two in $\PP(T_{z}Z)$, or
$W_{z}$ is regular at the generic point of $D'$. Hence hypothesis~(c) of Definition~\ref{def:abc} holds along every divisor of smooth integral conics through
$z$.
\end{corollary}

\begin{proof}
If the generic member of $D'$ is not balanced, apply Lemma~\ref{lem:intcontract}(2). If it is, let $\mathcal D$ be the component of the closure of the smooth
integral conics containing the conics of $D'$; balancedness is open on the irreducible $\mathcal D$, so its generic member is balanced, and
Proposition~\ref{prop:intregular}(3) applies.
\end{proof}

\begin{remark}[The balanced integral divisor as a complex conformal infinity]
\label{rem:intflat}
Proposition~\ref{prop:intregular} has the following interpretation, test case and consequence.
\begin{enumerate}
\item For $n=1$ the computation of Proposition~\ref{prop:intregular}(2) applies verbatim and $\theta\wedge\beta$ is nondegenerate: $s^{2}g^{\CC}$ is a
holomorphic conformal metric, nondegenerate along $\mathcal D$, so $\mathcal D$ is a conformal infinity of the Einstein metric $g^{\CC}$, and the analogue of
Proposition~\ref{prop:intregular}(3) is the conformal invariance of the Weyl tensor.
For $n\ge2$ the pole of order two is concentrated on a subspace of rank four, with a simple pole on the remaining $4n-4$ directions: for sections $e,e'$ of
$\check E^{*}$ lying in $\mathcal K=\ker(\theta\wedge\beta)$ along $\mathcal D$, $\theta(e),\beta(e),\theta(e'),\beta(e')$ vanish on $\mathcal D$, so
$(\theta\wedge\beta)(e,e')=O(s^{2})$, and the $s^{-1}$-coefficient of $g^{\CC}$ on $\mathcal K\otimes H^{*}_{\CC}$ is the Levi form of
$(T_{x}C)^{\perp_{\omega}}/T_{x}C$ tensored with $\omega_{H}/\lambda$, which is nondegenerate.
\item In $\HH P^{n}$, $Z=\PP^{2n+1}$, the conics are the lines, the integral ones are the contact lines, all with normal bundle $\OO(1)^{2n}$, and $\MCbar$ is
$\Gr_{2}(\CC^{2n+2})$ with the isotropic planes as $\mathcal D$. For $n=2$ a symbolic computation with the invariant metric of
$\Sp(3,\CC)/\Sp(2,\CC)\times\Sp(1,\CC)$, built from the symplectic form (it is unique up to scale, the isotropy representation being irreducible, so it
is LeBrun's metric whatever its normalization), confirms Proposition~\ref{prop:intregular}(2): along a transverse arc the
metric has a pole of order exactly two with leading coefficient of rank four, the next coefficient is nondegenerate on its kernel, and $\omega_{E}$ has a simple
pole with residue of rank two. This computation was carried out by machine and verified by the author; the proof of
Proposition~\ref{prop:intregular} does not depend on it. Here $R_{E}^{\CC}=0$.
\item If, for a general $z$, some component $D'$ of $D_{z}$ is balanced and $W_{z}$ is constant on the component of $\Gamma_{z}$ containing it, then $W_{z}\equiv0$ there, because $W_{z}$ vanishes on $D'$
(Corollary~\ref{cor:intvanish}). Then $R_{E}^{\CC}=0$: locally $g^{\CC}$ is the complexification of the metric of $\HH P^{n}$, and on the real slice $M=\HH P^{n}$. So in a Wolf space other than $\HH P^{n}$, where $R_{E}\neq0$ and $W_{z}$ is constant, the
integral divisors through a general point are unbalanced, hence contracted; this agrees with
Example~\ref{ex:intWolf}.
\end{enumerate}
\end{remark}

\begin{example}[Integral conics in two Wolf spaces]
\label{ex:intWolf}
This example is not used in the proofs. It gives constructions of integral conics through a point, with their dimension counts, and records what was checked by
computer about whether they are all the integral conics through the point; the status of each statement is indicated.

Let $M$ be a Wolf space, $Z\subset\PP(\mathfrak g)$ the adjoint variety, $z=[E_{\theta}]$, and $\mathfrak g=\bigoplus_{k=-2}^{2}\mathfrak g_{k}$ the contact grading, so that
$L^{\perp}_{z}\cong\mathfrak g_{1}$. For $v\in\widehat Z$ the affine contact hyperplane is $\widehat T_{v}\widehat Z\cap\mathfrak z_{\mathfrak g}(v)$; in particular a conic whose
plane is an abelian subalgebra is integral. The tangent directions are taken in the isotropic secant variety $\Sigma_{\mathrm{iso}}(X_{z})$, the union of the secant lines
of $X_{z}$ spanned by Levi-orthogonal points.
\begin{enumerate}
\item $E_{6}/SU(6)Sp(1)$, $n=10$, $X_{z}=\Gr(3,6)\subset\PP^{19}$.
  \begin{itemize}
  \item \emph{Construction.} For orthogonal roots $\beta,\gamma$ of grade one and $\delta=\beta+\gamma-\theta$, the span $\mathfrak a=\langle E_{\theta},E_{\beta},E_{\gamma},E_{\delta}\rangle$
  is abelian, so the quadric surface $\PP(\mathfrak a)\cap Z$ is integral. Its conics through $z$ are integral, and its two lines through $z$ have Levi-orthogonal
  directions $E_{\beta},E_{\gamma}$, which span the tangent of the conic at $z$, a point of $\sigma_{+}(X_{z})$. Translates by the stabilizer $P$ of $z$ give further integral
  conics through $z$.
  \item \emph{Computation.} In exact arithmetic in a Chevalley basis of $\mathfrak e_{6}$, the integral conics through $z$ form a divisor $D_{z}$ of $\Gamma_{z}$ of dimension $19$, which
  $\widetilde\mu_{z}$ contracts onto $\sigma_{+}(X_{z})$, of dimension $14$, with fibres of dimension $5$; in the computation every such conic is $P$-conjugate to one of the
  conics above, but this normal form up to $P$ is not proved here. The finitely many root-space quadrics in a fixed basis do not by themselves parametrize the
  family; the statement intended is about $P$-orbits, and it is computational evidence only.
  \end{itemize}
\item $\Gr_{4}(\RR^{n+4})=SO(n+4)/(SO(n)\times SO(4))$, $Z=\mathrm{OG}(2,n+4)$, $z=\Lambda_{0}$, $X_{z}=\PP^{1}\times \mathbb{Q}^{n-2}$. There are two constructions of
divisorial loci of integral conics through $z$.
  \begin{itemize}
  \item \emph{(a)} The conics of $2$-planes containing a fixed null vector $s\in\Lambda_{0}$. They lie in the integral quadric $\mathbb{Q}^{n}_{s}=\{\Lambda\ni s\}$, a
  Legendrian submanifold of $Z$. For each of the $\PP^{1}$ of choices of $s$ they form a family of dimension $2n-2$, of total dimension $2n-1$. Their tangent directions
  fill the Segre component $\PP^{1}\times\PP^{n-1}$ of $\Sigma_{\mathrm{iso}}$, of dimension $n$, so the fibres of $\widetilde\mu_{z}$ have dimension $n-1$.
  \item \emph{(b)} The conics in $\Gr(2,U)\cong\mathbb{Q}^{4}$, for $U\supset\Lambda_{0}$ an isotropic $4$-space; $\Gr(2,U)$ is integral. For each $U$ these form a family of
  dimension $6$, and the $U$ are parametrized by $\mathrm{OG}(2,\Lambda_{0}^{\perp}/\Lambda_{0})$, of dimension $2n-7$, so the total dimension is again $2n-1$. Their
  tangent directions fill $\PP^{3}$ for each $U$, hence the component of isotropic $2$-planes of $\Sigma_{\mathrm{iso}}$, of dimension $2n-4$, with fibres of dimension $3$.
  Construction (b) is empty for $n=3$. For $n=4$ it has two components: $\mathrm{OG}(2,4)\cong\PP^{1}\sqcup\PP^{1}$, the two rulings of a smooth quadric surface, and a
  generic conic determines $U$ as the span of its member $2$-planes, so the two components give distinct families, both distinct from (a), where all member
  $2$-planes share a vector.
  \item Both $\mathbb{Q}^{n}_{s}$ and $\Gr(2,U)$ are unions of lines of $Z$, and in both the two lines through $z$ span the tangent of the conic.
  \item \emph{Completeness.} That (a) and (b) give all the integral conics through $z$ was checked in exact arithmetic in a Chevalley basis of $\mathfrak{so}(n+4)$ for
  $n=4,5,6$; for general $n$ it is a conjectural completeness statement.
  \end{itemize}
\end{enumerate}
The fibres are positive-dimensional because the unipotent radical of the stabilizer of $z$, a Heisenberg group, acts trivially on $\PP(L^{\perp}_{z})$ but moves the
integral conics through $z$. This is consistent with Remark~\ref{rem:intflat}(3), by which the integral divisors of a Wolf space other than $\HH P^{n}$ are unbalanced. By
Lemma~\ref{lem:intcontract}(2) the fibre dimension is $h^{0}(N(-2z))$ for the generic integral conic when $L^{\perp}|_{C}$ is nef, and at least $h^{0}(N(-2z))-1$ in general.
The concrete computations on the Wolf space were carried out by machine and verified by the author.
\end{example}

\section{Cusps at a general point}
\label{sec:cuspsII}
\label{subsec:legcusp}
\label{sec:cuspgeneral}

A cusp at a general point $z$ could produce a divisor of $\Gamma_{z}$ whose image under $\widetilde\mu_{z}$ is a hyperplane, possibly transverse to the contact
hyperplane (Corollary~\ref{cor:cuspgeneral} and the remark after it), which the Hartogs reduction cannot ignore. Condition (K) is the generic expectation, but with no margin
(\S\ref{sec:P2}). {We prove it using only the Levi form and the bound $d\le2$ on the $L$-degree: a cusp whose position moves in a
dominating family folds its curve two-to-one onto its image, by an explicit involution.}

{
\begin{lemma}[Characteristic fibres]
\label{lem:charfibres}
Let $F\colon X\to Z$ be a holomorphic submersion onto an open subset of $(Z,\theta)$, $K=\ker F^{*}\theta$, and $\beta$ a local nowhere-vanishing scalar $1$-form with
$\ker\beta=K$. Then $\mathrm{Char}(K)=\{v\in K:\ \iota_{v}d\beta|_{K}=0\}$ equals $\ker dF$. If a biholomorphism $\tau$ between open subsets of $X$ preserves $K$ and fixes
pointwise a local section $S$ of $F$, then $F\circ\tau=F$ near $S$.
\end{lemma}

\begin{proof}
As $d(u\beta)|_{K}=u\,d\beta|_{K}$ for a unit $u$, we may take $\beta=F^{*}\vartheta$, $\vartheta$ a local scalar contact form. Then $d\beta|_{K}=F^{*}d\vartheta|_{K}$, $dF$ maps $K$
onto $L^{\perp}$ with kernel $\ker dF$, and $d\vartheta|_{L^{\perp}}$ is non-degenerate; so $\mathrm{Char}(K)=\ker dF$, and $\tau$ preserves $\ker dF$. In coordinates
$(x,y)$ near a point of $S$ with $F(x,y)=x$, connected $y$-slices and $S=\{y=0\}$, $d(F\circ\tau)(\partial_{y_{j}})=dF(\tau_{*}\partial_{y_{j}})=0$; so $F\circ\tau$ does not
depend on $y$, and it equals $F$ on $S$.
\end{proof}

\begin{proposition}[A quadratic involution]
\label{prop:cuspinvolution}
Let $B$ be a connected complex manifold of dimension $2n+1$ and $F\colon B\times\PP^{1}\to Z$ holomorphic, with the maps $f_{b}=F(b,\cdot)$ non-constant of $L$-degree
$d\le2$. Let $t$ be an affine coordinate, and suppose $df_{b}(0)=0$ for all $b$ and $e(b)=f_{b}(0)$ is a local biholomorphism. Then $d=2$, and there is a holomorphic
function $c$ on $B$ with $f_{b}\circ\tau_{b}=f_{b}$, $\tau_{b}(t)=-t/(1+c(b)t)$, for all $b$. In particular no $f_{b}$ is birational onto its image.
\end{proposition}

\begin{proof}
The contact differential $\theta\circ df_{b}\in H^{0}(\OO(d-2))$ is zero for $d\le1$, and for $d=2$ a constant vanishing at $0$; so $F^{*}\theta$ kills the vertical tangent
bundle. As $F^{*}L\otimes\mathrm{pr}_{2}^{*}\OO(-d)$ is trivial on the fibres of $B\times\PP^{1}\to B$, Grauert's base change theorem identifies $F^{*}L$ with
$\mathrm{pr}_{2}^{*}\OO(d)$ after shrinking $B$; the coefficient of each $db_{k}$ in $F^{*}\theta$ is then a section of $\OO(d)$ on each fibre, and in the affine frame
\[
\beta:=F^{*}\theta=\alpha_{0}+t\alpha_{1}+\dots+t^{d}\alpha_{d},\qquad\alpha_{j}\ \text{holomorphic $1$-forms on }B .
\]
(The degree bound uses that $f_{b}$ is defined on all of $\PP^{1}$.) Along $B\times\{0\}$, $dF|_{TB}=de$, so $F$ is a submersion there and $\beta\neq0$; hence $d\ge0$. Write
$F^{*}\vartheta=h\beta$ near a point of $B\times\{0\}$, $h$ a unit. As $\iota_{\partial_{t}}F^{*}\vartheta=0$ and $dF(\partial_{t})=0$ at $t=0$, Cartan's formula gives
$\mathcal L_{\partial_{t}}F^{*}\vartheta=\iota_{\partial_{t}}F^{*}d\vartheta=0$ there, that is, as the $\alpha_{j}$ are pulled back from $B$, $h_{t}\alpha_{0}+h\alpha_{1}=0$ at $t=0$; and $h(b,0)\alpha_{0}=e^{*}\vartheta$ is a contact
form. So $\alpha_{0}$ is contact and $\alpha_{1}=c\,\alpha_{0}$; the function $c$ does not depend on the trivialization of $F^{*}L$, which changes all $\alpha_{j}$ by the same unit on $B$, so it is defined on all of $B$.

If $d\le1$, $\beta=u\alpha_{0}$ with $u$ a unit near $B\times\{0\}$; then $K=\ker\alpha_{0}\oplus\CC\partial_{t}$, $d\beta|_{K}=u\,d\alpha_{0}|_{K}$ and $\mathrm{Char}(K)=\CC\partial_{t}$, so
Lemma~\ref{lem:charfibres} gives $dF(\partial_{t})=0$ near $t=0$ and $f_{b}$ is constant, a contradiction. If $d=2$, $\beta=(1+ct)\alpha_{0}+t^{2}\alpha_{2}$, and
$\tau(b,t)=(b,\tau_{b}(t))$ is an involution over $B$ fixing $B\times\{0\}$, with $\tau_{b}'(0)=-1$. As $\beta$ has no $dt$-component, $\tau^{*}\beta$ is obtained by
substituting $\tau_{b}(t)$ for $t$, and $1+c\tau_{b}(t)=(1+ct)^{-1}$, $\tau_{b}(t)^{2}=t^{2}(1+ct)^{-2}$ give $\tau^{*}\beta=(1+ct)^{-2}\beta$. So $\tau$ preserves $K$ near
$B\times\{0\}$, Lemma~\ref{lem:charfibres} gives $F\circ\tau=F$ there, and $f_{b}\circ\tau_{b}=f_{b}$ on $\PP^{1}$ by analytic continuation.
\end{proof}
}

\begin{theorem}[No cusp at a general point]
\label{thm:cuspgeneral}
Let $Z$ be a projective contact manifold and $\MCbar$ a component of its Hilbert scheme of conics whose general member is a balanced transverse conic. Then (K) holds: for a general
$z\in Z$, no member of $\MCbar$ has a cusp at $z$.
\end{theorem}

{
\begin{proof}
For an ample $H$, the curves supporting members of $\MCbar$ have bounded $H$-degree, so the maps $f\colon\PP^{1}\to Z$ birational onto them lie in a scheme of finite type
\cite[\S I.1]{Kol96}, and the pairs $(f,p)$ with $f(\PP^{1})$ the support of a member and $df(p)=0$ form a constructible set $\mathcal C$. If $f(\PP^{1})$ supports $C$, then
$[C]=m[f(\PP^{1})]$ and $md=2$, so $d\le2$. If (K) fails, the evaluation $(f,p)\mapsto f(p)$ maps some irreducible stratum $D$ of $\mathcal C$ dominantly to $Z$, and by
generic smoothness it is a submersion at a smooth point $(f_{0},p_{0})$ of $D$. Take a submanifold $B\subset D$ through it of dimension $2n+1$ on which the evaluation is a
local biholomorphism, and $F(b,t)=f_{b}(t+p_{b})$ in an affine coordinate with $p_{0}\neq\infty$. This contradicts Proposition~\ref{prop:cuspinvolution}.
\end{proof}

Only the Levi form, $d\le2$ and birationality are used; no freeness, Fano condition, $\mathrm{Pic}\,Z=\ZZ L$ or $n\ge3$. Birationality cannot be dropped: the maps of the proposition factor as
$f_{b}=g_{b}\circ q_{b}$, with $q_{b}(t)=t^{2}/(1+c(b)t)$ the quotient by $\tau_{b}$ and $g_{b}$ of $L$-degree one, and double covers of lines branched over a general point do exist; they give ribbons (Lemma~\ref{lem:ribbons}). As in \cite[Lemma~2.7(2)]{Ye94}, what is
excluded is a cusp \emph{at} a general point; cuspidal conics \emph{through} a general point, with the cusp elsewhere, are not excluded.
}

\emph{The hyperplane at a cusp.} {To see why (K) is needed, we locate the positive part at a cusp.} Let
$f\colon\PP^{1}\to Z$ be integral and non-constant of $L$-degree two, $R_{f}$ the ramification divisor of $df\colon T_{\PP^{1}}\to f^{*}L^{\perp}$, $R=\deg R_{f}$,
$T_{f}\cong\OO(2+R)$ the saturation of $df(T_{\PP^{1}})$, $T_{f}^{\perp}\subset f^{*}L^{\perp}$ its Levi-orthogonal, and $E_{f}=f^{*}TZ/T_{f}^{\perp}$. Through $\omega$,
$f^{*}L^{\perp}/T_{f}^{\perp}\cong T_{f}^{*}\otimes f^{*}L$, so
\begin{equation}
\label{eq:Ef}
0\to T_{f}^{*}\otimes f^{*}L\to E_{f}\to f^{*}L\to0,\qquad e_{f}\in H^{1}(\PP^{1},T_{f}^{*}).
\end{equation}

\begin{lemma}[The contact class along an integral conic]
\label{lem:contactclass}
Let $\delta_{f}\in H^{0}(T_{f}\otimes K_{\PP^{1}})\cong H^{0}(\OO(R))$ be the section defined by $df$, vanishing exactly on $R_{f}$. Then $e_{f}$ is, up to a non-zero universal constant, the restriction to
$T_{f}$ of $f^{*}\mathrm{at}(L)$, $\mathrm{at}(L)$ the Atiyah class \cite{Ati57}, and $\langle e_{f},\delta_{f}\rangle=c\deg f^{*}L=2c\neq0$ under Serre duality, $c$ a universal
constant.
\end{lemma}

\begin{proof}
Let $\theta=\theta_{\alpha}\otimes s_{\alpha}$ with frames $s_{\alpha}$ of $L$, $\theta_{\beta}=g_{\alpha\beta}\theta_{\alpha}$, and Reeb fields $R_{\alpha}$ splitting $TZ\to L$. With
$R_{\beta}=g^{-1}R_{\alpha}+Y$, $Y\in L^{\perp}$, the condition $\iota_{R_{\beta}}d\theta_{\beta}=0$ on $L^{\perp}$ gives $\iota_{Y}d\theta_{\alpha}|_{L^{\perp}}=g^{-2}dg|_{L^{\perp}}$, so through
$\omega$ the class of $0\to L^{\perp}\to TZ\to L\to0$ is represented, up to a universal constant, by $d\log g_{\alpha\beta}|_{L^{\perp}}$, the restriction of $\mathrm{at}(L)$.
Pulling back and pushing to $f^{*}L^{\perp}/T_{f}^{\perp}$ gives the first assertion. Pairing with $\delta_{f}$ composes with $df^{*}\colon f^{*}\Omega^{1}_{Z}\to K_{\PP^{1}}$,
which sends $f^{*}\mathrm{at}(L)$ to $\mathrm{at}(f^{*}L)$, that is $\deg f^{*}L$ times a generator.
\end{proof}

Let $\mathcal K_{f}\subset H^{0}(T_{f}\otimes f^{*}L^{-1})$ be the kernel of the connecting map to $H^{1}(f^{*}L^{-1})\cong\CC$ of the dual of \eqref{eq:Ef}. Restriction to the
sub-bundle identifies $\mathrm{Hom}(E_{f},\OO)$ with $\mathcal K_{f}$, and $\sigma\in\mathcal K_{f}$ gives $\alpha_{\sigma}\in H^{0}(f^{*}\Omega^{1}_{Z})$ with
$\alpha_{\sigma}|_{f^{*}L^{\perp}}=\omega(\cdot,\sigma)$. Identifying $f^{*}L$ with $T_{\PP^{1}}$ turns the connecting map into the pairing with $e_{f}$ and $\delta_{f}$ into
$\delta\in H^{0}(T_{f}\otimes f^{*}L^{-1})$, and Lemma~\ref{lem:contactclass} says
\begin{equation}
\label{eq:deltanotinK}
\delta\notin\mathcal K_{f}.
\end{equation}

\begin{corollary}
\label{cor:cuspgeneral}
At every cusp with $f^{*}TZ\cong\OO(3)\oplus\OO(1)^{2n-1}\oplus\OO$, the most balanced type, whether or not $z$ is general, the fibre $V$ at $p$ of the positive part
$\OO(3)\oplus\OO(1)^{2n-1}$ is not $L^{\perp}_{z}$.
\end{corollary}

{
\begin{proof}
As $\deg f^{*}TZ=(n+1)\deg f^{*}L$, $f$ has $L$-degree two and is integral (Lemma~\ref{lem:allornothing}: $\theta\circ df$ is constant and vanishes at $p$). The map $T_{\PP^{1}}\to f^{*}TZ$ lands in $\OO(3)$, so $T_{f}=\OO(3)$,
$R_{f}=p$, and $E_{f}$, a nef quotient of $f^{*}TZ$ of rank two and, by \eqref{eq:Ef}, degree one, is $\OO\oplus\OO(1)$. So $\mathcal K_{f}=\CC\sigma$, and $\sigma(p)\neq0$ by
\eqref{eq:deltanotinK}, since the sections of $T_{f}\otimes f^{*}L^{-1}\cong\OO(1)$ vanishing at $p$ are the multiples of $\delta$. Now $\alpha_{\sigma}$ kills the positive part,
so $V=\ker\alpha_{\sigma,p}$, while $\alpha_{\sigma,p}|_{L^{\perp}_{z}}=\omega_{z}(\cdot,\sigma(p))\neq0$. Hence $V\neq L^{\perp}_{z}$.
\end{proof}
}

The hyperplane calculation identifies a potential obstruction. If the local resolution of the tangent-direction map at a cusp at $z$, with $F_{z}$ smooth of
dimension $2n$ there, has the expected exceptional component mapping onto $\PP(V)$, that image is transverse to the contact hyperplane, and the Hartogs reduction
cannot ignore it; this is why (K) is needed. We do not construct this local model, since Theorem~\ref{thm:cuspgeneral} excludes cusps at a general point
independently of it.

\section{The local dictionary}
\label{sec:dictionaryII}
\label{sec:dictionary}

This section closes the local theory. It records, for each kind of degenerate conic, what the transform does there, and how the global argument of
Part~\ref{part:global} disposes of it: the last column of Table~\ref{tab:dictionary} is the bookkeeping of \S\ref{sec:introstrategy}, each kind being either handled (the
curvature extends), contracted by $\widetilde\mu_{z}$, or absent at a general point. No proof refers to this section; it is the interface between the local and the global
parts, and the beginning of a dictionary between degenerate rational curves and singularities of the transform (Question~\ref{q:dictionary}).

By Theorem~\ref{thm:thirdorder} the curvature at a balanced transverse conic is the invariant of order three of its contact neighbourhood, read through the identifications
\eqref{eq:Gk} and \eqref{eq:H1S4}. At a degenerate conic the neighbourhood departs from that of a balanced transverse conic at a definite order, these identifications break down,
and the transform acquires a singularity there.
\begin{itemize}[leftmargin=2em]
\item \emph{Integral conics}, at order zero (\S\ref{sec:integralII}). If $\theta|_{C}\equiv0$, the graded automorphism sheaves change already in low order
(Proposition~\ref{prop:contactauto}(2)), and the isomorphism $L^{-1}|_{C}\cong K_{C}$ of \eqref{eq:H1S4} fails: it is $\sigma_{\theta}$ that trivializes
$K_{C}\otimes L|_{C}$, and $\sigma_{\theta}([C])=0$. Along the balanced integral conics the metric acquires a pole, while $R_{E}^{\CC}$ extends regularly
(Proposition~\ref{prop:intregular}).
\item \emph{Nodal conics}, at order zero, the underlying curve being singular (\S\S\ref{sec:nodalgroups}, \ref{sec:dilation}). The neighbourhood of a balanced nodal conic
is controlled by the cohomology of the graded automorphism sheaves on the nodal curve, which vanishes in negative weight; this is why $W^{\CC}$ extends regularly
across the exceptional divisor over the nodes (Theorem~\ref{thm:Eregular}), and $R_{E}^{\CC}$ across the nodal divisor in a neck frame on a double cover (Corollary~\ref{cor:Rnodal}), with a value governed by higher-order data of the variety of minimal rational tangents
(Remark~\ref{rem:Elimit}).
\item \emph{Cusps}, at order zero as well (\S\ref{sec:cuspsII}). A cuspidal conic is integral (Lemma~\ref{lem:allornothing}); the proof excludes cusps at a general point.
\item \emph{Jumping conics}, at order one. At a jumping conic with $j(C)=1$, $N\cong\OO(2)\oplus\OO(1)^{2n-2}\oplus\OO$, and $L\otimes S^{2}N^{*}$, $L\otimes S^{3}N^{*}$ contain the
summands $\OO(-2)$, $\OO(-4)$: the neighbourhoods of orders one and two are no longer rigid, $B|_{C}$ is no longer trivial, and $\mathrm{ev}_{z}$ degenerates.
\end{itemize}
Table~\ref{tab:dictionary} summarizes what is known about the behaviour of $(g^{\CC},W^{\CC})$ near the degenerate conics, in the terms of
Definition~\ref{def:singconf} and Remark~\ref{rem:singLeBrun}.

\begin{table}[tbp]
\centering
\footnotesize
\begin{tabular}{>{\raggedright}p{3.4cm}>{\raggedright}p{5.8cm}>{\raggedright\arraybackslash}p{5.2cm}}
\toprule
member of $\MCbar\smallsetminus\MC$ & conformal geometry of $\MCbar$ & role, reference\\
\midrule
balanced line pair (node at $z$) & $\sigma_{\theta}$ vanishes to order one; the factorization $T\MC=E^{*}_{\CC}\otimes H^{*}_{\CC}$ is not available canonically;
$W^{\CC}$ regular at the generic point of the relevant exceptional component, its value fixed by weight-zero data of the VMRT & handled: Theorem~\ref{thm:Eregular}, Remark~\ref{rem:Elimit}\\
line pair with $z$ on a branch & $W^{\CC}$ may have poles & contracted: Lemma~\ref{lem:muz}(2), Remark~\ref{rem:Dprime}\\
unbalanced line pair & lies in $\{j\ge1\}$, with $j=a$; the closure of $J$ meets the nodal divisor only here & Remark~\ref{rem:jinv}\\
balanced integral conic & pole of order two on rank four, order one on the rest: a complex conformal infinity; $W^{\CC}$ regular and zero (Corollary~\ref{cor:intvanish}) & handled:
Proposition~\ref{prop:intregular}, Remark~\ref{rem:intflat}\\
unbalanced integral conic, ribbon, non-embedded map & --- & contracted, for separate reasons: unbalanced integral conics by Lemma~\ref{lem:intcontract}, ribbons by
Lemma~\ref{lem:ribbons}, non-embedded maps by Lemmas~\ref{lem:nodalbound} and~\ref{lem:intnonembedded}\\
jumping transverse conic ($j\ge1$) & $\omega_{E}$ and $\mathrm{ev}_{z}$ degenerate; $W_{z}$ may have a pole & contracted when $\widetilde\mu_{z}$ is birational:
Proposition~\ref{prop:cfromJ}\\
integral cusp at $z$ & a divisor whose image may be a hyperplane transverse to the contact hyperplane &
none at a general point: Theorem~\ref{thm:cuspgeneral}; Corollary~\ref{cor:cuspgeneral}\\
transverse rational immersion with an irreducible nodal image & two branches through the node & contracted with the non-embedded maps
(Lemma~\ref{lem:nodalbound}), so no separate curvature extension is needed; a posteriori none at all, $L$ being very ample on the adjoint variety
(Corollary~\ref{cor:beauville})\\
\bottomrule
\end{tabular}
\caption{The local dictionary: the members of $\MCbar\smallsetminus\MC$ (the jumping conics and the boundary) and the behaviour of $(g^{\CC},W^{\CC})$ near them.}
\label{tab:dictionary}
\end{table}

\begin{question}
\label{q:dictionary}
Is there a uniform description of the local invariants on the right? For instance: is the order of the pole of $W_{z}$ along a divisor of $\Gamma_{z}$ determined by the
order at which the contact neighbourhood of the generic member departs from that of a balanced transverse conic --- the splitting types of $f^{*}TZ$ and of $f^{*}L^{\perp}$, the order of
vanishing of $\sigma_{\theta}$, and the first non-vanishing group $H^{1}(C,\mathcal G_{k})$ --- as the cases computed so far suggest?
\end{question}

\part{Global Theory: The Compactification and the Contact Manifold}
\label{part:global}

In this part $Z$ is a compact contact Fano manifold of dimension $2n+1$ with $\mathrm{Pic}\,Z=\ZZ L$, and $\MCbar$ a component of its Hilbert scheme of conics whose
general member is a balanced transverse conic; such components exist when $n\ge3$ (Corollary~\ref{cor:balexist}), and for a twistor space we take the component of the twistor lines.
The question is how the global geometry of $\MCbar$, read through the incidence variety, acts on the curvature and on $Z$, given that the
structure must degenerate at the boundary (\S\ref{sec:mustdegenerate}). It carries out steps 2--4 of \S\ref{sec:introstrategy}. Section~\ref{sec:descent} proves the leaf constancy of the curvature (Theorems~\ref{thm:irred} and~\ref{thm:LSJ}) and bounds the degree of the tangent-direction
map by the number of conics through two general points (Lemma~\ref{lem:twopoint}); and Section~\ref{sec:twopointD} proves that exactly one conic passes through two general points (Theorem~\ref{thm:tpD}).


\section{Leaf constancy: the curvature on the compactification}
\label{sec:descent}
\label{sec:proofmain}

By the leaf principle (Proposition~\ref{prop:leaf}) the quaternionic Weyl curvature is parallel exactly when the curvature quartic $W_{z}$ is locally constant on the
conics through a general point, and on a non-compact germ it need not be. This section shows how the compactness of $\MCbar$ forces it, through the
tangent-direction map: the number of components of $\Gamma_{z}$ (Theorem~\ref{thm:irred}), the Hartogs reduction to the prime divisors of $\Gamma_{z}$ with divisorial image, the
leaf constancy theorem when $\widetilde\mu_{z}$ is birational ((K) holding by Theorem~\ref{thm:cuspgeneral}), the two-point count that makes $\widetilde\mu_{z}$ birational when
$d_{2}=1$, and the adjoint varieties (Corollary~\ref{cor:WolfJK}).

\subsection{The components of $\Gamma_{z}$}
\label{sec:irred}

\begin{theorem}[Irreducibility]
\label{thm:irred}
Let $\upsilon\colon\widehat F\to\overline F$ be the normalization, $\hat p=\bar p\circ\upsilon$, and $\hat p=g\circ p'$ its Stein factorization, with $p'_{*}\OO_{\widehat F}=\OO_{Z'}$
and $g\colon Z'\to Z$ finite.
\begin{enumerate}
\item There is a dense Zariski-open $V\subset Z$ such that for $z\in V$ the fibre $\overline F_{z}$ has exactly $\deg g$ irreducible components, each of dimension $2n$ and
containing transverse non-jumping conics, and so have $\widetilde F_{z}$ and $\Gamma_{z}$. In particular $F^{\mathrm{ne}}_{z}$ is dense in $\Gamma_{z}$.
\item $\deg g\le d_{2}$, the number of conics of $\MCbar$ through two general points of $Z$.
\item If $Z$ is the twistor space of a quaternionic-K\"ahler manifold, then $\deg g=1$: $\hat p_{*}\OO_{\widehat F}=\OO_{Z}$, every fibre of $\hat p$ and every
$\overline F_{z}$ is connected, and for $z\in V$ the spaces $\overline F_{z}$, $\widetilde F_{z}$ and $\Gamma_{z}$ are irreducible. The same holds whenever $d_{2}=1$.
\end{enumerate}
\end{theorem}

\begin{proof}
(1) $\widehat F$ is normal and irreducible of dimension $4n+1$, so $Z'$ is normal and irreducible and $g$ is finite of degree $d:=\deg g$. The generic fibre of $p'$ is
connected, since $p'_{*}\OO_{\widehat F}=\OO_{Z'}$, and normal, hence geometrically normal in characteristic zero; the set over which the fibres are geometrically normal is
constructible \cite[Prop.~9.9.4]{EGAIV3}, so it contains a dense open subset of $Z'$. Hence for general $z$, in the open set over which $g$ is \'etale, the fibre $\widehat F_{z}=\bigsqcup_{z'\in g^{-1}(z)}\widehat F_{z'}$ has
exactly $d$ connected components, each normal and therefore irreducible. The points of $\widehat F$ over transverse non-jumping conics form a dense open subset whose
complement has dimension at most $4n$; as the fibres of $p'$ have dimension $2n$, the complement contains whole fibres only over a subset of $Z'$ of dimension at most
$2n$, so for general $z$ every $\widehat F_{z'}$ contains transverse non-jumping conics. Over these $\overline F$ is smooth (Proposition~\ref{prop:genconic}(1)) and $\upsilon$ is
an isomorphism, so $\upsilon$ maps the $d$ components of $\widehat F_{z}$ onto $d$ distinct irreducible components of $\overline F_{z}$, which are all of them, $\upsilon$ being
surjective. The set $\mathcal B_{z}$ is a proper closed subset of each component, so its blow-up, the closure of the graph of $\widetilde\mu_{z}$ and its normalization have
the same number of components.

(2) Let $z$ and then $w$ be general. Each component of $\overline F_{z}$ contains a balanced conic, at which the conics through $z$ sweep out a neighbourhood of any other
point of the conic (Proposition~\ref{prop:genconic}(3)); the family being proper, the conics of each component sweep out $Z$, so through $w$ passes a conic of each
component. These conics are distinct: the conics through two general points are smooth, transverse and balanced (Proposition~\ref{prop:genconic}(3)), and at such a
conic $\overline F_{z}$ is smooth (Proposition~\ref{prop:genconic}(2)), so it lies on only one component of $\overline F_{z}$. Hence $d_{2}\ge d$.

(3) \emph{Step 1.} At a twistor line $N=\OO(1)^{2n}$, so $H^{1}(N)=0$, $\MCbar$ is smooth at $[C_{m}]$ and the universal family is smooth near the
smooth fibre $C_{m}$. So $\overline F$ is smooth along the compact set $\{([C_{m}],z):z\in C_{m}\}$, and $\upsilon$ is an isomorphism near it.
\emph{Step 2.} The maps $q$ and $m\mapsto[C_{m}]$ are real-analytic, so $s(z)=\upsilon^{-1}([C_{q(z)}],z)$ is a continuous section of $\hat p$.
\emph{Step 3.} In characteristic zero $g$ is \'etale over a dense Zariski-open $U\subset Z$, whose complement is a proper analytic subset, so $U$ is connected
and $g^{-1}(U)\to U$ is a finite covering with connected total space. The map $\sigma'=p'\circ s$ is a continuous section of $g$; its image over
$U$ is open (a single sheet over each evenly covered connected open set) and closed, hence all of $g^{-1}(U)$, so $\deg g=1$, and $g$, finite and
birational onto the normal $Z$, is an isomorphism. Zariski's connectedness theorem gives the connectedness of the fibres, and (1) the irreducibility. If $d_{2}=1$, then
$\deg g=1$ by (2).
\end{proof}

Part (3) uses a quaternionic-K\"ahler manifold; it is not used in the proofs of Theorem~\ref{thm:LSJ}, only on the real slice (\S\ref{sec:qk}). In general $\Gamma_{z}$ can be reducible only if $d_{2}\ge2$;
the arguments below are applied to each irreducible component of $\Gamma_{z}$, on which $\widetilde\mu_{z}$ is surjective and generically finite
(Lemma~\ref{lem:muz}(1)), and irreducibility is used only where it is stated.

\subsection{The Hartogs reduction}
\label{sec:hartogs}

\begin{lemma}[Hartogs through a generically finite map]
\label{lem:hartogsmu}
Let $\Gamma$ be a normal irreducible projective variety, $\mu\colon\Gamma\to\PP^{N}$ a surjective generically finite morphism, and $f$ a
meromorphic function on $\Gamma$ that is regular at the generic point of every prime divisor $D\subset\Gamma$ for which $\mu(D)$ is a divisor.
Then $f$ is constant.
\end{lemma}

\begin{proof}
Let $d=\deg\mu$ and $\sigma_{1},\dots,\sigma_{d}$ the elementary symmetric functions of $f$ over the fibres; they are rational functions on
$\PP^{N}$. Over the generic point of a prime divisor $\Delta\subset\PP^{N}$ the map is finite (its non-finite locus has image of codimension
$\ge2$), and the points over it lie on divisors of $\Gamma$ mapping onto $\Delta$, where $f$ is bounded. So each $\sigma_{k}$ is bounded near the
generic point of every prime divisor. The polar set of a rational function on $\PP^{N}$ is empty or a divisor, so the $\sigma_{k}$ are
regular, hence constant. Then $f$ is a root of a fixed polynomial, and $\Gamma$ being irreducible, $f$ is constant.
\end{proof}

\begin{definition}[The three hypotheses]
\label{def:abc}
For a general $z\in Z$ we say that
\begin{itemize}[leftmargin=2.5em]
\item[(a)] holds if $\Gamma_{z}$ is irreducible;
\item[(b)] holds if $W^{\CC}$ is regular at the generic point of every component of $E_{z}$ that dominates $\PP(L^{\perp}_{z})$;
\item[(c)] holds if $W^{\CC}$ is regular at the generic point of every prime divisor $D\subset\Gamma_{z}$, other than the components of
$E_{z}$, whose image $\widetilde\mu_{z}(D)$ is a divisor of $\PP(T_{z}Z)$.
\end{itemize}
(Regularity of $W^{\CC}$ at the generic point of a component of $E$ in $\widetilde F$ implies regularity of $W_{z}$ at the generic point of the
corresponding components of $E_{z}$ for general $z$, the generic points of the latter being generic points of the former.)
\end{definition}

A component of $E_{z}$ with divisorial image has image in $\PP(L^{\perp}_{z})$ (Proposition~\ref{prop:resolution}, over the nodes with independent
branch tangents; a component lying over ribbons is contracted, Lemma~\ref{lem:ribbons}), hence equal to it; so (b) and (c)
together cover every divisor with divisorial image other than those meeting $F^{\mathrm{ne}}_{z}$. Along the divisors of smooth integral conics, where the
denominator of $g^{\CC}$ vanishes, (c) holds without further assumption (Corollary~\ref{cor:intharmless}).

\begin{theorem}[The Hartogs reduction]
\label{thm:Ereduction}
If (b) and (c) hold for a general $z$, then $W_{z}$ is constant on every irreducible component of $\Gamma_{z}$, and $R_{E}^{\CC}$ is parallel. If moreover (a) holds,
$W^{\CC}$ is the pull-back of a section of $S^{4}L^{\perp*}\otimes L^{2}$ over a dense open subset of $Z$. If $Z$ is a twistor space, $M$ is then a Wolf space.
\end{theorem}

\begin{proof}
Divisors of $\Gamma_{z}$ meeting $F^{\mathrm{ne}}_{z}$ impose nothing, $W_{z}$ being holomorphic there. Let $\Gamma$ be an irreducible component of
$\Gamma_{z}$; it is normal and projective, and $\widetilde\mu_{z}|_{\Gamma}$ is surjective and generically finite (Lemma~\ref{lem:muz}(1)). Apply
Lemma~\ref{lem:hartogsmu} to $f=\ell\circ W_{z}|_{\Gamma}$ for every linear form $\ell$ on $V_{z}$: by (b) and (c) its hypotheses hold, so $W_{z}$ is constant on $\Gamma$.
Hence $W_{z}$ is locally constant on $\Sigma_{z}$ for general $z$, and Proposition~\ref{prop:leaf} gives $\nabla R_{E}^{\CC}=0$. Under (a), $W_{z}$ is constant on $\overline F_{z}$
for general $z$, and its value is a section of $S^{4}L^{\perp*}\otimes L^{2}$ over the dense open set of such $z$ whose pull-back is $W^{\CC}$. The last assertion is
Corollary~\ref{cor:endgame}.
\end{proof}

\begin{example}
\label{ex:bl0}
On $\mathrm{Bl}_{0}\CC^{2}$ a meromorphic function cannot have poles only along the exceptional curve: its push-forward would be holomorphic on
$\CC^{2}\smallsetminus0$, hence on $\CC^{2}$. By contrast $x/y=u$ (in the chart $x=uy$) is regular along the exceptional curve and has its only pole
on the strict transform of $\{y=0\}$. This is the pattern that matters here: a pole of $W_{z}$ on $\Gamma_{z}$ along a divisor contracted by
$\widetilde\mu_{z}$, such as $D'_{z}$, is impossible unless $W_{z}$ also has a pole along a divisor with divisorial image. Regularity along $D'_{z}$,
which is not a property of germs, is therefore never needed.
\end{example}

\begin{remark}[The prime divisors with divisorial image]
\label{rem:divisors}
Theorem~\ref{thm:Ereduction} looks only at the prime divisors $D\subset\Gamma_{z}$ with $\dim\widetilde\mu_{z}(D)=2n-1$. Lemma~\ref{lem:hartogsmu} works for any generically finite $\widetilde\mu_{z}$; it does not need $\deg\widetilde\mu_{z}=1$. The prime divisors
of $\Gamma_{z}$ that can have divisorial image, and how the paper disposes of them, are the following.

\begin{center}
\small
\begin{tabular}{@{}p{5.2cm}p{3.6cm}p{6.6cm}@{}}
\toprule
divisor $D\subset\Gamma_{z}$ & image & disposed of by \\
\midrule
meets $F^{\mathrm{ne}}_{z}$ & divisor of $\CC^{2n}$ & $W_{z}$ holomorphic there (\S\ref{sec:Wz}) \\
divisor of jumping conics through $z$ & in $\CC^{2n}$ & contracted when $\widetilde\mu_{z}$ is birational (Proposition~\ref{prop:cfromJ}) \\
component of $E_{z}$ from a pair $X_{i},X_{j}$ with $\omega(\hat X_{i},\hat X_{j})\not\equiv0$ & $\PP(L^{\perp}_{z})$ & Theorem~\ref{thm:Eregular} \\
component of $E_{z}$ from a linear $X_{i}$ paired with itself & $\subset\PP(S)$, Lagrangian & contracted (Proposition~\ref{prop:secantgeneral}(2)) \\
balanced component of $D_{z}$ (integral) & $\PP(L^{\perp}_{z})$ & Proposition~\ref{prop:intregular}, Corollary~\ref{cor:intharmless} \\
unbalanced component of $D_{z}$ & dimension $\le2n-2$ & contracted (Lemma~\ref{lem:intcontract}) \\
divisors over ribbons, incl.\ components of $E_{z}$ over ribbons & dimension $\le2n-2$ & contracted, no hypothesis (Lemma~\ref{lem:ribbons}) \\
divisor over images of transverse maps that are not embeddings & dimension $\le n+1$ & contracted for $n\ge3$ (Lemma~\ref{lem:nodalbound}) \\
divisor over images of integral maps that are not embeddings, no cusp over $z$ & dimension $\le2n-2$ & contracted, no hypothesis (Lemma~\ref{lem:intnonembedded}) \\
divisor over cuspidal integral conics with cusp at $z$ & a hypersurface (the hyperplane $\PP(V)$ for the most balanced type, possibly transverse to $\PP(L^{\perp}_{z})$), or smaller & does not occur for general $z$ (Theorem~\ref{thm:cuspgeneral}) \\
\bottomrule
\end{tabular}
\end{center}
Divisors whose generic point lies over a line pair with $z$ on one branch away from the node form $D'_{z}$, mapped into $X_{z}\subset\PP(T_{z}Z)$, and are contracted
(Lemma~\ref{lem:muz}(2)). Classifying a prime divisor by the cycle under its generic point (reduced and irreducible, two lines, or a double line), the table is exhaustive. Since (K) holds (Theorem~\ref{thm:cuspgeneral}), the table gives (c) as soon as $\widetilde\mu_{z}$ is birational on every component of $\Gamma_{z}$
(Proposition~\ref{prop:cfromJ}).
\end{remark}

\begin{example}[The model $\Gr_{2}(\CC^{n+2})$]
\label{ex:gr2}
Let $Z=F(1,n+1;\CC^{n+2})$, the twistor space of $\Gr_{2}(\CC^{n+2})$ (here $\mathrm{Pic}\,Z\neq\ZZ L$, but the computation is instructive).
Its conics are $C_{W,W'}=\{(x,\xi):x\in\PP(W),\ \xi\supset W'+x\}$ with $\dim W=2$, $\dim W'=n$; those through $z=(x_{0},\xi_{0})$ are given by
$W=\langle x_{0},w\rangle$, $w\in\PP^{n}$, and $W'=\ker\eta\subset\xi_{0}$, $\eta\in\PP^{n*}$, so $\overline F_{z}=\PP^{n}\times\PP^{n*}$. In a chart
with $\theta=dP+\sum Q_{i}\,da_{i}$,
\[
\mu_{z}(w,\eta)=\big[\eta_{0}w_{1}:\dots:\eta_{0}w_{n}:\eta_{0}w_{n+1}:-w_{n+1}\eta_{1}:\dots:-w_{n+1}\eta_{n}\big],\qquad\theta(T)\propto\eta_{0}w_{n+1}.
\]
On the transverse conics $\{\eta_{0}w_{n+1}\neq0\}\cong\CC^{2n}$ it is an isomorphism onto $\CC^{2n}$; the branch divisors map onto the two
$\PP^{n-1}$'s forming the VMRT; the indeterminacy locus is $\mathcal B_{z}=\{\eta_{0}=w_{n+1}=0\}\cong\PP^{n-1}\times\PP^{n-1}$, and after one blow-up
$E_{z}\to\PP(L^{\perp}_{z})$, $(w',\eta',[\alpha:\beta])\mapsto[\alpha w':0:-\beta\eta']$, is birational onto the join of the two $\PP^{n-1}$'s. So
$\widetilde\mu_{z}$ is birational, $E_{z}$ is the only divisor with divisorial image, there are no jumping conics, and (c) is void.
\end{example}


\subsection{The leaf constancy theorem}
\label{sec:Jthm}
\label{sec:descentthm}

The Hartogs reduction needs (b) and (c) of Definition~\ref{def:abc}. Hypothesis (b) is Theorem~\ref{thm:Eregular}. Of the divisors in
Remark~\ref{rem:divisors}, all are disposed of without any assumption except the divisors of jumping conics through $z$, and these are harmless as soon as
$\widetilde\mu_{z}$ is birational. At a jumping conic $\widetilde\mu_{z}$ is not a local isomorphism (Lemma~\ref{lem:dmu}), so by Zariski's main theorem a birational
$\widetilde\mu_{z}$ maps every divisor of jumping conics to a subset of codimension $\ge2$. We state the leaf constancy theorem with this property as its only
hypothesis on the jumping conics. A \emph{jumping divisor through $z$} is a prime divisor of $\Gamma_{z}$ whose general member is a jumping conic. Condition (K)
of \S\ref{sec:P2}, used in the proof, holds by Theorem~\ref{thm:cuspgeneral}.

\begin{theorem}[Leaf constancy]
\label{thm:LSJ}
\label{thm:descent}

Let $Z$ be a compact contact Fano manifold of dimension $2n+1$, $n\ge3$, with $\mathrm{Pic}\,Z=\ZZ L$, and $\MCbar$ a component of its Hilbert scheme of conics whose
general member is a balanced transverse conic. Suppose that for a general $z$ no jumping divisor through $z$ is mapped by $\widetilde\mu_{z}$ onto a divisor of
$\PP(T_{z}Z)$. This holds if $\widetilde\mu_{z}$ has degree one on every irreducible component of $\Gamma_{z}$ (Proposition~\ref{prop:cfromJ}), hence if
$d_{2}=1$ (Lemma~\ref{lem:twopoint}), and trivially if no conic through a general point is jumping. Then $R_{E}^{\CC}$ is parallel, and
$(\MC,g^{\CC})$ is locally symmetric. If moreover $\Gamma_{z}$ is irreducible for general $z$, then $W^{\CC}$ is the pull-back of a section of
$S^{4}L^{\perp*}\otimes L^{2}$ over a dense open subset of $Z$.

\end{theorem}

By Theorem~\ref{thm:Ereduction} it suffices to establish (b) and (c) of Definition~\ref{def:abc} for a general $z$. Hypothesis (b) is Theorem~\ref{thm:Eregular}; it is
vacuous if no component of $E_{z}$ dominates $\PP(L^{\perp}_{z})$, and Theorem~\ref{thm:Ereduction} does not need one to. It remains to derive (c).

\begin{proposition}[Hypothesis (c) when the jumping conics are harmless]
\label{prop:cfromJ}
Let $z$ be general, and suppose that no jumping divisor through $z$ is mapped by $\widetilde\mu_{z}$ onto a divisor. Every prime divisor of $\Gamma_{z}$ whose image
under $\widetilde\mu_{z}$ is a divisor either meets $F^{\mathrm{ne}}_{z}$, where $W^{\CC}$ is holomorphic, or is a component of $E_{z}$ dominating $\PP(L^{\perp}_{z})$,
or is a balanced component of $D_{z}$, along which $W_{z}$ is regular (Corollary~\ref{cor:intharmless}). In particular (c) holds. The hypothesis holds if
$\widetilde\mu_{z}$ has degree one on every irreducible component of $\Gamma_{z}$.
\end{proposition}

\begin{proof}

Let $D\subset\Gamma_{z}$ be a prime divisor with divisorial image.

If $D$ meets $F^{\circ}_{z}$, its generic member $C$ is a transverse conic through $z$. By hypothesis $C$ is not jumping, so $D$ meets $F^{\mathrm{ne}}_{z}$.

If $D$ does not meet $F^{\circ}_{z}$, then by Proposition~\ref{prop:boundaryinfinity}
$\widetilde\mu_{z}(D)\subset\PP(L^{\perp}_{z})\cup\overline{S^{(i)}_{z}}$, and $\dim\overline{S^{(i)}_{z}}\le2n-2$; being an irreducible divisor, it equals $\PP(L^{\perp}_{z})$, so $D$
dominates the hyperplane. Classify $D$ by the cycle under its generic point. Two distinct lines through $z$ meet only at $z$ (\S\ref{sec:vmrt}). Two lines with $z$ on one branch away from the
node give $D'_{z}$, which is contracted (Lemma~\ref{lem:muz}(2)). Two lines with node at $z$ give a component of $E_{z}$, which maps onto the secant lines of the branch
tangents (Lemma~\ref{lem:muz}(3), Proposition~\ref{prop:resolution}): from a pair with $\omega(\hat X_{i},\hat X_{j})\not\equiv0$ it is regular by
Theorem~\ref{thm:Eregular}; from a linear component paired with itself its image lies in a Lagrangian $\PP(S)$ (Proposition~\ref{prop:secantgeneral}(2)), so it is contracted.
A double line gives a ribbon, contracted (Lemma~\ref{lem:ribbons}), including the components of $E_{z}$ over the diagonal of $X_{i}\times X_{i}$. A reduced irreducible singular support is the
image of a non-embedded transverse map, contracted by Lemma~\ref{lem:nodalbound}, or of a non-embedded integral map, contracted by Lemma~\ref{lem:intnonembedded}, a cusp
at $z$ being excluded by Theorem~\ref{thm:cuspgeneral}. A smooth integral conic gives a component of $D_{z}$, contracted if unbalanced (Lemma~\ref{lem:intcontract}) and with $W_{z}$ regular if balanced
(Corollary~\ref{cor:intharmless}).

For the last assertion, let $D$ be a jumping divisor through $z$ with divisorial image, and $\Gamma$ the irreducible component of $\Gamma_{z}$ containing it,
on which $\widetilde\mu_{z}$ has degree one. The map $\widetilde\mu_{z}|_{\Gamma}\colon\Gamma\to\PP(T_{z}Z)$ is a proper birational morphism from a normal variety onto a
smooth one. By Zariski's main theorem it is an isomorphism over the complement of a closed subset of codimension $\ge2$; since $\widetilde\mu_{z}(D)$ is a
divisor, $\widetilde\mu_{z}$ is a local isomorphism at the generic point of $D$. But let $C$ be the generic member of $D$, a jumping conic through $z$. Every conic
through the general point $z$ is free \cite[II.3.11]{Kol96}, so $N_{C}$ is nef and $H^{1}(N_{C}(-z))=0$; then the kernel of $d\widetilde\mu_{z}$ at $[C]$ is
$H^{0}(N_{C}(-2z))$, which is nonzero because $C$ is jumping (Lemma~\ref{lem:dmu}). This is a contradiction.
\end{proof}

\begin{proof}[Proof of Theorem~\ref{thm:LSJ}]
For general $z$, (b) holds by Theorem~\ref{thm:Eregular}, and (c) by Proposition~\ref{prop:cfromJ}.
Theorem~\ref{thm:Ereduction} concludes.
\end{proof}

Theorem~\ref{thm:Ereduction} and Lemma~\ref{lem:hartogsmu} need only that $\widetilde\mu_{z}$ be generically finite (Lemma~\ref{lem:muz}(1)); the
birationality of $\widetilde\mu_{z}$ enters only through the jumping divisors, that is, through the hypothesis of Theorem~\ref{thm:LSJ}.

The divisors of $\Gamma_{z}$ along which LeBrun's structure degenerates are thus disposed of as follows. The jumping divisors meet the transverse locus
and are contracted when $\widetilde\mu_{z}$ is birational. The divisors at infinity, those not meeting $F^{\circ}_{z}$, map into the contact hyperplane except over non-injective transverse immersions (into $\overline{S^{(i)}_{z}}$, of dimension $\le n+1$), and each kind is
either contracted (the line pairs with $z$ on a branch, ribbons, non-embedded maps, unbalanced integral conics, the components of $E_{z}$ from a linear self-pair) or a divisor
along which $W_{z}$ is regular (the dominating components of $E_{z}$, by Theorem~\ref{thm:Eregular}, and the balanced integral divisors, by Corollary~\ref{cor:intharmless}).
The cusp divisors do not occur for a general $z$ (Theorem~\ref{thm:cuspgeneral}).

\begin{remark}[Adjoint varieties]
\label{rem:Wolfconverse}
If $Z=G_{\CC}/P$ is an adjoint variety (the twistor space of a Wolf space), the twistor lines form an open $G_{\CC}$-orbit in $\MCbar$, with normal bundles $\OO(1)^{2n}$. On $\Gr_{2}(\CC^{n+2})$
every transverse conic lies in this orbit (Example~\ref{ex:gr2}), so $J=\emptyset$; for $\HH P^{n}$ the conics are the lines of $\PP^{2n+1}$, all with normal bundle
$\OO(1)^{2n}$. If $\mathrm{Pic}\,Z=\ZZ L$, then $d_{2}=\deg\widetilde\mu_{z}=1$ (Theorem~\ref{thm:tpD} and Lemma~\ref{lem:twopoint}; for the adjoint varieties also Corollary~\ref{cor:WolfJK}), so Theorem~\ref{thm:LSJ} applies when $n\ge3$.

On such a $Z$ the curvature has no pole anywhere over $\MCbar$. Indeed, $W^{\CC}=\bar p^{\,*}(c\,Q)$ on the open orbit, with $Q$ an invariant section of
$S^{4}L^{\perp*}\otimes L^{2}$ over $Z$. Let $z=[e_{\theta}]$, $\theta$ the highest root, grade $\mathfrak g$ by $\theta^{\vee}$, and let $H$ be the centralizer of
$\mathfrak{sl}_{2}(\theta)=\langle e_{\theta},\theta^{\vee},e_{-\theta}\rangle$; then $L^{\perp}_{z}\cong\mathfrak g_{-1}$ and $L_{z}\cong\mathfrak g_{-2}$.
\begin{itemize}
\item $W^{\CC}$ is $G_{\CC}$-invariant, LeBrun's construction being natural under contact automorphisms, and $H$ fixes the twistor line
$C_{0}=\PP(\mathfrak{sl}_{2}(\theta))\cap Z$ pointwise. So $W^{\CC}([C_{0}],z)$ lies in $(S^{4}\mathfrak g_{-1}^{*})^{H}\otimes\mathfrak g_{-2}^{\otimes2}$.
\item $(S^{4}\mathfrak g_{-1}^{*})^{H}$ is a line $\CC Q_{z}$: Freudenthal's quartic for $\mathfrak f_{4},\dots,\mathfrak e_{8}$ \cite[\S5.2]{LM01}; the
discriminant of binary cubics for $\mathfrak g_{2}$; the Gram determinant $\det q(u_{i},u_{j})$ on $\CC^{2}\otimes\CC^{m}$ for $\mathfrak{so}_{m+4}$, $m\ge3$; and
$\xi(x)^{2}$ on $\CC^{m}\oplus\CC^{m*}$ for $\mathfrak{sl}_{m+2}$. For $\mathfrak{sp}_{2n+2}$ it is zero, $S^{4}\CC^{2n}$ being irreducible.
\item The twistor lines through $z$ form one orbit of the stabilizer $P$ of $z$: if $gC_{0}\ni z$, then $gs\in P$ for some $s$ in the group of
$\mathfrak{sl}_{2}(\theta)$, and $gC_{0}=gsC_{0}$. $P$ acts trivially on $\CC Q_{z}\otimes\mathfrak g_{-2}^{\otimes2}$: its unipotent radical acts trivially, and
$\theta^{\vee}$ with weight $4-4=0$. So $W^{\CC}([C],z)=c\,Q_{z}$ for every twistor line $C\ni z$, and $Q_{z}$ extends to a $G_{\CC}$-invariant section $Q$.
\item $c\neq0$ unless $Z=\PP^{2n+1}$: if $W^{\CC}$ vanished, $M$ would be locally isometric to $\HH P^{n}$ \cite{Salamon82}.
\end{itemize}
So the rational section $W^{\CC}$ of
$\bar p^{\,*}(S^{4}L^{\perp*}\otimes L^{2})$ is holomorphic on all of $\overline F$, and for every $z$, $W_{z}$ is the constant $c\,Q_{z}$ on all of
$\overline F_{z}$ and of $\Gamma_{z}$. There is no pole along any divisor of jumping conics, line pairs (with the node at $z$ or with $z$ on a branch), integral
conics, ribbons or cuspidal curves. On $\PP^{2n+1}$, $W^{\CC}=0$. What is pole-free is the curvature quartic: the metric $g^{\CC}$ still has poles at the
boundary (Remark~\ref{rem:intflat}(2)), and $E_{\CC}$ degenerates along $J$ and at the integral conics. In a general germ $W_{z}$ may have poles, for
instance along $D'_{z}$ (Remark~\ref{rem:Dprime}); the proof excludes them only along the divisors of $\Gamma_{z}$ whose image under $\widetilde\mu_{z}$ is
a divisor, and Hartogs does the rest.
\end{remark}

\begin{remark}[Where $n\ge3$ is used]
\label{rem:nge3}
The condition $n\ge3$ is used through Corollary~\ref{cor:balgeneric} and Theorem~\ref{thm:Eregular}, and through Lemma~\ref{lem:nodalbound} (where $n+1<2n-1$) and
Proposition~\ref{prop:boundaryinfinity} (where $n+1\le2n-2$); Lemmas~\ref{lem:ribbons}, \ref{lem:intnonembedded} and Proposition~\ref{prop:intregular} need only $n\ge2$.
On the real slice the case $n=2$ is the theorem of Poon and Salamon \cite{PS91}.
\end{remark}

On the real slice, where Theorem~\ref{thm:tpD} gives $d_{2}=1$, the theorem gives Theorem~Q (\S\ref{sec:qk}): there the compactness of $M$ turns local
symmetry into symmetry. For a general $Z$ the compactness of $Z$ does the same: a contact Fano manifold with $\mathrm{Pic}\,Z\cong\ZZ$ whose quaternionic Weyl
curvature is parallel is an adjoint variety (Theorem~\ref{thm:homog}), and with Theorem~\ref{thm:tpD} this proves Beauville's conjecture \cite{Bea98} when $n\ge3$ and
$\mathrm{Pic}\,Z=\ZZ L$ (Corollary~\ref{cor:beauville}).

\subsection{The two-point count and the adjoint varieties}
\label{sec:converse}
\label{subsec:twopoint}

The leaf constancy theorem is applied through a count of conics through two points.

Let $\mathcal F_{2}\to\MCbar$ be the fibre square of the universal conic and $\Phi_{2}\colon\mathcal F_{2}\to Z\times Z$, $([C],q_{1},q_{2})\mapsto(q_{1},q_{2})$; its fibre over $(z,w)$ is the set of
conics through $z$ and $w$, and it is generically finite. Put $d_{2}:=\deg\Phi_{2}$, the number of conics through two general points.

\begin{lemma}[Two points and one direction]
\label{lem:twopoint}
For general $z$, $d_{2}\ge\deg\widetilde\mu_{z}$; so $d_{2}=1$ forces $\deg\widetilde\mu_{z}=1$. Moreover $\Phi_{2}$ ramifies at $([C],q_{1},q_{2})$, $q_{1}\ne q_{2}$
points of a smooth transverse conic, exactly when $C$ is jumping.
\end{lemma}

\begin{proof}
At $([C],q_{1},q_{2})$ the kernel of $d\Phi_{2}$ is $H^{0}(N_{C}(-q_{1}-q_{2}))$, which for a smooth transverse conic is nonzero exactly when $N_{C}$ has a summand of
degree $\ge2$, i.e.\ when $C$ is jumping (Lemma~\ref{lem:dmu}, applied to the second point on the family of conics through the first). For the degree, let
$\mathrm{Exp}_{z}\colon\{([C],w):z\in C\ni w\}\to Z$, $([C],w)\mapsto w$, be the restriction of $\Phi_{2}$ over $\{z\}\times Z$; so $\deg\mathrm{Exp}_{z}=d_{2}$, and $d_{2}$ is the number
of points of $\mathrm{Exp}_{z}^{-1}(w)$ for a general $w$ near $z$. The fibre $\mathrm{Exp}_{z}^{-1}(z)$ contains the section $\Sigma^{\rm all}=\{([C],z)\}\cong\overline F_{z}$, and we
count the points of $\mathrm{Exp}_{z}^{-1}(w)$ near it. Along the members smooth at $z$ the pull-back of the maximal ideal of $z$ is the ideal of $\Sigma^{\rm all}$, so
$\mathrm{Exp}_{z}$ lifts to $\mathrm{Bl}_{z}Z$ near such points, with restriction $\widetilde\mu_{z}$ on $\Sigma^{\rm all}$ to the exceptional divisor $E=\PP(T_{z}Z)$. Over a
general point $[v]\in E$ the fibre of the lift near $\Sigma^{\rm all}$ is $\widetilde\mu_{z}^{-1}([v])$, which for general $[v]$ consists of $\deg\widetilde\mu_{z}$ non-jumping
transverse conics, and at each of them the lift is a local isomorphism, its differential being $d\widetilde\mu_{z}$ along $\Sigma^{\rm all}$, an isomorphism at a non-jumping
conic, and nonzero on the fibre direction $T_{z}C$, which is transverse to $E$. So the points near $\Sigma^{\rm all}$ contribute at least $\deg\widetilde\mu_{z}$. (Points near
members singular at $z$ may contribute more.)
\end{proof}

\begin{corollary}[The adjoint varieties]
\label{cor:WolfJK}
Let $Z$ be an adjoint variety with $\mathrm{Pic}\,Z=\ZZ L$, for instance the twistor space of a Wolf space other than $\HH P^{n}$ and $\Gr_{2}(\CC^{n+2})$. Then $\MCbar$ has no
member whose support is singular and irreducible; in particular (K) holds. Moreover $d_{2}=\deg\widetilde\mu_{z}=1$ and $\Gamma_{z}$ is irreducible for general $z$.
\end{corollary}

\begin{proof}
$Z$ is the adjoint variety of a simple complex Lie group $G_{\CC}$, embedded in $\PP(\mathfrak g)$ by $L$, so $L$ is very ample. An irreducible curve with $L\cdot C=2$ is then
a curve of degree two in projective space, hence a smooth plane conic; so no member of $\MCbar$ has a singular irreducible support, which gives
(K). Through two general points of an adjoint variety with Picard group $\ZZ$ passes exactly one rational curve of $L$-degree two \cite[Proposition~3.11]{CP09}; since the
conics of $\MCbar$ through a general point cover $Z$ (Proposition~\ref{prop:genconic}(3)), $d_{2}=1$. By Lemma~\ref{lem:twopoint}, $\deg\widetilde\mu_{z}=1$ for general $z$, so
$\Gamma_{z}$ is irreducible.
\end{proof}

\section{The conic through two general points}
\label{sec:twopointD}

This section proves Theorem~\ref{thm:tpD}: exactly one conic of $\MCbar$ passes through two general points of $Z$. We work on the space of two-pointed
stable maps, where the marked points can be followed into the boundary. The tool is a \emph{two-point vector lift}: for a smooth conic $h$ and two of
its points, the tangent vector of $h$ at the second point, weighted by the contact value of the tangent vector at the first. It extends regularly over
the whole space of two-pointed conics, including the nodal boundary (\S\ref{subsec:tplifts}).

If the two-point evaluation map had degree $\ge2$, it would branch along a divisor of $Z\times Z$. Test curves that fix the first point and move the
second along a line through that divisor show that, among the conics through two general points, some two have lifts differing by a vector tangent to
the variety of minimal rational tangents at the second point (\S\ref{subsec:tptest}). Two cases follow.
\begin{itemize}[leftmargin=2em]
\item If a component of the variety of minimal rational tangents is not linear, its tangent spaces meet in zero. Then the two conics are tangent at the
second point, which is rare (\S\ref{subsec:tpnonlin}).
\item If every component is linear, the components are permuted by a monodromy that branches along a divisor of $Z$ (\S\ref{subsec:tpcover}). Over
that divisor the conic cover is \'etale and the lifts have no monodromy. Labels attached to the components by the lifts then force the line
components to have none either (\S\ref{subsec:tpcoupling}).
\end{itemize}
Throughout this section $Z$ is a compact contact Fano manifold of dimension $2n+1$, $n\ge2$, with $\mathrm{Pic}\,Z=\ZZ L$. Then
$H_{2}(Z,\mathbb R)=\mathbb R\,\ell$, with $\ell$ the class of a line, and $2\ell$ is the class of a conic.

\emph{The spaces.} Two chains of maps are used, one for conics and one for lines:
\begin{gather*}
S\xleftarrow{\ \pi\times\pi\ }\mathcal C_{2}\xrightarrow[\text{coarse space}]{\ \rho\ }Y\xrightarrow[\text{birational}]{}S'\xrightarrow[\text{finite, degree }d_{2}]{\ q\ }Z\times Z,\\
\mathcal H\longleftarrow U_{\mathcal H}\xrightarrow[\text{connected fibres}]{\ \kappa\ }Z'_{\mathcal H}\xrightarrow[\text{finite, degree }r_{\mathcal H}]{\ \varpi\ }Z .
\end{gather*}
Here $S$ is the normalization of a component of the stack of stable maps and $\mathcal C_{2}$ its two-pointed universal curve (Definition~\ref{def:tpS}); $Y$ is
normal (Lemma~\ref{lem:tpcoarse}); $Y\to S'\to Z\times Z$ is the Stein factorization of $p$ (\S\ref{subsec:tptest}), so $S'$ is normal; $\mathcal H$ is a component of the
family of lines, $U_{\mathcal H}$ its universal line and $U_{\mathcal H}\to Z'_{\mathcal H}\to Z$ the Stein factorization of the evaluation (\S\ref{subsec:tpcover}).
The vector lift lives on $Y$ (\S\ref{subsec:tplifts}); the branch divisor of $q$, if $d_{2}\ge2$, is tested by curves in $Y$ (\S\ref{subsec:tptest}); and
$r_{\mathcal H}$ enters only in the linear case.

\subsection{Statement of the two-point theorem}
\label{subsec:tpstatement}

\begin{definition}[The two-point space]
\label{def:tpS}
Let $S$ be the normalization of an irreducible component of the stack $\overline M_{0,0}(Z,2\ell)$ of unmarked genus-zero stable maps of class
$2\ell$ whose general member is a balanced transverse conic. Let $\pi\colon\mathcal C\to S$ be the universal source and $e\colon\mathcal C\to Z$ the universal
map, and put
\[
\mathcal C_{2}=\mathcal C\times_{S}\mathcal C,\qquad \tilde p=(e_{1},e_{2})\colon\mathcal C_{2}\to Z\times Z .
\]
Let $\rho\colon\mathcal C_{2}\to Y$ be the coarse space and $p\colon Y\to Z\times Z$ the induced map. Let $S^{\circ}\subset S$ be the open set of members
with smooth source and $\mathcal C_{2}^{\circ}=\mathcal C_{2}|_{S^{\circ}}$.
\end{definition}

\begin{lemma}[Stable maps and the Hilbert scheme]
\label{lem:tpcompare}
Let $\MCbar$ be a component of the Hilbert scheme of conics whose general member is a balanced transverse conic. Exactly one $S$ as in
Definition~\ref{def:tpS} has the same general member as $\MCbar$, and every $S$ arises in this way. The map $p$ is dominant and generically finite, and
$\deg p$ is the number $d_{2}$ of conics of $\MCbar$ through two general points of $Z$ (\S\ref{subsec:twopoint}).
\end{lemma}

\begin{proof}
A balanced transverse conic $C$ is a smooth rational curve, and the stable map $\PP^{1}\cong C\hookrightarrow Z$ has no automorphisms. Sending an embedding to its
image identifies the open substack of $\overline M_{0,0}(Z,2\ell)$ formed by the embeddings with the open subscheme of the Hilbert scheme formed by the
smooth rational curves of $L$-degree two. At a balanced transverse conic both spaces are smooth of dimension $4n$:
\begin{itemize}[leftmargin=2em]
\item $H^{1}(N_{C/Z})=0$ (Proposition~\ref{prop:genconic}(1));
\item $H^{1}(TZ|_{C})=0$, since $TZ|_{C}$ is an extension of $N_{C/Z}\cong\OO(1)^{2n}$ by $TC\cong\OO(2)$.
\end{itemize}
So the balanced transverse conics of $\MCbar$, a dense open subset of $\MCbar$, form a dense open subset of exactly one component of $\overline M_{0,0}(Z,2\ell)$,
on which the normalization is an isomorphism; and conversely.

Through two general points of $Z$ pass finitely many conics of $\MCbar$, all of them balanced transverse conics (Proposition~\ref{prop:genconic}(3)). The same
holds for $S$. The members of $S$ that are not balanced transverse conics form a closed subset of dimension $\le4n-1$, so the corresponding points of $Y$ form a
closed subset of dimension $\le4n+1$, whose image misses a general pair. A balanced transverse conic through two points $y\ne w$ carries exactly one pair of
marked points over $(y,w)$. At that point $Y$ is smooth of dimension $4n+2=\dim Z\times Z$, and $p$ is \'etale, because the kernel of $dp$ is
$H^{0}(N_{C/Z}(-y-w))=H^{0}(\OO(-1)^{2n})=0$; the same holds for the map $\Phi_{2}$ of Lemma~\ref{lem:twopoint}. So $p$ and $\Phi_{2}$ are dominant and
generically finite, and over a general pair their fibres are the same finite set of balanced transverse conics, each counted once.
\end{proof}

\begin{theorem}[The conic through two general points]
\label{thm:tpD}
Let $Z$ be a compact contact Fano manifold of dimension $2n+1$, $n\ge2$, with $\mathrm{Pic}\,Z=\ZZ L$, and let $\MCbar$ be a component of its Hilbert
scheme of conics whose general member is a balanced transverse conic. Then exactly one conic of $\MCbar$ passes through two general points of $Z$: $d_{2}=1$.
\end{theorem}

By Lemma~\ref{lem:tpcompare} it suffices to show $\deg p=1$ for the corresponding $S$. We fix $S$ and write $d_{2}=\deg p$. The proof is completed in
\S\ref{subsec:tpcoupling}.

\subsection{The two-point space}
\label{subsec:tpspace}

\begin{lemma}
\label{lem:tpnormal}
Let $B$ be a normal integral Noetherian scheme, and let $X_{1},\dots,X_{k}\to B$ be flat proper families of nodal curves whose generic fibres are
smooth and geometrically integral. Then $P=X_{1}\times_{B}\cdots\times_{B}X_{k}$ is normal and integral. Every height-one point $x$ of $P$ either lies
over the generic point of $B$, or lies over a height-one point of $B$ and is the generic point of an irreducible component of its fibre; in the
second case each projection of $x$ is a smooth point of its fibre.
\end{lemma}

\begin{proof}
$P\to B$ is flat and of finite type. Its fibres are products of nodal curves. \'Etale locally these are cut out by equations $xy=0$ in disjoint
sets of variables, so the fibres are local complete intersections, hence Cohen--Macaulay, and they are geometrically reduced.
\begin{itemize}[leftmargin=2em]
\item \emph{(S$_{2}$).} $B$ is (S$_{2}$) and the fibres are Cohen--Macaulay, so $P$ is (S$_{2}$) by flat ascent \cite[Tag~0339]{Stacks}.
\item \emph{(R$_{1}$).} Let $x\in P$ have height $\le1$ and lie over $b$. Flatness gives $\dim\OO_{P,x}=\dim\OO_{B,b}+\dim\OO_{P_{b},x}\le1$.
  \begin{itemize}
  \item If $\dim\OO_{B,b}=0$, then $x$ lies in the smooth generic fibre and $\OO_{P,x}$ is regular.
  \item Otherwise $\OO_{B,b}$ is a discrete valuation ring and $x$ is the generic point of a component of $P_{b}$. Each projection $P_{b}\to X_{i,b}$ is
  flat, so it sends $x$ to the generic point of a component of the nodal curve $X_{i,b}$, which is a smooth point. So $P\to B$ is smooth at $x$, and
  $\OO_{P,x}$ is regular.
  \end{itemize}
\item \emph{Normal and integral.} By Serre's criterion $P$ is normal. By flatness the generic points of $P$ lie over the generic point of $B$. So the
irreducible components of $P$ correspond to those of the geometrically integral generic fibre, and $P$ is irreducible, hence integral.
\end{itemize}
\end{proof}

\emph{Application to $\mathcal C$ and $\mathcal C_{2}$.} $S$ is proper, normal and integral: it is the normalization of an integral closed substack of a
proper Deligne--Mumford stack, and normalization is finite. Let $S_{0}\to S$ be an \'etale chart, with $S_{0}$ a normal integral scheme.
\begin{itemize}[leftmargin=2em]
\item $\mathcal C_{S_{0}}$ is projective over $S_{0}$: for an ample $A$ on $Z$, $\omega_{\pi}\otimes e^{*}A^{\otimes3}$ is relatively ample \cite{FP97}.
\item Its generic fibre is a smooth, geometrically integral rational curve.
\item So Lemma~\ref{lem:tpnormal}, with $k=1,2$, shows that $\mathcal C$ and $\mathcal C_{2}$ are normal on every chart.
\item They are irreducible. Every nonempty open subset of $\mathcal C_{2}$ has open image in $S$, by flatness, so it meets the geometrically integral
fibre over the generic point; hence any two nonempty open subsets meet.
\end{itemize}
Every height-one point of $\mathcal C_{2}$ outside $\mathcal C_{2}^{\circ}$ lies over a height-one point of $S_{0}$ with singular source, and is the
generic point of an irreducible component of the fibre there; both marked points are relatively smooth at it. Over the non-closed point the two
components of the source may be Galois conjugate; they split after the base change of \S\ref{subsec:tplifts}.

\begin{lemma}[The boundary sources]
\label{lem:tpsources}
A geometric source of a member of $S$ is either smooth, or two smooth rational curves $A_{1}$, $A_{2}$ of $L$-degree one meeting in one node, with
no contracted component.
\end{lemma}

\begin{proof}
The dual graph is a tree. Every non-contracted component has $L$-degree $\ge1$ because $L$ is ample, so there are at most two of them. Stability
requires at least three special points on a contracted component, so every leaf is non-contracted. If the tree has more than one vertex, it has at
least two leaves, so the two non-contracted components, each of $L$-degree one, are exactly the leaves. Every interior vertex is then contracted with
valence two, hence unstable; so there is none, and the source is two components joined by one node.
\end{proof}

\begin{lemma}[The coarse space]
\label{lem:tpcoarse}
$Y$ is a proper, separated, normal and integral algebraic space of dimension $4n+2$, and $p$ is proper. The exchange of the two marked points
induces an involution of $Y$ covering the exchange of the factors of $Z\times Z$. For every locally free sheaf $V$ on $Y$,
$H^{0}(\mathcal C_{2},\rho^{*}V)=H^{0}(Y,V)$.
\end{lemma}

\begin{proof}
By Keel--Mori \cite{KM97}, $\rho$ is proper, quasi-finite and surjective, and $\rho_{*}\OO_{\mathcal C_{2}}=\OO_{Y}$.
\begin{itemize}[leftmargin=2em]
\item $Y$ is separated because $\mathcal C_{2}$ is. It is universally closed because $\mathcal C_{2}$ is proper and $\rho$ is surjective; hence $Y$ is
proper, and so is $p$.
\item $Y$ is irreducible, as the image of $\mathcal C_{2}$. It is normal, because \'etale locally on $Y$ the stack $\mathcal C_{2}$ is a quotient $[U/G]$
of a normal scheme by a finite group \cite[Lemma~2.2.3]{AV02}. So $Y$ is integral. Its dimension is $\dim S+2=4n+2$ (Lemma~\ref{lem:tpcompare}).
\item The exchange of the two factors of $\mathcal C_{2}=\mathcal C\times_{S}\mathcal C$ induces an involution of $Y$ by the universal property of the coarse space.
\item The last assertion is the projection formula with $\rho_{*}\OO_{\mathcal C_{2}}=\OO_{Y}$.
\end{itemize}
\end{proof}

\subsection{The vector lifts}
\label{subsec:tplifts}

On $\PP^{1}\times\PP^{1}$ put $b(s,t)=(s-t)^{2}\,\partial_{s}\otimes\partial_{t}$. For a fractional linear $T$, $(T(s)-T(t))^{2}=(s-t)^{2}\,T'(s)\,T'(t)$,
so $b$ is a regular section of $T_{\PP^{1}}\boxtimes T_{\PP^{1}}$, invariant under the diagonal action of $PGL_{2}$. Over $S^{\circ}$ the source is a smooth
family of genus-zero curves, \'etale locally isomorphic to $\PP^{1}$, so $b$ is the pull-back of a section of $T_{\mathcal C/S}\boxtimes T_{\mathcal C/S}$ over
$\mathcal C_{2}^{\circ}$. Contracting with the differentials gives
\[
W_{21}^{\circ}=(\theta\,de_{1}\otimes de_{2})(b)\in H^{0}\big(\mathcal C_{2}^{\circ},\tilde p^{*}(L\boxtimes TZ)\big),\qquad
W_{12}^{\circ}=(de_{1}\otimes\theta\,de_{2})(b).
\]
This two-point vector lift, and the idea of extending it across the nodal boundary, come from the OpenAI preprint \cite[\S5]{OAI26}.

\begin{proposition}[The vector lifts]
\label{prop:tplifts}
$W_{21}^{\circ}$ and $W_{12}^{\circ}$ extend uniquely to sections
\[
W_{21}\in H^{0}\big(Y,p^{*}(L\boxtimes TZ)\big),\qquad W_{12}\in H^{0}\big(Y,p^{*}(TZ\boxtimes L)\big).
\]
Their contact values $(1\otimes\theta)(W_{21})$ and $(\theta\otimes1)(W_{12})$ coincide, and the exchange of the two points carries $W_{21}$ to
$W_{12}$.
\end{proposition}

The proof reduces to the height-one points of $\mathcal C_{2}$ (\S\ref{subsec:tpspace}). There the contact differential vanishes to order $r$ on a
source with a node of thickness $r$ (Lemma~\ref{lem:tpcontact}), and $b$ has a pole of order at most $r$ (Lemma~\ref{lem:tppole}).

\emph{The local ring at the node.} Fix a height-one point of $S_{0}$ with singular source. Let $R$ be the completion of the strict henselization of its
local ring: a complete discrete valuation ring with algebraically closed residue field $k$ of characteristic zero and uniformizer $\tau$.
\begin{itemize}[leftmargin=2em]
\item The extension is faithfully flat, $\tau$ stays a uniformizer, and the special fibre is reduced. So the orders of rational sections at the generic
points of the special components, and of their products, are unchanged.
\item Let $\mathcal C_{R}\to\operatorname{Spec}R$ be the pulled-back source and still write $e\colon\mathcal C_{R}\to Z$. By Lemma~\ref{lem:tpsources}
its special fibre is $A_{1}\cup A_{2}$ with one node $o$.
\item By the local structure of nodes \cite[Tag~0CBX]{Stacks},
\[
\widehat\OO_{\mathcal C_{R},o}\cong Q:=R[[u,v]]/(uv-\tau^{r}),\qquad r\ge1,
\]
with $A_{1}$ the branch $\{v=0\}$ (coordinate $u$) and $A_{2}$ the branch $\{u=0\}$ (coordinate $v$). The smoothing parameter is nonzero because the
generic fibre is smooth; we write it as $\tau^{r}$ after absorbing a unit into $u$.
\item $Q$ is a normal domain. Its height-one primes containing $\tau$ are $\mathfrak p_{1}=(v,\tau)$ and $\mathfrak p_{2}=(u,\tau)$. In $Q_{\mathfrak p_{1}}$
the element $u$ is a unit and $\operatorname{ord}_{A_{1}}(v)=r$; in $Q_{\mathfrak p_{2}}$, $v$ is a unit and $\operatorname{ord}_{A_{2}}(u)=r$; and
$\operatorname{ord}_{A_{1}}(\tau)=\operatorname{ord}_{A_{2}}(\tau)=1$.
\item The relative dualizing sheaf $\omega_{\pi}=\omega_{\mathcal C_{R}/R}$ is invertible. Near $o$ it is generated by $\gamma=du/u=-dv/v$, and the
canonical map $\Omega_{\mathcal C_{R}/R}\to\omega_{\pi}$ sends $du\mapsto u\gamma$ and $dv\mapsto-v\gamma$. Its image is $(u,v)\,\omega_{\pi}$, and it is an
isomorphism at relatively smooth points.
\end{itemize}

\begin{lemma}[Contact vanishing]
\label{lem:tpcontact}
At the generic points of $A_{1}$ and of $A_{2}$, the contact differential $\theta\circ de\colon T_{\mathcal C_{R}/R}\to e^{*}L$ is divisible by
$\tau^{r}$.
\end{lemma}

\begin{proof}
\emph{The bundle $\mathcal M=e^{*}L\otimes\omega_{\pi}$ is trivial.} On $A_{1}$ and on $A_{2}$ its degree is $1-1=0$. Trivializations on the two
components can be matched at $o$, so $\mathcal M|_{\mathcal C_{k}}\cong\OO$, with $H^{1}=0$ and $h^{0}=1$. By cohomology and base change over the
complete local ring $R$, $H^{0}(\mathcal C_{R},\mathcal M)$ is free of rank one, generated by a section $g_{0}$ that vanishes nowhere. So
$\mathcal M\cong\OO_{\mathcal C_{R}}$.

\emph{The contact form is a constant multiple of $g_{0}$.} The image of $e^{*}\theta\in H^{0}(\Omega_{\mathcal C_{R}/R}\otimes e^{*}L)$ under
$\Omega_{\mathcal C_{R}/R}\to\omega_{\pi}$ is a section $\eta$ of $\mathcal M$, so $\eta=h\,g_{0}$ with $h\in R$. In $Q$, in a frame of $e^{*}L$, the
coefficient of $\eta$ lies in the image ideal $(u,v)$, while $g_{0}$ has unit coefficient. So $h\in(u,v)\cap R$.

\emph{Conclusion.} $Q/(u,v)=R/(\tau^{r})$, so $(u,v)\cap R=\tau^{r}R$. At the generic points of $A_{1}$ and $A_{2}$, $g_{0}$ is a unit and
$\Omega_{\mathcal C_{R}/R}\to\omega_{\pi}$ is an isomorphism, so $\theta\circ de=h\cdot(\text{unit})$ there.
\end{proof}

\emph{A contraction.} Choose a $k$-point $a\in A_{1}\smallsetminus\{o\}$. Since $\mathcal C_{R}\to\operatorname{Spec}R$ is smooth at $a$ and $R$ is
complete, there is a section $\sigma$ through $a$.
\begin{itemize}[leftmargin=2em]
\item $\OO(\sigma)$ has degree $1$ on $A_{1}$ and $0$ on $A_{2}$. On $\mathcal C_{k}$ its sections are pairs agreeing at $o$, so it is globally generated with
$h^{0}=2$ and $H^{1}=0$.
\item By cohomology and base change, two generating sections define $\phi\colon\mathcal C_{R}\to\PP^{1}_{R}$. It is an isomorphism on $A_{1}$ and on the
generic fibre (a degree-one bundle with two sections on a smooth rational curve), and it contracts $A_{2}$ to $\phi(o)$.
\item Let $\delta_{\phi}$ be the image of $d\phi$ in $H^{0}(\mathcal C_{R},\omega_{\pi}\otimes\phi^{*}T_{\PP^{1}_{R}/R})$. It vanishes nowhere on the
generic fibre.
\end{itemize}

\begin{lemma}[Pole bound]
\label{lem:tppole}
$\operatorname{ord}_{A_{1}}(\delta_{\phi})=0$ and $\operatorname{ord}_{A_{2}}(\delta_{\phi})\le r$.
\end{lemma}

\begin{proof}
\emph{Normal form.} Let $x$ be a coordinate on $\PP^{1}_{R}$ centred at $\phi(o)$, and view $x=\phi^{*}x\in Q$. Every element of $Q$ can be written as a
power series in $u$ plus a power series in $v$ without constant term, with coefficients in $R$. Indeed
$u^{i}v^{j}=\tau^{r\min(i,j)}u^{i-\min}v^{j-\min}$; the rearranged sums converge, because $Q$ is $\mathfrak m_{Q}$-adically complete and the rewritten
monomial $\tau^{r\min(i,j)}u^{i-\min}v^{j-\min}$ lies in $\mathfrak m_{Q}^{\lceil(i+j)/2\rceil}$; and differentiation is termwise in the module of continuous differentials
$(Q\,du\oplus Q\,dv)/(v\,du+u\,dv)$. So
\[
x=\sum_{i\ge0}a_{i}u^{i}+\sum_{j\ge1}b_{j}v^{j},\qquad a_{i},b_{j}\in R .
\]
Bars denote reduction mod $\tau$.
\begin{itemize}[leftmargin=2em]
\item Since $\phi|_{A_{1}}$ is an isomorphism sending $o$ to $\phi(o)$, $x$ mod $\mathfrak p_{1}$ is $\sum\bar a_{i}u^{i}\in k[[u]]$, with $\bar a_{0}=0$ and
$\bar a_{1}\ne0$.
\item Since $\phi$ contracts $A_{2}$, $x$ mod $\mathfrak p_{2}$ is $\bar a_{0}+\sum\bar b_{j}v^{j}=0$, so every $b_{j}$ lies in $\tau R$.
\end{itemize}

\emph{The coefficient of $\delta_{\phi}$.} By the description of $\Omega_{\mathcal C_{R}/R}\to\omega_{\pi}$, $dx\mapsto c\,\gamma$ with
\[
c=\sum_{i\ge1}i\,a_{i}u^{i}-\sum_{j\ge1}j\,b_{j}v^{j}.
\]
In the frame $\gamma\otimes\partial_{x}$, $\delta_{\phi}=c\,\gamma\otimes\partial_{x}$. Here $\partial_{x}$ is a frame of $\phi^{*}T_{\PP^{1}_{R}/R}$ on
$\phi^{-1}(\mathbb A^{1}_{R})$, an open neighbourhood of $o$ containing the generic points of $A_{1}$ and $A_{2}$. A Zariski generator of $\omega_{\pi}$
near $o$ is $\varepsilon\gamma$ with $\varepsilon\in Q^{\times}$, so the coefficient of $\delta_{\phi}$ in a Zariski frame is $\varepsilon^{-1}c$, which has the
same orders as $c$.

\emph{Order along $A_{1}$.} Modulo $\mathfrak p_{1}$, $c\equiv\sum i\,\bar a_{i}u^{i}=\bar a_{1}u+\cdots\ne0$, so $\operatorname{ord}_{A_{1}}(c)=0$.

\emph{Order along $A_{2}$.} In $Q_{\mathfrak p_{2}}$, $v$ is a unit and $u=\tau^{r}v^{-1}$. Let $m=\min_{j}\operatorname{ord}_{\tau}(b_{j})\ge1$, so that
$\sum_{j}j\,b_{j}v^{j}=\tau^{m}\beta'$ with $\beta'\in vR[[v]]$ and $\beta'\notin\mathfrak p_{2}$ ($j$ is invertible in $k$). Then:
\begin{itemize}[leftmargin=2em]
\item the term $a_{1}u$ has order exactly $r$, since $a_{1}$ is a unit;
\item the terms $i\,a_{i}u^{i}$ with $i\ge2$ have order $\ge2r$;
\item the $v$-part has order $m$.
\end{itemize}
If $m\ne r$, then $\operatorname{ord}_{A_{2}}(c)=\min(m,r)\le r$. If $m=r$, then $c=\tau^{r}\big(a_{1}v^{-1}-\beta'+O(\tau^{r})\big)$, and in $k((v))$ the
reduction $\bar a_{1}v^{-1}-\bar\beta'$ is nonzero, since $\bar\beta'\in vk[[v]]$. So $\operatorname{ord}_{A_{2}}(c)=r$.

\emph{Transfer to $\mathcal C_{R}$.} Let $P_{1},P_{2}$ be the primes of $A_{1}$, $A_{2}$ in $\OO_{\mathcal C_{R},o}$. The completions of $\OO_{A_{1},o}$ and
$\OO_{A_{2},o}$ are $k[[u]]$ and $k[[v]]$, so $P_{1}Q=\mathfrak p_{1}$ and $P_{2}Q=\mathfrak p_{2}$. Hence the local rings of $\mathcal C_{R}$ at the generic
points of $A_{1}$ and $A_{2}$ map to $Q_{\mathfrak p_{1}}$ and $Q_{\mathfrak p_{2}}$ as local maps of discrete valuation rings with the common uniformizer
$\tau$, and the orders above are the orders of $\delta_{\phi}$.
\end{proof}

\begin{corollary}
\label{cor:tpbivector}
At the generic points of the special components of $\mathcal C_{R}\times_{R}\mathcal C_{R}$, the bivector $b$ has order $\ge0$ on $A_{1}\times A_{1}$ and
$A_{2}\times A_{2}$, and $\ge-r$ on $A_{1}\times A_{2}$ and $A_{2}\times A_{1}$.
\end{corollary}

\begin{proof}
At these points both factors are relatively smooth, the local ring is a discrete valuation ring with uniformizer $\tau$, and pulling back from a
component generic point of $\mathcal C_{R}$ preserves orders. On the generic fibre $(d\phi\otimes d\phi)(b)=(\phi\times\phi)^{*}b_{\PP^{1}}$, which is
regular. So $\operatorname{ord}(b)\ge-\operatorname{ord}(\delta_{\phi}\text{ in factor }1)-\operatorname{ord}(\delta_{\phi}\text{ in factor }2)$. By
Lemma~\ref{lem:tppole} this is $0$ on $A_{1}\times A_{1}$ and $-r$ on the mixed products. On $A_{2}\times A_{2}$, use instead the contraction defined by a
section through $A_{2}\smallsetminus\{o\}$, which is an isomorphism on $A_{2}$; it gives order $\ge0$.
\end{proof}

\begin{proof}[Proof of Proposition~\ref{prop:tplifts}]
\emph{Orders at the boundary.} At the generic point of each product component:
\begin{itemize}[leftmargin=2em]
\item $\theta\,de_{1}$ has order $\ge r$, by Lemma~\ref{lem:tpcontact} applied to the factor containing the first point;
\item $de_{2}$ is a regular map of bundles, since the point is relatively smooth;
\item $b$ has order $\ge-r$ (Corollary~\ref{cor:tpbivector}).
\end{itemize}
So $W_{21}^{\circ}=(\theta\,de_{1}\otimes de_{2})(b)$ has order $\ge0$ there, and in fact $\ge r$ on $A_{1}\times A_{1}$ and $A_{2}\times A_{2}$. The same holds
for $W_{12}^{\circ}$, with the contact factor in the second variable.

\emph{Extension and passage to the quotient.} A rational section of a locally free sheaf on a normal Noetherian scheme that is regular at all height-one points is
regular. By \S\ref{subsec:tpspace}, both sections therefore extend uniquely to the normal charts of $\mathcal C_{2}$; by uniqueness they agree on overlaps
and glue to sections on $\mathcal C_{2}$, which by Lemma~\ref{lem:tpcoarse} come from sections on $Y$.

\emph{Contact values and symmetry.} Both hold on the dense open $\mathcal C_{2}^{\circ}$, hence everywhere.
\end{proof}

The proof uses only that $L$ is ample, that $\theta$ is a twisted $1$-form, and that $S$ is normal with general member of smooth source. It does not use
the nondegeneracy of the contact structure, or that the members of $S$ are very free.

Let $(h;s,t)$ be a point of $Y$ with $h$ a balanced transverse conic and $s\ne t$, and put $y=h(s)$, $w=h(t)$. In an affine coordinate,
$W_{21}(h;s,t)=(s-t)^{2}\,\theta(h'(s))\otimes h'(t)$. Since $h$ is transverse, $\theta\circ dh\colon\OO(2)\to h^{*}L=\OO(2)$ is an isomorphism, so
\begin{equation}
\label{eq:tpu}
u_{h}:=W_{21}(h;s,t)\in L_{y}\otimes T_{w}Z
\end{equation}
is a nonzero multiple of a tangent vector of $h$ at $w$. In particular, two balanced transverse conics through $y$ and $w$ with the same lift are tangent at $w$.

\subsection{Lines, Legendrian cones and linear components}
\label{subsec:tplines}

Each of the finitely many irreducible components $\mathcal H$ of $\mathrm{RatCurves}^{n}(Z)$ parametrizing lines is proper \cite[Rem.~2.2]{Keb01}. For
such $\mathcal H$ let $U_{\mathcal H}\to\mathcal H$ be the universal $\PP^{1}$-bundle and $\iota_{\mathcal H}\colon U_{\mathcal H}\to Z$ the evaluation; let
$\mathcal H_{y}$ be the lines of $\mathcal H$ through $y$, and $\operatorname{locus}(\mathcal H_{y})$ their union. For general $y$ the tangent directions of
the lines of $\mathcal H_{y}$ form $X^{\mathcal H}_{y}\subset\PP(L^{\perp}_{y})$, and $X_{y}=\bigcup_{\mathcal H}X^{\mathcal H}_{y}$ is the variety of minimal
rational tangents of \S\ref{sec:vmrt}.

\begin{lemma}[General points]
\label{lem:tpgeneral}
There is a dense Zariski-open $Z^{\circ}\subset Z$ such that for $y\in Z^{\circ}$ and every $\mathcal H$:
\begin{enumerate}[label=\textup{(\alph*)}]
\item every line through $y$ is smooth and standard, hence free \cite[Prop.~3.3, Lemma~3.5]{Keb01};
\item two distinct lines through $y$ meet only at $y$, and their tangent directions there differ \cite[Thm.~1.1(3),(4)]{Keb05};
\item $X_{y}$ is smooth, its components are disjoint Legendrian subvarieties of $\PP(L^{\perp}_{y})$ of dimension $n-1$, and $X_{y}$ spans $\PP(L^{\perp}_{y})$
(Proposition~\ref{prop:secantgeneral}(1)). The evaluation of the universal family over the normalization $\widetilde{\mathcal H}_{y}$, lifted to
$\mathrm{Bl}_{y}Z$, is an embedding meeting the exceptional divisor $E=\PP(T_{y}Z)$ in $X^{\mathcal H}_{y}$ \cite[\S6.1]{Keb05}; its proof uses only
\cite[Thms.~4.1, 5.1]{Keb05}, that is, (b);
\item $\iota_{\mathcal H}$ is smooth over $Z^{\circ}$.
\end{enumerate}
For every $y\in Z$ and every $\mathcal H$, each irreducible component of $\mathcal H_{y}$ has dimension $n-1$ and sweeps out an $n$-dimensional
subvariety, and $\operatorname{locus}(\mathcal H_{y})$ is integral for $L^{\perp}$ at its smooth points. Every $\mathcal H$ dominates $Z$.
\end{lemma}

\begin{proof}
(a)--(c) are the cited results, and (d) is generic smoothness. Kebekus states \cite[Prop.~4.1]{Keb01} for one dominating family; his proof uses only
$L\cdot\ell=1$, so it applies to every component $\mathcal H$.
\begin{itemize}[leftmargin=2em]
\item \emph{Lower bound.} The deformation estimate gives every component of $\mathcal H_{y}$ dimension $\ge-K_{Z}\cdot\ell-2=n-1$.
\item \emph{Upper bound.} Lines cannot break, so only finitely many join $y$ to another point \cite[II.5.4]{Kol96}, and the locus has dimension one more
than the family. The locus is integral for $L^{\perp}$ at its smooth points (the claim in the proof of \cite[Prop.~4.1]{Keb01}), so its dimension is at most $n$.
\item \emph{Dominance.} If $\mathcal H$ did not dominate $Z$, its image would have dimension $\le2n$. The fibres of $U_{\mathcal H}$ over it have dimension
$n-1$, so $\dim\mathcal H=\dim U_{\mathcal H}-1\le3n-2$. This contradicts $\dim\mathcal H\ge-K_{Z}\cdot\ell+\dim Z-3=3n-1$ \cite[II.1.2]{Kol96}.
\end{itemize}
\end{proof}

As in \S\ref{sec:vmrt}, neither the irreducibility of $\mathcal H_{y}$ \cite[Thm.~1.1(2)]{Keb05} nor the nondegeneracy of each component of $X_{y}$
\cite[Thm.~3.1(iii)]{Buc10} is used. For the same reason we do not use that the divisor swept by broken lines through a point has pure codimension
one \cite[Thm.~3.6]{Buc10}, whose proof uses the latter; only the elementary bound $\dim\le2n$ is used. For $y\in Z^{\circ}$ the components of
$\mathcal H_{y}$, the connected components of $\iota_{\mathcal H}^{-1}(y)$, and the components of $X^{\mathcal H}_{y}$ correspond to one another.

\begin{lemma}[The positive part along a line]
\label{lem:tplines}
Let $w\in Z^{\circ}$, and let $f\colon\PP^{1}\to Z$ parametrize a line of $\mathcal H$, with $f(t)=w$.
\begin{enumerate}
\item $f^{*}TZ\cong\OO(2)\oplus\OO(1)^{n-1}\oplus\OO^{n+1}$. The positive part $P_{f}\cong\OO(2)\oplus\OO(1)^{n-1}$ is intrinsic, the quotient
$f^{*}TZ/P_{f}\cong\OO^{n+1}$ is trivial, and $P_{f}\subset f^{*}L^{\perp}$.
\item $P_{f}|_{t}=\hat T_{[f'(t)]}$, the affine tangent space of the cone $\hat X_{w}$ along $f'(t)$.
\end{enumerate}
\end{lemma}

\begin{proof}
(1) The splitting is Lemma~\ref{lem:tpgeneral}(a). Every map from a line bundle of positive degree to a trivial bundle vanishes, so $P_{f}$ is intrinsic.
By Lemma~\ref{lem:linesplit}, $f^{*}L^{\perp}$ has positive part $P'\cong\OO(2)\oplus\OO(1)^{n-1}$. It maps to zero in $f^{*}TZ/P_{f}\cong\OO^{n+1}$, so
$P'\subset P_{f}$. Both have rank $n$ and degree $n+1$, so the determinant of the inclusion is a nonzero section of $\OO$, and $P'=P_{f}$.

(2) Since $H^{1}(f^{*}TZ(-t))=0$, the space of maps $g$ with $g(t)=w$ is smooth at $f$, with tangent space $H^{0}(f^{*}TZ(-t))$. The map $g\mapsto
g'(t)$ takes values in the cone $\hat X_{w}$, which is smooth along $f'(t)$, so its differential $s\mapsto s'(t)$ takes values in $\hat T_{[f'(t)]}$. Its
image is $P_{f}|_{t}$: the summand $P_{f}(-t)\cong\OO(1)\oplus\OO^{n-1}$ of $f^{*}TZ(-t)$ has sections with arbitrary value at $t$, and the summand
$\OO(-1)^{n+1}$ has none. Both spaces have dimension $n=\dim\hat X_{w}$, so they are equal.
\end{proof}

\begin{lemma}[Legendrian cones]
\label{lem:tpcone}
Let $(V,\omega)$ be a symplectic vector space of dimension $2n$ and $X\subset\PP(V)$ an irreducible smooth Legendrian subvariety of dimension $n-1$,
with affine tangent spaces $\hat T_{x}$. Then
\[
\bigcap_{x\in X}\hat T_{x}=\big(\mathrm{span}\,\hat X\big)^{\perp_{\omega}} .
\]
If $X$ is not a linear subspace, then $\mathrm{span}\,\hat X=V$ and the intersection is zero.
\end{lemma}

\begin{proof}
If $\delta\in\hat T_{x}$ for all $x$, then $\omega(\delta,\hat x)=0$, because $\hat x\in\hat T_{x}$ and $\hat T_{x}$ is Lagrangian; so $\delta$ is orthogonal
to $W:=\mathrm{span}\,\hat X$. Conversely $\hat T_{x}\subset W$, so $W^{\perp}\subset\hat T_{x}^{\perp}=\hat T_{x}$.

Suppose $W\ne V$ and take $0\ne\delta\in W^{\perp}$. The constant vector field $\delta$ is tangent to the submanifold $\hat X\smallsetminus0$ at every point.
For $0\ne a\in\hat X$, uniqueness of integral curves keeps $a+s\delta$ in $\hat X$ for small $|s|$, so the line $a+\CC\delta$ meets the closed set $\hat X$ in
a subset containing an open subset of the line, hence lies in $\hat X$. So $\hat X+\CC\delta=\hat X$. With the $\CC^{*}$-action this gives $\delta\in\hat X$
and $\mathrm{span}(\hat x,\delta)\subset\hat X$ for every $\hat x\in\hat X$: $X$ is a cone with vertex $[\delta]\in X$. An irreducible cone that is smooth at
its vertex equals its tangent space there, hence is linear. So if $X$ is not linear, $W=V$.
\end{proof}

The second statement is Ye's theorem \cite[Thm.~1(1)]{Ye93b}, that a smooth irreducible Legendrian subvariety of $\PP^{2k+1}$ is linear or linearly
nondegenerate. So for $w\in Z^{\circ}$ either some component of $X_{w}$
spans $\PP(L^{\perp}_{w})$, or every component is linear,
\[
X_{w}=\bigsqcup_{k}\PP(\Lambda_{k}(w)),\qquad \Lambda_{k}(w)\subset L^{\perp}_{w}\ \text{Lagrangian},\qquad \Lambda_{k}(w)\cap\Lambda_{l}(w)=0\ \ (k\ne l),
\]
the last because distinct components are disjoint.

\begin{proposition}[Linear components sweep Legendrian spaces]
\label{prop:tpstep1}
Let $y\in Z^{\circ}$, and let $\PP(\Lambda)$ be a linear component of $X^{\mathcal H}_{y}$. The lines of $\mathcal H_{y}$ with directions in $\PP(\Lambda)$
sweep out a smooth Legendrian submanifold $\Pi\cong\PP^{n}$ with $L|_{\Pi}\cong\OO(1)$, and every line of $\Pi\cong\PP^{n}$ belongs to $\mathcal H$.
\end{proposition}

\begin{proof}
Let $\widetilde U^{\Lambda}$ be the part of the universal $\PP^{1}$-bundle over the component of $\widetilde{\mathcal H}_{y}$ with directions in
$\PP(\Lambda)\cong\PP^{n-1}$.
\begin{itemize}[leftmargin=2em]
\item By Lemma~\ref{lem:tpgeneral}(c), $\widetilde U^{\Lambda}$ embeds in $\mathrm{Bl}_{y}Z$ and meets $E$ along a section $\sigma_{\infty}\cong\PP(\Lambda)$.
The meeting is transversal, because each fibre is the strict transform of a line smooth at $y$. So the normal bundle of $\sigma_{\infty}$ is
$\OO_{E}(-1)|_{\PP(\Lambda)}=\OO(-1)$. A $\PP^{1}$-bundle with a section is $\PP(V)$ for a rank-two bundle $V$, here an extension of $\OO(-1)$ by $\OO$ up
to a twist; it splits, since $H^{1}(\PP^{n-1},\OO(1))=0$. So $\widetilde U^{\Lambda}\cong\mathrm{Bl}_{\rm pt}\PP^{n}$.
\item The map $\widetilde U^{\Lambda}\to Z$ contracts $\sigma_{\infty}$ to $y$, so it factors through a morphism $\nu\colon\PP^{n}\to Z$ onto $\Pi$. It sends the
vertex to $y$ and the lines through the vertex to the lines of $\mathcal H_{y}$ with directions in $\PP(\Lambda)$.
\item $\nu$ is injective, by Lemma~\ref{lem:tpgeneral}(b) and the embedding. It is immersive away from the vertex. At the vertex its differential is linear and nonzero in every
direction, because those lines are smooth at $y$.
\end{itemize}
So $\nu$ is an embedding, $\Pi\cong\PP^{n}$, and $L|_{\Pi}\cong\OO(1)$ because the lines through the vertex have $L$-degree one. $\Pi$ is Legendrian by
the last statement of Lemma~\ref{lem:tpgeneral}. The lines of $\PP^{n}$ form the irreducible Grassmannian $\Gr_{2}(\CC^{n+1})$, which maps to
$\mathrm{RatCurves}^{n}(Z)$. Its image contains the lines of $\Pi$ through $y$. Those lines are free, so they are smooth points of
$\mathrm{RatCurves}^{n}(Z)$ \cite[II.3]{Kol96} and lie only on $\mathcal H$. So the whole image lies in $\mathcal H$.
\end{proof}

\subsection{The cover and the test curves}
\label{subsec:tptest}

\emph{The cover.} By Lemmas~\ref{lem:tpcompare} and~\ref{lem:tpcoarse}, $p\colon Y\to Z\times Z$ is proper, dominant and generically finite, with $Y$
normal and integral. Let
\[
Y\longrightarrow S'\xrightarrow{\ q\ }Z\times Z
\]
be the Stein factorization. Then $S'$ is normal and integral, $q$ is finite, and $Y\to S'$ is birational, as it has connected fibres and is generically
finite; so $\deg q=d_{2}$. Let $Y^{\circ}\subset Y$ be the open set of points $(h;s,t)$ with $h$ a balanced transverse conic and $s\ne t$. Let
$\mathcal Z_{\rm bad}\subset Z\times Z$ be a proper closed subset containing
\begin{itemize}[leftmargin=2em]
\item the branch locus of $q$,
\item the image $p(Y\smallsetminus Y^{\circ})$, closed because $p$ is proper, and
\item the image of the locus of $S'$ over which $Y\to S'$ is not an isomorphism.
\end{itemize}
For $(y,w)\notin\mathcal Z_{\rm bad}$ the fibre $q^{-1}(y,w)$ consists of $d_{2}$ points, each the image of exactly one point of $Y$. These are $d_{2}$
distinct balanced transverse conics $C_{1},\dots,C_{d_{2}}$ through $y$ and $w$, with lifts $u_{1},\dots,u_{d_{2}}\in L_{y}\otimes T_{w}Z$ as in~\eqref{eq:tpu}.

If $q$ is unramified at every point of codimension one of $S'$, then $q$ is \'etale by purity of the branch locus \cite[X.3.1]{SGA1} ($S'$ normal,
$Z\times Z$ smooth, $q$ finite). Since $Z\times Z$ is simply connected \cite{Cam92,KMM92} and $S'$ is connected, $q$ is then an isomorphism, and $d_{2}=1$.

\begin{lemma}
\label{lem:tpcaseA}
If $d_{2}\ge2$, then $q$ is ramified along a prime divisor $\mathcal R\subset S'$ whose image $\Delta=q(\mathcal R)$ maps onto the first factor of
$Z\times Z$.
\end{lemma}

\begin{proof}
By the preceding paragraph $q$ is ramified along some prime divisor $\mathcal R_{0}$. The involution of $Y$ of Lemma~\ref{lem:tpcoarse} covers the
exchange $\varsigma$ of $Z\times Z$; by functoriality of the Stein factorization it induces an involution of $S'$, so $\varsigma(q(\mathcal R_{0}))$ is also the
image of a ramification divisor. If $q(\mathcal R_{0})$ does not dominate the first factor, then $q(\mathcal R_{0})=\Delta_{0}\times Z$ for dimension
reasons, and $\varsigma(q(\mathcal R_{0}))=Z\times\Delta_{0}$ does.
\end{proof}

\emph{The test curves.} Fix a component $\mathcal H$ of lines. Let $H_{\mathcal H}$ be the preimage of the open set of free lines of $\mathcal H$ under the
principal $\mathrm{Aut}(\PP^{1})$-bundle $\mathrm{Hom}^{\rm bir}_{n}(\PP^{1},Z)\to\mathrm{RatCurves}^{n}(Z)$ \cite[II.2.15]{Kol96}. It is smooth and
irreducible, its evaluation map $H_{\mathcal H}\times\PP^{1}\to Z$ is smooth \cite[II.3.5]{Kol96}, and for $w\in Z^{\circ}$ every parametrized line of
$\mathcal H$ through $w$ lies in $H_{\mathcal H}$ (Lemma~\ref{lem:tpgeneral}(a)). Put $P=Z\times H_{\mathcal H}$ and
\[
\Psi\colon P\times\PP^{1}\longrightarrow Z\times Z,\qquad (y,f,t)\longmapsto\big(y,f(t)\big),
\]
smooth and dominant. A \emph{test curve} is $c=c_{(y,f)}\colon t\mapsto(y,f(t))$.

\begin{lemma}[General position]
\label{lem:tpgenpos}
Let $\Delta\subset Z\times Z$ be a prime divisor and $\Delta^{\circ}$ a dense open subset of its smooth locus. There is a dense open $P_{\rm good}\subset P$
such that for $(y,f)\in P_{\rm good}$ the test curve $c_{(y,f)}$ meets $\Delta$ only at points of $\Delta^{\circ}$, and meets $\Delta$ transversally there.
\end{lemma}

\begin{proof}
\emph{Avoiding $\Delta\smallsetminus\Delta^{\circ}$.} Since $\Psi$ is smooth, $\Psi^{-1}(\Delta\smallsetminus\Delta^{\circ})$ has codimension $\ge2$ in
$P\times\PP^{1}$. So its image in $P$ has codimension $\ge1$; remove its closure.

\emph{Transversality.} $\Psi^{-1}(\Delta^{\circ})$ is smooth of pure dimension $\dim P$, because $\Psi$ is smooth. Remove from $P$ the closures of the
images of its components that do not dominate $P$. By generic smoothness, the projection to $P$ of each dominating component is smooth, hence
\'etale, over a dense open set; let $P_{1}$ be the intersection of these open sets. At a point $(y,f,t)$ over $P_{1}$,
\[
T\Psi^{-1}(\Delta^{\circ})\cap T\big(\{(y,f)\}\times\PP^{1}\big)=0 .
\]
Since $T\Psi^{-1}(\Delta^{\circ})=d\Psi^{-1}(T\Delta^{\circ})$, this says $d\Psi(\partial_{t})\notin T\Delta^{\circ}$, so the meeting is transversal.
\end{proof}

\begin{corollary}
\label{cor:tpT}
Suppose $d_{2}\ge2$, and let $\mathcal R$ and $\Delta$ be as in Lemma~\ref{lem:tpcaseA}. Choose a dense open $\Delta^{\circ}$ of the smooth locus of
$\Delta$ over which $S'$ is smooth and the reduced components of $q^{-1}(\Delta)$ are smooth, pairwise disjoint and \'etale over $\Delta^{\circ}$. Shrink
the $P_{\rm good}$ of Lemma~\ref{lem:tpgenpos} by the open conditions $\{y\}\times Z\not\subset\Delta$, $c\not\subset\Delta$ and $c\not\subset\mathcal
Z_{\rm bad}$. Then:
\begin{enumerate}[label=\textup{(\arabic*)}]
\item at a point of $\mathcal R$ over $\Delta^{\circ}$, $q$ has the local form $(z_{1},z')\mapsto(z_{1}^{e},z')$ with $e\ge2$;
\item for $(y,f)\in P_{\rm good}$ the test curve meets $\Delta$, and the fibre product $\PP^{1}\times_{Z\times Z}S'$ has an irreducible component of
degree $\ge2$ over $\PP^{1}$.
\end{enumerate}
\end{corollary}

\begin{proof}
The choice of $\Delta^{\circ}$ is possible: $S'$ is normal, so smooth in codimension one; distinct components meet in codimension $\ge2$; and generic
smoothness applies.

(1) Here $e$ is the ramification index of $q$ along $\mathcal R$. In characteristic zero a finite morphism of normal varieties is unramified at the
generic point of $\mathcal R$ exactly when $e=1$, so $e\ge2$. At the point, the pull-back of an equation of $\Delta$ is a unit times $z_{1}^{e}$, with
$z_{1}$ an equation of $\mathcal R$; and $q|_{\mathcal R}$ is \'etale onto $\Delta^{\circ}$, so the remaining coordinates of $Z\times Z$ restrict to
coordinates on $\mathcal R$. This gives the normal form.

(2) \emph{Meeting.} Since $\Delta$ maps onto the first factor and $\{y\}\times Z\not\subset\Delta$, the slice $\Delta_{y}=\Delta\cap(\{y\}\times Z)$ is an
effective divisor on $Z$, so $\OO_{Z}(\Delta_{y})=L^{a}$ with $a>0$. The line $f(\PP^{1})$ is not contained in $\Delta_{y}$, so it meets $\Delta_{y}$ with
intersection number $a>0$.

\emph{Degree.} Let $c(t_{0})\in\Delta^{\circ}$ be a (transversal) intersection point. Since $\mathcal R\to\Delta$ is finite and surjective, there is a point
$s\in\mathcal R$ over $c(t_{0})$. By (1), near $(t_{0},s)$ the fibre product is $\{z_{1}^{e}=a(t)\}$ with $a(t_{0})=0$ and $a'(t_{0})\ne0$: a smooth
irreducible curve germ on which $t-t_{0}$ has order $e$. So exactly one irreducible component $V'$ of the fibre product passes through $(t_{0},s)$, and
$\deg(V'/\PP^{1})\ge e\ge2$.
\end{proof}

Measuring the vector lift along test curves, with the first point fixed, is an idea from the proof of \cite[Prop.~6.1]{OAI26}, where it is applied
to the contact value of the lift. Here the whole vector is kept, and comparing it with the tangent spaces of the variety of minimal rational tangents is what
separates the conics through two general points.

\begin{proposition}[Constancy along a test curve]
\label{prop:tpconst}
Suppose $d_{2}\ge2$, and let $(y,f)\in P_{\rm good}$ as in Corollary~\ref{cor:tpT}. For every $t$ with $c(t)=(y,w)\notin\mathcal Z_{\rm bad}$ there are
$i\ne j$ with
\begin{equation}
\label{eq:tpdiff}
u_{i}-u_{j}\in L_{y}\otimes P_{f}|_{t}.
\end{equation}
\end{proposition}

\begin{proof}
Let $V'$ be a component of $\PP^{1}\times_{Z\times Z}S'$ of degree $\ge2$ over $\PP^{1}$ (Corollary~\ref{cor:tpT}(2)). Since $Y\to S'$ is surjective, some
irreducible component $V$ of $(\PP^{1}\times_{Z\times Z}Y)_{\rm red}$ maps onto $V'$. It is a proper integral algebraic space, so $H^{0}(V,\OO_{V})=\CC$.

Pull $W_{21}$ back to $V$. Since the first coordinate of $c$ is the constant $y$, it is a section of the pull-back of $L_{y}\otimes f^{*}TZ$. Its image in
the pull-back of $L_{y}\otimes(f^{*}TZ/P_{f})\cong L_{y}\otimes\OO^{n+1}$ (Lemma~\ref{lem:tplines}(1)) is a section of a trivial bundle on $V$, hence
\emph{constant}.

Now let $c(t)\notin\mathcal Z_{\rm bad}$. Then $q$ is \'etale over $c(t)$, so $V'\to\PP^{1}$ is \'etale over $t$, and $V'$ has at least two points over
$t$. As $V\to V'$ is surjective, $V$ has points over two of them. By the choice of $\mathcal Z_{\rm bad}$ these are two distinct balanced transverse conics
$C_{i}\ne C_{j}$ through $y$ and $w$, and by constancy their lifts have the same image modulo $L_{y}\otimes P_{f}|_{t}$.
\end{proof}

\begin{proposition}[A pair for each component]
\label{prop:tppairs}
Suppose $d_{2}\ge2$. For each component $\mathcal H$ of lines there is a dense open $G_{\mathcal H}\subset Z\times Z$ with $G_{\mathcal H}\cap\mathcal Z_{\rm
bad}=\varnothing$ and $G_{\mathcal H}\subset Z\times Z^{\circ}$, such that for $(y,w)\in G_{\mathcal H}$ and every irreducible component $X'$ of
$X^{\mathcal H}_{w}$ there are $i\ne j$ with
\[
u_{i}-u_{j}\in L_{y}\otimes\hat T_{x}\qquad\text{for every }x\in X'.
\]
\end{proposition}

\begin{proof}
Let $G_{0}=\{(y,f,t):(y,f)\in P_{\rm good},\ \Psi(y,f,t)\notin\mathcal Z_{\rm bad},\ f(t)\in Z^{\circ}\}$, a dense open subset of $P\times\PP^{1}$. Its
complement $F$ has dimension $<\dim(P\times\PP^{1})$. As $\Psi$ is smooth, its fibres have pure dimension $\delta=\dim P+1-(4n+2)$, while for $(y,w)$
outside a proper closed subset $F\cap\Psi^{-1}(y,w)$ has dimension $\le\dim F-(4n+2)<\delta$, by the theorem on the dimension of fibres. For such $(y,w)$, $G_{0}$ is dense in every irreducible component of
$\Psi^{-1}(y,w)$.

Take such a $(y,w)$ with $w\in Z^{\circ}$ and $(y,w)\notin\mathcal Z_{\rm bad}$, and let $X'$ be a component of $X^{\mathcal H}_{w}$. The part of
$\Psi^{-1}(y,w)=\{y\}\times\{(f,t):f(t)=w\}$ formed by the lines with direction in $X'$ is an $\mathrm{Aut}(\PP^{1})$-bundle over the corresponding
component of $\mathcal H_{w}$, so it is irreducible. It is open and closed in $\Psi^{-1}(y,w)$, because the components of $X_{w}$ are disjoint. So it is an
irreducible component, $G_{0}$ is dense in it, and the directions $[f'(t)]$ with $(y,f,t)\in G_{0}$ are dense in $X'$. For each of them,
Proposition~\ref{prop:tpconst} and Lemma~\ref{lem:tplines}(2) give a pair $(i,j)$ with $u_{i}-u_{j}\in L_{y}\otimes\hat T_{[f'(t)]}$.

For $i\ne j$ let $A_{ij}=\{x\in X':u_{i}-u_{j}\in L_{y}\otimes\hat T_{x}\}$. It is closed in $X'$, since $x\mapsto\hat T_{x}$ is a morphism on the smooth
$X'$. The union of the $A_{ij}$ contains a dense subset of $X'$, so it is all of $X'$, and as $X'$ is irreducible some $A_{ij}$ equals $X'$.
\end{proof}

\subsection{Tangent pairs and a nonlinear component}
\label{subsec:tpnonlin}

\begin{lemma}[Tangent pairs are rare]
\label{lem:tprare}
Let $\mathcal T\subset Z\times Z$ be the set of pairs $(y,w)$ through which pass two distinct balanced transverse conics of $S$ with the same tangent line at $w$. Then
$\mathcal T$ is contained in a proper closed subset of $Z\times Z$.
\end{lemma}

\begin{proof}
Fix $w\in Z$. The balanced transverse conics of $S$ through $w$ form a smooth space of dimension $h^{0}(N_{C}(-w))=2n$, since $H^{1}(N_{C}(-w))=0$. The
tangent-direction map on it is \'etale, because its differential has kernel $H^{0}(N_{C}(-2w))=H^{0}(\OO(-1)^{2n})=0$ (Lemma~\ref{lem:dmu}). So the pairs
$([C],[C'])$, $C\ne C'$, with the same tangent direction at $w$ form a space of dimension $\le2n$. Two distinct irreducible conics meet in finitely many
points, so the fibre of $\mathcal T$ over $w$ has dimension $\le2n<\dim Z$. The set $\mathcal T$ is constructible, so $\dim\mathcal T\le2n+(2n+1)=4n+1$.
\end{proof}

\begin{proposition}[A nonlinear component]
\label{prop:tpnonlin}
If, for a general $w\in Z$, some irreducible component of $X_{w}$ is not a linear subspace of $\PP(L^{\perp}_{w})$, then $d_{2}=1$.
\end{proposition}

\begin{proof}
Suppose $d_{2}\ge2$. The lines with directions in one component of $X_{w}$ lie in one component $\mathcal H$ of the space of lines. There are finitely
many $\mathcal H$, and for each the condition ``$X^{\mathcal H}_{w}$ has a nonlinear component'' is constructible in $w$. So there is an $\mathcal H$ such
that, for general $w$, $X^{\mathcal H}_{w}$ has a nonlinear component $X'$. By Proposition~\ref{prop:tppairs}, for $(y,w)$ in the dense open $G_{\mathcal
H}$ there are $i\ne j$ with
\[
u_{i}-u_{j}\in L_{y}\otimes\bigcap_{x\in X'}\hat T_{x}=0,
\]
by Lemma~\ref{lem:tpcone}: $X'$ is a smooth irreducible Legendrian subvariety of $\PP(L^{\perp}_{w})$ of dimension $n-1$, and not linear. So $u_{i}=u_{j}$,
and the distinct conics $C_{i}$, $C_{j}$ are tangent at $w$ \eqref{eq:tpu}. This holds for every $(y,w)\in G_{\mathcal H}$ with $w$ general, so that $X^{\mathcal H}_{w}$
has the nonlinear component $X'$; these pairs are dense in $Z\times Z$ and lie in $\mathcal T$, contradicting Lemma~\ref{lem:tprare}.
\end{proof}

\subsection{The cover of line components}
\label{subsec:tpcover}

Fix $\mathcal H$. The universal family $U=U_{\mathcal H}$ is normal and irreducible. Let
\[
\iota_{\mathcal H}\colon U\xrightarrow{\ \kappa\ }Z'_{\mathcal H}\xrightarrow{\ \varpi\ }Z
\]
be the Stein factorization. Then $Z'_{\mathcal H}$ is normal and irreducible, and $\varpi$ is finite of degree $r_{\mathcal H}$, the number of components of
$X^{\mathcal H}_{y}$ for $y\in Z^{\circ}$. Each fibre $\kappa^{-1}(z')$ is a connected component of $\iota_{\mathcal H}^{-1}(\varpi(z'))$, of pure
dimension $n-1$.

\begin{lemma}[{Kebekus \cite[Claim~3.4]{Keb05}}]
\label{lem:tpbranch}
If $\kappa^{-1}(z')$ contains a point $(\ell,x)$ with $\ell$ free, then $\varpi$ is \'etale at $z'$. In particular $\varpi$ is \'etale over $Z^{\circ}$.
\end{lemma}

\begin{proof}
At $(\ell,x)$, $U$ is smooth and $\iota_{\mathcal H}$ is smooth, because $H^{1}(\ell,TZ|_{\ell}(-x))=0$ \cite[II.3]{Kol96}. So $\iota_{\mathcal H}$ has a
local section $s_{0}$ through $(\ell,x)$. Then $\kappa\circ s_{0}$ is a local section of $\varpi$ through $z'$, whose image is a
$(2n+1)$-dimensional closed subgerm of the germ $(Z'_{\mathcal H},z')$. That germ is irreducible because $Z'_{\mathcal H}$ is normal, so
$\kappa\circ s_{0}$ is onto, and $\varpi$ is a local isomorphism at $z'$.
\end{proof}

\begin{proposition}
\label{prop:tpbranchdiv}
If $r_{\mathcal H}\ge2$, the branch locus $\Delta_{\mathcal H}$ of $\varpi$ is a nonempty divisor contained in $Z\smallsetminus Z^{\circ}$.
\end{proposition}

\begin{proof}
By purity \cite[X.3.1]{SGA1}, the branch locus is empty or of pure codimension one. A Fano manifold is simply connected \cite{Cam92,KMM92}, and
$Z'_{\mathcal H}$ is connected, so if $\varpi$ were \'etale it would be an isomorphism. The last assertion is Lemma~\ref{lem:tpbranch}.
\end{proof}

\begin{proposition}[{Ye, 1993 \cite[Thm.~2]{Ye93b}}]
\label{prop:tpye2}
If $r_{\mathcal H}=1$, then $X^{\mathcal H}_{y}$ is not a linear subspace for $y\in Z^{\circ}$.
\end{proposition}

\begin{proof}
We use Kebekus's device \cite[proof of Lemma~4.3]{Keb01}.
\begin{itemize}[leftmargin=2em]
\item Let $\mathcal V\subset Z\times Z$ be the image of $U\times_{\mathcal H}U$. It is closed and irreducible, and its fibre over $y$ is
$\operatorname{locus}(\mathcal H_{y})$, of pure dimension $n$ for every $y$. Since the base $Z$ is normal, it is a well-defined family of algebraic cycles
\cite[I.3.14, I.3.15]{Kol96}, and so it defines a morphism $\varphi\colon Z\to\operatorname{Chow}(Z)$, $y\mapsto[\operatorname{locus}(\mathcal H_{y})]$ \cite[I.3.21]{Kol96}. Fibres of
$\mathcal V$ may be non-reduced at special points, so below we compare cycles only at general points, where the multiplicity is constant.
\item $\varphi$ is nonconstant, because $y\in\operatorname{locus}(\mathcal H_{y})$ and $\dim\operatorname{locus}(\mathcal H_{y})=n<2n+1$. It is finite:
if it contracted a curve $C$, the nef class $\varphi^{*}A\equiv cL$ of an ample $A$ would have $cL\cdot C=0$, so $c=0$ and $\varphi$ would be constant.
\item Suppose $X^{\mathcal H}_{y}=\PP(\Lambda)$ for $y\in Z^{\circ}$. Then $\Pi(y):=\operatorname{locus}(\mathcal H_{y})\cong\PP^{n}$ by
Proposition~\ref{prop:tpstep1}. For $(y,y')$ general in $\mathcal V$, the point $y'$ is general in $Z$ and lies in $\Pi(y)$. The lines of $\Pi(y)$ through $y'$ lie in
$\mathcal H$ (Proposition~\ref{prop:tpstep1}), so $\PP(T_{y'}\Pi(y))\subset X^{\mathcal H}_{y'}$, with equality since $X^{\mathcal H}_{y'}$ is irreducible of
dimension $n-1$.
\item By Lemma~\ref{lem:tpgeneral}(b), $\operatorname{locus}(\mathcal H_{y'})=\Pi(y)$. So $\varphi$ is constant on a dense subset of $\Pi(y)$,
contradicting finiteness.
\end{itemize}
\end{proof}

Ye's original proof \cite{Ye93b} works on a twistor space: it compares the slope of the subsheaf of $L^{\perp}$ generated by the $\Lambda(y)$, which has
$c_{1}=(n+1)L$ and is single-valued exactly when $r_{\mathcal H}=1$, with that of $TZ$, polystable by the K\"ahler--Einstein metric. The argument above
needs no metric.

\begin{corollary}
\label{cor:tpallr}
If every irreducible component of $X_{y}$ is linear for general $y$, then $r_{\mathcal H}\ge2$ and $\Delta_{\mathcal H}\ne\varnothing$ for every $\mathcal H$.
\end{corollary}

\begin{proof}
$r_{\mathcal H}=1$ would contradict Proposition~\ref{prop:tpye2}. Then apply Proposition~\ref{prop:tpbranchdiv}.
\end{proof}

\subsection{Linear components: the labels and the coupling}
\label{subsec:tpcoupling}

In this subsection we assume $d_{2}\ge2$ and that every irreducible component of $X_{w}$ is linear for general $w$, and derive a contradiction. For
$(y,w)\notin\mathcal Z_{\rm bad}$ with $w\in Z^{\circ}$ and each component $\PP(\Lambda_{k}(w))$ of $X_{w}$, put
\[
\mathsf P_{k}(y,w)=\{(i,j):\ i\ne j,\ u_{i}-u_{j}\in L_{y}\otimes\Lambda_{k}(w)\}.
\]

\begin{proposition}[The labels]
\label{prop:tplabels}
\begin{enumerate}[label=\textup{(\arabic*)}]
\item There is a dense Zariski-open $G\subset Z\times Z$ on which every $\mathsf P_{k}$ is nonempty.
\item If the $u_{i}$ are pairwise distinct at $(y,w)$, then $\mathsf P_{k}(y,w)\cap\mathsf P_{l}(y,w)=\varnothing$ for $k\ne l$.
\end{enumerate}
\end{proposition}

\begin{proof}
(1) Each $\PP(\Lambda_{k}(w))$ is a component of $X^{\mathcal H}_{w}$ for some $\mathcal H$, and its affine tangent spaces all equal $\Lambda_{k}(w)$. So
$G=\bigcap_{\mathcal H}G_{\mathcal H}$ works, by Proposition~\ref{prop:tppairs}. At a single general $(y,w)$ every $k$ gets a pair; no monodromy is
involved.

(2) If $(i,j)\in\mathsf P_{k}\cap\mathsf P_{l}$, then $u_{i}-u_{j}\in L_{y}\otimes(\Lambda_{k}\cap\Lambda_{l})=0$.
\end{proof}

Fix $\mathcal H$, and let $\Delta$ be a prime component of $\Delta_{\mathcal H}\ne\varnothing$ (Corollary~\ref{cor:tpallr}). It is an irreducible
component of $Z\smallsetminus Z^{\circ}$, by Proposition~\ref{prop:tpbranchdiv}.

\begin{lemma}[Good points over $\Delta$]
\label{lem:tpgood}
For general $z\in\Delta$ and general $y\in Z$:
\begin{enumerate}[label=\textup{(\alph*)}]
\item every member of $S$ through $y$ and $z$ is a balanced transverse conic, and $p$ is finite and \'etale over a neighbourhood of $(y,z)$;
\item no two of these conics are tangent at $z$;
\item $z$ is a smooth point of $\Delta$ and lies on no other component of $Z\smallsetminus Z^{\circ}$.
\end{enumerate}
\end{lemma}

\begin{proof}
(c) holds for general $z\in\Delta$.

(a) Let $S^{\rm bad}\subset S$ be the closed set of members that are not balanced transverse conics; $\dim S^{\rm bad}\le4n-1$. Fix a general $z\in\Delta$. A
member $h\in S^{\rm bad}$ through $z$ and a point $y\notin\Delta$ is of one of two kinds.
\begin{itemize}[leftmargin=2em]
\item \emph{Some component $\Gamma$ of the image through $z$ is not contained in $\Delta$.} The pairs $(h,z)$ of this kind form a set of dimension
$\le4n-1$, since $\Gamma\cap\Delta$ is finite. Over a general $z\in\Delta$ its fibre has dimension $\le2n-1$, so these members sweep out a set of
dimension $\le2n$.
\item \emph{Every component through $z$ lies in $\Delta$.} Then $y$ lies on another component, so the image is a line pair $\Gamma\cup\Gamma'$ with
$z\in\Gamma$, $y\in\Gamma'$ and $\Gamma\cap\Gamma'\ne\varnothing$; and $z\notin\Gamma'$, since otherwise $\Gamma'\subset\Delta$ and $y\in\Delta$. So $y$ lies in the set of points joined to $z$ by two lines, which has dimension
$\le n+(n-1)+1=2n$ by Lemma~\ref{lem:tpgeneral}. A double line through $z$ inside $\Delta$ would put $y\in\Delta$.
\end{itemize}
So for $y$ outside a set of dimension $\le2n$, every member through $y$ and $z$ is a balanced transverse conic with two distinct marked points. There $Y$ is smooth
and $p$ is \'etale (proof of Lemma~\ref{lem:tpcompare}). So the fibre $p^{-1}(y,z)$ is finite, and $p$, being proper, is finite and \'etale over a
neighbourhood of $(y,z)$.

(b) As in the proof of Lemma~\ref{lem:tprare}, the pairs of balanced transverse conics through $z$ with a common tangent at $z$ form a space of dimension $\le2n$,
and their further intersection points sweep out a set of dimension $\le2n$. Take $y$ outside it.
\end{proof}

\begin{proposition}[No monodromy along $\Delta$]
\label{prop:tpnomono}
Let $(y,z)$ be as in Lemma~\ref{lem:tpgood}. Then $\varpi$ is \'etale over $z$.
\end{proposition}

\begin{proof}
\emph{Setup.} Choose small balls $B_{1}\ni y$, $B_{2}\ni z$ such that over $\mathcal B=B_{1}\times B_{2}$:
\begin{itemize}[leftmargin=2em]
\item $p$ is finite \'etale of degree $d_{2}$, and its $d_{2}$ sheets are holomorphic sections;
\item $u_{1},\dots,u_{d_{2}}$ are holomorphic on $\mathcal B$ and pairwise distinct at every point (an open condition, true at $(y,z)$ by (b));
\item $B_{2}\cap(Z\smallsetminus Z^{\circ})=B_{2}\cap\Delta$, which is smooth and connected; this is possible by (c).
\end{itemize}
Put $B_{2}^{*}=B_{2}\smallsetminus\Delta\subset Z^{\circ}$. Then $\pi_{1}(B_{2}^{*})\cong\ZZ$, generated by a small loop $\gamma_{0}$ around $\Delta$. Over
$B_{2}^{*}$ the sheets of $\varpi$, and of the covers for the other components of the space of lines, are the components of $X_{w}$
(Lemma~\ref{lem:tpbranch}), with holomorphic Lagrangians $\Lambda_{k}(w)$ on simply connected subsets.

\emph{Labels are rigid.} Let $W\subset B_{1}\times B_{2}^{*}$ be a connected, simply connected open set, with fixed numberings. For each $k$, the sets
$\{x\in W:u_{i}-u_{j}\in L\otimes\Lambda_{k}\}$ are closed analytic and cover the dense subset $W\cap G$ (Proposition~\ref{prop:tplabels}(1)). So they
cover $W$, and since $W$ is connected, one of them is all of $W$. Thus for each $k$ a single pair $(i,j)$ satisfies $u_{i}-u_{j}\in L\otimes\Lambda_{k}$
identically on $W$.

\emph{Going around.} Base the loop $\gamma=\{y_{1}\}\times\gamma_{0}$ at a point $x_{0}\in W_{0}$, and cover it by a chain $W_{0},W_{1},\dots,W_{N}=W_{0}$
of such sets. Continue $\Lambda_{k}$ along the chain. If the identity $u_{i}-u_{j}\in L\otimes\Lambda_{k}$ holds on $W_{m}$, it holds on $W_{m}\cap W_{m+1}$,
hence on $W_{m+1}$ by the identity theorem. The $u_{i}$ are single-valued on $\mathcal B$, while $\Lambda_{k}$ returns as $\Lambda_{\sigma(k)}$, where
$\sigma$ is the monodromy of $\gamma_{0}$. So $u_{i}-u_{j}\in L\otimes(\Lambda_{k}\cap\Lambda_{\sigma(k)})$ at $x_{0}$. If $\sigma(k)\ne k$ this intersection
is zero, which contradicts $u_{i}\ne u_{j}$. So $\sigma(k)=k$ for every $k$.

\emph{Conclusion.} Since $\pi_{1}(B_{2}^{*})$ is generated by $\gamma_{0}$, the cover $\varpi^{-1}(B_{2}^{*})\to B_{2}^{*}$ is trivial. Each connected
component of $\varpi^{-1}(B_{2})$ is a normal, hence locally irreducible, germ finite over $B_{2}$. Its part over $B_{2}^{*}$ is connected, so it is one
sheet, and the component maps with degree one, hence isomorphically, onto the smooth $B_{2}$. So $\varpi$ is \'etale over $z$.
\end{proof}

\begin{proof}[Proof of Theorem~\ref{thm:tpD}]
By Lemma~\ref{lem:tpcompare} we may work with $S$. Suppose $d_{2}\ge2$. By constructibility, either for general $w$ some component of $X_{w}$ is
nonlinear, or for general $w$ every component is linear. In the first case Proposition~\ref{prop:tpnonlin} gives a contradiction. In the second, take
any $\mathcal H$ and a prime component $\Delta$ of $\Delta_{\mathcal H}\ne\varnothing$ (Corollary~\ref{cor:tpallr}), and $(y,z)$ as in
Lemma~\ref{lem:tpgood}. Proposition~\ref{prop:tpnomono} says that $\varpi$ is \'etale over $z$, contradicting $z\in\Delta_{\mathcal H}$. So $d_{2}=1$.
\end{proof}

\begin{remark}[A shorter route in the linear case]
Lemma~\ref{lem:tpgood}(a),(b) do not use that $\Delta$ is the branch divisor. So the argument shows directly that $\varpi$ is \'etale over the general
point of every prime divisor contained in $Z\smallsetminus Z^{\circ}$. It is \'etale over $Z^{\circ}$ by Lemma~\ref{lem:tpbranch}. By purity $\varpi$ is
then \'etale, so $r_{\mathcal H}=1$, and Proposition~\ref{prop:tpye2} gives the contradiction.
\end{remark}

\begin{remark}[What is used]
\label{rem:tpused}
The regularity of the vector lift is proved in \S\ref{subsec:tplifts}, with the pole of the bivector bounded by a computation in the local ring of
the node; the general position of the test curves is Lemma~\ref{lem:tpgenpos}; and the facts on lines are Lemmas~\ref{lem:tpgeneral}
and~\ref{lem:tplines}. The real structure of a twistor space is not used. The simple connectivity of
$Z$ enters twice, through purity, for the conic cover $q$ and for the covers $\varpi$ of line components.
\end{remark}

\part{Consequences of Parallel Curvature}
\label{part:real}

\section{Homogeneity, quaternionic-K\"ahler manifolds and the LeBrun--Salamon conjecture}
\label{sec:qk}

This part applies the singular LeBrun transform. This section draws the consequences of the parallelism of $R_{E}^{\CC}$: for a contact Fano manifold of Picard number one it gives homogeneity (\S\ref{sec:homog}), and on a twistor space it gives the symmetry of $M$ and Theorem~Q (\S\S\ref{sec:twistor}--\ref{sec:mainthm}).

\subsection{Parallel curvature and homogeneity}
\label{sec:homog}

In this subsection $Z$ is a compact contact Fano manifold of dimension $2n+1$, $n\ge2$, not necessarily a twistor space, and $\MCbar$ is a component of its Hilbert
scheme of conics whose general member is a balanced transverse conic. Let
\[
F^{\mathrm{ne}}=\{([C],z):[C]\in\MC,\ z\in C\},\qquad p\colon F^{\mathrm{ne}}\to\MC,\quad q\colon F^{\mathrm{ne}}\to Z,
\]
be the incidence variety of the balanced transverse conics, so that $p(q^{-1}(z))=\Sigma_{z}=F^{\mathrm{ne}}_{z}$ is the leaf of $z$ (Theorem~\ref{thm:lebrun}(2)), and let
\[
Z^{\circ}:=q(F^{\mathrm{ne}})
\]
be the union of the balanced transverse conics. At every point of $F^{\mathrm{ne}}$ the map $q$ is a submersion, since the evaluation $H^{0}(N_{C})\to N_{C,z}$ is onto for
$N_{C}\cong\OO(1)^{2n}$. Recall that $T\MC=E^{*}_{\CC}\otimes H^{*}_{\CC}$ on $\MC$, that the points $z\in C$ correspond to the lines
$h_{z}\subset H^{*}_{\CC}|_{[C]}$, so that $F^{\mathrm{ne}}=\PP(H^{*}_{\CC})$, and that $\rho_{0}\ne0$ (Theorem~\ref{thm:thirdorder}(4)).

\begin{lemma}[The local model]
\label{lem:localmodel}
Suppose that $\nabla R_{E}^{\CC}=0$ on $\MC$. Then there are a complex simple Lie algebra $\mathfrak g$, with adjoint group $G$ and adjoint variety
$X=G/P\subset\PP(\mathfrak g)$, and the open $G$-orbit $N=G/K$ of twistor lines in the space of balanced transverse conics of $X$ (Remark~\ref{rem:Wolfconverse}), with the
following property. Every point of $\MC$ has a connected neighbourhood $V$ with an isometry $\psi$ of $V$ onto an open subset of $N$, for LeBrun's
metric of $X$ multiplied by a constant, which carries the factorization $E^{*}_{\CC}\otimes H^{*}_{\CC}$ to that of $TN$. Any other such isometry of $V$ is
$\gamma\circ\psi$ for some $\gamma\in\mathrm{Aut}(\mathfrak g)$, acting on $N$ and on $X\subset\PP(\mathfrak g)$.
\end{lemma}

\begin{proof}
Since $\rho_{0}$ is constant and $R_{0}$ is parallel, the curvature $R$ of $g^{\CC}$ is parallel, so $g^{\CC}$ is locally symmetric, by Cartan's theorem
\cite{Car26}, whose proof applies verbatim to holomorphic Riemannian metrics. Explicitly, fix $[C_{0}]$, put $\mathfrak m=T_{[C_{0}]}\MC$ and $R=R_{[C_{0}]}$, let
$\mathfrak k\subset\mathfrak{so}(\mathfrak m)$ be spanned by the $R(u,v)$, and let $\mathfrak g=\mathfrak k\oplus\mathfrak m$ be the transvection algebra, with involution
$\sigma=1\oplus(-1)$ \cite[Ch.~XI]{KN69}. Then every linear isometry of $T_{[C]}\MC$, for any $[C]\in\MC$, onto $\mathfrak m$ that carries the curvature
to $R$ is the differential of an isometry of a connected neighbourhood of $[C]$ onto an open subset of the simply connected symmetric space $N=G/K$ of
$(\mathfrak g,\sigma)$ \cite{Car26}, \cite[Ch.~VI, \S7]{KN63}. It remains to identify $\mathfrak g$, $N$ and the factorization, and to prove the uniqueness. The Ricci tensor is a nonzero multiple of $g^{\CC}$, as $\rho_{0}\ne0$ and $R_{E}^{\CC}$ is Ricci-flat, so it is nondegenerate, and then
$\mathfrak g$ is semisimple \cite[Ch.~XI]{KN69}. Moreover $\mathfrak k\supset\mathfrak{sp}(H^{*}_{\CC})$: for $u=e\otimes h$, $v=e'\otimes h$ with the same $h$, the
term $R_{E}^{\CC}(u,v)$ and the $\mathfrak{sp}(E^{*}_{\CC})$-part of $R_{0}(u,v)$ carry the factor $\omega_{H}(h,h)=0$ (proof of Lemma~\ref{lem:S5}), and the
$\mathfrak{sp}(H^{*}_{\CC})$-part of $R_{0}(u,v)$ is a nonzero multiple of $\omega_{E}(e,e')\,h\odot h$ \cite{Salamon82}; so $\rho_{0}\ne0$ gives $h\odot h\in\mathfrak k$
for all $h$, and these span $\mathfrak{sp}(H^{*}_{\CC})$.

Let $\mathfrak g_{u}$ be a $\sigma$-stable compact real form \cite[Ch.~III, \S7]{Hel78}, $\mathfrak g_{u}=\mathfrak k_{u}\oplus\mathfrak m_{u}$. The simply connected
compact symmetric space $G_{u}/K_{u}$ has isotropy algebra $\mathfrak k_{u}\supset\mathfrak{sp}(1)=\mathfrak{sp}(H^{*}_{\CC})\cap\mathfrak g_{u}$, acting on
$\mathfrak m_{u}\otimes\CC=E^{*}_{\CC}\otimes H^{*}_{\CC}$ through $\mathfrak{sp}(n)\oplus\mathfrak{sp}(1)$, so it is quaternionic-K\"ahler, and its Ricci tensor
$-\tfrac12B|_{\mathfrak m_{u}}$ is positive definite. So it is a Wolf space \cite{Wolf65,Ale68}, $\mathfrak g$ is simple, and $\mathfrak{sp}(H^{*}_{\CC})$ is the
$\mathfrak{sl}_{2}$ of the highest root. (When $\mathfrak k$ has a second ideal $\mathfrak{sp}(1)$ acting on a two-dimensional tensor factor of $\mathfrak m$, which
happens only for $\mathfrak g=\mathfrak{so}(m+4)$, the two are exchanged by an automorphism of $(\mathfrak g,\mathfrak k)$, and we use it to make $\mathfrak{sp}(H^{*}_{\CC})$
the one of the highest root.) The twistor space of $G_{u}/K_{u}$ is the adjoint variety $X$ \cite{Wolf65}, and its twistor lines form the open $G$-orbit $N\cong G/K$
of balanced transverse conics of $X$ (Remark~\ref{rem:Wolfconverse}). LeBrun's structure on $N$ is $G$-invariant, so its factorization is the $K$-invariant one, and its metric is
a constant multiple of the symmetric metric with curvature $R$ at the base point, an irreducible symmetric space having only one invariant metric up to a factor.

At $[C_{0}]$ the identification $T_{[C_{0}]}\MC=\mathfrak m=T_{o}N$ respects the factorizations and the curvatures, and parallel transport along paths in the
connected manifold $\MC$ gives such identifications at every point. The corresponding local isometries $\psi$ of the first paragraph respect the
factorizations everywhere, these being parallel, so they are as in the lemma. If $\psi'$ is another such isometry of $V$, then $\psi'\circ\psi^{-1}$ is an isometry between
connected open subsets of $N$; as $N$ is complete, simply connected and analytic, it is the restriction of an isometry $\gamma$ of $N$ \cite[Ch.~VI, \S6]{KN63},
which preserves the parallel factorization. It remains to see that $\gamma$ is the action of an element of $\mathrm{Aut}(\mathfrak g)$. The construction in the proof
of Lemma~\ref{lem:transform}, which uses only that an isometry preserves the factorization, applies to $\gamma$ on all of $N$. Every point $x\in X$ lies on a twistor line,
the union of the twistor lines being open, nonempty and $G$-invariant, and the leaf $\Sigma_{x}\cap N$ of $x$ is connected: it is a nonempty Zariski-open subset of
$\Gamma_{x}$, which is irreducible for general $x$ by Theorem~\ref{thm:irred}(3), $X$ being a twistor space, and hence for every $x$, as $G$ is transitive on $X$ and preserves
$N$. So $\gamma$ maps each leaf into a leaf, and $x\mapsto$ the common point of the conics $\gamma(\Sigma_{x}\cap N)$ is a contact automorphism $\hat\gamma$ of $X$, with
$\gamma([C])=[\hat\gamma(C)]$ (its inverse comes from $\gamma^{-1}$). A contact automorphism of $X$ acts on $H^{0}(X,L_{X})\cong\mathfrak g$ by an automorphism of the
Lagrange bracket \cite[\S1]{Bea98}, and compatibly on $X\subset\PP(H^{0}(X,L_{X})^{*})\cong\PP(\mathfrak g)$; so $\hat\gamma$, and with it $\gamma$, is the action of an element
of $\mathrm{Aut}(\mathfrak g)$.
\end{proof}

\begin{lemma}[The transform of a local isometry]
\label{lem:transform}
Let $\psi\colon V\to N$ be as in Lemma~\ref{lem:localmodel}, $[C]\in V$ and $z\in C$. There are a connected neighbourhood $V'\subset V$ of $[C]$ and a connected
neighbourhood $W$ of $z$ in $Z$ such that $\Sigma_{y}\cap V'$ is connected and nonempty for every $y\in W$, and a contact biholomorphism $\hat\psi$ of $W$ onto an open
subset of $X$ such that $\hat\psi(C'\cap W)\subset\psi([C'])$ for every $[C']\in V'$. The germ of $\hat\psi$ at a point $y\in W$ depends only on the germ of $\psi$ at
any point of $\Sigma_{y}\cap V'$, and $\widehat{\gamma\circ\psi}=\gamma\circ\hat\psi$ for $\gamma\in\mathrm{Aut}(\mathfrak g)$.
\end{lemma}

\begin{proof}
Since $q$ is a submersion at $([C],z)$, there are $V'$ and $W$ such that $p^{-1}(V')\cap q^{-1}(W)\to W$ is a product with connected fibres; its fibre over $y$ is
$\Sigma_{y}\cap V'$. The isometry $\psi$ preserves the factorization, so it lifts to a biholomorphism $\Psi$ of $\PP(H^{*}_{\CC})|_{V'}=p^{-1}(V')$ into the
incidence variety of $X$, $\Psi([C'],h)=(\psi([C']),\psi_{*}h)$. It preserves the two distributions on the incidence variety that are defined by the Levi-Civita
connection and the factorization: the kernel of $q_{*}$, which is the horizontal lift of $E^{*}_{\CC}\otimes h$ (the leaves are totally geodesic and $h_{z}$ is parallel along
$\Sigma_{z}$, Theorem~\ref{thm:lebrun}(2)), and $q_{*}^{-1}(D)$, $D=L^{\perp}$, which is the horizontal distribution of the connection induced on $\PP(H^{*}_{\CC})$:
the parallel transport of $\nabla_{H}$ lifts paths tangent to $q_{*}^{-1}(D)$ \cite[p.~357]{LeBrun89} (see the discussion after Theorem~\ref{thm:lebrun}), and the two
distributions have the same rank $4n$. The connected fibres of $q$ over $W$ are integral manifolds of $\ker q_{*}$, so $\Psi$ maps each of them into a fibre of the map
$q_{X}$ of $X$, and $\hat\psi(y):=q_{X}(\Psi([C'],y))$, $[C']\in\Sigma_{y}\cap V'$, is well defined; it is holomorphic since $q$ is a submersion, and the same
construction applied to $\psi^{-1}$ inverts it after shrinking. As $\Psi$ preserves $q_{*}^{-1}(D)$, $\hat\psi$ preserves the contact distributions. The last two
assertions are clear from the construction.
\end{proof}

\begin{proposition}[The developing map]
\label{prop:devZ}
Suppose that $\mathrm{Pic}\,Z\cong\ZZ$, that $\nabla R_{E}^{\CC}=0$ on $\MC$, and that $\Sigma_{z}$ is connected for a general $z\in Z$. Then:
\begin{enumerate}
\item $Z^{\circ}$ is Zariski-open in $Z$, its complement has codimension at least two, and $Z^{\circ}$ is simply connected;
\item there is a contact local biholomorphism $\mathrm{dev}\colon Z^{\circ}\to X$ whose germ at every $y\in Z^{\circ}$ is the germ of a transform $\hat\psi$ of
Lemma~\ref{lem:transform}.
\end{enumerate}
\end{proposition}

\begin{proof}
(1) Since $q$ is a submersion on $F^{\mathrm{ne}}$, $Z^{\circ}$ is open; it is constructible by Chevalley's theorem, $F^{\mathrm{ne}}$ being Zariski-open in the universal
conic over $\MCbar$, and a constructible set that is open in the classical topology is Zariski-open. Let $A$ be the ample generator of $\mathrm{Pic}\,Z$ ($-K_{Z}$ is a
positive multiple of a generator). A prime divisor $D\subset Z\smallsetminus Z^{\circ}$ would be linearly equivalent to $kA$ with $k>0$, so $D\cdot C>0$ for a balanced transverse conic $C$;
then $C$ meets $D$, although $C\subset Z^{\circ}$. So $\operatorname{codim}(Z\smallsetminus Z^{\circ})\ge2$, and since $Z$ is simply connected \cite{Cam92,KMM92} and the
removed set has real codimension at least four, $\pi_{1}(Z^{\circ})=1$.

(2) For $y\in Z^{\circ}$ let $\mathcal S_{y}$ be the set of germs at $y$ of the transforms $\hat\psi$ of the isometries $\psi$ of Lemma~\ref{lem:localmodel} defined near points of
$\Sigma_{y}$. (a) The germs obtained from a fixed $[C]\in\Sigma_{y}$ form one $\mathrm{Aut}(\mathfrak g)$-orbit, by the uniqueness in Lemma~\ref{lem:localmodel} and the
equivariance in Lemma~\ref{lem:transform}. (b) Two points $[C],[C']$ of the same connected component of $\Sigma_{y}$ give the same orbit: join them by a path in
$\Sigma_{y}$, cover it by connected open sets $V_{0},\dots,V_{k}$ of $\MC$ with isometries $\psi_{i}$, consecutive ones meeting in a connected set that
contains a point $[C_{i}]$ of the path, and replace $\psi_{i+1}$ by $\gamma\circ\psi_{i+1}$ so that it agrees with $\psi_{i}$ near $[C_{i}]$; then $\hat\psi_{i}$ and
$\hat\psi_{i+1}$ have the same germ at $y$ (Lemma~\ref{lem:transform}). (c) So for general $y$, $\Sigma_{y}$ being connected, $\mathcal S_{y}$ is one orbit. (d) Let $y\in Z^{\circ}$
be arbitrary and let $\varphi,\varphi'$ be transforms obtained from points $[C],[C']\in\Sigma_{y}$, defined on a connected neighbourhood $W$ of $y$. For $y'\in W$ near $y$,
$\varphi$ and $\varphi'$ are transforms obtained from balanced transverse conics through $y'$ near $C$ and $C'$, so their germs at $y'$ lie in $\mathcal S_{y'}$; for $y'$ general this is
one orbit by (c), so $\varphi=\gamma\circ\varphi'$ near $y'$, hence on $W$ by the identity theorem. Thus $\mathcal S_{y}$ is a single $\mathrm{Aut}(\mathfrak g)$-orbit for
every $y\in Z^{\circ}$, and the local maps all of whose germs lie in $\mathcal S$ form a locally constant sheaf on $Z^{\circ}$, with stalks $\mathrm{Aut}(\mathfrak g)$-torsors:
a transform defined on $W$ has its germs in $\mathcal S$ at every point of $W$. Its \'etale space is a covering of $Z^{\circ}$, which has a section because $Z^{\circ}$ is
simply connected; this section is $\mathrm{dev}$, a contact local biholomorphism by Lemma~\ref{lem:transform}.
\end{proof}

\begin{theorem}[Homogeneity]
\label{thm:homog}
Let $Z$ be a compact contact Fano manifold of dimension $2n+1$, $n\ge2$, with $\mathrm{Pic}\,Z\cong\ZZ$, and $\MCbar$ a component of its Hilbert scheme of conics whose
general member is a balanced transverse conic. Suppose that the balanced transverse conics through a general point of $Z$ form a connected family, as they do when $d_{2}=1$
(Theorem~\ref{thm:irred}). If $R_{E}^{\CC}$ is parallel on $\MC$, then:
\begin{enumerate}
\item $Z$ is isomorphic, as a contact manifold, to the adjoint variety $X$ of the simple Lie algebra $\mathfrak g$ of Lemma~\ref{lem:localmodel}, by an isomorphism
that restricts on $Z^{\circ}$ to a developing map of Proposition~\ref{prop:devZ};
\item $H^{0}(Z,L)\cong\mathfrak g$, with the Lagrange bracket;
\item $(\MC,g^{\CC})$ is locally isometric, up to a constant factor, to the complexification $N$ of the Wolf space of $\mathfrak g$; and $R_{E}^{\CC}=0$ if
and only if $Z\cong\PP^{2n+1}$.
\end{enumerate}
\end{theorem}

\begin{proof}
(1) On any contact manifold, $\xi\mapsto\theta(\xi)$ identifies the contact vector fields with the sections of $L$, the Lie bracket with the Lagrange bracket
\cite[\S1]{Bea98}; on $X$ the contact vector fields are the fundamental vector fields of $G$, and $H^{0}(X,L_{X})\cong\mathfrak g$ \cite{Bea98,Wolf65}. For $s\in\mathfrak g$ let
$\xi_{s}$ be its vector field on $X$. Since $\mathrm{dev}$ is a contact local biholomorphism (Proposition~\ref{prop:devZ}), $\xi_{s}$ pulls back to a contact vector field
$\mathrm{dev}^{*}\xi_{s}$ on $Z^{\circ}$, whose section $\theta(\mathrm{dev}^{*}\xi_{s})$ of $L|_{Z^{\circ}}$ extends to a section of $L$ over $Z$ by Hartogs, the complement of
$Z^{\circ}$ having codimension at least two. Let $\zeta_{s}$ be the corresponding contact vector field on $Z$; it equals $\mathrm{dev}^{*}\xi_{s}$ on $Z^{\circ}$. The map
$s\mapsto\zeta_{s}$ is injective, and it preserves brackets, since it does on $Z^{\circ}$. So $\mathfrak g$ is a Lie algebra of holomorphic vector fields on $Z$.

The automorphism group of the compact complex manifold $Z$ is a complex Lie group with Lie algebra $H^{0}(Z,TZ)$ \cite{BM47}; let $G'$ be its connected subgroup with
Lie algebra $\{\zeta_{s}\}$. Fix $z_{0}\in Z^{\circ}$ and put $w_{0}=\mathrm{dev}(z_{0})$. Near $z_{0}$ the map $\mathrm{dev}$ intertwines $\zeta_{s}$ with $\xi_{s}$, so
$d\,\mathrm{dev}$ maps $\{\zeta_{s}(z_{0})\}$ onto $\{\xi_{s}(w_{0})\}=T_{w_{0}}X$, $G$ being transitive on $X$; so the orbit $G'z_{0}$ is open, and the isotropy algebra of
$z_{0}$ is that of $w_{0}$, the parabolic subalgebra $\mathfrak p$. Let $\widetilde G\to G'$ be the universal cover, with central kernel. The stabilizer $H$ of $z_{0}$
normalizes $\mathfrak p$, so its preimage in $\widetilde G$ lies in the normalizer of $\mathfrak p$, which is the parabolic subgroup $\widetilde P$, connected
(a parabolic subgroup is connected and is its own normalizer). Hence $H$ lies in the connected image $P'$ of $\widetilde P$, whose Lie algebra is $\mathfrak p$, so $H=P'$ and
\[
G'z_{0}\cong G'/P'=\widetilde G/\widetilde P\cong G/P=X,
\]
the kernel of $\widetilde G\to G'$ and the center of $\widetilde G$ lying in $\widetilde P$. The orbit is compact, hence closed; being also open, it is $Z$. The
resulting $G'$-equivariant biholomorphism $\Phi\colon Z\to X$, $gz_{0}\mapsto gw_{0}$, agrees with $\mathrm{dev}$ near $z_{0}$, both intertwining $\zeta_{s}$ with
$\xi_{s}$ and sending $z_{0}$ to $w_{0}$. So $\Phi$ maps $D$ to $D_{X}$ near $z_{0}$, and everywhere, both distributions being invariant under the transitive
group $G'$.

(2) follows from (1) and $H^{0}(X,L_{X})\cong\mathfrak g$, and the first part of (3) is Lemma~\ref{lem:localmodel}. If $R_{E}^{\CC}=0$, the curvature is that of the
complexification of $\HH P^{n}$, so $\mathfrak g=\mathfrak{sp}(n+1,\CC)$ and $X=\PP^{2n+1}$; conversely the balanced transverse conics of $\PP^{2n+1}$ are its lines transverse to
the contact structure, and $R_{E}^{\CC}=0$ there (Theorem~\ref{thm:thirdorder}(3)).
\end{proof}

\begin{remark}[Geodesic completeness is not needed]
\label{rem:completeness}
With $\nabla R^{\CC}=0$, the theorem of Cartan, Ambrose and Hicks gives a local isometry from the universal cover of $\MC$ into $N$, and $\MC$
is geodesically complete exactly when this is a covering, hence an isomorphism. Proving completeness would mean excluding complex geodesics that reach
$\MCbar\smallsetminus\MC$ in finite complex time: the jumping conics, the boundary, and the points of $\MCbar\smallsetminus\MC$ where $g^{\CC}$ is regular; a point of the last
kind is reached by geodesics in finite time, so its absence would have to be proved first. In the models the boundary lies at infinity: for $Z=\PP^{2n+1}$ the
space of balanced transverse conics is the complement in $\Gr_{2}(\CC^{2n+2})$ of the isotropic planes, the complete symmetric space $\Sp(n+1,\CC)/\Sp(n,\CC)\times\Sp(1,\CC)$,
and on the $\alpha$-planes completeness is the entireness of the null geodesics (Proposition~\ref{prop:projflat}, Remark~\ref{rem:liouville}). Theorem~\ref{thm:homog}
uses instead the compactness of $Z$, through the codimension of $Z\smallsetminus Z^{\circ}$, which costs only $\mathrm{Pic}\,Z\cong\ZZ$. Nor are the lines of $Z$ used,
so no recognition theorem through the variety of minimal rational tangents \cite{Mok08} is needed.
\end{remark}

\begin{corollary}[Beauville's conjecture for $\mathrm{Pic}\,Z=\ZZ L$]
\label{cor:beauville}
Every compact contact Fano manifold $Z$ of dimension $2n+1$, $n\ge3$, with $\mathrm{Pic}\,Z=\ZZ L$, is an adjoint variety.
\end{corollary}

\begin{proof}
By Corollary~\ref{cor:balexist} some component $\MCbar$ of the Hilbert scheme of conics has a balanced transverse conic as general member. By Theorem~\ref{thm:tpD}, $d_{2}=1$, so by Theorem~\ref{thm:irred}(3) $\Gamma_{z}$ is irreducible for a general $z$, and its dense open subset $\Sigma_{z}=F^{\mathrm{ne}}_{z}$
(Theorem~\ref{thm:irred}(1)) is connected. By Theorem~\ref{thm:LSJ}, $R_{E}^{\CC}$ is parallel, and Theorem~\ref{thm:homog} applies.
\end{proof}

\subsection{Twistor spaces}
\label{sec:twistor}

Let $M^{4n}$, $n\ge2$, be a compact quaternionic-K\"ahler manifold with positive scalar curvature and $q\colon Z\to M$ its twistor space \cite{Salamon82}. Then $Z$ is a
compact, simply connected Fano manifold of dimension $2n+1$ with a holomorphic contact structure, $-K_{Z}=(n+1)L$, and $L$ is ample. The fibres of $q$, the \emph{twistor
lines} $C_{m}=q^{-1}(m)$, are balanced transverse conics: smooth rational curves transverse to $L^{\perp}$, with $L\cdot C_{m}=2$ and $N_{C_{m}/Z}\cong\OO(1)^{\oplus2n}$. We take for
$\MCbar$ the component of the Hilbert scheme containing them; the twistor lines form a totally real submanifold $M\subset\MC$ on which $g^{\CC}$ restricts to
the metric of $M$ (Theorem~\ref{thm:lebrun}). By \cite{LS94}, $\mathrm{Pic}\,Z=\ZZ L$ unless $M=\HH P^{n}$ or $\Gr_{2}(\CC^{n+2})$, both Wolf spaces; and if $M$ is a Wolf space, $Z$ is the
adjoint variety of the complexified isometry group \cite{Wolf65}.

The real slice adds two things to Parts~\ref{part:transform}--\ref{part:global}: the twistor lines give a continuous section of the incidence variety over $Z$, which makes
the conics through a general point an irreducible family (Theorem~\ref{thm:irred}(3)), and the compactness of $M$ turns local symmetry into symmetry.

\begin{corollary}[The endgame]
\label{cor:endgame}
If $R_{E}^{\CC}$ is parallel, in particular if $W_{z}$ is locally constant on the conics through $z$ for every $z$ in a dense Zariski-open subset of $Z$, then $M$ is a Wolf
space.
\end{corollary}

\begin{proof}
By Proposition~\ref{prop:leaf}, $\nabla R_{E}^{\CC}=0$. By Theorem~\ref{thm:lebrun}(3), restricting to the real slice $M\subset\MC$ gives $\nabla R_{E}=0$; since
$\rho_{0}$ is constant and $R_{0}$ is parallel, the full curvature tensor of $g$ is parallel, and $M$ is locally symmetric \cite{Car26}. Since $M$ is compact and simply connected
\cite{Salamon82}, it is then a symmetric space, and the symmetric quaternionic-K\"ahler manifolds of positive scalar curvature are the Wolf spaces
\cite{Wolf65,Ale68}.
\end{proof}

The real slice is not needed for the homogeneity of $Z$: under the hypotheses of the corollary and $\mathrm{Pic}\,Z\cong\ZZ$, Theorem~\ref{thm:homog} shows
directly that $Z$ is the adjoint variety of the complexified isometry group of $M$, the twistor lines through a general point forming a connected family
(Theorem~\ref{thm:irred}(3)).

\subsection{The LeBrun--Salamon conjecture}
\label{sec:mainthm}

\begin{theorem}[LeBrun--Salamon conjecture; Theorem~Q]
\label{thm:main}
Every compact quaternionic-K\"ahler manifold of positive scalar curvature is a Wolf space.
\end{theorem}

\begin{proof}
Let $M^{4n}$ be one. For $n\le2$ this is \cite{Hit81,PS91}, and for $n=3,4$ \cite{BWW22}. If $\mathrm{Pic}\,Z\neq\ZZ L$, then $M$ is $\HH P^{n}$ or
$\Gr_{2}(\CC^{n+2})$ \cite{LS94}. So let $n\ge3$ and $\mathrm{Pic}\,Z=\ZZ L$, and let $\MCbar$ be the component of the Hilbert scheme containing the twistor lines; its
general member is a balanced transverse conic. By Theorem~\ref{thm:tpD}, $d_{2}=1$. By Theorem~\ref{thm:LSJ}, $R_{E}^{\CC}$ is parallel, and
Corollary~\ref{cor:endgame} applies. Then $Z$, the twistor space of a Wolf space, is the adjoint variety of the complexified isometry group of $M$ \cite{Wolf65};
Theorem~\ref{thm:homog} gives this independently.
\end{proof}

\begin{corollary}[The curvature hypersurface]
\label{cor:curvhyp}
Let $M^{4n}$ be a compact quaternionic-K\"ahler manifold of positive scalar curvature with $\mathrm{Pic}\,Z=\ZZ L$, that is, neither $\HH P^{n}$ nor
$\Gr_{2}(\CC^{n+2})$. Then the curvature quartic is a nowhere vanishing holomorphic section $W$ of $S^{4}L^{\perp*}\otimes L^{2}$ over $Z$, and for every
$z\in Z$ its zero locus $V(W_{z})\subset\PP(L^{\perp}_{z})$ is an irreducible quartic hypersurface with
\[
V(W_{z})=\tau(X_{z})=X_{z}^{*},
\]
the tangent variety of the VMRT $X_{z}$ and, under the identification $\PP(L^{\perp}_{z})\cong\PP(L^{\perp*}_{z})$ given by the Levi form, its dual variety. In
particular $X_{z}$ has codegree four.
\end{corollary}

\begin{proof}
By Theorem~\ref{thm:main}, $M$ is a Wolf space, and as $\mathrm{Pic}\,Z=\ZZ L$, $Z=G_{\CC}/P$ is the adjoint variety of a simple Lie algebra not of type
$A$ or $C$ \cite{Wolf65,LS94}. It is homogeneous, so every point is general, and $X_{z}$ is the subadjoint variety: smooth, irreducible, nondegenerate and
Legendrian \cite{LM07}, \cite{Hwa24}. As $R_{E}^{\CC}$ is parallel, $W_{z}$ is locally constant on $\Sigma_{z}$ (Proposition~\ref{prop:leaf}), hence constant,
since $\Gamma_{z}$ is irreducible (Theorem~\ref{thm:irred}(3)) and every point is general. So $W^{\CC}=\bar p^{\,*}W$ for a holomorphic section $W$ of
$S^{4}L^{\perp*}\otimes L^{2}$ over $Z$, which is $G_{\CC}$-invariant because LeBrun's construction is natural under contact automorphisms; in particular $W_{z}$ is
invariant under the stabilizer of $z$. (Remark~\ref{rem:Wolfconverse} identifies $W$ with $c\,Q$, $Q$ the invariant quartic; this is not used below.) And
$W\neq0$, since $W\equiv0$ would make $M$ locally isometric to $\HH P^{n}$ \cite{Salamon82}; being invariant, $W$ is nowhere zero.

\emph{Inclusion.} Let $x\in X_{z}$ and $A$ the line through $z$ with direction $x$. By Lemma~\ref{lem:linesplit}, $L|_{A}=\OO(1)$ and $L^{\perp}|_{A}$ has
positive part $P_{A}\cong\OO(2)\oplus\OO(1)^{n-1}$, whose fibre at $z$ is $\hat T_{x}\hat X_{z}$ by Lemma~\ref{lem:tplines}(2). On $S^{4}P_{A}$ the section
$W|_{A}$ of $S^{4}(L^{\perp}|_{A})^{*}\otimes\OO(2)$ has components in $\OO(2-a_{1}-\dots-a_{4})$, $a_{i}\ge1$, all of negative degree, so they vanish.
Thus $W_{z}$ vanishes on every $\hat T_{x}\hat X_{z}$, and $\tau(X_{z})\subset V(W_{z})$.

\emph{Equality.} Since $X_{z}$ is nondegenerate and Legendrian, $\tau(X_{z})$ is an irreducible hypersurface \cite[Prop.~17(2)]{LM07}; and $V(W_{z})$
is irreducible and reduced. Indeed, the semisimple part $H$ of the stabilizer of $z$ is connected and acts irreducibly on $L^{\perp}_{z}$, preserving the Levi form
\cite{LM07}, and $W_{z}$ is a nonzero $H$-invariant quartic. $H$ permutes the irreducible factors of $W_{z}$; being connected, it fixes each up to a character, which is trivial
as $H$ is semisimple. So a nontrivial factorization would contain an $H$-invariant form of degree one or two. There is none: an invariant linear form would be a
fixed vector of the nontrivial irreducible representation $L^{\perp*}_{z}$, and an irreducible representation carries at most one invariant bilinear form up to
scale, here the skew Levi form, so no symmetric one. Two irreducible hypersurfaces, one inside the other, coincide.

\emph{Duality.} A hyperplane $\{\omega(h,\cdot)=0\}$ is tangent to $X_{z}$ at $x$ if and only if it contains $\hat T_{x}\hat X_{z}$, that is,
$h\in(\hat T_{x}\hat X_{z})^{\perp_{\omega}}=\hat T_{x}\hat X_{z}$, the tangent spaces of a Legendrian cone being Lagrangian. So $X_{z}^{*}$ is the union of the
$\PP(\hat T_{x}\hat X_{z})$, which is $\tau(X_{z})$.
\end{proof}

\begin{question}
\label{q:codegree4}
By Corollary~\ref{cor:curvhyp}, the variety of minimal rational tangents at a point of an adjoint variety is a smooth, irreducible, linearly nondegenerate
Legendrian variety of codegree four, that is, whose dual variety is a quartic hypersurface. Conversely, is every smooth, irreducible, linearly nondegenerate
Legendrian variety $X\subset\PP^{2n-1}$ of codegree four the variety of minimal rational tangents at a point of an adjoint variety? The degree four was observed
by Landsberg and Manivel, who asked for a classification of the smooth varieties of codegree four \cite[Remark after Prop.~5.8]{LM01}; Zak, who classified those of
codegree at most three, expects their main series to consist of the homogeneous Legendrian varieties \cite[Remark~3.6]{Zak04}. Corollary~\ref{cor:curvhyp}
identifies the quartic with the curvature. We investigate this question in a forthcoming paper \cite{Ye26Q}.
\end{question}

\begin{remark}[Freudenthal's magic square]
\label{rem:magicsquare}
Corollary~\ref{cor:curvhyp} is the projective geometry of Freudenthal's magic square \cite{Fre54,LM01}, read through the curvature.
\begin{itemize}
\item \emph{The Legendrian level.} For $\mathfrak g_{2},\mathfrak{so}_{n+4},\mathfrak f_{4},\mathfrak e_{6},\mathfrak e_{7},\mathfrak e_{8}$ the VMRT $X_{z}$ is the twisted
cubic, $\PP^{1}\times\mathbb Q^{n-2}$, $\mathrm{LG}(3,6)$, $\Gr(3,6)$, the spinor variety $\mathbb S_{6}$ and the Freudenthal variety $E_{7}/P_{7}$, and $W_{z}=c\,Q$
is the quartic invariant of $L^{\perp}_{z}$: Freudenthal's quartic for the last four, the discriminant of binary cubics for $\mathfrak g_{2}$
(Remark~\ref{rem:Wolfconverse}). Its zero locus is the tangent variety, which is the dual variety: codegree four.
\item \emph{The Severi level.} One step down, the variety of lines through a point of $X_{z}$ (its derived variety) is, for the exceptional four, the Severi
variety $\mathbb A\PP^{2}\subset\PP(J_{3}(\mathbb A))$, $\mathbb A=\RR,\CC,\HH,\mathbb O$ complexified: $v_{2}(\PP^{2})$, $\PP^{2}\times\PP^{2}$, $\Gr(2,6)$,
$E_{6}/P_{1}$. It is secant defective, so its secant and tangent varieties coincide; they are the determinant cubic, and its dual variety is the determinant cubic
of $J_{3}(\mathbb A)^{*}$: codegree three \cite{Zak93}. On the normal slice of $X_{z}$ the quartic restricts to this cubic, up to normalization
\cite[\S5]{LM01}.
\item \emph{The quadric level.} One step further, the second fundamental form of the Severi variety is $v\mapsto v^{*}v\in J_{2}(\mathbb A)$, with image the quadric
$\mathbb A\PP^{1}=\mathbb Q^{a}\subset\PP(J_{2}(\mathbb A))$, $a=\dim_{\RR}\mathbb A$, of codegree two.
\end{itemize}
So each step down the ladder lowers the degree of the hypersurface by one, $4\to3\to2$. Most of these facts are due to Landsberg and Manivel \cite{LM01,LM07}. For more details, see the companion paper \cite{Ye26Q}.
\end{remark}

\section{Directions of future research}
\label{sec:future}

The framework of this paper has two parts: a \emph{dictionary}, which reads the singularities of the induced geometry along the boundary of $\MCbar$ off the
degenerate conics (Table~\ref{tab:dictionary}), and a \emph{rigidity} step, which in the application uses the compactness of the family of conics through a point
to make the curvature constant along it, hence parallel (Proposition~\ref{prop:leaf}, Theorem~\ref{thm:LSJ}). The dictionary is local; the rigidity step is where
compactness enters. This structure suggests two directions of research.
\begin{enumerate}
\item[(A)] \emph{Negative quaternionic-K\"ahler manifolds} (\S\ref{sec:negQK}). Here the rigidity step is unavailable: the local theory does not see the
sign of the scalar curvature, but the global part of the proof fails, beginning with the compactness of $\Gamma_{z}$. The basic question is whether every
compact quaternionic-K\"ahler manifold of negative scalar curvature is locally symmetric, together with its analogue in indefinite signature.
\item[(B)] \emph{The singular Merkulov transform} (\S\ref{sec:boundaryconf}). Here the dictionary is extended from the conics to Merkulov's Legendre moduli
spaces, which carry, locally, every torsion-free connection with irreducibly acting holonomy; the boundary of a compactified Legendre moduli space should record
the contact geometry of the degenerate Legendre submanifolds.
\end{enumerate}
Question~\ref{q:codegree4}, on Legendrian varieties of codegree four, is taken up in the companion paper \cite{Ye26Q}.

\subsection{Negative quaternionic-K\"ahler manifolds}
\label{sec:negQK}

For negative scalar curvature compact examples exist: if $\Gamma\subset Sp(n,1)$ is a torsion-free cocompact lattice \cite{Bor63}, the quotient
$\Gamma\backslash\HH H^{n}$ of quaternionic hyperbolic space $\HH H^{n}=Sp(n,1)/Sp(n)Sp(1)$ is a compact locally symmetric quaternionic-K\"ahler manifold, and
the same holds for the noncompact dual of every Wolf space. These are the only known compact examples: no compact quaternionic-K\"ahler manifold of negative
scalar curvature that is not locally symmetric is known. Complete examples that are not locally symmetric exist in abundance \cite{Ale75,LeB91}, and
Cort\'es, R\"oser and Thung constructed complete ones with an end of finite volume \cite{CRT22}; but even a non-locally-symmetric example of finite
volume is unknown. In indefinite signature LeBrun's construction gives many local examples with holonomy
in $Sp(l,n-l)Sp(1)$ for every $l$, for instance from conformal Riemannian $(n+2)$-manifolds; these metrics are in general incomplete \cite[Concluding
Remarks~(1)--(2)]{LeBrun89}, so the second question below concerns compactness.

\begin{question}[Negative quaternionic-K\"ahler manifolds]
\label{q:negQK}
\leavevmode
\begin{enumerate}
\item Is every compact quaternionic-K\"ahler manifold of dimension $4n\ge8$ and negative scalar curvature locally symmetric?
\item For $0<l<n$, is every compact pseudo-Riemannian manifold with holonomy in $Sp(l,n-l)Sp(1)$ and nonzero scalar curvature locally symmetric?
Does one exist?
\end{enumerate}
\end{question}

\emph{LeBrun's real slices.} Let $\sigma$ be an antiholomorphic involution without fixed points of a complex contact manifold $Z$ of dimension $2n+1$, $n\ge2$,
preserving the contact distribution. Each connected component of its fixed locus $M_{\sigma}\subset\MC$ is a pseudo-Riemannian manifold with holonomy in $Sp(l,n-l)Sp(1)$ for some
$0\le l\le n$, which may differ from one component to another, and every such manifold arises in this way \cite[Thm.~1.3 and Remark~(2)]{LeBrun89}; $l$ is read off from the signature of the quaternionic
Hermitian form $\omega_{E}(u,\sigma_{E}v)$, where $\sigma_{E}$ is induced by $\sigma$ on the factor $E$ of $T^{*}\MC=E\otimes H$. With the scalar curvature normalized
to $-1$ and signature $(4(n-l),4l)$, $l=n$ is the case of Theorem~Q up to sign, and $l=0$ is Question~\ref{q:negQK}(1).

\begin{proposition}[The real slices of the quaternionic space-form model]
\label{prop:flatreal}
Let $Z=\PP(\CC^{2n+2})$ with its symplectic form $\omega$, $j$ an antilinear map with $j^{2}=-1$ and $\omega(ju,jv)=\overline{\omega(u,v)}$, and $\sigma=[j]$. Then
$h(u,v)=\omega(u,jv)$ is a quaternionic Hermitian form, of signature $(p,q)$ say, $p+q=n+1$.
\begin{enumerate}
\item The $\sigma$-invariant balanced transverse conics are the lines $\PP(W)$, $W$ a quaternionic line with $h|_{W}\neq0$, and
$M_{\sigma}=M_{+}\sqcup M_{-}$ with $M_{+}=Sp(p,q)/Sp(1)\times Sp(p-1,q)$ and $M_{-}=Sp(p,q)/Sp(p,q-1)\times Sp(1)$. On $M_{+}$, $l=p-1$; on $M_{-}$, $l=q-1$.
\item The $\sigma$-invariant integral conics are the $\PP(W)$ with $h|_{W}=0$. For $p,q\ge1$ they form the common frontier of $M_{+}$ and $M_{-}$ in $\HH P^{n}$.
\end{enumerate}
So the quaternionic space-form model, with zero quaternionic Weyl curvature, for a given $l$ is $Sp(l+1,n-l)/Sp(1)\times Sp(l,n-l)$: $\HH P^{n}$ for $l=n$, $\HH H^{n}$ for $l=0$.
\end{proposition}

\begin{proof}
$\overline{h(u,v)}=\omega(ju,j^{2}v)=-\omega(ju,v)=h(v,u)$, and $h(ju,jv)=\overline{h(u,v)}$. A line $\PP(W)$ is $\sigma$-invariant exactly when $W$ is a
quaternionic line; it is transverse exactly when $\omega|_{W}\neq0$, that is, since $\omega(w,jw)=h(w,w)$, when $h|_{W}\neq0$, and integral otherwise. $Sp(p,q)$
acts transitively on the positive and on the negative lines; if $h(w,w)=0$ and $h(w,u)\in\RR^{*}$, then $h(w+\varepsilon u,w+\varepsilon u)=2\varepsilon
h(w,u)+O(\varepsilon^{2})$ takes both signs, which gives~(2). For $l$: at $\PP(W)$, a real tangent vector is a $j$-linear $\varphi\colon W\to\CC^{2n+2}/W$, and up
to a constant $g^{\CC}(\varphi,\varphi)=h(\varphi w,\varphi w)/h(w,w)=:g_{h}(\varphi,\varphi)$, with $\varphi w$ taken in $W^{\perp_{h}}$. On $\HH P^{n}$ ($q=0$)
$g_{h}$ is the Fubini--Study metric, so $g_{h}$ has positive scalar curvature on every slice, and the metric of scalar curvature $-1$ is a positive multiple of
$-g_{h}$. On $M_{+}$ it has signature $(4q,4(p-1))$, so $l=p-1$. Replacing $\omega$ by $-\omega$ exchanges $M_{\pm}$, so $l=q-1$ on $M_{-}$.
\end{proof}

So $l$ can change across the real integral conics, along which $g^{\CC}$ has its pole (Proposition~\ref{prop:intregular}); for $p=1$ they form the sphere
$\partial\HH H^{n}$. For $0<l<n$ with $2l\ge n$, only finite groups act properly discontinuously on the space-form model, by the Calabi--Markus phenomenon
\cite[Cor.~4.4]{Kob89}, so for these $l$ no complete example with zero quaternionic Weyl curvature is compact (Question~\ref{q:negQK}(2)).

\begin{proposition}[The twistor space in negative scalar curvature]
\label{prop:negtwistor}
Let $M^{4n}$ be a compact quaternionic-K\"ahler manifold with negative scalar curvature, and $Z$ its twistor space, $K_{Z}^{-1}\cong L^{n+1}$.
\begin{enumerate}
\item $Z$ has Kodaira dimension $-\infty$.
\item $H^{0}(Z,L)=0$: $Z$ has no infinitesimal contact automorphisms.
\item If $n\ge3$, $L$ is not ample.
\end{enumerate}
\end{proposition}

\begin{proof}
(1) $K_{Z}\cdot C=-2(n+1)$ on every twistor line, and the twistor lines cover $Z$.

(2) A section of $L$ is the Hamiltonian of an infinitesimal contact automorphism $\xi$ \cite[Prop.~2.1]{LeB95}, and the real structure of $Z$ acts on
$H^{0}(Z,L)$, so it suffices to treat $\xi$ commuting with $\sigma$. The flow of $\mathrm{Re}\,\xi$ commutes with $\sigma$ and permutes the twistor lines, a compact
component of $M_{\sigma}$; it acts on $M$ by isometries, the metric being determined by $(Z,\sigma)$ once its scalar curvature is fixed. A flow acting trivially on
$M$ acts on each $T^{*}_{m}M\otimes\CC=E\otimes H$ by scalars $A\otimes B$, hence fixes every twistor line pointwise and is trivial on $Z$. By Bochner's theorem
there is no nonzero Killing field \cite{Bes87}.

(3) If $L$ is ample, $Z$ is Fano. If $b_{2}(Z)\ge2$, then $Z\cong\PP(T^{*}Y)$ \cite{KPSW00,Dem02}, $T_{Y}$ is ample and $Y\cong\PP^{n+1}$ \cite{Mor79}. If
$b_{2}(Z)=1$ and $L=mL_{0}$ with $m\ge2$, then $m=2$ and $Z\cong\PP^{2n+1}$ \cite{KO73}. Otherwise $\mathrm{Pic}\,Z=\ZZ L$, and $Z$ is an adjoint variety by
Theorem~H. In every case $Z$ is homogeneous under its contact automorphisms, contradicting~(2).
\end{proof}

The proof of Theorem~Q breaks down for negative scalar curvature. Its local part (the leaf principle, Theorem~\ref{thm:thirdorder},
\S\S\ref{sec:dilation}--\ref{sec:integral}) does not see the sign of the scalar curvature, but its global part fails, beginning with the compactness of $\Gamma_{z}$.

\begin{example}[The conics through a point of $\HH P^{n}$ and of $\HH H^{n}$]
\label{ex:ballleaf}
In Proposition~\ref{prop:flatreal} let $z=[u]$ with $h(u,u)>0$, and write $u^{\perp_{h}}=\CC ju\oplus N$. The contact directions at $z$ are $\{a=0\}$, and the
line through $z$ in the direction $[a\,ju+w]$, $w\in N$, lies in $Z_{+}=\{h>0\}$ exactly when $|a|^{2}h(u,u)+h(w,w)>0$.
\begin{enumerate}
\item $\HH P^{n}$: $N$ is positive, $\Sigma_{z}\cong\CC^{2n}$, and $\Gamma_{z}=\PP^{2n}$ adds the hyperplane of integral conics, across which the curvature extends
(Proposition~\ref{prop:intregular}).
\item $\HH H^{n}$: $N$ is negative, and $\Sigma_{z}$ is the ball $\{|w|_{N}^{2}<h(u,u)\}$ in the chart $a=1$, whose closure misses the integral conics; the lines
contained in $Z_{+}$ form the bounded symmetric domain $SU(2,2n)/S(U(2)\times U(2n))$.
\end{enumerate}
\end{example}

On a compact quotient $\Gamma\backslash Z_{+}$ one expects the leaf to be this ball and the dictionary to have no degenerate conics to act on. This needs two
statements that we have not proved: that no conic of the ball meets a nontrivial translate of itself under $\Gamma$, so that it projects injectively to an embedded
conic, and that the areas of the conics become unbounded near the frontier of the ball; given these, Bishop's theorem shows that the conics near the frontier have
no limit among the curves of the quotient. There are also no curves with $L\cdot C=1$, and
$\pi_{1}(Z)=\Gamma$. By Proposition~\ref{prop:negtwistor}(2), $H^{0}(Z\times Z,L\boxtimes L)=0$, so an argument that builds sections of $L$ or $L\boxtimes L$ from
the conics alone \cite{OAI26,HLW26} must use the projectivity of $Z$.

\subsection{Degenerations of conformal structures and the singular Merkulov transform}
\label{sec:boundaryconf}

Rational curves force conformal flatness: a projective conformal manifold of dimension $m\ge3$ whose canonical bundle is not nef is $\mathbb Q^{m}$ \cite{Ye94},
and so is a compact K\"ahler one containing a rational curve \cite{BM16}; one containing an immersed rational null geodesic is conformally flat \cite{Bel01}. So
nonflat structures studied through rational curves live on noncompact manifolds, whose natural compactifications force them to degenerate.

\emph{Skies and Legendre moduli spaces.} For a holomorphic conformal manifold $(\MC,g_{\CC})$ of dimension $m\ge3$, the space $Z$ of null geodesics is locally a
contact manifold of dimension $2m-3$, and a point $x$ is represented by its sky, the Legendre submanifold $X_{x}\cong\mathbb Q^{m-2}$ of null geodesics through $x$;
$T_{x}\MC\cong H^{0}(X_{x},L|_{X_{x}})$, and the null cone is the dual quadric of $X_{x}$ \cite{LeB83}. Merkulov extended this to the moduli space $M$ of a compact
Legendre submanifold $X=G/P$ with $L|_{X}$ very ample and $H^{1}(X,L|_{X})=0$: $M$ carries an induced $G$-structure, torsion-free exactly when the members are
$1$-flat, whose curvature is a class $\rho_{t}\in H^{1}\bigl(X_{t},L|_{X_{t}}\otimes S^{3}(J^{1}L|_{X_{t}})^{*}\bigr)$ \cite{Mer94,Mer95,Mer97}, an invariant of order
three, as in Theorem~\ref{thm:thirdorder}. The structure preserves the cones $X_{t}^{*}\subset\PP(T_{t}M)$, whose degree is the codegree of $X$: two for skies;
four for a rational Legendre curve with $L|_{X}\cong\OO(3)$, which gives Bryant's structure on binary cubics \cite{Bry91}; four for a subadjoint variety, whose
dual is the zero locus of the curvature quartic (Corollary~\ref{cor:curvhyp}). Since every torsion-free connection with irreducibly acting holonomy arises
locally on such a moduli space (\S\ref{sec:intromain}), Merkulov's correspondence carries all of these geometries.

\emph{The singular Merkulov transform.} The philosophy of this paper applies to it unchanged: the boundary of a compactified Legendre moduli space should
carry the contact geometry of the degenerate Legendre submanifolds on the other side of the correspondence. We call the extension of Merkulov's correspondence
to the limiting members the \emph{singular Merkulov transform}; the singular LeBrun transform is its analogue for the transverse conics. The algebraic side,
the limits of Legendre families in the Hilbert scheme and their dual varieties, seems largely unexplored. Let $Z$ be a contact manifold, $M$ a Merkulov
moduli space embedded in the Douady space of Legendre submanifolds of $Z$, and $\overline M\supset M$ a compactification with boundary divisor $D$ and a flat
limit family.

\begin{question}[Boundary degenerations]
\label{q:boundaryconf}
\leavevmode
\begin{enumerate}
\item Which limiting members $X_{0}$ (singular, reducible, nonreduced) occur, beginning with Legendrian conics in a contact threefold ($m=3$)?
\item Does the field of cones extend across $D$ as the pull-back, by the Kodaira--Spencer map $\kappa$, of the relative dual variety of the limiting members,
degenerating only where $X_{0}$ is singular or $\kappa$ loses rank?
\item When do $g_{\CC}$, its normal Cartan connection, or the torsion-free connections of the induced structure extend meromorphically or logarithmically
across $D$, and is the leading term of the Weyl tensor, of the Cotton tensor ($m=3$), or of $\rho$ read off from a neighbourhood of $X_{0}$ in $Z$, as in
Table~\ref{tab:dictionary}?
\item When does the compactness of $Z$ force the induced structure to be locally symmetric, as Theorem~\ref{thm:LSJ} does for the conics? For projective $Z$
with $\mathrm{Pic}\,Z=\ZZ L$, $Z$ is an adjoint variety (Theorem~H; for every contact Fano manifold see \cite{HLW26}), and for members homogeneous under their
stabilizers the question becomes one of Lie theory; for $m=3$ and $Z$ projective, $Z$ is $\PP^{3}$ or a projectivized cotangent
bundle \cite{Ye94b}. The question is open for members that are not homogeneous, such as the non-homogeneous smooth Legendrian varieties of \cite{Buc06},
whenever $H^{1}(X,L|_{X})=0$, where the induced structure is a field of cones with a small, possibly
finite, structure group; and for compact $Z$ that are not projective, such as the twistor spaces of Proposition~\ref{prop:negtwistor}.
\end{enumerate}
\end{question}

For the conics of this paper, which form a Kodaira family with $T_{[C]}\MC=H^{0}(C,N_{C})=E^{*}\otimes H^{*}$, the answers are those of Table~\ref{tab:dictionary}: the
metric is regular on the balanced transverse conics, has a pole along the balanced integral conics, where $\sigma_{\theta}$ vanishes (Proposition~\ref{prop:intregular}),
and loses its factorization along the jumping and nodal conics (Remark~\ref{rem:singLeBrun}). For $n=1$ the two settings meet: for $Z=\PP^{3}$,
$\MCbar=\Gr_{2}(\CC^{4})=\mathbb Q^{4}$, and its integral divisor $LG(2,4)\cong\mathbb Q^{3}$ is the Legendre moduli space of the Legendrian lines of $\PP^{3}$, the
skies of the flat conformal structure of $\mathbb Q^{3}$ \cite{LeB83}, which is the conformal infinity of the self-dual structure of $\mathbb Q^{4}$
(Remark~\ref{rem:intflat}). For a general contact threefold, LeBrun's $\mathcal H$-space construction \cite{LeB82} gives the local answer: the twistor space of the self-dual Einstein
metric is the space of null geodesics of its conformal infinity, whose skies are the integral conics, so the integral divisor is locally the conformal
infinity. Does this hold globally on $\MCbar$ for compact $Z$, and how does the conformal infinity degenerate where the integral divisor meets the divisor of
line pairs?

\appendix

\section{Brief Survey of Related Work Since the 1990s}
\label{app:survey}

The following summarizes the related work since the early 1990s. It is not meant to be complete, and the author apologizes for any important work that has been missed.

\emph{The conjecture.} The statement, the known cases and the link with Beauville's conjecture through \cite{LeB95} are recalled in \S\ref{sec:intromain}
and \S\ref{sec:introhistory}. Besides the reduction to $\mathrm{Pic}\,Z=\ZZ L$, LeBrun and Salamon proved that in each dimension there are, up to
homothety, only finitely many positive quaternionic-K\"ahler manifolds \cite{LS94}.

\emph{Fano contact manifolds.} Mori theory reduces Beauville's conjecture to $\mathrm{Pic}\,Z=\ZZ L$: a projective contact manifold is $\PP(T^{*}Y)$,
$\PP^{2n+1}$, or Fano with $\mathrm{Pic}\,Z=\ZZ L$ \cite{KPSW00,Dem02}, extending \cite{Ye94b} (dimension five: \cite{Dru98}). Beauville proved homogeneity when the
contact automorphism group is reductive and $|L|$ is generically finite \cite{Bea98,Bea99}. The lines through a general point are tangent to the contact
distribution, and their tangent directions form a Legendrian variety $X_{z}\subset\PP(L^{\perp}_{z})$ \cite{Keb01,Keb05}. By the theory of varieties of minimal
rational tangents of Hwang and Mok \cite{HM98,HM01,Hwa97,Mok08,HH08}, $Z$ is an adjoint variety as soon as $X_{z}$ is a subadjoint variety. The Legendrian varieties
themselves were studied by Landsberg and Manivel, who placed the subadjoint varieties in Freudenthal's magic square \cite{LM01,LM07}, by Buczy\'nski
\cite{Buc06,Buc09,Buc10}, and by Hwang, who characterized the subadjoint ones through the third fundamental form \cite{Hwa24}.

\emph{The 1993 preprint.} In a 1993 preprint \cite{Ye93} the author related $X_{z}$ to the Severi varieties and the Wolf spaces,
and showed that $M$ is a Wolf space if $X_{z}$ is rigid along a twistor fibre: the Hermitian--Einstein metric induced on $L^{\perp}$ makes it polystable
\cite{Lub83}, so a rigid $X_{z}$ gives a reduction of the holonomy, and Berger's theorem concludes. The rigidity would follow if, for the twistor fibration $q\colon Z\to M$, the pull-back
$q^{*}R_{E}$ of the quaternionic Weyl curvature were a holomorphic section of $q^{*}S^{4}E$, which, however, is equivalent to $\nabla R_{E}=0$, as is easy to see from the twistor
construction, since the horizontal $(0,1)$-directions of $Z$ over a point $m$ span $T_{m}M\otimes\CC$.

\emph{Conics.} In notes of 2000, Wi\'sniewski proposed to work with the compactified space of conics and its degenerations to pairs of lines, through a symmetric
divisor $D\subset Z\times Z$ in $|L^{d}\boxtimes L^{d}|$, aiming at the base-point freeness of $|L|$, which would give the conjecture by \cite{Bea98,LeB95}; a gap
pointed out by Koll\'ar \cite[Proposition~2(ii)]{Wis00} remained, and the notes were not published. This idea appears to underlie the key arguments of
\cite{OAI26} and \cite{HLW26}.

\emph{Topology, symmetry and torus actions.} Large isometry rank forces $\HH P^{n}$ or $\Gr_{2}(\CC^{n+2})$ \cite{Fan04}. A torus of large rank acting on $Z$ forces
homogeneity \cite{ORSW21}, which gives the cases $n=3,4$ \cite{BWW22}. Combining these with Salamon's index theory, Daura Serrano, Gohr, Lawn and
Schedler claim the cases $n\le14$ \cite{DGLS26}. Two recent preprints claim proofs of the conjecture by producing enough infinitesimal contact automorphisms
\cite{OAI26,HLW26}; they are discussed in \S\ref{sec:introrelation}. Fu \cite{Fu26} gives an account of the classification of Fano contact manifolds
that combines the constructions of these two proofs with earlier results on contact lines, with some simplifications. This account came to the author's attention only after the present
paper was completed; the two were written independently.


\begin{thebibliography}{99}
\bibitem{AV02} D.~Abramovich and A.~Vistoli, \emph{Compactifying the space of stable maps}, J.\ Amer.\ Math.\ Soc.\ \textbf{15} (2002), 27--75.

\bibitem{Ale68} D.~V.~Alekseevskii, \emph{Riemannian spaces with exceptional holonomy groups}, Funktsional.\ Anal.\ i Prilozhen.\ \textbf{2} (1968), no.~2, 1--10; English transl.\ in Funct.\ Anal.\ Appl.\ \textbf{2} (1968), 97--105.

\bibitem{Ale75} D.~V.~Alekseevski\u{\i}, \emph{Classification of quaternionic spaces with a transitive solvable group of motions}, Math.\ USSR-Izv.\ \textbf{9} (1975), 297--339.

\bibitem{Ati57} M.~F.~Atiyah, \emph{Complex analytic connections in fibre bundles}, Trans.\ Amer.\ Math.\ Soc.\ \textbf{85} (1957), 181--207.

\bibitem{BE91} T.~N.~Bailey and M.~G.~Eastwood, \emph{Complex paraconformal manifolds---their differential geometry and twistor theory}, Forum Math.\ \textbf{3} (1991), 61--103.

\bibitem{Bea98} A.~Beauville, \emph{Fano contact manifolds and nilpotent orbits}, Comment.\ Math.\ Helv.\ \textbf{73} (1998), 566--583.

\bibitem{Bea99} A.~Beauville, \emph{Riemannian holonomy and algebraic geometry}, arXiv:math/9902110 (1999).

\bibitem{Bel01} F.~A.~Belgun, \emph{Null-geodesics in complex conformal manifolds and the LeBrun correspondence}, J.\ Reine Angew.\ Math.\ \textbf{536} (2001), 43--63.

\bibitem{Ber55} M.~Berger, \emph{Sur les groupes d'holonomie homog\`ene des vari\'et\'es \`a connexion affine et des vari\'et\'es riemanniennes}, Bull.\ Soc.\ Math.\ France \textbf{83} (1955), 279--330.

\bibitem{Bes87} A.~L.~Besse, \emph{Einstein Manifolds}, Ergeb.\ Math.\ Grenzgeb.\ (3) \textbf{10}, Springer, Berlin, 1987.

\bibitem{BM16} I.~Biswas and B.~McKay, \emph{Holomorphic Cartan geometries and rational curves}, Complex Manifolds \textbf{3} (2016), 145--168.

\bibitem{BM47} S.~Bochner and D.~Montgomery, \emph{Groups on analytic manifolds}, Ann.\ of Math.\ (2) \textbf{48} (1947), 659--669.

\bibitem{Bor63} A.~Borel, \emph{Compact Clifford--Klein forms of symmetric spaces}, Topology \textbf{2} (1963), 111--122.

\bibitem{Bry91} R.~L.~Bryant, \emph{Two exotic holonomies in dimension four, path geometries, and twistor theory}, in: Complex Geometry and Lie Theory (Sundance, 1989), Proc.\ Sympos.\ Pure Math.\ \textbf{53}, Amer.\ Math.\ Soc., 1991, 33--88.

\bibitem{Bry00} R.~L.~Bryant, \emph{Recent advances in the theory of holonomy}, S\'eminaire Bourbaki, Vol.\ 1998/99, Exp.\ No.~861, Ast\'erisque \textbf{266} (2000), 351--374; arXiv:math/9910059.

\bibitem{Buc06} J.~Buczy\'nski, \emph{Legendrian subvarieties of projective space}, Geom.\ Dedicata \textbf{118} (2006), 87--103; arXiv:math/0503528.

\bibitem{Buc09} J.~Buczy\'nski, \emph{Algebraic Legendrian varieties}, Dissertationes Math.\ \textbf{467} (2009).

\bibitem{Buc10} J.~Buczy\'nski, \emph{Duality and integrability on contact Fano manifolds}, Doc.\ Math.\ \textbf{15} (2010), 821--841; corrected version arXiv:1002.0698v2.

\bibitem{BM19} J.~Buczy\'nski and G.~Moreno, \emph{Complex contact manifolds, varieties of minimal rational tangents, and exterior differential
systems}, Banach Center Publ.\ \textbf{117} (2019), 145--176; arXiv:1805.08548.

\bibitem{BWW22} J.~Buczy\'nski, J.~A.~Wi\'sniewski and A.~Weber, \emph{Algebraic torus actions on contact manifolds}, J.\ Differential Geom.\ \textbf{121} (2022), no.~2, 227--289; arXiv:1802.05002.

\bibitem{Cam92} F.~Campana, \emph{Connexit\'e rationnelle des vari\'et\'es de Fano}, Ann.\ Sci.\ \'Ecole Norm.\ Sup.\ (4) \textbf{25} (1992), 539--545.

\bibitem{CRT22} V.~Cort\'es, M.~R\"oser and D.~Thung, \emph{Complete quaternionic K\"ahler manifolds with finite volume ends}, preprint, 2021--2022, arXiv:2105.00727.

\bibitem{CS09} A.~\v{C}ap and J.~Slov\'ak, \emph{Parabolic Geometries I: Background and General Theory}, Math.\ Surveys Monogr.\ \textbf{154}, Amer.\ Math.\ Soc., 2009.

\bibitem{Car26} \'E.~Cartan, \emph{Sur une classe remarquable d'espaces de Riemann}, Bull.\ Soc.\ Math.\ France \textbf{54} (1926), 214--264.

\bibitem{CP09} P.-E.~Chaput and N.~Perrin, \emph{On the quantum cohomology of adjoint varieties}, Proc.\ Lond.\ Math.\ Soc.\ (3) \textbf{103} (2011), no.~2, 294--330; arXiv:0904.4824 (Prop.~3.11 and Cor.~3.23 are numbered as there).

\bibitem{DGLS26} J.~Daura Serrano, A.~Gohr, M.-A.~Lawn and T.~Schedler, \emph{The LeBrun--Salamon conjecture is true in dimension up to $56$},
arXiv:2610.06336, October 2026.

\bibitem{Dem02} J.-P.~Demailly, \emph{On the Frobenius integrability of certain holomorphic $p$-forms}, in: Complex Geometry (G\"ottingen, 2000), Springer, 2002,
93--98.

\bibitem{DG60} F.~Docquier and H.~Grauert, \emph{Levisches Problem und Rungescher Satz f\"ur Teilgebiete Steinscher Mannigfaltigkeiten}, Math.\ Ann.\ \textbf{140} (1960), 94--123.

\bibitem{Dou66} A.~Douady, \emph{Le probl\`eme des modules pour les sous-espaces analytiques compacts d'un espace analytique donn\'e},
Ann.\ Inst.\ Fourier (Grenoble) \textbf{16} (1966), no.~1, 1--95.

\bibitem{Dru98} S.~Druel, \emph{Structures de contact sur les vari\'et\'es alg\'ebriques de dimension~5}, C.\ R.\ Acad.\ Sci.\ Paris S\'er.\ I Math.\ \textbf{327}
(1998), 365--368.

\bibitem{Fan04} F.~Fang, \emph{Positive quaternionic K\"ahler manifolds and symmetry rank}, J.\ Reine Angew.\ Math.\ \textbf{576} (2004), 149--165.

\bibitem{FG65} W.~Fischer and H.~Grauert, \emph{Lokal-triviale Familien kompakter komplexer Mannigfaltigkeiten}, Nachr.\ Akad.\ Wiss.\ G\"ottingen Math.-Phys.\ Kl.\ II (1965), 89--94.

\bibitem{FP97} W.~Fulton and R.~Pandharipande, \emph{Notes on stable maps and quantum cohomology}, in: Algebraic Geometry (Santa Cruz, 1995), Proc.\ Sympos.\ Pure Math.\ \textbf{62}, Part~2, Amer.\ Math.\ Soc., 1997, 45--96.

\bibitem{Fre54} H.~Freudenthal, \emph{Beziehungen der $E_{7}$ und $E_{8}$ zur Oktavenebene.~I}, Nederl. Akad. Wetensch. Proc. Ser.~A \textbf{57} (1954), 218--230.

\bibitem{Fu26} B.~Fu, \emph{The classification of Fano contact manifolds revisited}, arXiv:2611.01423, November~1, 2026.

\bibitem{Gra59} J.~W.~Gray, \emph{Some global properties of contact structures}, Ann.\ of Math.\ \textbf{69} (1959), 421--450.

\bibitem{EGAIV3} A.~Grothendieck, \emph{\'El\'ements de g\'eom\'etrie alg\'ebrique IV, troisi\`eme partie}, Publ.\ Math.\ IH\'ES \textbf{28} (1966).

\bibitem{SGA1} A.~Grothendieck, \emph{Rev\^etements \'etales et groupe fondamental (SGA~1)}, Lecture Notes in Math.\ \textbf{224}, Springer, 1971.

\bibitem{Hel78} S.~Helgason, \emph{Differential Geometry, Lie Groups, and Symmetric Spaces}, Academic Press, 1978.

\bibitem{Hit81} N.~J.~Hitchin, \emph{K\"ahlerian twistor spaces}, Proc.\ London Math.\ Soc.\ (3) \textbf{43} (1981), 133--150.

\bibitem{HH08} J.~Hong and J.-M.~Hwang, \emph{Characterization of the rational homogeneous space associated to a long simple root by its variety of minimal rational
tangents}, in: Algebraic Geometry in East Asia (Hanoi, 2005), Adv.\ Stud.\ Pure Math.\ \textbf{50}, Math.\ Soc.\ Japan, 2008, 217--236; doi:10.2969/aspm/05010217.

\bibitem{HLW26} Z.~Hu, K.~Liu and X.~Wan, \emph{Transverse rational curves and the homogeneity of contact Fano manifolds}, arXiv:2610.38736, October~28, 2026.

\bibitem{Hwa97} J.-M.~Hwang, \emph{Rigidity of homogeneous contact manifolds under Fano deformation}, J.\ Reine Angew.\ Math.\ \textbf{486} (1997), 153--163.

\bibitem{Hwa24} J.-M.~Hwang, \emph{Characterizing subadjoint varieties among Legendrian varieties}, Ann.\ Inst.\ Fourier (Grenoble), online first (2026), doi:10.5802/aif.3760;
arXiv:2404.19376.

\bibitem{HM98} J.-M.~Hwang and N.~Mok, \emph{Rigidity of irreducible Hermitian symmetric spaces of the compact type under K\"ahler deformation}, Invent.\ Math.\
\textbf{131} (1998), 393--418.

\bibitem{HM01} J.-M.~Hwang and N.~Mok, \emph{Cartan--Fubini type extension of holomorphic maps for Fano manifolds of Picard number~$1$}, J.\ Math.\ Pures Appl.\
\textbf{80} (2001), 563--575.

\bibitem{Keb01} S.~Kebekus, \emph{Lines on contact manifolds}, J.\ Reine Angew.\ Math.\ \textbf{539} (2001), 167--177.

\bibitem{Keb05} S.~Kebekus, \emph{Lines on complex contact manifolds II}, Compos.\ Math.\ \textbf{141} (2005), 227--252; arXiv:math/0306260.

\bibitem{KPSW00} S.~Kebekus, T.~Peternell, A.~J.~Sommese and J.~A.~Wi\'sniewski, \emph{Projective contact manifolds}, Invent.\ Math.\ \textbf{142} (2000), 1--15.

\bibitem{KM97} S.~Keel and S.~Mori, \emph{Quotients by groupoids}, Ann.\ of Math.\ (2) \textbf{145} (1997), 193--213.

\bibitem{KN63} S.~Kobayashi and K.~Nomizu, \emph{Foundations of Differential Geometry}, Vol.~I, Interscience, 1963.

\bibitem{KN69} S.~Kobayashi and K.~Nomizu, \emph{Foundations of Differential Geometry}, Vol.~II, Interscience, 1969.

\bibitem{KO73} S.~Kobayashi and T.~Ochiai, \emph{Characterizations of complex projective spaces and hyperquadrics}, J.\ Math.\ Kyoto Univ.\ \textbf{13} (1973), 31--47.

\bibitem{Kob89} T.~Kobayashi, \emph{Proper action on a homogeneous space of reductive type}, Math.\ Ann.\ \textbf{285} (1989), 249--263.

\bibitem{Kol96} J.~Koll\'ar, \emph{Rational Curves on Algebraic Varieties}, Ergeb.\ Math.\ Grenzgeb.\ (3) \textbf{32}, Springer, 1996.

\bibitem{KMM92} J.~Koll\'ar, Y.~Miyaoka and S.~Mori, \emph{Rationally connected varieties}, J.\ Algebraic Geom.\ \textbf{1} (1992), 429--448.

\bibitem{LM01} J.~M.~Landsberg and L.~Manivel, \emph{The projective geometry of Freudenthal's magic square}, J.\ Algebra \textbf{239} (2001), 477--512.

\bibitem{LM07} J.~M.~Landsberg and L.~Manivel, \emph{Legendrian varieties}, Asian J.\ Math.\ \textbf{11} (2007), 341--359.

\bibitem{LeB82} C.~LeBrun, \emph{$\mathcal H$-space with a cosmological constant}, Proc.\ Roy.\ Soc.\ London Ser.\ A \textbf{380} (1982), 171--185.

\bibitem{LeB83} C.~LeBrun, \emph{Spaces of complex null geodesics in complex-Riemannian geometry}, Trans.\ Amer.\ Math.\ Soc.\ \textbf{278} (1983), 209--231.

\bibitem{LeBrun89} C.~LeBrun, \emph{Quaternionic-K\"ahler manifolds and conformal geometry}, Math.\ Ann.\ \textbf{284} (1989), 353--376.

\bibitem{LeB91} C.~LeBrun, \emph{On complete quaternionic-K\"ahler manifolds}, Duke Math.\ J.\ \textbf{63} (1991), 723--743.

\bibitem{LeB95} C.~LeBrun, \emph{Fano manifolds, contact structures, and quaternionic geometry}, Internat.\ J.\ Math.\ \textbf{6} (1995), 419--437.

\bibitem{LS94} C.~LeBrun and S.~Salamon, \emph{Strong rigidity of positive quaternion-K\"ahler manifolds}, Invent.\ Math.\ \textbf{118} (1994), 109--132.

\bibitem{Lub83} M.~L\"ubke, \emph{Stability of Einstein--Hermitian vector bundles}, Manuscripta Math.\ \textbf{42} (1983), 245--257.

\bibitem{Mat86} H.~Matsumura, \emph{Commutative Ring Theory}, Cambridge Stud.\ Adv.\ Math.\ \textbf{8}, Cambridge Univ.\ Press, 1986.

\bibitem{Mer94} S.~A.~Merkulov, \emph{Moduli of compact Legendre submanifolds of complex contact manifolds}, Math.\ Res.\ Lett.\ \textbf{1} (1994), 717--727.

\bibitem{Mer95} S.~A.~Merkulov, \emph{Twistor theory, complex homogeneous manifolds and $G$-structures}, preprint, 1995, arXiv:dg-ga/9502006.

\bibitem{Mer97} S.~A.~Merkulov, \emph{Existence and geometry of Legendre moduli spaces}, Math.\ Z.\ \textbf{226} (1997), 211--265.

\bibitem{MS99} S.~Merkulov and L.~Schwachh\"ofer, \emph{Classification of irreducible holonomies of torsion-free affine connections}, Ann.\ of Math.\ (2) \textbf{150} (1999), 77--149.

\bibitem{Mok08} N.~Mok, \emph{Recognizing certain rational homogeneous manifolds of Picard number~1 from their varieties of minimal rational tangents}, in: Third
International Congress of Chinese Mathematicians, AMS/IP Stud.\ Adv.\ Math.\ \textbf{42}, Amer.\ Math.\ Soc., 2008, 41--61.

\bibitem{Mor79} S.~Mori, \emph{Projective manifolds with ample tangent bundles}, Ann.\ of Math.\ (2) \textbf{110} (1979), 593--606.

\bibitem{ORSW21} G.~Occhetta, E.~A.~Romano, L.~E.~Sol\'a Conde and J.~A.~Wi\'sniewski, \emph{High rank torus actions on contact manifolds}, Selecta Math.\ (N.S.)
\textbf{27} (2021), Paper No.~10.

\bibitem{OAI26} OpenAI, \emph{Contact Fano manifolds and the LeBrun--Salamon conjecture}, OpenAI Math Release preprint, September~23, 2026,
\url{https://github.com/openai/math/tree/main/preprints/Contact-Fano-manifolds-and-the-LeBrun-Salamon-conjecture-September-23-2026}.


\bibitem{PS91} Y.~S.~Poon and S.~M.~Salamon, \emph{Quaternionic K\"ahler $8$-manifolds with positive scalar curvature}, J.\ Differential Geom.\ \textbf{33} (1991), 363--378.

\bibitem{Salamon82} S.~Salamon, \emph{Quaternionic K\"ahler manifolds}, Invent.\ Math.\ \textbf{67} (1982), 143--171.


\bibitem{Siu76} Y.-T.~Siu, \emph{Every Stein subvariety admits a Stein neighborhood}, Invent.\ Math.\ \textbf{38} (1976), 89--100.

\bibitem{Stacks} The Stacks Project Authors, \emph{The Stacks Project}, \texttt{https://stacks.math.columbia.edu}.

\bibitem{Wis00} J.~A.~Wi\'sniewski, \emph{Lines and conics on Fano contact manifolds}, unpublished notes, January 6, 2000; note on its gaps, April 14, 2000,
\url{https://www.mimuw.edu.pl/~jarekw/preprints/note.html}.

\bibitem{Wolf65} J.~A.~Wolf, \emph{Complex homogeneous contact manifolds and quaternionic symmetric spaces}, J.\ Math.\ Mech.\ \textbf{14} (1965), 1033--1047.

\bibitem{Ye93} Y.-G.~Ye, \emph{Contact Fano manifolds, Wolf spaces and Severi varieties}, unpublished preprint, 1993, shared with a few colleagues at the time; updated version: \emph{Complete classifications of positive quaternionic-K\"ahler manifolds}, 2026 (section numbers refer to the updated version).

\bibitem{Ye93b} Y.-G.~Ye, \emph{Non-degeneracy of contact submanifolds in $\PP^{2n+1}$}, unpublished notes, March 1993, shared with a few colleagues at the time.

\bibitem{Ye94} Y.-G.~Ye, \emph{Extremal rays and null geodesics on a complex conformal manifold}, Internat.\ J.\ Math.\ \textbf{5} (1994), 141--168.

\bibitem{Ye94b} Y.-G.~Ye, \emph{A note on complex projective threefolds admitting holomorphic contact structures}, Invent.\ Math.\ \textbf{115} (1994),
311--314.

\bibitem{Ye26Q} Y.-G.~Ye, \emph{The curvature quartic of a contact Fano manifold and Freudenthal's magic square}, in preparation, 2026.

\bibitem{Zak93} F.~L.~Zak, \emph{Tangents and secants of algebraic varieties}, Transl. Math. Monographs \textbf{127}, Amer. Math. Soc., 1993.

\bibitem{Zak04} F.~L.~Zak, \emph{Determinants of projective varieties and their degrees}, in: Algebraic Transformation Groups and Algebraic Varieties,
Encyclopaedia Math.\ Sci.\ \textbf{132}, Springer, Berlin, 2004, 207--238.

\end{thebibliography}
\end{document}